\documentclass[reqno, 11pt, oneside]{amsart}

\usepackage[utf8]{inputenc}
\usepackage{amssymb} 
\usepackage{mathrsfs}  
\usepackage{amsmath}  
\usepackage{euscript} 
\usepackage{amsthm} 
\usepackage{mathrsfs} 
\usepackage{dsfont}  
\usepackage{amsfonts}     
\usepackage{marginnote}
\usepackage{latexsym} 
\usepackage{amsopn} 
\usepackage{geometry}
\usepackage{amscd} 
\usepackage{mathrsfs}
\usepackage{caption} 
\usepackage{aurical}  
\usepackage[english]{babel}
\usepackage[utf8]{inputenc}
\usepackage{esint} 
\usepackage{graphicx}
\usepackage[dvipsnames]{xcolor}
\usepackage{multicol}
\usepackage{nomencl}
\makeglossary
\makenomenclature  
\usepackage{refcount}
\usepackage{float}
\usepackage[titletoc]{appendix}
\usepackage{accents}
\usepackage{scalerel,stackengine}
\usepackage{tikz-cd}
\usepackage{footmisc}
\usepackage{setspace}
\usepackage{verbatim}

\stackMath
\newcommand\rwhat[1]{%
	\savestack{\tmpbox}{\stretchto{%
			\scaleto{%
				\scalerel*[\widthof{\ensuremath{#1}}]{\kern-.6pt\bigwedge\kern-.6pt}%
				{\rule[-\textheight/2]{1ex}{\textheight}}
			}{\textheight}%
		}{0.5ex}}%
	\stackon[1pt]{#1}{\tmpbox}%
}
\newlength{\dhatheight}

\newcommand{\set}[1]{{\left\{#1\right\}}}

\newcommand{\floornorm}[1]{\lfloor #1 \rfloor}

\newcommand{\eps}{\varepsilon}

\newcommand{\s}{{\sigma}}

\newcommand{\floorNorm}[1]{\left| \! \left\lfloor #1 \right\rfloor \! \right|}
\newcommand{\Floornorm}[1]{\bigg\lfloor #1 \bigg\rfloor}

\renewcommand{\Im}{\operatorname{Im}}

\newcommand{\id}{\operatorname{Id}}

\newcommand{\co}[1]{\textit{#1}}

\usepackage{hyperref}
\hypersetup{
colorlinks=true,
linkcolor=blue,
citecolor=LimeGreen,
urlcolor=blue,
linktoc=all
}

\newtheorem{thm}{Theorem}[section]
\newtheorem*{thm*}{Theorem}
\newtheorem{prop}[thm]{Proposition}
\newtheorem{lemma}[thm]{Lemma}

\newtheorem{cor}[thm]{Corollary}
\newtheorem*{cor*}{Corollary}
\newtheorem{rmk}[thm]{Remark}
\newtheorem{ex}{Example}

\newtheorem{defn}[thm]{Definition}
\newtheorem{notation}[thm]{Notation}

\numberwithin{equation}{section}

\newcommand{\om}{{\omega}}

\newcommand{\C}{{\mathbb C}}

\newcommand{\N}{{\mathbb N}}

\newcommand{\R}{{\mathbb R}}

\newcommand{\T}{{\mathbb T}}
\newcommand{\Z}{{\mathbb Z}}

\newcommand{\cA}{{\mathcal A}}

\newcommand{\tc}{{\mathtt{c}}}

\newcommand{\tm}{{\mathtt{m}}}

\usepackage{bm}

\newcommand{\im}{{\rm i}}
\newcommand{\jap}[1]{\langle #1 \rangle}

\newcommand{\nnorm}[1]{{\left\vert\kern-0.25ex\left\vert\kern-0.25ex\left\vert #1 
		\right\vert\kern-0.25ex\right\vert\kern-0.25ex\right\vert}}

\usepackage{framed,enumitem}

\newcommand{\norm}{{|\hspace{-0.075cm}|\hspace{-0.075cm}|}} 
\newcommand{\ldbracket}{[\![}
\newcommand{\rdbracket}{]\!]}
\newcommand{\normS}[1]{\ldbracket #1 \rdbracket}

\usepackage{mathtools}
\mathtoolsset{showonlyrefs}

\newcounter{paraga}[subsubsection]
\renewcommand{\theparaga}{{\bf\arabic{paraga}.}}
\newcommand{\paraga}{\medskip \addtocounter{paraga}{1}
	\noindent{\theparaga\ } }

\makeatletter

\def\l@subsection{\@tocline{2}{0pt}{2.5pc}{5pc}{}}
\def\l@subsubsection{\@tocline{3}{0pt}{4.5pc}{5pc}{}}
\renewcommand\tocchapter[3]{%
	\indentlabel{\@ifnotempty{#2}{\ignorespaces#2.\quad}}#3%
}

\author{Lorenzo Baroni}
\address{\scriptsize{Dipartimento di Matematica e Fisica, Universit\`a degli Studi di Roma Tre, 
		Largo San L. Murialdo 1, 00146 Roma}}
\email{lorenzo.baroni@uniroma3.it}

\author{Jessica Elisa Massetti}
\address{\scriptsize{Dipartimento di Matematica, Universit\`a degli Studi di Roma ``Tor Vergata", 
		Via della Ricerca Scientifica 1, 00133 Roma}}
\email{massetti@mat.uniroma2.it}

\author{Michela Procesi}
\address{\scriptsize{Dipartimento di Matematica e Fisica, Universit\`a degli Studi di Roma Tre, 
		Largo San L. Murialdo 1, 00146 Roma}}
\email{michela.procesi@uniroma3.it}

\date{\today}

\begin{document}

\title[Liouville-Arnold integrability in infinite dimension]{Foliations by Lagrangian tori and Liouville-Arnold integrability in infinite dimension}
	
\begin{abstract}
	We provide a first generalization of the Liouville--Arnold theorem to the infinite-dimensional setting, with particular emphasis on integrable Hamiltonian systems on infinite lattices and PDEs on the circle. We work on the phase space $\ell^\infty_w\times \T^\infty$, where $w$ is a sequence of weights that need not decay, and consider families of Hamiltonians of the form $y_k+P_k(y,\varphi)$ in involution. In this setting, under suitable regularity and smallness conditions, we prove that the common level sets form a foliation by Lagrangian tori and explicitly construct analytic action-angle variables. As a first application, we consider commuting families of Hamiltonians on infinite lattices near an elliptic fixed point. We then study the KdV equation near zero in weighted $\ell^\infty$ spaces. In the open region where all spectral gaps are open, we prove Liouville--Arnold integrability providing a foliation by invariant Lagrangian tori.
\end{abstract}
	
	\maketitle
	\setcounter{tocdepth}{2}
\tableofcontents

\section{Introduction}

In the present article we provide a first generalization of Liouville-Arnold theorem on foliations by invariant Lagrangian tori for completely integrable systems to the infinite dimensional setting, with particular emphasis on systems on infinite lattices and PDEs on the circle. Although several results concerning the complete integrability of infinite-dimensional Hamiltonian systems have been established (see e.g., \cite{Zakh-Fadd71} for seminal results on the KdV, \cite{Ablowitz81} for the IST method, or \cite{kappeler2003kdv} for a modern Birkhoff coordinates approach), none of them provides a decomposition of the phase space as rigid as the one provided by the Liouville-Arnold theorem in the finite-dimensional case. This is largely due to the choice of topology for the phase space, a crucial issue that will be discussed in detail later on.

In this context, the present work aims to fill this gap. To the best of our knowledge, the following results provide the first instance of analytic invariant foliations by infinite-dimensional tori for infinite-dimensional Hamiltonian systems, alongside the first construction of action-angle coordinates that are analytic in an open neighborhood of the phase space (see Theorems \ref{completethm} and \ref{COMPL}). These results establish a more faithful analogue of the Liouville-Arnold theorem in the infinite-dimensional setting. 

\medskip
Specifically, this work addresses the following questions:
\begin{itemize}[leftmargin=*]
	\item \textbf{Foliation into infinite-dimensional tori (Theorem \ref{FFoliation}):} \emph{When do the level sets of a family of analytic Hamiltonians induce a regular foliation of an open subset of the phase space into infinite-dimensional invariant tori?}
	
	\item \textbf{Action-Angle variables (Theorem \ref{SYMPL}):} \emph{Assuming the family of Hamiltonians is in involution, what additional analytical regularity assumptions are required to construct a set of action-angle coordinates adapted to the foliation?}
	
	\item \textbf{Liouville-Arnold integrability (Theorem \ref{completethm}):}
	\emph{When the integrals of motion induce both a regular foliation and a set of action-angle coordinates, is it possible to establish an infinite-dimensional analogue of the Liouville-Arnold theorem for a given analytic Hamiltonian commuting with the family?}	
\end{itemize}

We give positive answers to the above questions under appropriate hypotheses on the  topology of the phase space, on the commuting Hamiltonians and their vector fields. In particular, in order to have a submanifold structure, we need our tori to be non compact. As a consequence our phase space will be modeled on weighted $\ell^\infty$ spaces.  As a drawback, the argument of Liouville that identifies regular compact level sets with tori fails to work. To bypass this problem, we restrict to a specific class of first integrals, which are ``close" to the linear actions, so that we can use perturbative arguments. Of course this setting can be generalized, but we believe that our result represents a significant step forward, because it contains most of the technical difficulties at the analytic level.

Of course, if on the one hand these results are interesting per se as they bring to the infinite dimensional setting a more intrinsic and geometric view point, on the other hand it is important to verify wether systems that are known to possess a complete set of constants of motion fit our framework. 

We wish to transfer many finite-dimensional arguments to our setting and this require rather strong regularity assumptions. However, stronger regularity conditions imply less applicability; this requires a delicate balance. The general frame will be majorant analytic functions of infinitely many angles and actions with uniformly bounded analyticity radius (see Definition \ref{matcalA}). We highlight that this is a much weaker condition with respect to \cite{procmont} and~\cite{cgpro}, and it is a major source of difficulties.  

With this in mind, it is natural to consider either lattice systems with an elliptic fixed point at zero, or even more specifically, Lax integrable PDEs on the circle. 

\begin{itemize}[leftmargin=*]
	\item \textbf{Applications to the Korteweg-deVries equation (Theorem \ref{LA_KdV}):} \emph{we prove that there exists an open and dense set of sufficiently small and regular initial data, whose evolution is supported on a leaf of a local analytic foliation by invariant Lagrangian tori. Moreover, in such ``tubular neighborhoods" we construct symplectic action-angle variables.}
	\end{itemize}

We consider KdV as a paradigmatic example (the same discussion applies verbatim to NLS).
Indeed, classical results show that the finite-gap solutions of KdV form a dense subset of $L^2$ and lie on finite-dimensional invariant tori. Moreover, all solutions lie on compact invariant tori that are limits of finite-gap tori.
This structure, however, is incompatible with a geometric foliation.
In our setting, the union of the finite-gap tori fails to be dense due to the nonseparability of the $\ell^\infty$-based phase space.
In our proof, we use the same kind of argument in the sense that our starting point is the map from the initial data to the spectral gaps.
However, we work on an open neighbourhood in which all the spectral gaps are open.
The first issue to address is to prove that this map is analytic on that neighbourhood.
A parallel approach, pursued in \cite{kappeler2003kdv}, is to construct Birkhoff variables for KdV.
As a by-product, one obtains a full description of the invariant tori.
The Birkhoff map is known to be analytic on the Sobolev spaces $H^s(\T)$ for all $s \ge 0$ (see \cite{KMT2005} for the extension to $H^{-1}(\T)$ in the case of KdV).
While there are results on this map in $\ell^\infty$-based spaces of very low regularity (see \cite{KM}), proving our result by showing that the Birkhoff map is analytic on these spaces seems far from trivial.

\vspace{.5 cm}
In the following we shall contextualize the above questions to our setting and give a more precise description of the results.  
	
\subsection{Foliations by infinite tori and action-angle coordinates}\label{chfol}
Let us first rapidly review our strategy in the finite dimensional frame. 
Let $\mathbb{T}^n:=\mathbb{R}^n/2\pi\mathbb{Z}^n$ be the $n$-dimensional torus, $n\geq 1$, and consider its cotangent bundle $T^*\mathbb{T}^n \simeq\mathbb{R}^n \times \mathbb{T}^n$.  Let $\mathscr{U}$ be an open neighborhood of $\set{0}\times \mathbb{T}^n$, and let $F: \mathscr{U} \to \R^n$ be a real analytic function. 

In the following we shall say that \co{the level sets of $F$ give rise to a real analytic foliation of a neighborhood of $\set{0}\times \mathbb{T}^n$ by $n$-dimensional tori} if there exist
 $c\in\R^n$ and $r>0$ such that $\mathscr{N}_0 := F^{-1}(\mathscr{U}_{r,c}(\R^n)) \subset \mathscr{U}$\footnote{Given a normed space $(X,|\cdot|_X)$, the open ball of radius $r$ centered at the origin is 
 \begin{equation*}
     \mathscr{U}_r(X):=\{x\in X : |x|_X<r\}.
 \end{equation*}
Its closure is denoted by $\overline{\mathscr{U}}_r(X)$.  For $X\subset \C^\Z$ and $x_0 \in X$, the open ball centered at $x_0$ is:
	\begin{equation*}
	\mathscr{U}_{r,x_0}(X) := x_0 + \mathscr{U}_{r}(X) = \{x\in X : |x - x_0|_X < r \}.
	\end{equation*}
	The closure is denoted by $\overline{\mathscr{U}}_{r,x_0}(X)$.} is an open neighborhood of $\set{0}\times \mathbb{T}^n$ and there exists a bi-analytic diffeomorphism $\chi: \mathscr{N}_0 \to \mathscr{U}_{r,c}(\R^n)\times \mathbb{T}^n$ such that the following diagram commutes

\begin{center}
			\begin{tikzcd}[column sep=large]
				\mathscr{N}_0 \arrow[d, "\textstyle F"'] \arrow[r, "\textstyle \chi"] & \mathscr{U}_{r,c}(\R^n)\times \mathbb{T}^n \arrow[ld, "\textstyle p_1"] \\
				\mathscr{U}_{r,c}(\R^n) & 
			\end{tikzcd}
		\end{center}
where $p_1$ denotes the projection onto the first factor. In particular, $F$ also defines a trivial analytic fiber bundle over 
$\mathscr{U}_{r,c}(\R^n)$ with typical fiber analytically diffeomorphic to $\mathbb{T}^n$.

Let us define
\begin{equation*}
\Gamma_0:\mathbb{R}^n \times \mathbb{T}^n \to \mathbb{R}^n \, , \quad (y,\varphi) \mapsto y.
\end{equation*}
In this line of thought, it is clear that the level sets of $\Gamma_0$ give rise to a global real analytic foliation with leaves $\Gamma_0^{-1}(\gamma) \cong \mathbb{T}^n$, $\forall \gamma \in \mathbb{R}^n$, with $\chi = \id$. By a straightforward application of a contraction lemma one can prove the robustness of such a trivial foliation under sufficiently regular perturbations. More precisely, let 
\begin{equation*}
P: \mathscr{U}_r(\mathbb{R}^n)\times \mathbb{T}^n \to \R^n
\end{equation*}
be a real analytic map, sufficiently small with respect to some appropriate norm, 
and define
\begin{equation*}
\Gamma_P(y,\varphi)=y+P(y,\varphi).
\end{equation*}
Then, the level sets of $\Gamma_P$ give rise to a real analytic foliation of a neighborhood of $\set{0}\times \mathbb{T}^n$ by $n$-dimensional tori, for an open set of ``sufficiently small'' values of $\Gamma_P$.

In Theorem \ref{FFoliation}, we extend this result to an infinite-dimensional framework. To this end, we must first specify the setting carefully (see Section \ref{fsetting}), we choose a  topology that allows us to perform differential calculus and to deal with functions in normed spaces where suitable Cauchy estimates for the differentials can be provided. More precisely we work in a phase space which is a product of an infinite torus with the uniform topology, denoted $\T^\infty$, and a weighted $\ell^\infty$ space denoted by $\ell^\infty_w$. This result is obtained through a fixed point argument; the main difficulty is to introduce the correct functional setting of majorant analytic functions from $\ell^\infty_w\times \T^\infty$ to $ \ell^\infty_w$ admitting a holomorphic extension in the angles in the uniform domain $|\Im\varphi_j| < s$. Note that this is a much weaker condition with respect to the previous literature (e.g. \cite{Corsi20211231,procmont,cgpro}), and a source of serious technical challenges throughout the paper.  


\subsection*{Action-Angle coordinates in infinite dimension}
The second question we address is whether there exists a symplectic straightening of the foliation by infinite tori provided by Theorem \ref{FFoliation}, possibly under a further smallness condition on the perturbation in a stronger norm (see Definition \ref{definitionnorm}) and assuming that the components of the map $\Gamma_P$ Poisson-commute.

To this end, in Section \ref{SecRHam} we develop a functional framework in which we introduce an infinite-dimensional Poisson structure together with the notion of a symplectomorphism. The straightening map (see Theorem \ref{SYMPL}) induces a new coordinate system $(I, \theta)$ in which $I$ depends on the original variables $(y, \varphi)$ only through $\Gamma_P$. These coordinates, which we refer to as \emph{action-angle} coordinates, are adapted to the foliation. The coordinates $I = (I_k)_{k \in \mathbb{N}}$ label the leaves, while $\theta = (\theta_k)_{k \in \mathbb{N}}$ provides coordinates on each leaf. To construct these coordinates, we extend the classical generating function method (see \cite{arnold1978mathematical} for a more detailed discussion) to an infinite-dimensional setting by adapting the standard Hamiltonian machinery to our context. In infinite dimensions, one must carefully check that all the algebraic manipulations are justified by quantitative estimates. It is here that our functional setting, in particular the notions of majorant analyticity and regularity (see Definitions \ref{RF}, \ref{definitionnorm} and \ref{regchange}), plays an essential role.

Theorem \ref{completethm} is the dynamical counterpart of this result and states that any analytic Hamiltonian that Poisson-commutes with every component of $\Gamma_P$ depends only on the actions. To the best of our knowledge, this is the first result to give a geometric decomposition of phase space analogous to that of finite-dimensional completely integrable Hamiltonian systems. 

To conclude this subsection, we give the informal statement of Theorem \ref{completethm}.

\begin{thm*}[Liouville--Arnold integrability on $\ell^\infty_w \times \T^\infty$, informal statement]
Let $r>0$ and let $H$ be a majorant analytic Hamiltonian $H:\mathscr{U}_r(\ell^\infty_w)\times \T^\infty \to \R$. Assume that $H$ admits a family of first integrals $\Gamma_P^{(k)}=y_k+P_k(y,\varphi)$ in  involution. If the $P^{(k)}$'s satisfy an appropriate uniform bound and a summably analytic condition on the associated family of Hamiltonian vector fields, then the following holds:
	\begin{itemize}
	\item[\textnormal{(i)}] The common level sets of the $\Gamma^{(k)}_P$ form an analytic foliation by invariant Lagrangian tori $\T^\infty$ on an open neighbourhood of $\{0\}\times \T^\infty$.
	\item[\textnormal{(ii)}] There exist analytic action-angle coordinates $(I,\theta)$ adapted to this foliation. In these coordinates $H=h(I)$ and the equations of motion reduce to
	\[
	\dot{I}_k=0,\qquad
	\dot{\theta}_k=\partial_{I_k}h(I),
	\qquad k\in\N.
	\]
\end{itemize}
\end{thm*}

\subsection{Lattice systems}

There are examples of integrable Hamiltonian systems on lattices that admit hierarchies of commuting Hamiltonians. For the periodic Toda lattice (see \cite{Toda89}) with any fixed number of particles, in \cite{Henrici_2008} the authors construct global Birkhoff coordinates on each symplectic leaf, in which the Toda Hamiltonian depends only on the actions.

Here we consider majorant analytic Hamiltonians $H$ on the phase space
\begin{equation*}
	\ell^\infty_{\eta}(\C):=\{ \, u\in \C^\Z \, : \, \sup_{k \in \mathbb{N}} \eta_k|u_k|<\infty \,  \}, \qquad \eta_k>0,
\end{equation*}
admitting a family of sufficiently regular first integrals in involution with respct to the standard Poisson bracket. We also assume that they have the form
\begin{equation}\label{form}
\Gamma^{(k)}(u)=|u_k|^2+f^{(k)}(u), \qquad f^{(k)}(u)=O(|u|^3),
\qquad k\in\N.
\end{equation}
To relate this setting to the results of the preceding subsections, fix a sequence $\xi=(\xi_k)_{k\in\N}$ such that
\[
0<\inf_{k\in\N}\eta_k^2\xi_k
\leq \sup_{k\in\N}\eta_k^2\xi_k<\infty,
\]
and consider the torus $\mathcal{T}_\xi=\{u:|u_k|^2=\xi_k,\ k\in\N\}$.
Setting $w_k=\eta_k^2$, the symplectic polar coordinates
\[
u_k=\sqrt{\xi_k+y_k}\,e^{\im\varphi_k},
\qquad k\in\N,
\]
identify a tubular neighborhood of $\mathcal{T}_\xi$ with
$\mathscr{U}_r(\ell^\infty_w)\times\T^\infty$ for
$0<r<\inf_k w_k\xi_k$. Under this identification, denoted by
$\Phi_\xi$, the integrals take the form
\begin{equation}\label{phiintro}
\Gamma^{(k)}\circ\Phi_\xi-\xi_k
=y_k+P_k(y,\varphi),
\qquad P_k=f^{(k)}\circ\Phi_\xi.
\end{equation}
The notion of majorant analyticity we adopt ensures that the maps in \eqref{phiintro} are majorant analytic and in particular they admit a holomorphic extension on a complex domain whose analyticity radius is uniform with respect to the angles, therey fitting the framework of the previously established abstract result. This property is particularly relevant for lattice systems. In Definition~\ref{matcalA}, the Fourier weight is $e^{s|\nu|_1}$, where $|\nu|_1=\sum_j|\nu_j|$. Thus the Fourier mode $e^{\im(\varphi_{k+1}-\varphi_k)}$, for instance, carries the same weight $e^{2s}$ at every lattice site and is translation invariant. On the contrary all the literature consider the much stronger weight  $e^{s|\nu|_*}$, where $|\nu|_*=\sum_j h_j|\nu_j|$ and $h_j\to\infty$ as $j\to\infty$ (cf.~\cite{Corsi20211231,procmont,cgpro}). This condition would simplify the proof in action-angles, but then it would not be possible to apply it to lattice systems, unless one imposes some really unnatural strong hypothesis on the system itself. 

Let us rephrase informally the integrability result for Hamiltonian lattice systems. For the precise statement, see Theorem \ref{COMPL}.

\begin{thm*}[Liouville--Arnold integrability on an infinite lattice, informal statement]
	Let $H$ be majorant analytic on a ball of radius $R$ in $\ell^\infty_\eta(\C)$ and assume it admits a family of commuting first integrals $\Gamma^{(k)}$ of the form \eqref{form}. Assume in particular that the corrections $f^{(k)}$'s satisfy an appropriate uniform bound  and a summably analytic condition on the family of Hamiltonian vector fields. Then there exist $0<\xi_-\leq\xi_+<R^2$ such that, for every $\xi$ satisfying
	\[
	\xi_-\leq\eta_k^2\xi_k\leq\xi_+,
	\qquad k\in\N,
	\]
	the following hold in a neighborhood of $\mathcal{T}_\xi$.
	\begin{itemize}
		\item[\textnormal{(i)}] The common level sets of the $\Gamma^{(k)}$ form an analytic foliation by invariant Lagrangian tori $\T^\infty$ locally modeled on $\ell^\infty$.
		\item[\textnormal{(ii)}] There exist analytic action-angle coordinates $(I,\theta)$ adapted to this foliation. In these coordinates $H=h(I)$ and the equations of motion reduce to
		\[
		\dot{I}_k=0,\qquad
		\dot{\theta}_k=\partial_{I_k}h(I),
		\qquad k\in\N.
		\]
	\end{itemize}
\end{thm*}

This result has to be compared with the infinite-dimensional Vey theorem of Kuksin and Perelman~\cite{KuksinPerelman2010} and the normal form theorem of Bambusi and Stolovitch~\cite{BaSt} for summably analytic families of commuting vector fields. Both results concern analytic Birkhoff nromal forms in a neighborhood of an elliptic equilibrium in a Hilbert space. For the construction of action-angle coordinates in Theorem \ref{COMPL}, however, it suffices to verify the hypotheses of Theorem~\ref{completethm} directly in symplectic polar coordinates on a tubular neighborhood of $\mathcal{T}_\xi$. This is possible because the condition $\inf_k\eta_k|u_k|>0$ defines an open subset of $\ell^\infty_\eta(\C)$. In a Hilbert sequence space, finitely supported sequences are dense, so every neighborhood meets singularities of the symplectic polar coordinates.Therefore, in the weighted $\ell^\infty$ topology we can, in principle, establish the foliation near $\mathcal{T}_\xi$ without requiring the first integrals to extend analytically to a ball about the origin. We use precisely this observation in the application to KdV below.

\subsection{Application to KdV}
We consider the Korteweg-de Vries (KdV) equation with periodic boundary conditions for a real-valued function $u$ with zero mean:
\begin{equation}\label{kdv}
u_t + u_{xxx} - 6uu_x = 0, \quad x \in \R /\Z, \quad \int_0^1 u(x,t) dx = 0.
\end{equation}
The condition $\int_0^1 u(x,0) dx = 0$ is preserved by the KdV flow. Indeed, since $6uu_x - u_{xxx} = (3u^2 - u_{xx})_x$, the total mass is a conserved quantity. By multiplying both sides of the equation by $u$ and integrating, one can also find the conserved Hamiltonian:
\begin{equation}\label{ham_u}
H(u) = \int_0^1 \left( \frac{u_x^2}{2} + u^3 \right) dx.
\end{equation}
In \cite{Gardner71}, Gardner showed that $L^2_0(\R/\Z, \R)$\footnote{We denote by $L^2_0(\R/\Z, \R)$ the subspace of $L^2(\R/\Z, \R)$ of functions with zero mean.} carries a Poisson structure defined by the bracket
\begin{equation}\label{poissongardner}
\{F, G\} := \int_{0}^1 \nabla_{L^2} F \cdot \partial_x \nabla_{L^2} G \, dx,
\end{equation}
where $\nabla_{L^2}$ denotes the $L^2$-gradient, making KdV an infinite-dimensional Hamiltonian system with Hamiltonian vector field $V_H = \partial_x \nabla_{L^2} H$.

KdV is Lax integrable in the sense that the spectrum of the Schr\"odinger with potential $u$ is conserved by the flow. It is well known that a good set of constants of motion is the family of the squared gap lengths, and there are several integrablity results on this PDE in large part with initial data in Sobolev-Hilbert spaces $H^s$ (see, for e.g., \cite{kappeler2003kdv}). 

Our goal is to prove the Liouville-Arnold integrability of the KdV equation by showing that the map constructed from the complete list of the squared gaps $\gamma_n^2$ expressed in symplectic polar coordinates, satisfies the requirements of the Complete Integrability Theorem \ref{completethm}. 

Note that the symplectic form induced on the space of Fourier coefficients by the Gardner Poisson brackets is not the standard Darboux one, but is given by $\im\sum_{k \in \mathbb{N}} k^{-1} du_k \wedge du_{-k}$. To recover the canonical structure, we first perform the change of variables
\begin{equation}
\quad   z_\ell=\langle \ell \rangle^{-1/2}u_\ell, \qquad \jap{\ell}:=\max\{1,|\ell|\}.
\end{equation}
Clearly, since $u$ is real-valued and has zero mean, we have $z_{-\ell}=\bar{z}_\ell$ for all $\ell \in \Z$, and $z_0=0$. In these coordinates, the symplectic form assumes the standard Darboux form $\om = \im \sum_{\ell \in \mathbb{N}} dz_\ell \wedge dz_{-\ell}$, and the KdV equation can be rewritten as an infinite system of coupled oscillators:
\begin{equation}\label{lattice_kdv_intro}
\dot{z}_k = \im k^3 z_k + 3\im k \sum_{\substack{\ell+m=k \\ \ell,m \neq 0}} \sqrt{|\ell m|} \, z_\ell z_m, \quad k \in \mathbb{Z} \setminus \{0\},
\end{equation}
where the corresponding Hamiltonian is given by:
\begin{equation}
H_{KdV}(z) = \sum_{k \in \mathbb{N}} k^3 |z_k|^2 + \sum_{\substack{k+\ell+m=0 \\ k,\ell,m \neq 0}} \sqrt{|k \ell m|} \, z_k z_\ell z_m,
\end{equation}
and the phase space is
\begin{equation*}
	\ell^\infty_{a,p}:=\{ \, z\in \ell^2(\C) \, : \, \sup_{k \in \mathbb{N}} e^{a|k|}\jap{k}^{p}|z_k|<\infty \,  \}, \qquad p>2, a \ge 0.
\end{equation*}

\begin{thm}[Liouville-Arnold Integrability of KdV on $\T$]\label{LA_KdV}
	Let $p>2$ and $a \geq0$. Consider the equation \eqref{lattice_kdv_intro}. There exist $0<\xi_-\leq\xi_+$ such that, for any sequence $\xi = (\xi_k)_{k \in \mathbb{N}}$ of positive real numbers satisfying
	\begin{equation}\label{grazie}
		\xi_-\leq \inf_{k \in \mathbb{N}} e^{2a|k|}\jap{k}^{2p} \xi_k \leq \sup_{k \in \mathbb{N}} e^{2a|k|}\jap{k}^{2p} \xi_k \leq \xi_+,
	\end{equation}
	the following hold in a neighborhood of the torus $\mathcal{T}_\xi = \{ z \in \ell^\infty_{a,p} : |z_k|^2 = \xi_k \}$.
	\begin{itemize}
		\item[\textnormal{(i)}] \textbf{Foliation by Infinite Invariant Tori:} The phase space is locally foliated by analytic invariant tori analytically diffeomorphic to $\T^\infty$. Each torus is a level set of the spectral gaps and supports an almost-periodic Hamiltonian flow.
		\item[\textnormal{(ii)}] \textbf{Action-Angle Coordinates:} There exists an analytic symplectomorphism 
		$$\Psi: z \mapsto (I, \theta)$$
		that introduces action-angle coordinates.
		\item[\textnormal{(iii)}] \textbf{Integration by Quadratures:} In these coordinates, the Hamiltonian $H_{KdV} \circ \Psi^{-1}$ depends only on the actions $I$, and the Hamilton's equations reduce to the linear flow:
		\begin{equation}
			\begin{cases}
				\dot{I}_k = 0 \\
				\dot{\theta}_k = \omega_k(I),
			\end{cases} \quad k \in \mathbb{N},
		\end{equation}
		where the frequencies are given by $\omega_k(I) = \partial_{I_k} (H_{KdV} \circ \Psi^{-1})$.
	\end{itemize}
\end{thm}

Our strategy is to construct the map $z \to (\gamma^2_n(z))_{n\in\N}$ via perturbative argument. A first challenge would be to prove analyticity of the above map at zero, as is well known in $\ell^2$-based spaces. In fact, the condition \eqref{grazie} allows us to bypass this problem completely and study the squared gaps directly in symplectic polar coordinates.

In this setting, the required regularity is the so called \textquotedblleft summable analyticity". More precisely, we impose a smallness condition on an appropriated norm of the sum of the Hamiltonian vector fields associated to the squared gap lengths.


\subsection*{Acknowledgments} 
J.E. Massetti has been supported by the INdAM GNAMPA Project “Growth of Sobolev norms and energy cascades in dispersive Hamiltonian PDE ” - CUP E53C25002010001 and acknowledges also the support of the Department of Excellence grant MatMod@TOV (2023- 27), awarded to the Department of Mathematics at University of Rome Tor Vergata.

\section{The phase space and properties of majorant analytic functions and conjugacies}\label{fsetting}

\subsection{The phase space $\ell^\infty_w\times\T^\infty$}
Let us start by introducing the main objects that define the phase space.

\subsubsection{The infinite torus} We consider the infinite-dimensional complexified torus\footnote{Throughout this paper, \(\N\) denotes the set of positive natural numbers \(\{1,2,\dots\}\), while \(\N_0\) denotes the set of natural numbers \(\{0,1,2,\dots\}\).}
\begin{equation}
    \T^\infty_\C:=\bigg\{ \, \varphi=(\varphi_j)_{j \in \N}, \, \varphi_j \in \C/2\pi\Z \, : \, \sup_{j \in \N}|\mathrm{Im}(\varphi_j)|<\infty \, \bigg\}
\end{equation}
endowed with the topology induced by the uniform distance
\begin{equation}
d_{\T^\infty_\C}(\varphi,\varphi'):=\sup_{j \in \N}\inf_{n \in \Z}|\varphi_j-\varphi'_j+2\pi n|.
\end{equation}
With this topology, \(\T^\infty_\C\) has a natural structure of a complex-analytic abelian Banach--Lie group modeled on \(\ell^\infty(\N,\C)\). We remark that \(\T^\infty_\C\) is naturally identified with the quotient Banach--Lie group \(\ell^\infty(\N,\C)/\ell^\infty(\N,2\pi\Z)\), and the quotient map \(\ell^\infty(\N,\C) \to \T^\infty_\C\) is the universal analytic covering map onto \(\T^\infty_\C\).

We define the \emph{infinite torus} \(\T^\infty\) as the topological subgroup of \(\T^\infty_\C\) consisting of those \(\varphi\) such that \(\mathrm{Im}(\varphi_j)=0\) for all \(j\in\N\). Clearly, \(\T^\infty\) is a real-analytic abelian Banach--Lie group modeled on \(\ell^\infty(\N,\R)\). Indeed, \(\T^\infty\) is analytically diffeomorphic to \(\ell^\infty(\N,\R)/\ell^\infty(\N,2\pi\Z)\), and the quotient map
\[
\pi_\infty:\ell^\infty(\N,\R) \to \T^\infty
\]
is the universal analytic covering map onto \(\T^\infty\).

Finally, given \(s>0\), we denote by \(\T^\infty_s\) the thickened infinite torus, i.e., the complex neighborhood of \(\T^\infty\) in \(\T^\infty_\C\) defined by the condition \(\sup_{j\in\N}|\mathrm{Im}(\varphi_j)|<s\).

In this paper, we shall consider real-analytic functions on \(\T^\infty\) that admit a holomorphic extension to \(\T^\infty_s\) for some \(s>0\). By the above discussion, for each such \(s\), such functions correspond to holomorphic functions on the domain 
\[
\ell^\infty_s:=\bigg\{ \, x \in \ell^\infty(\N,\C) \, : \, \sup_{j\in\N}|\mathrm{Im}(x_j)|<s \, \bigg\}
\] 
that are invariant under translations by \(\ell^\infty(\N,2\pi\Z)\).

We emphasize that there are other possible infinite-dimensional generalizations of the torus $\T^n$, with $n\in\N$, that one could consider. For example, if one equips the infinite product of circles with the product topology, there would be no way to introduce a topological manifold structure on the resulting infinite torus. Roughly speaking, the cylinders form a basis for the product topology, and non-contractible loops can be found in each cylinder. This contrasts with the local simple connectedness required for any topological manifold.

\subsubsection{The phase space} Given a sequence of positive real numbers $w=(w_j)_{j \in \mathbb{N}}$, we define the real weighted Banach space $\ell^\infty_w$ of conjugate momenta as
\begin{equation*}
\ell^{\infty}_w:=\bigg\{y \in \mathbb{R}^\mathbb{N} \ : \ |y|_{\ell^{\infty}_w}:=\sup_{j \in \mathbb{N}}w_j |y_j|<\infty \bigg\}.
\end{equation*}
In the following, we shall denote by $\ell^\infty_\C(w)$ the complex analogue of $\ell^\infty_w$.

In this paper, we study the integrability of infinite-dimensional Hamiltonian systems on the phase space
\begin{equation}
    \ell^\infty_w\times \T^\infty.
\end{equation}

 It is worth stressing that the integrability result obtained in this paper does not depend on the specific choice of the weight sequence $w$. From a strict geometric viewpoint, defining the weak symplectic structure
 \begin{equation}
 \sum_{k \in \N} dy_k \wedge d\varphi_k
 \end{equation}
 would require that $\ell^\infty_w$ continuously embeds into $T^*_\varphi \T^\infty=(\ell^\infty(\N,\R))^*$. In that case, our phase space would be a subbundle of the cotangent bundle of the infinite torus $\T^\infty$. For this embedding to hold, the sequence formed by the inverse weights must be summable, meaning that 
 \begin{equation}
 \tc(w):=\sum_{k \in \N} w_k^{-1} < \infty.
 \end{equation}
 However, our approach bypasses this restrictive geometric constraint by relying on a robust Poisson structure. The Hamiltonian structure on $\ell^\infty_w\times \T^\infty$ can be introduced for any weight $w$ by defining the Poisson bracket
 \begin{equation}
     \{f,g\}=\sum_{k \in \N}\frac{\partial f}{\partial y_k}\frac{\partial g}{\partial \varphi_k}-\frac{\partial f}{\partial \varphi_k}\frac{\partial g}{\partial y_k}
 \end{equation}
 on a suitable subalgebra of the algebra of smooth functions (see Subsection \ref{RHPS}), which we call the subalgebra of \emph{regular Hamiltonians}, and subsequently defining the Hamiltonian vector field $V^f$ associated with any function $f$ in this subalgebra as the unique vector field that satisfies the relation $\mathscr{L}_{V^f}g=\{g,f\}$, where $\mathscr{L}_{V^f}g$ denotes the Lie derivative of $g$ along~$V^f$.

\subsection{Space of functions}\label{Sof}
Let us now fix a weight $w$, namely a sequence of positive real numbers $w=(w_j)_{j \in \mathbb{N}}$, once and for all. In this subsection, we introduce and discuss some relevant properties of a class of real analytic functions on $\ell^\infty_w \times \mathbb{T}^\infty$. 

From now on we adopt the following notations:
\begin{equation*}
    \ell^\infty:=\ell^\infty(\N,\R), \qquad \ell^\infty_\C:=\ell^\infty(\N,\C).
\end{equation*}

 We consider analytic functions on infinite tori which are represented by trigonometric series on the universal covering $\ell^\infty$ of $\mathbb{T}^\infty$, in the sense explained in the previous subsection. In the following, we shall adopt the notation $\varphi$ not to denote elements of $\mathbb{T}^\infty$, but instead to refer to elements of the universal covering $\ell^\infty$. 

We consider formal power-trigonometric series of the type
\begin{equation}\label{fourier0}
	f(y,\varphi)=\sum_{\alpha \in \mathbb{N}_0^{(\mathbb{N})}}\sum_{\nu \in \mathbb{Z}^{(\mathbb{N})}}f_{\alpha,\nu}y^\alpha e^{\im \nu \cdot \varphi}, \qquad y \in \ell^\infty_w, \,  \varphi \in \ell^\infty,
	\end{equation}
with $f_{\alpha,\nu} \in \mathbb{C}$, $f_{\alpha,\nu}= \bar f_{\alpha,-\nu}$ for all $(\nu, \alpha) \in \mathbb{Z}^{(\mathbb{N})}\times \mathbb{N}_0^{(\mathbb{N})}$. Given $X,Y \subseteq \Z$, $X^{(Y)}$ denotes the set of sequences with finite support
	\begin{equation*}
	X^{(Y)}:=\bigg\{(\lambda_k)_{k \in Y} \ : \ \lambda_k \in X \ \forall k \in Y, \ |\lambda|:=\sum_{k \in Y}|\lambda_{k}|<\infty \bigg\}.
	\end{equation*}
    In \eqref{fourier0} we have used the standard multi-indices notation for a list $z=(z_k)_{k \in Y}\in \C^Y$ and $\lambda \in X^{(Y)}$:
	\begin{equation*}
	z^\lambda := \prod_{k \in Y}z_k^{\lambda_k}.
	\end{equation*}
To ease the notation we set
\begin{equation}
    \sum_{\alpha \in \mathbb{N}_0^{(\mathbb{N})}}\sum_{\nu \in \mathbb{Z}^{(\mathbb{N})}}=\sum_{\alpha,\nu}
\end{equation}
and we denote the monomials also by
\begin{equation}
\tm_{\alpha,\nu}(y,\varphi)=y^\alpha e^{\im \nu\cdot \varphi}.
\end{equation}
Finally, given $r,s>0$, we define
\begin{equation}
D_{r,s}:=\mathscr{U}_{r}(\ell^\infty_\C(w))\times \ell^\infty_s.
\end{equation}
    
\begin{defn}[The space of power-trigonometric series $\mathcal{A}$]\label{matcalA}
	Given $r,s>0$, we define $\mathcal{A}_{r,s}$ as the space of formal power-trigonometric series $f$ of the form \eqref{fourier0} such that the norm
	\begin{equation}\label{normR}
	|f|_{r,s}:= \sum_{\alpha,\nu}|f_{\alpha,\nu}|\sup_{(y,\varphi) \in D_{r,s}}|\tm_{\alpha,\nu}(y,\varphi)|
	\end{equation}
    is finite. 
\end{defn}
Note that, if we denote by $v_{w,r}\in \ell^{\infty}_w$ the sequence defined by $(v_{w,r})_j=rw_j^{-1}$, for all $j \in \mathbb{N}$, then the $|\cdot|_{r,s}$-norm of $f \in \mathcal{A}_{r,s}$ can be expressed as
    \begin{equation}
	|f|_{r,s}= \sum_{\alpha, \nu}|f_{\alpha,\nu}|v_{w,r}^\alpha e^{s|\nu|}.
	\end{equation}
Finally, we denote by $\mathcal{A}$ the space of power-trigonometric series defined by
	\begin{equation*}
	\mathcal{A}:=\bigcup_{r,s>0}\mathcal{A}_{r,s}. \ 
	\end{equation*} 

\begin{rmk}[Properties of $|\cdot|_{r,s}$]\label{algebra} It is straightforward to see that the weighted norm defined above satisfies the following properties
\begin{enumerate}
    \item Monotonicity: the norm $|\cdot|_{r,s}$ is monotone increasing with respect to the indices $r,s > 0$.
    \item Algebra: for any $f,g\in\cA_{r,s}$ we have that $$|fg|_{r,s} \le |f|_{r,s}|g|_{r,s}.$$
\end{enumerate}
    \end{rmk}

Note that, elements in $\cA$ are in fact the Taylor-Fourier series of majorant analytic functions. More precisely, given $f \in \mathcal{A}$, we define its \emph{majorant} $f_M$ as the series 
	\begin{equation*}
	f_M(y,\varphi):=\sum_{\alpha, \nu}|f_{\alpha,\nu}|y^\alpha e^{\im \nu \cdot \varphi}.
	\end{equation*}
We shall say that $f$ defines a \emph{majorant real analytic} function if there exists $r>0$ such that
	\begin{equation*}
	\mathscr{U}_r(\ell^\infty_w)\times \ell^\infty \to \mathbb{R} \ , \qquad (y,\varphi) \mapsto f_M(y,\varphi),
	\end{equation*}
	is real analytic. The following holds (see appendix \ref{appanalytic} for the proof).

\begin{prop}[Majorant analyticity and Taylor-Fourier representation]\label{manalytic}
       If $f \in \mathcal{A}_{r,s}$ for some $r,s>0$, then the function
	\begin{equation*}
	\mathscr{U}_r(\ell^\infty_w)\times \ell^\infty \to \mathbb{R} \ , \qquad (y,\varphi) \mapsto f(y,\varphi) \ ,
	\end{equation*}
	is majorant real analytic. Moreover, for any $\nu\in\Z^{(\N)}$, the map
	\begin{equation*}
	f_\nu:\mathscr{U}_r(\ell^\infty_w) \to \mathbb{C} \ , \qquad f_{\nu}(y)=\sum_{\alpha}f_{\alpha,\nu}y^\alpha ,
	\end{equation*}
  satisfies
	\begin{equation}\label{ansia1}
	f_{\nu}(y)= \int_{\mathbb{T}^\infty}f(y,\varphi) e^{- \im \nu \cdot \varphi} \ d\varphi := \lim_{N \to \infty} \frac{1}{(2\pi)^{N}}\int_{\mathbb{T}^{N}} f(y,\varphi) e^{-\im \nu\cdot\varphi} \,
	d \varphi_{1} \ldots d \varphi_{N}   
	\end{equation}
	for any $y \in \mathscr{U}_r(\ell^\infty_w)$.
   \end{prop}
   For the proof, see appendix \ref{appanalytic}. 
   
   In particular, any \(f \in \mathcal{A}_{r,s}\) has the property that its majorant, and hence \(f\) itself, extends holomorphically to \(D_{r,s}\). Since the extension is invariant under translations of \(\varphi\) by \(\ell^\infty(\N,2\pi\Z)\), it descends to a holomorphic function on \(\mathscr{U}_r(\ell^\infty_w)\times \T^\infty_s\).

\bigskip
Let us now discuss the differentiability of functions in $\mathcal{A}$. For any $f \in \mathcal{A}$, and every $n \in \N$, we define 
\begin{equation}
    d^nf:=\sum_{\alpha,\nu}f_{\alpha,\nu}\,d^n\,\tm_{\alpha,\nu},
\end{equation}
where $d^n\,\tm_{\alpha,\nu}$ denotes the $n$-th Fréchet differential of the monomial $\tm_{\alpha,\nu}$. Note that, if $f \in \mathcal{A}_{r,s}$, then $d^nf$ is a map defined on $\mathscr{U}_r(\ell^\infty_w) \times \ell^\infty$ taking values in the space of continuous $n$-linear operators $(\ell^\infty_w \times \ell^\infty)^n \to \R$. We denote by $d^n_yf$ and by $d^n_\varphi f$ the restriction of $d^nf$ to the subspace $\ell^\infty_w \times \{0\}$ and $\{0\}\times \ell^\infty$, respectively.
   
\begin{lemma}[Cauchy estimates]\label{precauchy}
	For all \(r,s>0\) and any \(f \in \mathcal{A}_{r,s}\), the following estimates hold for all \(\rho \in (0,r)\), \(\sigma \in (0,s)\), and for all $n \in \N$:
	\begin{align}
		\sum_{\alpha,\nu}|f_{\alpha,\nu}|\sup_{(y,\varphi) \in D_{r-\rho,s}}|d^n_y \tm_{\alpha,\nu}(y,\varphi)|^{\mathrm{op}}
		&\leq \left( \frac{r}{r-\rho} \right)^n \frac{n^n e^{-n}}{\rho^n} |f|_{r,s}, \label{cauchyy} \\[1ex]
		\sum_{\alpha,\nu}|f_{\alpha,\nu}|\sup_{(y,\varphi) \in D_{r,s-\sigma}}|d^n_\varphi\tm_{\alpha,\nu}(y,\varphi)|^{\mathrm{op}}
		&\leq \frac{n^n e^{-n}}{\sigma^n} |f|_{r,s}, \label{cauchyphi} \\[1ex]
	\sum_{\alpha,\nu}|f_{\alpha,\nu}|\sup_{(y,\varphi) \in D_{r-\rho,s-\sigma}}|d^n \tm_{\alpha,\nu}(y,\varphi)|^{\mathrm{op}}
		&\leq \left( \frac{r}{(r-\rho)\rho} + \frac{1}{\sigma} \right)^n n^n e^{-n} |f|_{r,s}. \label{cauchyunified}
	\end{align}
\end{lemma}

\begin{proof}
	Write
	\[
	\tm_{\alpha,\nu}(y,\varphi)=Y_\alpha(y)\Phi_\nu(\varphi),
	\qquad
	Y_\alpha(y)=y^\alpha,
	\qquad
	\Phi_\nu(\varphi)=e^{\im\nu\cdot\varphi}.
	\]
        For \(d_y^n\), we have
	\[
	d_y^n \tm_{\alpha,\nu}(y,\varphi) = d^n Y_\alpha(y) \Phi_\nu(\varphi).
	\]
	Thus
	\[
	\sum_{\alpha,\nu}|f_{\alpha,\nu}|\sup_{(y,\varphi)\in D_{r-\rho,s}}|d_y^n \tm_{\alpha,\nu}|^{\mathrm{op}}
	=
	\sum_{\alpha,\nu}|f_{\alpha,\nu}|\sup_{y \in \mathscr{U}_{r-\rho}(\ell^\infty_\C(w))}|d^n Y_\alpha(y)|^{\mathrm{op}} e^{s|\nu|}.
	\]
	By the elementary Cauchy estimates for monomials,
	\[
	\sup_{y \in \mathscr{U}_{r-\rho}(\ell^\infty_\C(w))}|d^n Y_\alpha(y)|^{\mathrm{op}}
	\le
	\left(\frac{r}{r-\rho}\right)^n \frac{n^n e^{-n}}{\rho^n} v_{w,r}^\alpha.
	\]
	Therefore \eqref{cauchyy} follows. Similarly, for \(d_\varphi^n\),
	\[
	d_\varphi^n \tm_{\alpha,\nu}(y,\varphi) = Y_\alpha(y) d^n \Phi_\nu(\varphi),
	\]
	and
	\[
	\sup_{y \in \mathscr{U}_{r}(\ell^\infty_\C(w))}|Y_\alpha(y)| \le v_{w,r}^\alpha,
	\qquad
	\sup_{\varphi \in \ell^\infty_{s-\sigma}}|d^n \Phi_\nu(\varphi)|^{\mathrm{op}}
	\le
	\frac{n^n e^{-n}}{\sigma^n} e^{s|\nu|}.
	\]
	Hence \eqref{cauchyphi} follows. For \(h_i=(u_i,v_i)\in \ell^\infty_w\times\ell^\infty\), \(i=1,\dots,n\), the Leibniz rule gives
	\[
	d^n \tm_{\alpha,\nu}(y,\varphi)[h_1,\dots,h_n]
	=
	\sum_{S\subseteq\{1,\dots,n\}}
	d^{|S|}Y_\alpha(y)[(u_i)_{i\in S}]\,
	d^{n-|S|}\Phi_\nu(\varphi)[(v_j)_{j\notin S}].
	\]
	Taking the operator norm and using the triangle inequality,
	\[
	|d^n \tm_{\alpha,\nu}(y,\varphi)|^{\mathrm{op}}
	\le
	\sum_{S\subseteq\{1,\dots,n\}}
	|d^{|S|}Y_\alpha(y)|^{\mathrm{op}}\,
	|d^{n-|S|}\Phi_\nu(\varphi)|^{\mathrm{op}}.
	\]
	The right-hand side depends only on \(|S|\). Grouping the subsets by their cardinality \(p=|S|\), and noting that there are exactly \(\binom{n}{p}\) subsets of \(\{1,\dots,n\}\) of cardinality \(p\), we obtain
	\[
	|d^n \tm_{\alpha,\nu}(y,\varphi)|^{\mathrm{op}}
	\le
	\sum_{p=0}^n
	\binom{n}{p}
	|d^p Y_\alpha(y)|^{\mathrm{op}}\,
	|d^{n-p}\Phi_\nu(\varphi)|^{\mathrm{op}}.
	\]
	We now bound the sum over \(\alpha,\nu\). For the \(Y_\alpha\) part, as above,
	\[
	\sup_{y \in \mathscr{U}_{r-\rho}(\ell^\infty_\C(w))}|d^p Y_\alpha(y)|^{\mathrm{op}}
	\le
	\left(\frac{r}{(r-\rho)\rho}\right)^p p^p e^{-p}\, v_{w,r}^\alpha.
	\]
	For the \(\Phi_\nu\) part,
	\[
	\sup_{\varphi \in \ell^\infty_{s-\sigma}}|d^{n-p}\Phi_\nu(\varphi)|^{\mathrm{op}}
	\le
	\left(\frac{1}{\sigma}\right)^{n-p} (n-p)^{n-p} e^{-(n-p)} e^{s|\nu|}.
	\]
	Therefore
	\[
	\begin{aligned}
		&\sum_{\alpha,\nu}|f_{\alpha,\nu}|
		\sup_{D_{r-\rho,s-\sigma}}|d^n \tm_{\alpha,\nu}|^{\mathrm{op}} \\
		&\le
		\sum_{p=0}^n \binom{n}{p}
		\left(\frac{r}{(r-\rho)\rho}\right)^p
		\left(\frac{1}{\sigma}\right)^{n-p}
		p^p (n-p)^{n-p} e^{-n}
		\sum_{\alpha,\nu}|f_{\alpha,\nu}| v_{w,r}^\alpha e^{s|\nu|} \\
		&\le
		n^n e^{-n}
		\sum_{p=0}^n \binom{n}{p}
		\left(\frac{r}{(r-\rho)\rho}\right)^p
		\left(\frac{1}{\sigma}\right)^{n-p}
		|f|_{r,s} \\
		&=
		\left(\frac{r}{(r-\rho)\rho}+\frac{1}{\sigma}\right)^n
		n^n e^{-n} |f|_{r,s}.
	\end{aligned}
	\]
	This proves \eqref{cauchyunified}. Finally, the corresponding estimates for \(d^n_y f\), \(d^n_\varphi f\), and \(d^n f\) follow by linearity and the triangle inequality for the operator norm.
\end{proof}

From the estimates in Lemma \ref{precauchy} one can easily prove that the series
	\[
	\sum_{\alpha,\nu} f_{\alpha,\nu}\, d^n \tm_{\alpha,\nu}
	\]
converges totally, and hence uniformly, on a neighborhood of any point in $D_{r,s}$ with respect to the operator norm. Consequently, term-by-term differentiation is justified, and the sum coincides with the \(n\)-th Fréchet differential \(d^n f\) of the function \(f\). The same argument applies to the partial differentials $d_y^nf$ and $d^n_\varphi f$. Therefore, Lemma \ref{precauchy} also implies that any $f \in \mathcal{A}$ is infinitely differentiable and, for any $n \in \N$, the $n^\mathrm{th}$-differential is a majorant analytic map taking values in the space of bounded $n$-linear real valued operators on $(\ell^\infty_w\times \ell^\infty)^n$.

\subsection{Changes of coordinates}\label{Coc}

Given $r,s>0$, we denote by $(\mathcal{A}^{\infty,w}_{r,s},|\cdot|^{\infty,w}_{r,s})$ the Banach space of functions 
	\[ N:\mathscr{U}_{r}(\ell^\infty_w) \times \ell^\infty \to \ell^\infty_w, \quad N(y,\varphi)=(N^{(m)}(y,\varphi))_{m \in \mathbb{N}}, \] 
	such that, for any $m \in \mathbb{N}$, $N^{(m)}\in \mathcal{A}_{r,s}$, with 
	\begin{equation*}
	|N|^{\infty,w}_{r,s}:=\sup_{m \in \mathbb{N}}w_m|N^{(m)}|_{r,s}<\infty.
	\end{equation*}
	Accordingly, $(\mathcal{A}^{\infty}_{r,s},|\cdot|^{\infty}_{r,s})$ denotes the Banach space of functions 
	\[ \Omega:\mathscr{U}_{r}(\ell^\infty_w) \times \ell^\infty \to \ell^\infty, \quad \Omega(y,\varphi)=(\Omega^{(m)}(y,\varphi))_{m \in \mathbb{N}}, \] 
	such that, for any $m \in \mathbb{N}$, $\Omega^{(m)}\in \mathcal{A}_{r,s}$, with
	\begin{equation*}
	|\Omega|^{\infty}_{r,s}:=\sup_{m \in \mathbb{N}}|\Omega^{(m)}|_{r,s}<\infty.
	\end{equation*}

Given $r,s>0$, let $\mathcal{C}_{r,s}(\mathcal{A})$ be the affine space passing through the identity and directed by $\mathcal{A}^{\infty,w}_{r,s}\times\mathcal{A}^{\infty}_{r,s}$, namely
\begin{equation*}
\mathcal{C}_{r,s}(\mathcal{A}):= \left\{\Phi:\mathscr{U}_{r}(\ell^{\infty}_w)\times \ell^\infty \to \ell^\infty_w \times \ell^\infty \, : \, \Phi=\id+(N,\Omega), \, N \in \mathcal{A}^{\infty,w}_{r,s}, \ \Omega \in \mathcal{A}^{\infty}_{r,s} \right\}.
\end{equation*}
For any $\rho,\sigma\geq0$, let $\mathcal{C}^{\rho,\sigma}_{r,s}(\mathcal{A})$ be the subset of those maps $\Phi=\id+(N,\Omega)\in\mathcal{C}_{r,s}(\mathcal{A})$ such that 
\begin{equation*}
|N|^{\infty,w}_{r,s}\leq \rho \qquad \text{and} \qquad |\Omega|^{\infty}_{r,s}\leq \sigma \ .
\end{equation*}
In this case, we say that $\Phi$ is close to identity. Finally, we denote by $\mathcal{C}(\mathcal{A})$ the space of maps 
\begin{equation*}
\mathcal{C}(\mathcal{A}):=\bigcup_{r,s>0}\mathcal{C}_{r,s}(\mathcal{A}) \ .
\end{equation*}

\begin{rmk}\label{closetoidanal}
	Functions belonging to either $\mathcal{A}^{\infty,w}_{r,s}$ or $\mathcal{A}^{\infty}_{r,s}$ are majorant real analytic, meaning that, if one defines the majorant $F_M$ of a given $F $ belonging to $ \mathcal{A}^{\infty,w}_{r,s}$ ($\text{respecively to} \, \mathcal{A}^{\infty}_{r,s}$) by setting \[ F_M(y,\varphi)=(F_M^{(m)}(y,\varphi))_{m \in \mathbb{N}} \ , \] then $F_M$ admits an extension to an analytic function on $D_{r,s}$ taking values in $\ell^{\infty}_{w,C}$ ({respectively in}  $\ell^\infty_\C)$. Moreover, the $n^{\text{th}}$ differential is given by 
	\begin{equation*}
	d^nF(y,\varphi)[h^{(1)},\dots,h^{(n)}]=(d^nF^{(m)}(y,\varphi)[h^{(1)},\dots,h^{(n)}])_{m \in \mathbb{N}}.
	\end{equation*}
    The proof relies on Lemma \ref{uniformanalytic} and on the fact that the map 
    \[
    \ell^\infty_\C(w) \ni (y_j)_{j \in \mathbb{N}} \mapsto (w_jy_j)_{j \in \mathbb{N}}\in \ell^\infty_\C
    \]
    is globally bi-analytic. As a consequence, any $\Phi \in \mathcal{C}(\mathcal{A})$ is majorant real analytic.
\end{rmk}

\begin{defn}[Changes of coordinates]\label{mchange}
	A function $\Phi \in \mathcal{C}_{r_0,s_0}(\mathcal{A})$ is a change of coordinates if there exist $\rho \in [0,r_0/2)$ and $\sigma \in [0,s_0/2)$ such that:
	\begin{enumerate}
		\item[(a)] for all $r \in (0,r_0]$,
		\begin{equation*}
		\Phi: \mathscr{U}_{r}(\ell^{\infty}_w)\times \ell^\infty \to \mathscr{U}_{r+\rho}(\ell^\infty_w) \times \ell^\infty
		\end{equation*}
		with the property that, for all $s\in(0,s_0]$, and every $f \in \mathcal{A}_{r+\rho,s+\sigma}$, the pullback $\Phi^*f$ belongs to $\mathcal{A}_{r,s}$ with
		\begin{equation*}\label{mchangea}
		|\Phi^*|_{r,s}^{\rho,\sigma}:=\sup_{f\neq 0}\frac{|\Phi^*f|_{r,s}}{|f|_{r+\rho,s+\sigma}}<\infty;
		\end{equation*}
		
		\item[(b)] there exists $\Psi \in \mathcal{C}_{r_0-\rho,s_0-\sigma}(\mathcal{A})$ such that, for all $r \in (0,r_0-\rho]$,
		\begin{equation*}
		\Psi: \mathscr{U}_{r}(\ell^{\infty}_w)\times \ell^\infty \to \mathscr{U}_{r+\rho}(\ell^\infty_w) \times \ell^\infty
		\end{equation*}
		and, for all $s\in(0,s_0-\sigma]$, and every $f \in \mathcal{A}_{r+\rho,s+\sigma}$, the pullback $\Psi^*f$ belongs to $\mathcal{A}_{r,s}$ with
		\begin{equation*}\label{mchangeb}
		|\Psi^*|_{r,s}^{\rho,\sigma}:=\sup_{f\neq 0}\frac{|\Psi^*f|_{r,s}}{|f|_{r+\rho,s+\sigma}}<\infty,
		\end{equation*}
		which is the inverse of $\Phi$ in the following sense:
		\begin{equation*}
		\Phi \circ \Psi =\id_{\mathscr{U}_{r}(\ell^\infty_w) \times \ell^\infty}, \quad \forall r\in (0,r_0-\rho], \qquad \Psi \circ \Phi =\id_{\mathscr{U}_{r}(\ell^\infty_w) \times \ell^\infty}, \quad \forall r\in (0, r_0-2\rho].
		\end{equation*}
	\end{enumerate}
\end{defn}

\begin{rmk}\label{changeanal}
By Remark \ref{closetoidanal}, any change of coordinates $\Phi$ is a bi-analytic diffeomorphism.
\end{rmk}

Note that, given $r_0>0$ and $\rho \in [0,r_0)$, any map 
\[ \Phi:\mathscr{U}_{r_0}(\ell^\infty_w)\times \ell^\infty \to \ell^\infty_w\times \ell^\infty, \quad \Phi=\id+(N,\Omega), \]
satisfying condition $(a)$ of the above definition is such that $N\in\mathcal{A}^{\infty,w}_{r_0,s_0}$. Indeed, one has the bound:
\begin{equation*}
\begin{aligned}
\sup_{k \in \mathbb{N}}w_k|N^{(k)}|_{r_0,s_0}&=\sup_{k \in \mathbb{N}}w_k|\Phi^*y_k-\id^*y_k|_{r_0,s_0} \\
&\leq (r_0+\rho)\sup_{k \in \mathbb{N}}\frac{|\Phi^*y_k|_{r_0,s_0}}{w_k^{-1}(r_0+\rho)}+r_0 \\
&\leq (r_0+\rho)|\Phi^*|_{r_0,s_0}^{\rho,\sigma}+r_0.
\end{aligned}
\end{equation*}

The same argument does not apply to $\Omega$ because the maps $(y,\varphi) \mapsto \varphi_k$, $k \in \mathbb{N}$, are not periodic and, hence, do not belong to the class $\mathcal{A}$. One is tempted to follow the same reasoning using the family of functions $(y,\varphi) \mapsto e^{i\varphi_k}$, $k \in \mathbb{N}$, and subsequently attempt to provide an estimate for the logarithm. It can be shown that this approach works if $\Omega$ depends only on finitely many angles; however, we are neither able to extend the result to the general case nor to provide a counterexample. For this reason, we decided to consider the condition $\Omega \in \mathcal{A}^{\infty}_{r_0,s_0}$ as an assumption.

\begin{prop}[Estimate on the pullback of functions]\label{composition}
Let $r_0,s_0>0$ and $\rho,\sigma\geq0$, and let 
\[
\Phi=\mathrm{id}+(N,\Omega)\in \mathcal C^{\rho,\sigma}_{r_0,s_0}(\mathcal A).
\]
Then for every $r\in(0,r_0]$, $s\in(0,s_0]$, and every
\[
f\in \mathcal A_{r+\rho,\,s+\sigma},
\]
the pullback $\Phi^*f$ is well-defined as a formal power-trigonometric 
series and belongs to $\mathcal A_{r,s}$. Moreover,
\[
|\Phi^*f|_{r,s}\le |f|_{r+\rho,\,s+\sigma}.
\]
\end{prop}

\begin{proof}[Proof of Proposition \ref{composition}]
Fix $r\in(0,r_0]$, $s\in(0,s_0]$, and $f\in\mathcal A_{r+\rho,s+\sigma}$. 
We begin by expanding the pullback:
\[
\Phi^*f
=
f(y+N,\varphi+\Omega)
=
\sum_{\nu,\alpha} f_{\nu\alpha} (y+N)^\alpha 
e^{\mathrm i \nu\cdot\varphi} e^{\mathrm i \nu\cdot\Omega}.
\]
For each pair $(\nu,\alpha)$, write
\[
(y+N)^\alpha e^{\mathrm i \nu\cdot\Omega}
=
\sum_{\lambda,\beta} a_{\lambda\beta}^{(\nu,\alpha)} 
y^\beta e^{\mathrm i \lambda\cdot\varphi}.
\]
Substituting this into the previous expression yields
\[
\Phi^*f
=
\sum_{\nu,\alpha}\sum_{\lambda,\beta}
f_{\nu\alpha} a_{\lambda\beta}^{(\nu,\alpha)}
y^\beta e^{\mathrm i(\lambda+\nu)\cdot\varphi}.
\]

We now show that this double series converges absolutely in 
$\mathcal A_{r,s}$. By the triangle inequality and 
$|\lambda+\nu|\le |\lambda|+|\nu|$,
\[
\begin{aligned}
&\sum_{\nu,\alpha}\sum_{\lambda,\beta}
|f_{\nu\alpha}||a_{\lambda\beta}^{(\nu,\alpha)}|
v_{w,r}^{\beta} e^{s|\lambda+\nu|}
\\
&\le
\sum_{\nu,\alpha} |f_{\nu\alpha}| e^{s|\nu|}
\sum_{\lambda,\beta} |a_{\lambda\beta}^{(\nu,\alpha)}|
v_{w,r}^{\beta} e^{s|\lambda|}.
\end{aligned}
\]
The inner sum is the norm of 
$(y+N)^\alpha e^{\mathrm i \nu\cdot\Omega}$ in $\mathcal A_{r,s}$. 
By submultiplicativity,
\[
\bigl|(y+N)^\alpha e^{\mathrm i \nu\cdot\Omega}\bigr|_{r,s}
\le
\bigl|(y+N)^\alpha\bigr|_{r,s}
\bigl|e^{\mathrm i \nu\cdot\Omega}\bigr|_{r,s}.
\]

Since $r\le r_0$ and $s\le s_0$, monotonicity gives
\[
|N|^{\infty,w}_{r,s}\le \rho,\qquad |\Omega|^{\infty}_{r,s}\le \sigma.
\]
For the polynomial factor, for each component $j$,
\[
|y_j|_{r,s}= r w_j^{-1},\qquad
|N^{(j)}|_{r,s}\le w_j^{-1}|N|^{\infty,w}_{r,s}\le \rho w_j^{-1}.
\]
Hence
\[
|y_j+N^{(j)}|_{r,s}\le (r+\rho)w_j^{-1},
\]
and by submultiplicativity,
\[
\bigl|(y+N)^\alpha\bigr|_{r,s}
\le \prod_j \bigl((r+\rho)w_j^{-1}\bigr)^{\alpha_j}
= v_{w,r+\rho}^{\alpha}.
\]

For the exponential factor, expanding in series and using 
submultiplicativity,
\[
\bigl|e^{\mathrm i \nu\cdot\Omega}\bigr|_{r,s}
\le e^{|\nu|\,|\Omega|^{\infty}_{r,s}}
\le e^{|\nu|\sigma}.
\]

Therefore,
\[
\sum_{\lambda,\beta} |a_{\lambda\beta}^{(\nu,\alpha)}|
v_{w,r}^{\beta} e^{s|\lambda|}
\le
v_{w,r+\rho}^{\alpha} e^{|\nu|\sigma}.
\]
Substituting this into the double series gives
\[
\begin{aligned}
&\sum_{\nu,\alpha}\sum_{\lambda,\beta}
|f_{\nu\alpha}||a_{\lambda\beta}^{(\nu,\alpha)}|
v_{w,r}^{\beta} e^{s|\lambda+\nu|}
\\
&\le
\sum_{\nu,\alpha} |f_{\nu\alpha}| e^{s|\nu|} v_{w,r+\rho}^{\alpha} 
e^{|\nu|\sigma}
\\
&=
\sum_{\nu,\alpha} |f_{\nu\alpha}| v_{w,r+\rho}^{\alpha} 
e^{(s+\sigma)|\nu|}
=
|f|_{r+\rho,s+\sigma}<\infty.
\end{aligned}
\]
Thus the double series converges absolutely. We may therefore reorder 
and group terms to write $\Phi^*f$ in canonical form
\[
\Phi^*f
=
\sum_{\mu,\beta} g_{\mu\beta} y^\beta e^{\mathrm i \mu\cdot\varphi},
\]
where
\[
g_{\mu\beta}
=
\sum_{\nu,\alpha} f_{\nu\alpha} a_{\mu-\nu,\beta}^{(\nu,\alpha)}.
\]
Finally, using the triangle inequality and the absolute convergence 
just established,
\[
|\Phi^*f|_{r,s}
=
\sum_{\mu,\beta} |g_{\mu\beta}| v_{w,r}^{\beta} e^{s|\mu|}
\le
\sum_{\nu,\alpha}\sum_{\lambda,\beta}
|f_{\nu\alpha}||a_{\lambda\beta}^{(\nu,\alpha)}|
v_{w,r}^{\beta} e^{s|\lambda+\nu|}
\le
|f|_{r+\rho,s+\sigma}.
\]
This proves both $\Phi^*f\in\mathcal A_{r,s}$ and the desired estimate.
\end{proof}

In particular, given $r_0,s_0>0$ and $\rho \in [0,r_0/2)$, $\sigma \in [0,s_0/2)$, every $\Phi \in \mathcal{C}^{\rho,\sigma}_{r_0,s_0}(\mathcal{A})$ that admits an inverse $\Psi \in \mathcal{C}^{\rho,\sigma}_{r_0-\rho,s_0-\sigma}(\mathcal{A})$, in the sense of Definition \ref{mchange}, is a close to identity change of coordinates.

\section{Robustness of foliations by infinite tori}\label{sectionstability}
\noindent In this section we study the robustness of foliations by infinite tori under real analytic perturbations. The result is established for an open neighborhood of $\{0\} \times \mathbb{T}^\infty$ within the Banach manifold $\ell^\infty_w \times \mathbb{T}^\infty$. 

\begin{notation}\label{gammap}
	Given $P \in \mathcal{A}^{\infty,w}_{r,s}$, we denote by $\Gamma_P$ the map 
	\begin{equation*}
	\Gamma_P:\mathscr{U}_{r}(\ell^\infty_w)\times \mathbb{T}^\infty \to \ell^\infty_w, \quad \Gamma_P(y,\pi_\infty(\varphi))=y+P(y,\varphi), \quad \forall y \in \mathscr{U}_r(\ell^\infty_w), \, \varphi \in \ell^\infty.
	\end{equation*}
    We recall that $\pi_\infty:\ell^\infty \to \T^\infty$ is the universal covering.
	Accordingly, $\Gamma_0$ denotes the projection onto the first factor.
\end{notation}

Our goal is to show that if $P$ is sufficiently small, namely if $\Gamma_P$ is a perturbation of the projection onto the first factor, then its level sets are infinite tori and give rise to an analytic foliation of a neighborhood of $\{0\}\times \mathbb{T}^\infty$. More precisely, we have the following:

\begin{thm}[Foliation by infinite tori]\label{FFoliation}
	Given $r_0>0$, if $P \in \mathscr{U}_{r_0/4}(\mathcal{A}^{\infty,w}_{r_0,s_0})$ for some $s_0>0$, then for all $r\leq r_0-2|P|^{\infty,w}_{r_0,s_0}$ one has that:
	\begin{itemize}
		\item[(a)] for all $\gamma \in \mathscr{U}_r(\ell^\infty_w)$, the level set $\Gamma_P^{-1}(\gamma)$ is the graph of a majorant analytic function from $\mathbb{T}^\infty$ to $\ell^\infty_w$ and hence it is analytically diffeomorphic to $\mathbb{T}^\infty$;
		\item[(b)] $\mathscr{N}_{0,r} := \Gamma_P^{-1}(\mathscr{U}_r(\ell^\infty_w))$ is a neighborhood of $\{0\}\times \mathbb{T}^\infty$, and there is a bi-analytic map $\chi$ such that the following diagram commutes:
		\begin{equation}\label{diagr}
			\begin{tikzcd}[column sep=large]
				\mathscr{N}_{0,r} \arrow[d, " \textstyle \Gamma_P"'] \arrow[r, "\textstyle \chi"] & \mathscr{U}_r(\ell^\infty_w)\times \mathbb{T}^\infty \arrow[ld, "\textstyle \Gamma_0"] \\
				\mathscr{U}_r(\ell^\infty_w) & 
			\end{tikzcd}
		\end{equation}
	\end{itemize}
\end{thm}

It is worth noting that Theorem \ref{FFoliation} extends to a broader class of maps of the form 
\begin{equation*}
\Gamma(y,\varphi)=f(y)+P(y,\varphi),
\end{equation*}
provided that $f$ is analytic and has an invertible differential. In fact, by performing a local Taylor expansion of $f$ around a point $y_0$, the map $\Gamma$ can be recast as 
$$f(y_0) + df(y_0)[y-y_0] + \mathcal{O}(\|y-y_0\|_{\ell^\infty_w}^2)+P(y,\varphi).$$ 
Therefore, the statement still holds  provided one sets $\tilde{\gamma} = y_0+df(y_0)^{-1}[\gamma - f(y_0)]$ and requires the smallness condition on the redefined perturbation 
$$\tilde{P}(y,\varphi)=df(y_0)^{-1}(\mathcal{O}(\|y-y_0\|^2_{\ell^\infty_w})+P(y,\varphi)).$$

Theorem \ref{FFoliation} can be rephrased as a fixed point problem in the space of conjugacies. More precisely, we show that given $\gamma$ sufficiently small, the equation
\begin{equation}
   y+P(y,\varphi)=\gamma, \quad \forall \varphi \in \T^\infty,
\end{equation}
is solved by parametrizing $y$ as $y(\gamma,\varphi)=\gamma-A(\gamma,\varphi)$, where $A$ is a real analytic function satisfying
\begin{equation}\label{A}
P(\gamma-A(\gamma,\varphi))=A(\gamma,\varphi), \quad \forall \varphi \in \T^\infty.
\end{equation}
This amounts to proving that, given $F \in \mathcal{A}^{\infty,w}_{r,s}$ sufficiently small, and $\rho<r$, the map
\begin{equation}
    \mathscr{F}: \mathscr{U}_{\rho}(\mathcal{A}^{\infty,w}_{r-\rho,s}) \ni G \mapsto \Phi^*_{G}F \in \mathcal{A}^{\infty,w}_{r-\rho,s} \,, \quad \text{where} \quad  \Phi_G: (y,\varphi)\mapsto (y+G(y,\varphi),\varphi)\,,
\end{equation}
admits a fixed point. In fact, the following statement holds.

\begin{thm}\label{fix}
Given $r,s>0$, let $F \in \mathcal{A}^{\infty,w}_{r,s}$ and define $\rho:=2|F|^{\infty,w}_{r,s}$. Then
    \begin{enumerate}[label=(\roman*), ref=(\roman*)]
        \item\label{fix-1} if $\rho<r\,$, $\mathscr{F}$ admits a unique fixed point $\hat{G} \in \overline{\mathscr{U}}_{\rho/2}(\mathcal{A}^{\infty,w}_{r-\rho})\,$.
        \item\label{fix-2} if $\rho<r/2\,$, $\Phi_{-F} \in \mathcal{C}_{r,s}^{\rho,0}(\mathcal{A})$ is a change of coordinates with inverse $\Phi_{\hat{G}}\in \mathcal{C}^{\rho,0}_{r-\rho,s}(\mathcal{A})\,$.
        \end{enumerate}
\end{thm}

The proof of Theorem \ref{fix} is deferred to the end of the section; we first prove the main result.

\begin{proof}[Proof of Theorem \ref{FFoliation}]
By setting $F=-P$, it follows from \ref{fix-1} of Theorem \ref{fix} that there exists $A\in\overline{\mathscr{U}}_{\rho/2}(\mathcal{A}^{\infty,w}_{r_0-\rho})$, where $\rho:=2|P|^{\infty,w}_{r_0,s_0}$, such that $\Phi^*_{-A}P=A$. In particular, for all $\gamma \in \mathscr{U}_{r_0-\rho}(\ell^\infty_w)$, one has that
\begin{equation}
    \Gamma^{-1}_P(\gamma):=\{\, (y,\pi_\infty(\varphi))\in\mathscr{U}_{r_0}(\ell^\infty_w)\times \T^\infty \, : \, y+P(y,\varphi)=\gamma  \,\}
\end{equation}
is the graph of the majorant analytic function $\T^\infty \ni \varphi \mapsto y_\gamma(\varphi)=\gamma-A(\gamma,\varphi)$. This proves the first point of Theorem \ref{FFoliation}. 

Moreover, for all $r \in (0,r_0-\rho]$, the open subset $\mathscr{N}_{0,r}:=\Gamma^{-1}_P(\mathscr{U}_r(\ell^\infty_w))$ is equal to $q \circ \Phi_{-A}(\mathscr{U}_r(\ell^\infty_w)\times \ell^\infty)$, where $q: \ell^\infty_w \times \ell^\infty \to \ell^\infty_w \times \mathbb{T}^\infty$ is the map defined by $q(y, \varphi) = (y, \pi_\infty(\varphi))$, and is a neighborhood of $\{0\}\times \T^\infty$. Indeed, by item \ref{fix-2} of Theorem~\ref{fix}, $\Phi_P\in\mathcal{C}^{\rho,0}_{r_0,s_0}$ is a change of coordinates with inverse $\Phi_{-A}\in\mathcal{C}^{\rho,0}_{r_0-\rho,s_0}$. Therefore,
\begin{align}\label{nbdhof0T}
		\{0\}\times \ell^\infty &\subset \mathscr{U}_{r_0-2\rho}(\ell^\infty_w)\times \ell^\infty 
        \\&= \Phi_{-A}\circ\Phi_{P}(\mathscr{U}_{r_0-2\rho}(\ell^\infty_w)\times \ell^\infty)
        \\& \subset \Phi_{-A}(\mathscr{U}_{r_0-\rho}(\ell^\infty_w)\times \ell^\infty)\,.
		\end{align}
Finally, by the surjectivity and the analyticity of $q$, it follows from Remark \ref{changeanal} that, for all $r \in (0,r_0-\rho]$, the map $\chi:\mathscr{N}_{0,r} \to \mathscr{U}_{r_0-2\rho}(\ell^\infty_w)\times \T^\infty$ defined by $\chi \circ q = q \circ \Phi_P$ is a bi-analytic diffeomorphism that makes the diagram \eqref{diagr} commute.
\end{proof}

\begin{rmk}
    
Theorem \ref{FFoliation} provides the essential geometric framework for extending the \emph{Liouville-Arnold} theorem to an infinite-dimensional setting by interpreting $\Gamma_P$ as a momentum map. Nevertheless, a fundamental distinction from Arnold's original derivation (see, e.g., \cite{arnold1978mathematical}) lies in the characterization of the level sets of the momentum map: while the classical (finite-dimensional) approach recovers the tori as the quotient of the group action defined by the joint flow (an algebraic-geometric construction) with respect to the stabilizer, our method is fundamentally analytic, describing the tori as graphs of majorant functions. It should be noted, however, that this result is perturbative in nature, as it requires the momentum map to be a small perturbation of the projection onto the first factor, thus restricting the class of systems to which it applies. A primary focus of the next section is the construction of an infinite-dimensional version of \emph{action-angle} coordinates adapted to a foliation by tori, such as the one investigated in this section.
\end{rmk}

The rest of the section is devoted to the proof of the technical Theorem \ref{fix}. 

\subsection{Proof of Theorem \ref{fix}}
\begin{prop}\label{gat}
	Given $r,s>0$, let $F \in \mathcal{A}^{\infty,w}_{r,s}$. Then, for any $\rho \in (0,r)$, the map
	\begin{equation*}
	\mathscr{F}:\mathscr{U}_{\rho}(\mathcal{A}^{\infty,w}_{r-\rho,s}) \ni G \mapsto \Phi^*_GF \in \mathcal{A}^{\infty,w}_{r-\rho,s}
	\end{equation*}
	is Fréchet differentiable. Moreover, the following estimate holds:
	\begin{equation}\label{boundgat}
	|d\mathscr{F}(G)|_{\text{op}}\leq \frac{r}{r-\rho_G}\frac{1}{e\rho_G}|F|^{\infty,w}_{r,s},
	\end{equation}
	where $\rho_G:=\min\{\rho-|G|_{r-\rho,s}^{\infty,w}, r/2 \}$.
	\begin{proof}
		Let us show that the Fréchet differential $d\mathscr{F}(G)$ of $\mathscr{F}$ at $G\in\mathscr{U}_{\rho}(\mathcal{A}^{\infty,w}_{r-\rho,s})$ is the linear map that associates to each $H\in \mathcal{A}^{\infty,w}_{r-\rho,s}$ the function in $\mathcal{A}^{\infty,w}_{r-\rho,s}$ given by
		\begin{equation*}
		(\Phi^*_G \ \partial_yF)[H]:(y,\varphi) \mapsto \partial_y F (y+G(y,\varphi),\varphi)[H(y,\varphi)].
		\end{equation*}
		Let $H\in \mathcal{A}^{\infty,w}_{r-\rho,s}$ be such that $G+H \in \mathscr{U}_\rho(\mathcal{A}^{\infty,w}_{r-\rho,s})$. We must provide a suitable estimate for the following quantity:
		\begin{equation}\label{quantity}
		| \ \mathscr{F}(G+H)-\mathscr{F}(G)-(\Phi^*_G \ \partial_yF)[H] \ |^{\infty,w}_{r-\rho,s}.
		\end{equation}
		Note that, for $|t|\leq 1$, the equality
		\begin{equation*}
		\begin{aligned}
		&(\mathscr{F}^{(m)}(G+tH))(y,\varphi)-(\mathscr{F}^{(m)}(G))(y,\varphi)-t\partial_yF^{(m)}(y+G(y,\varphi),\varphi)[H(y,\varphi)] \\
		&\quad =\sum_{\nu,\alpha}F^{(m)}_{\nu\alpha}\bigg((y+G(y,\varphi)+tH(y,\varphi))^\alpha-(y+G(y,\varphi))^\alpha \\
		&\qquad -t\sum_{k \in \mathbb{N}}H^{(k)}(y,\varphi)\alpha_k(y+G(y,\varphi))^{\alpha-e_k}\bigg)e^{\im \nu\cdot \varphi}
		\end{aligned}
		\end{equation*}
		holds pointwise for all $(y,\varphi) \in \mathscr{U}_{r-\rho}(\ell^\infty_w)\times \ell^\infty$. Defining the analytic function
		\begin{equation*}
		g_\alpha:t \mapsto g_\alpha(t):=(y+G(y,\varphi)+tH(y,\varphi))^\alpha,
		\end{equation*}
		one can write the above formula as
		\begin{equation*}
		\sum_{\nu,\alpha}F^{(m)}_{\nu\alpha}\left(\int_0^t g_\alpha''(\tau)(t-\tau)d\tau\right)e^{\im \nu\cdot \varphi},
		\end{equation*}
		which is equivalent to
		\begin{equation*}
		\sum_{\nu,\alpha}F^{(m)}_{\nu\alpha}\left(\int_{0}^t\sum_{k_1,k_2 \in \mathbb{N}}(H^{(k_1)}H^{(k_2)})(y,\varphi)\frac{\partial^2 y^\alpha}{\partial y_{k_1} \partial y_{k_2}}\circ \Phi_{G+\tau H}(y,\varphi) (t-\tau)d\tau\right)e^{\im \nu \cdot \varphi},
		\end{equation*}
		where, by abuse of notation, $y^\alpha$ denotes the majorant analytic function $(y,\varphi)\mapsto y^\alpha$.
		By the second item of Remark \ref{algebra} and the triangle inequality, the $|\cdot|_{r-\rho,s}$ norm of the above expression is bounded by
		\begin{equation}\label{imply}
		(|H|^{\infty,w}_{r-\rho,s})^2\sum_{\nu,\alpha}|f_{\alpha,\nu}^{(m)}|\left(\int_0^t\sum_{k_1,k_2 \in \mathbb{N}}w_{k_1}^{-1}w_{k_2}^{-1}\left|\frac{\partial^2 y^\alpha}{\partial y_{k_1} \partial y_{k_2}}\circ \Phi_{G+\tau H}\right|_{r-\rho,s} (t-\tau)d\tau\right)e^{s|\nu|}.
		\end{equation}
		Evaluating at $t=1$, Proposition \ref{composition} yields
		\begin{equation*}
		\eqref{imply}\leq (|H|^{\infty,w}_{r-\rho,s})^2\sum_{\nu,\alpha}|f_{\alpha,\nu}^{(m)}|\left(\int_0^1\sum_{k_1,k_2 \in \mathbb{N}}w_{k_1}^{-1}w_{k_2}^{-1}\left|\frac{\partial^2 y^\alpha}{\partial y_{k_1} \partial y_{k_2}}\right|_{r-\rho',s} (1-\tau)d\tau\right)e^{s|\nu|},
		\end{equation*}
		for any $\rho'$ such that $0<\rho'\leq\rho-|G|^{\infty,w}_{r-\rho,s}-|H|^{\infty,w}_{r-\rho,s}$. Furthermore, by Lemma \ref{precauchy}, one can estimate the above expression by
		\begin{equation*}
		\left(\frac{r}{r-\rho'}\right)^2\frac{2e^{-2}}{(\rho')^2}(|H|^{\infty,w}_{r-\rho,s})^2|F^{(m)}|_{r,s}.
		\end{equation*}
		Hence,
		\begin{equation*}
		(|H|^{\infty,w}_{r-\rho,s})^{-1} \times \eqref{quantity} \leq \left(\left(\frac{r}{r-\rho'}\right)^2\frac{2e^{-2}}{(\rho')^2}|F|^{\infty,w}_{r,s}\right)|H|^{\infty,w}_{r-\rho,s},
		\end{equation*}
		which tends to $0$ as $|H|^{\infty,w}_{r-\rho,s} \to 0$.
		
		It remains to show that $\mathcal{A}^{\infty,w}_{r-\rho,s} \ni H \mapsto (\Phi^*_G \ \partial_y F)[H] \in \mathcal{A}^{\infty,w}_{r-\rho,s}$ is a bounded operator. The algebra property (see Remark \ref{algebra}) and the triangle inequality imply:
		\begin{equation}\label{onehas}
		\begin{aligned}
		| \ (\Phi^*_G \ \partial_y F)[H] \ |^{\infty,w}_{r-\rho,s} &\leq \sup_{m \in \mathbb{N}}\sum_{k \in \mathbb{N}}|H^{(k)}|_{r-\rho,s}\left|\frac{\partial F^{(m)}}{\partial y_k}\circ \Phi_{G}\right|_{r-\rho,s} \\
		&\leq |H|^{\infty,w}_{r-\rho,s}\sup_{m \in \mathbb{N}}\sum_{k \in \mathbb{N}}w_k^{-1}\left|\frac{\partial F^{(m)}}{\partial y_k}\circ \Phi_{G}\right|_{r-\rho,s}.
		\end{aligned}
		\end{equation}
		Invoking again Proposition \ref{composition} and Lemma \ref{precauchy}, we obtain:
		\begin{equation*}
		\eqref{onehas} \leq |H|^{\infty,w}_{r-\rho,s}\sup_{m \in \mathbb{N}}\sum_{k \in \mathbb{N}}w_k^{-1}\left|\frac{\partial F^{(m)}}{\partial y_k}\right|_{r-\rho',s} \leq |H|^{\infty,w}_{r-\rho,s}\frac{r}{r-\rho'}\frac{1}{e\rho'}|F|^{\infty,w}_{r,s},
		\end{equation*}
		where $0<\rho'\leq\rho-|G|^{\infty,w}_{r-\rho,s}$. The function
		\begin{equation*}
		(0,r) \ni \rho' \mapsto \frac{1}{(r-\rho')\rho'}
		\end{equation*}
		is strictly convex and attains its minimum at $\rho'=r/2$. Hence, $\mathscr{F}$ is Fréchet differentiable and its differential satisfies:
		\begin{equation*}
		|d\mathscr{F}(G)|_{\text{op}}=\sup_{H \in \overline{\mathscr{U}}_1} | \ (\Phi^*_G \ \partial_yF)[H] \ |^{\infty,w}_{r-\rho,s} \leq \frac{r}{r-\rho_G}\frac{1}{e\rho_G}|F|^{\infty,w}_{r,s},
		\end{equation*}
		where $\rho_G:=\min\{\rho-|G|_{r-\rho,s}^{\infty,w}, r/2 \}$. This concludes the proof.
	\end{proof}
\end{prop}
 
\begin{proof}[Proof of Theorem \ref{fix}]
		Since $|G|^{\infty,w}_{r-\rho,s}\leq \rho/2 \leq \rho$, it follows from Proposition \ref{composition} that
		\begin{equation*}
		 |\Phi^*_GF|^{\infty,w}_{r-\rho,s} \leq |F|^{\infty,w}_{r,s} = \rho/2.
		\end{equation*}
		Then, for all $G \in \overline{\mathscr{U}}_{\rho/2}(\mathcal{A}^{\infty,w}_{r-\rho,s}) $, the sequence $(G_n)_{n \in \N}$ defined by
		\begin{equation}\label{seq}
		G_{n+1}:=\mathscr{F}(G_n), \qquad G_1:=G
		\end{equation}
		is contained in $\overline{\mathscr{U}}_{\rho/2}(\mathcal{A}^{\infty,w}_{r-\rho,s})$.
		
		Since $\mathscr{F}$ is Fréchet differentiable (see Proposition \ref{gat}), if we show that its Fréchet differential $d\mathscr{F}$ satisfies
		\begin{equation*}
		\sup_{G \in \overline{\mathscr{U}}_{\rho/2}(\mathcal{A}^{\infty,w}_{r-\rho,s})}|d\mathscr{F}(G)|_{op}\leq K<1,
		\end{equation*}
		by Lagrange's inequality, one has
		\begin{equation*}
		|\mathscr{F}(G)-\mathscr{F}(G')|^{\infty,w}_{r-\rho,s}\leq \sup_{G \in \overline{\mathscr{U}}_{\rho/2}(\mathcal{A}^{\infty,w}_{r-\rho,s})}|d\mathscr{F}(G)|_{op} |G-G'|^{\infty,w}_{r-\rho,s}\leq K|G-G'|^{\infty,w}_{r-\rho,s}.
		\end{equation*}
		Namely, $\mathscr{F}$ would turn out to be a contraction with a Lipschitz constant equal to $K<1$. We make use of \eqref{boundgat} for obtaining a uniform bound for the operator norm of the Fréchet differential of $\mathscr{F}$. First, note that
		\begin{equation*}
			\min_{G \in \overline{\mathscr{U}}_{\rho/2}(\mathcal{A}^{\infty,w}_{r-\rho,s})}\rho_G=\min\{\rho/2,r/2\}=\rho/2
		\end{equation*}
        for $\rho<r$. Then, \eqref{boundgat} gives
		\begin{equation*}
			\sup_{G \in \overline{\mathscr{U}}_{\rho/2}(\mathcal{A}^{\infty,w}_{r-\rho,s})}|d\mathscr{F}(G)|_{\text{op}} \leq \frac{1}{e}\frac{r}{r-\rho/2}\frac{2|F|^{\infty,w}_{r,s}}{\rho}<\frac{1}{e}\frac{r}{r-r/2}=\frac{2}{e}<1.
		\end{equation*}
		
		Therefore, there exists $\hat{G}\in \overline{\mathscr{U}}_{\rho/2}(\mathcal{A}^{\infty,w}_{r-\rho,s})$ such that, for all $G \in \overline{\mathscr{U}}_{\rho/2}(\mathcal{A}^{\infty,w}_{r-\rho,s})$, the sequence $(G_n)_{n \in \N}$ defined in \eqref{seq} satisfies
		\begin{equation*}
		\lim_{n \to \infty}|G_n-\hat{G}|^{\infty,w}_{r-\rho,s}=0.
		\end{equation*}
		Of course, $\hat{G}$ satisfies $\hat{G}=\Phi^*_{\hat{G}}F$. This proves item \ref{fix-1}.

		Since for all $r'\leq r$ and $s'\leq s$, $|F|^{\infty,w}_{r,s}\leq |F|_{r,s}^{\infty,w} = \rho/2$, one has $|y-F|_{r',s'}^{\infty,w}\leq r'+\rho/2 \leq r'+\rho$. Then $$\Phi_{-F}:\mathscr{U}_{r'}(\ell^\infty_w)\times \ell^\infty \to \mathscr{U}_{r'+\rho}(\ell^\infty_w)\times \ell^\infty\,, \quad \text{for all} \,  0<r'\leq r \,.$$ Moreover, by Proposition \ref{composition} with $N=-F$ and $\Omega=0$, for all $r' \in (0,r]$, $s' \in (0,s]$, and every $f \in \mathcal{A}_{r'+\rho,s'}$
		\begin{equation*}
		|\Phi^*_{-F}f|_{r',s'} \leq |f|_{r'+\rho,s'}.
		\end{equation*}
	Let us prove that $\Phi_{\hat{G}}$ satisfies all the properties stated in $(b)$ of Definition \ref{mchange}. Since for all $0<r' \leq r-\rho$ and $s' \leq s$, $|y+\hat{G}|^{\infty,w}_{r',s'} \leq r'+\rho/2 \leq r'+\rho$, one has $$ \Phi_{\hat{G}}:\mathscr{U}_{r'}(\ell^\infty_w)\times \ell^\infty \to \mathscr{U}_{r'+\rho}\times \ell^\infty\,.$$
		Proposition \ref{composition}, with $N=\hat{G}$ and $\Omega=0$ gives that, for all $r' \in (0,r-\rho]$, $s' \in (0,s]$, and every $f \in \mathcal{A}_{r'+\rho,s'}$
		\begin{equation*}
		|(\Phi_{\hat{G}})^*f|_{r',s'} \leq |f|_{r'+\rho,s'}.
		\end{equation*}
		The fact that $\Phi_{\hat{G}}$ is the inverse of $\Phi_{-F}$ follows from item \ref{fix-1} of Theorem \ref{fix}. The further condition $|F|^{\infty,w}_{r,s}<r/4$ ensures that the left inverse is also well defined.
	\end{proof}

\section{Poisson algebra of regular Hamiltonians on $\ell^\infty_w\times\T^\infty$ and symplectomorphisms}\label{SecRHam}
In the following we shall introduce the class of Hamiltonian we deal with together with the Poisson structure and related admissible symplectomorphisms.

\noindent Let $q: \ell^\infty_w \times \ell^\infty \to \ell^\infty_w \times \mathbb{T}^\infty$ be the map $q(y, \varphi) = (y, \pi_\infty(\varphi))$. We shall identify analytic vector fields defined on an open neighborhood $\mathscr{N}_0 \subset \ell^\infty_w \times \mathbb{T}^\infty$ of $\{0\} \times \mathbb{T}^\infty$ with analytic vector fields on the lifted neighborhood $q^{-1}(\mathscr{N}_0) \subset \ell^\infty_w \times \ell^\infty$ that are $2\pi$-periodic in each variable $\varphi_k$, $k \in \mathbb{N}$.

\begin{defn}[Space of analytic vector fields]\label{analvect}
Consider the real Banach space $\ell^\infty_w \times \ell^\infty$ and two real numbers $r,s>0$. We denote by $\mathfrak{X}_{r,s}$ the space of vector fields
\begin{equation*}
	V:\mathscr{U}_r(\ell^\infty_w)\times \ell^\infty \to \ell^\infty_w \times \ell^\infty, \quad V(y,\varphi)= (V_y(y,\varphi), V_\varphi(y,\varphi))
\end{equation*}
such that
\begin{equation}\label{vf}
V_y \in \mathcal{A}^{\infty,w}_{r,s} \quad \text{and} \quad V_\varphi \in \mathcal{A}^\infty_{r,s}.
\end{equation}
Accordingly, we denote by $\mathfrak{X}$ the space of vector fields 
$$ \mathfrak{X}:=\bigcup_{r,s>0} \mathfrak{X}_{r,s}.$$
\end{defn}

Note that, as a consequence of Theorem \ref{uniformanalytic}, any vector field in $V\in\mathfrak{X}$ is real analytic as a map from $\mathscr{U}_r(\ell^\infty_w)\times \ell^\infty$ to $\ell^\infty_w \times \ell^\infty$ and hence\footnote{See for e.g., \cite{PoTr} Appendix B.}, for every $(\underline{y},\underline{\varphi}) \in \mathscr{U}_r(\ell^\infty_w)\times\ell^\infty$, there exists $\eps>0$ and a unique solution curve
\begin{equation*}
	\phi_{V}(\underline{y},\underline{\varphi}):(-\eps,\eps)\to \mathscr{U}_r(\ell^\infty_w)\times\ell^\infty, \quad \begin{cases}
    \dfrac{d}{dt}\phi_V^t(\underline{y},\underline{\varphi})=V(\phi_V^t(\underline{y},\underline{\varphi})) \\
	    \phi_V^0(\underline{y},\underline{\varphi})=(\underline{y},\underline{\varphi}),
        \end{cases}
\end{equation*}
which is an analytic function of $t$ and $(\underline{y},\underline{\varphi})$.

\begin{defn}[Lie derivative]
Given a vector field $V \in \mathfrak{X}$, the Lie derivative of any $f \in \mathcal{A}$ in direction of $V$ is defined as
\begin{equation*}
	\mathscr{L}_Vf:=\left. \frac{d}{dt} \right|_{t=0}(\phi^t_V)^*f.
\end{equation*}
\end{defn}

\begin{prop}\label{proplie}
	The Lie derivative of $f \in \mathcal{A}$ in the direction of $V \in \mathfrak{X}$ can be expressed as
	\begin{equation}\label{liee}
	\mathscr{L}_Vf = \sum_{k \in \mathbb{N}} \left( \frac{\partial f}{\partial y_k} V_y^{(k)} + \frac{\partial f}{\partial \varphi_k} V_\varphi^{(k)} \right).
	\end{equation}
		Moreover, for any given $r,s>0$ and for all $V \in \mathfrak{X}_{r,s}$, the map $f \to \mathscr{L}_Vf$ defines a bounded linear operator
		\begin{equation*}
			\mathscr{L}_V:\mathcal{A}_{r,s}\to\mathcal{A}_{r-\rho,s-\sigma}
		\end{equation*}
		for all $\rho \in (0,r)$ and every $\sigma\in (0,s)$. More precisely, the following estimate holds
		\begin{equation}\label{est}
	|\mathscr{L}_Vf|_{r-\rho,s-\sigma}\leq\bigg(\frac{r}{e\rho(r-\rho)}|V_y|^{\infty,w}_{r,s}+\frac{1}{e\sigma}|V_\varphi|^{\infty}_{r,s}\bigg)|f|_{r,s}.
	\end{equation}
    \end{prop}

		\begin{proof}
			By the analyticity of $\phi_V$ with respect to $t$, one can easily verify that \eqref{liee} holds pointwise for any monomial $f = \tm_{\alpha,\nu}$, with $\alpha \in \mathbb{N}_0^{(\mathbb{N})}$ and $\nu \in \mathbb{Z}^{(\mathbb{N})}$. Therefore \eqref{liee} follows from the fact that, for all $(y,\varphi) \in \mathscr{U}_{r}(\ell^\infty_w) \times \ell^\infty$, the series
\begin{gather}
\sum_{\nu,\alpha}f_{\alpha,\nu}(y(t))^\alpha e^{\im\nu\cdot\varphi(t)} \label{serie1}  \\
\sum_{\nu,\alpha}f_{\alpha,\nu}\sum_{k \in \N}\alpha_kV_y^{(k)}(y(t),\varphi(t))(y(t))^{\alpha-e_k} e^{\im \nu\cdot \varphi(t)}+\im\nu_kV_\varphi^{(k)}(y(t),\varphi(t))(y(t))^{\alpha}e^{\im \nu \cdot \varphi(t)}, \label{serie2}
\end{gather}
where $(y(t), \varphi(t)) = \phi_V^t(y, \varphi)$, converge totally in a neighborhood of $t=0$. Finally, using the triangle inequality, algebraic properties, and Cauchy estimates, we have the following:
	\begin{equation*}
	\begin{aligned}
	|\mathscr{L}_Vf|_{r-\rho,s-\sigma} & \leq \sum_{k \in \N} \left( \bigg| \frac{\partial f}{\partial y_k}\bigg|_{r-\rho,s-\sigma} |V_y^{(k)}|_{r-\rho,s-\sigma} + \bigg|\frac{\partial f}{\partial \varphi_k}\bigg|_{r-\rho,s-\sigma} |V_\varphi^{(k)}|_{r-\rho,s-\sigma} \right) \\ & \leq |V_y|^{\infty,w}_{r,s}\sum_{k \in \N}w_k^{-1}\bigg| \frac{\partial f}{\partial y_k}\bigg|_{r-\rho,s-\sigma}+|V_\varphi|^{\infty}_{r,s}\sum_{k \in \N}\bigg| \frac{\partial f}{\partial \varphi_k}\bigg|_{r-\rho,s-\sigma} \\
	& \leq \bigg(\frac{r}{e\rho(r-\rho)}|V_y|^{\infty,w}_{r,s}+\frac{1}{e\sigma}|V_\varphi|^{\infty}_{r,s}\bigg)|f|_{r,s},
	\end{aligned}
	\end{equation*}
	which proves the bound in \eqref{est}.
	\end{proof}

 \subsection{Regular Hamiltonians and Poisson structure}\label{RHPS}
 \begin{defn}[Regular functions]\label{RF}
 	We say that a function $f$ is regular if there exist $r,s>0$ such that $f \in \mathcal{A}_{r,s}$ and 
 	\begin{equation*}
 	\|f\|_{r,s}:=\max\bigg\{|f|_{r,s}, \ \sup_{k\in \N}w_k\bigg|\frac{\partial f}{\partial \varphi_k}\bigg|_{r,s}, \ r\sup_{k\in \N}\bigg|\frac{\partial f}{\partial y_k}
 	\bigg|_{r,s}\bigg\}<\infty.
 	\end{equation*}
 	We denote such space by $ \ \mathcal{R} \ $. Moreover, for any $ \ r,s>0 \ $, $ \ \mathcal{R}_{r,s} \ $ denotes the subspace of $ \ \mathcal{R} \ $ of functions with $ \ \|f\|_{r,s}<\infty \ $.
 \end{defn}

 \begin{defn}[Poisson bracket]\label{PB}
	For any $f,g \in \mathcal{R}$, we define their Poisson bracket
	\begin{equation}\label{pb}
		\{f,g\}:=\sum_{k \in \N}\frac{\partial f}{\partial y_k}\frac{\partial g}{\partial \varphi_k }-\frac{\partial f}{\partial \varphi_k}\frac{\partial g}{\partial y_k}.
	\end{equation}
	\end{defn}

\begin{prop}\label{algebraregular}
For any $r,s>0$, $\mathcal{R}_{r,s}$ is a Banach algebra with respect to the product of functions. 
\begin{proof}
	Given $f,g \in \mathcal{R}_{r,s}$, by the triangular inequality and the algebra property enjoyed by the norm $|\cdot|_{r,s}$ (see Proposition \ref{algebra}), we have
	\begin{equation*}
		r\sup_{k \in \N}\bigg|\frac{\partial fg}{\partial y_k}\bigg|_{r,s}\leq r\sup_{k \in \N}\bigg|\frac{\partial f}{\partial y_k}\bigg|_{r,s}|g|_{r,s}+|f|_{r,s} \ r\sup_{k \in \N}\bigg|\frac{\partial g}{\partial y_k}\bigg|_{r,s}\leq 2\|f\|_{r,s}\|g\|_{r,s}.
	\end{equation*}
	The same estimate holds for the derivative with respect to $\varphi$. Thus, by Proposition \ref{algebra}, one obtains
	\begin{equation*}
		\|fg\|_{r,s}\leq 2\|f\|_{r,s}\|g\|_{r,s} \ ,
	\end{equation*}
	which gives the claim.
	\end{proof}
\end{prop}

\begin{prop}
	The space of regular functions $\mathcal{R}$ equipped with the product of functions and the bracket \eqref{pb} is a Poisson algebra. Namely, for all $r,s>0$ and every $f,g \in \mathcal{R}_{r,s}$, we have that $ \ \{f,g\} \in \mathcal{R}_{r-\rho,s-\sigma} \ $ and the following estimate holds 
	\begin{equation*}
		\|\{f,g\}\|_{r-\rho,s-\sigma}\leq \frac{2}{e}\bigg(\frac{r}{r-\rho}\frac{1}{\rho}+\frac{1}{r\sigma}\bigg) \|f\|_{r,s}\|g\|_{r,s} \ . 
	\end{equation*}
	for all $\rho \in (0,r)$ and $ \ \sigma \in (0,s) \ $.
	
	\begin{proof}
		Let $r,s>0$ and $f,g \in \mathcal{R}_{r,s}$. The following bounds hold for any $\rho \in (0,r)$, $\sigma\in (0,s)$, by triangular inequality and Proposition \ref{algebra}.
\begin{equation}\label{prima}
|\{f,g\}|_{r-\rho,s-\sigma} \leq \sum_{k \in \N}\bigg|\frac{\partial f}{\partial y_k}\bigg|_{r-\rho,s-\sigma}\bigg|\frac{\partial g}{\partial \varphi_k}\bigg|_{r-\rho,s-\sigma}+\bigg|\frac{\partial f}{\partial \varphi_k}\bigg|_{r-\rho,s-\sigma}\bigg|\frac{\partial g}{\partial y_k}\bigg|_{r-\rho,s-\sigma} \ ,
\end{equation}
\begin{equation}\label{seconda}
\begin{aligned}
\sup_{\ell \in \N}\bigg|\frac{\partial \{f,g\}}{\partial y_\ell}\bigg|_{r-\rho,s}\leq &\sup_{\ell \in \N}\sum_{k \in \N}\bigg|\frac{\partial^2 f}{\partial y_\ell \partial y_k}\bigg|_{r-\rho,s}\bigg|\frac{\partial g}{\partial \varphi_k}\bigg|_{r-\rho,s} \\
+&\sup_{\ell \in \N}\sum_{k \in \N}\bigg|\frac{\partial f}{\partial y_k}\bigg|_{r-\rho,s}\bigg|\frac{\partial^2 g}{\partial y_\ell\partial \varphi_k}\bigg|_{r-\rho,s} + f \leftrightarrow g \ ,
\end{aligned}
\end{equation}	
\begin{equation}\label{terza}
\begin{aligned}
\sup_{\ell \in \N}w_\ell\bigg|\frac{\partial \{f,g\}}{\partial \varphi_\ell}\bigg|_{r,s-\s}\leq &\sup_{\ell \in \N}w_\ell\sum_{k \in \N}\bigg|\frac{\partial^2 f}{\partial \varphi_\ell \partial y_k}\bigg|_{r-\rho,s-\sigma}\bigg|\frac{\partial g}{\partial \varphi_k}\bigg|_{r-\rho,s-\sigma}  \\
+&\sup_{\ell \in \N}w_\ell\sum_{k \in \N}\bigg|\frac{\partial f}{\partial y_k}\bigg|_{r-\rho,s-\sigma}\bigg|\frac{\partial^2 g}{\partial \varphi_\ell\partial \varphi_k}\bigg|_{r-\rho,s-\sigma} + f \leftrightarrow g \ .
\end{aligned}
\end{equation}
By using the regularity of $f$ and $g$, and Cauchy estimate on the derivatives with respect to the $\varphi_k$-variables, $k \in \N$, one has
\begin{equation*}
\eqref{prima} \leq \frac{2}{e}\frac{1}{r\sigma}\|f\|_{r,s}\|g\|_{r,s}.
\end{equation*}
By invoking Lemma \ref{precauchy}, we obtain the following estimates:
\begin{equation*}
\begin{aligned}
(r-\rho)\times \eqref{seconda}=& (r-\rho)\sup_{\ell \in \N}\sum_{k \in \N} w_k^{-1} \bigg|\frac{\partial^2 f}{\partial y_\ell \partial y_k}\bigg|_{r-\rho,s-\sigma} w_k \bigg| \frac{\partial g}{\partial \varphi_k}\bigg|_{r-\rho,s-\sigma}  \\
+&(r-\rho)\sup_{\ell \in \N}\sum_{k \in \N}\bigg|\frac{\partial f}{\partial y_k} \bigg|_{r-\rho,s-\sigma}\bigg|\frac{\partial^2 g}{\partial y_\ell\partial \varphi_k}\bigg|_{r-\rho,s-\sigma} + f \leftrightarrow g \\
& \leq \frac{2}{e}\bigg(\frac{1}{\rho}+\frac{1}{r\sigma}\bigg)\|g\|_{r,s}\|f\|_{r,s}
\end{aligned}
\end{equation*}
and
\begin{equation*}
\begin{aligned}
\eqref{terza}=& \sup_{\ell \in \N}w_\ell\sum_{k \in \N}w_{k}^{-1}\bigg|\frac{\partial^2 f}{\partial \varphi_\ell \partial y_k}\bigg|_{r-\rho,s-\sigma}w_{k}\bigg|\frac{\partial g}{\partial \varphi_k}\bigg|_{r-\rho,s-\sigma} + \\
&\sup_{\ell \in \N}w_\ell\sum_{k \in \N}\bigg|\frac{\partial f}{\partial y_k}\bigg|_{r-\rho,s-\sigma}\bigg|\frac{\partial^2 g}{\partial \varphi_\ell\partial \varphi_k}\bigg|_{r-\rho,s-\sigma} + f \leftrightarrow g \\
& \leq \frac{2}{e}\bigg(\frac{r}{r-\rho}\frac{1}{\rho}+\frac{1}{r\sigma}\bigg)\|g\|_{r,s}\|f\|_{r,s} \ .
\end{aligned}
\end{equation*}
Hence the claim follows.
\end{proof}
\end{prop}

Motivated by Definition \ref{PB}, given differentiable, real valued functions $f$ and $g$ defined on $\mathscr{U}_r(\ell^\infty_w)\times \ell^\infty$, for some $r>0$, we denote by $\{f,g\}$ any expression of the form
\begin{equation*}
	\sum_{k \in \N}\frac{\partial f}{\partial y_k}\frac{\partial g}{\partial \varphi_k }-\frac{\partial f}{\partial \varphi_k}\frac{\partial g}{\partial y_k}
\end{equation*}
whenever it is well defined. By the same estimates as in Proposition \ref{proplie} it follows that, for all $f \in \mathcal{R}$, the Poisson bracket defines an operator
\begin{equation}
    \{\cdot,f\}:\mathcal{A}\to\mathcal{A}.
\end{equation}
More precisely, given $ \ r,s>0 \ $, for any $H \in \mathcal{A}_{r,s}$ and $f \in \mathcal{R}_{r,s}$, the following estimate holds
\begin{equation*}
		|\{H,f\}|_{r-\rho,s-\sigma}\leq \frac{1}{e}\bigg(\frac{r}{r-\rho}\frac{1}{\rho}+\frac{1}{r\sigma}\bigg)\|f\|_{r,s}|H|_{r,s}
	\end{equation*}
for all $\rho \in (0,r)$ and $\sigma \in (0,s)$.

\begin{rmk}\label{hamvfrmk}
	To any regular function $ \ f \in \mathcal{R} \ $ there is an associated Hamiltonian vector field, namely the unique vector field $ \ V^f \in \mathfrak{X} \ $ such that
	\begin{equation*}
		\mathscr{L}_{V^f}-\{\cdot,f\}=0 \quad \text{on} \quad \mathcal{A}.
	\end{equation*}
	By a direct computation, one can easily verify that
	\begin{equation}\label{hamvfregular}
	(V_y^f)^{(k)}=-\frac{\partial f}{\partial \varphi_k}, \quad (V_\varphi^f)^{(k)}=\frac{\partial f}{\partial y_k} \ .
	\end{equation}
	In contrast with the finite-dimensional case, the Hamiltonian vector field associated with a general analytic function is not necessarily a well-defined map from the phase space into its tangent bundle. Indeed, the differential of a real analytic function defined on a Banach space $X$ is a map taking values in the dual space $X^*$. This problem is absent in finite dimension, where the dual of a vector space is isomorphic to the space itself, and the non-degeneracy of the symplectic form always ensures that the Hamiltonian vector field associated with any analytic Hamiltonian is a well-defined analytic section of the tangent bundle.
\end{rmk}
 
 \subsection{Regular changes of coordinates and symplectomorphisms}\label{Rcoc}
Following the approach of the previous section, we aim to define appropriate norms for controlling changes of coordinates that preserve the regularity of functions in the sense of Definition~\ref{RF}. Given a change of coordinates $\Phi$ and a regular function $f$, the Hamiltonian vector field associated with the pullback of $f$ by $\Phi$ satisfies the following relation:
\begin{equation*}
V^{\Phi^*f}=(-J(d\Phi)^tJ) V^f\circ \Phi,
\end{equation*}
where $(d\Phi)^t$ denotes the transpose of $d\Phi$ and $J:\R^\N \times \R^\N \to \R^\N \times \R^\N$ is the linear operator 
\begin{equation*}
    J=\begin{pmatrix}
    0 && \id_{\ell^\infty} \\ -\id_{\ell^\infty_w} && 0
    \end{pmatrix}.
\end{equation*}
In this line of thought, we seek to define norms that ensure the boundedness of
\begin{equation*}
    (d\Phi)^t:\ell^\infty \times \ell^\infty_w \to \ell^\infty \times \ell^\infty_w.
\end{equation*}

\begin{defn}\label{definitionnorm}
	Given $r,s>0$, we denote by $\mathcal{R}^{\infty,w}_{r,s}$ the space of functions $$N:\mathscr{U}_{r}(\ell^\infty_w) \times \ell^\infty \to \ell^\infty_w, \quad N(y,\varphi)=(F^{(m)}(y,\varphi))_{m \in \N},$$ such that, for any $m \in \N$, $N^{(m)}\in \mathcal{A}_{r,s}$, with
	\begin{equation*}
	\|N\|^{\infty,w}_{r,s}:=\max\bigg\{|N|_{r,s}^{\infty,w}, \ r\sup_{k \in \N}\sum_{m=1}^\infty \bigg|\frac{\partial N^{(m)}}{\partial y_k}\bigg|_{r,s} \ , \ \sup_{k \in \N}w_k\sum_{m=1}^\infty \bigg|\frac{\partial N^{(m)}}{\partial \varphi_k}\bigg|_{r,s} \bigg\}<\infty \ .
	\end{equation*}
	Accordingly, $ \ \mathcal{R}^{\infty}_{r,s} \ $ denotes the space of functions 
	$$\Omega:\mathscr{U}_{r}(\ell^\infty_w) \times \ell^\infty \to \ell^\infty, \quad \Omega(y,\varphi)=(\Omega^{(m)}(y,\varphi))_{m \in \N},$$ such that, for any $m \in \N$, $\Omega^{(m)}\in \mathcal{A}_{r,s}$, with
	\begin{equation*}
	\|\Omega\|^{\infty}_{r,s}:=\max\bigg\{|\Omega|_{r,s}^{\infty}, \ r\sup_{k \in \N}\sum_{m=1}^\infty w_m^{-1} \bigg|\frac{\partial \Omega^{(m)}}{\partial y_k}\bigg|_{r,s} \ , \ \sup_{k \in \N}w_k\sum_{m=1}^\infty w_m^{-1} \bigg|\frac{\partial \Omega^{(m)}}{\partial \varphi_k}\bigg|_{r,s} \bigg\}<\infty \ .
	\end{equation*}
\end{defn}

Given $ \  r,s>0 \  $, we denote by $ \  \mathcal{C}_{r,s}(\mathcal{R}) \ $ the affine space passing through the identity and directed by $\mathcal{R}^{\infty,w}_{r,s}\times\mathcal{R}^{\infty}_{r,s}$, namely
$$ \mathcal{C}_{r,s}(\mathcal{R}):= \{\Phi:\mathscr{U}_{r}(\ell^{\infty}_w)\times \ell^\infty \to \ell^\infty_w \times \ell^\infty \ : \ \Phi=\id+(N,\Omega) \ , \ N \in \mathcal{R}^{\infty,w}_{r,s}, \ \Omega \in \mathcal{R}^{\infty}_{r,s}\} \ .$$

\noindent For any $ \ \rho,\sigma>0 \ $, let $ \ \mathcal{C}^{\rho,\sigma}_{r,s}(\mathcal{R}) \ $ be the subset of those maps $ \ \Phi=\id+(N,\Omega)\in\mathcal{C}_{r,s}(\mathcal{R}) \ $ such that 
\begin{equation*}
|N|^{\infty,w}_{r,s}\leq \rho \qquad \text{and} \qquad |\Omega|^{\infty}_{r,s}\leq \sigma \ .
\end{equation*}
In this case, we say that $ \ \Phi \ $ is close to identity.

 Accordingly, we denote by $ \ \mathcal{C}(\mathcal{R}) \ $ the space of functions 
\begin{equation*}
\mathcal{C}(\mathcal{R}):=\bigcup_{r,s>0}\mathcal{C}_{r,s}(\mathcal{R}) \ .
\end{equation*}

\begin{rmk}\label{adj}
	A map $ \ \Phi \in \mathcal{C}(\mathcal{R})$ has the property that the operator $-J(d\Phi)^tJ$ induces a bounded map from $ \ \mathfrak{X}_{r,s} \ $ to itself.
\end{rmk}

\begin{prop}[Estimate on the pullback of a regular function]\label{zzz}
	Given $ \ r_0,s_0>0 \ $, let $ \ \rho \in [0,r_0) \ $ and $ \ \sigma \in [0,s_0) \ $. Consider a function $ \ \Psi=\id+(N,\Omega) \in \mathcal{C}^{\rho,\sigma/2}_{r_0-\rho,s_0-\sigma}(\mathcal{R}) \ $. Then, for all $ \ r \in (\rho,r_0] \ $, $ \ s \in (\sigma,s_0] \ $, and every $ \ f \in \mathcal{R}_{r,s} \ $, the pullback of $ \ f \ $ by $ \ \Psi \ $
	belongs to $ \ \mathcal{R}_{r-\rho,s-\sigma} \ $ with estimate
	\begin{equation*}
	\|\Psi^*f\|_{r-\rho,s-\sigma} \leq 2(1+\rho/r+\sigma/2)\|f\|_{r,s}.
	\end{equation*}
	\begin{proof}
		The hypotheses of this proposition ensure that all the assumptions needed for applying Proposition \ref{composition} are fulfilled. Then, $|\Psi^*f|_{r-\rho,s-\sigma}\leq 2|f|_{r,s}$. Let us perform estimates on the derivatives of $\Psi^*f$. Let $x,\theta$ be the coordinates defined by $x=y+N(y,\varphi)$, $\theta=\varphi+\Omega(y,\varphi)$. By the triangular inequality applied to the formula given by the chain rule, the algebra property proved in Proposition \ref{algebra} and the estimate for the pullback of a majorant analytic function stated in Proposition \ref{composition}, one has
		\begin{equation*}
		\begin{aligned}
		r\sup_{k \in \N} \bigg|\frac{\partial \ \Psi^*f}{\partial y_k}\bigg|_{r-\rho,s-\sigma} 
		& \leq 
		r\sup_{k \in \N}
		\sum_{m=1}^{\infty}
		\bigg| \Psi^*\frac{\partial f}{\partial x_m} \bigg|_{r-\rho,s-\sigma} \bigg|\delta_{mk}+\frac{\partial N^{(m)}}{\partial y_k}\bigg|_{r-\rho,s-\sigma}\\
		&+r\sup_{k \in \N}
		\sum_{m=1}^{\infty}\bigg|\Psi^*\frac{\partial f}{\partial \theta_m}\bigg|_{r-\rho,s-\sigma}
		\bigg|\frac{\partial \Omega^{(m)}}{\partial y_k}\bigg|_{r-\rho,s-\sigma} \\
		& \leq 2r\sup_{k \in \N}\bigg| \frac{\partial f}{\partial x_k} \bigg|_{r,s} \\
		&+ 2r \sup_{k \in \N}\sum_{m=1}^\infty
		\bigg| \frac{\partial f}{\partial x_m} \bigg|_{r,s} \bigg|\frac{\partial N^{(m)}}{\partial y_k}\bigg|_{r-\rho,s-\sigma}+
		\bigg|\frac{\partial f}{\partial \theta_m}\bigg|_{r,s}
		\bigg|\frac{\partial \Omega^{(m)}}{\partial y_k}\bigg|_{r-\rho,s-\sigma}.
		\end{aligned}
		\end{equation*}
		By the regularity of $f$, and applying the assumptions, the above quantity is bounded by
		\begin{gather*}
		2\|f\|_{r,s}\bigg(1+	\sup_{k \in \N}\sum_{m=1}^\infty
		\bigg|\frac{\partial N^{(m)}}{\partial y_k}\bigg|_{r-\rho,s-\sigma}+
		r\sup_{k \in \N}\sum_{m=1}^\infty w_m^{-1}\bigg|\frac{\partial \Omega^{(m)}}{\partial y_k}\bigg|_{r-\rho,s-\sigma}\bigg) \\
		\leq 2(1+\rho/r+\sigma/2) \|f\|_{r,s} \ .
		\end{gather*}
		At this stage it only remains to give a suitable bound for the term of the norm which involves the derivatives with respect to the periodic variables. 	By the same argument of before we get
		\begin{equation*}
		\begin{aligned}
		\sup_{k \in \N} w_k \bigg|\frac{\partial \ \Psi^*f}{\partial \varphi_k}\bigg|_{r-\rho,s-\sigma} 
		& \leq \sup_{k \in \N}w_k
		\sum_{m=1}^{\infty}\bigg|\Psi^*\frac{\partial f}{\partial x_m}\bigg|_{r-\rho,s-\sigma}
		\bigg|\frac{\partial N^{(m)}}{\partial \varphi_k}\bigg|_{r-\rho,s-\sigma} \\
		&+\sup_{k \in \N} w_k
		\sum_{m=1}^{\infty}
		\bigg| \Psi^*\frac{\partial f}{\partial \theta_m} \bigg|_{r-\rho,s-\sigma} \bigg|\delta_{mk}+\frac{\partial \Omega^{(m)}}{\partial \varphi_k}\bigg|_{r-\rho,s-\sigma} \\
		& \leq 2\sup_{k \in \N}w_k\bigg| \frac{\partial f}{\partial \varphi_k} \bigg|_{r,s} \\
		&+ 2 \sup_{k \in \N}\sum_{m=1}^\infty
		\bigg| \frac{\partial f}{\partial x_m} \bigg|_{r,s} \bigg|\frac{\partial N^{(m)}}{\partial \varphi_k}\bigg|_{r-\rho,s-\sigma}+
		\bigg|\frac{\partial f}{\partial \theta_m}\bigg|_{r,s}
		\bigg|\frac{\partial \Omega^{(m)}}{\partial \varphi_k}\bigg|_{r-\rho,s-\sigma} \ ,
		\end{aligned}
		\end{equation*}
		which is bounded by
		\begin{gather*}
		2\|f\|_{r,s}\bigg(1+\frac{1}{r}\sup_{k \in \N}\sum_{m=1}^\infty
		\bigg|\frac{\partial N^{(m)}}{\partial \varphi_k}\bigg|_{r-\rho,s-\sigma}+
		\sup_{k \in \N}\sum_{m=1}^\infty w_m^{-1}\bigg|\frac{\partial \Omega^{(m)}}{\partial \varphi_k}\bigg|_{r-\rho,s-\sigma}\bigg) \\
		\leq 2(1+\rho/r+\sigma/2) \|f\|_{r,s} \ .
		\end{gather*}
		The claim is proved.
	\end{proof}
\end{prop}

Let us introduce a class of change of coordinates which preserves the regularity of functions belonging to the class $ \ \mathcal{R} \ $.

\begin{defn}[Regular changes of coordinates]\label{regchange}
	A function \( \Phi \in \mathcal{C}_{r_0,s_0}(\mathcal{R}) \) is called a regular change of coordinates if it is a change of coordinates such that there exist \( \rho \in [0,r_0/2) \) and \( \sigma \in [0,s_0/2) \) for which:
	\begin{enumerate}
		\item[(a)] for all \( r \leq r_0 \),
		\begin{equation*}
		\Phi: \mathscr{U}_{r}(\ell^{\infty}_w)\times \ell^\infty \to \mathscr{U}_{r+\rho}(\ell^\infty_w) \times \ell^\infty,
		\end{equation*}
		and, for all \( r \in (\rho,r_0+\rho] \), \( s\in(\sigma,s_0+\sigma] \), and every \( f \in \mathcal{R}_{r,s} \), \( \Phi^*f \in \mathcal{R}_{r-\rho,s-\sigma} \), with
		\begin{equation*}\label{mchangea}
		\|\Phi^*\|_{r,s}^{\rho,\sigma}:=\sup_{f\neq 0}\frac{\|\Phi^*f\|_{r-\rho,s-\sigma}}{\|f\|_{r,s}}<\infty;
		\end{equation*}

		\item[(b)] its inverse \( \Psi \in \mathcal{C}_{r_0-\rho,s_0-\sigma}(\mathcal{R}) \) satisfies
		\begin{equation*}
		\Psi: \mathscr{U}_{r}(\ell^{\infty}_w)\times \ell^\infty \to \mathscr{U}_{r+\rho}(\ell^\infty_w) \times \ell^\infty, \quad \forall r\leq r_0-\rho,
		\end{equation*}
		and, for all \( r \in (\rho,r_0] \), \( s\in(\sigma,s_0] \), and every \( f \in \mathcal{R}_{r,s} \), \( \Psi^*f \in \mathcal{R}_{r-\rho,s-\sigma} \), with
		\begin{equation*}\label{mchangeb}
		\|\Psi^*\|_{r,s}^{\rho,\sigma}:=\sup_{f\neq 0}\frac{\|\Psi^*f\|_{r-\rho,s-\sigma}}{\|f\|_{r,s}}<\infty.
		\end{equation*}
	\end{enumerate}
\end{defn}

\begin{rmk}\label{regrmk}
Note that, given $ \ r_0,s_0>0 \ $, $ \ \rho \in [0,r_0) \ $ and $ \ \sigma \in [0,s_0)$, any invertible map $ \ \Phi \in \mathcal{C}^{\rho,\sigma/2}_{r_0,s_0}(\mathcal{R}) \ $ with inverse $ \ \Psi \in \mathcal{C}^{\rho,\sigma/2}_{r_0-\rho,s_0-\sigma}(\mathcal{R}) \ $ is a (close to identity) regular change of coordinates, as follows by Proposition \ref{zzz}.
\end{rmk}

\begin{defn}[Symplectomorphisms]\label{sympl}
	A regular change of coordinates $ \ \Phi \in \mathcal{C}_{r_0,s_0}(\mathcal{R}) \ $ is a symplectomorphism if, for all $ \ r \in (\rho, r_0 + \rho] \ $ and $ \ s \in (\sigma, s_0 + \sigma] \ $, with $ \ \sigma \ $ and $ \ \rho \ $ as in Definition \ref{regchange}, the following equalities hold in $ \ \mathscr{U}_{r-\rho}(\ell^\infty_w) \times \ell^\infty \ :$
	\begin{equation}
	\begin{aligned}
	(i) \quad &\Phi^*\{H,f\}= \{\Phi^*H , \Phi^*f \} \ , \qquad && \text{for all } H\in\mathcal{A}_{r,s} \ , \ f\in \mathcal{R}_{r,s} \ , \\
	(ii) \quad &\Phi^*\{H,\varphi_m\}=\{\Phi^*H,\Phi^* \varphi_m\} \ , \quad && \text{for all } H\in\mathcal{A}_{r,s} \ , \ m \in \mathbb{N} \ .
	\end{aligned}
	\end{equation}
\end{defn}

\begin{rmk}\label{symplrmk}
	Note that, if $ \ \Phi \in \mathcal{C}_{r_0,s_0}(\mathcal{R}) \ $ is a symplectomorphism, then its inverse $ \ \Psi\in \mathcal{C}_{r_0-\rho,s_0-\sigma}(\mathcal{R}) \ $ satisfies $ \ (i) \ $ and $ \ (ii) \ $ of the above definition for all $ \ r \in (\rho,r_0] \ $ and all $ \ s \in (\sigma,s_0] \ $.
\end{rmk}

The following proposition establishes a sufficient condition for a regular map to be a symplectomorphism. Specifically, we show that, in analogy with the finite-dimensional case, the symplectic nature of a transformation can be determined by checking whether it preserves the fundamental Poisson brackets of the coordinate components.

\begin{prop}\label{fundamentalpb}
	Let $ \ \Phi = \id + (N, \Omega) : \mathscr{U}_{r_0}(\ell^{\infty}_{w}) \times \ell^\infty \to \ell^\infty_w \times \ell^\infty \ $ be a regular change of coordinates with $ |N|^{\infty,w}_{r_0,s_0}< \rho \ $. Let us denote by $ \ \hat{N} \ $ and $ \ \hat{\Omega} \ $ the functions $\hat{N}(y,\varphi) = y + N(y,\varphi)$ and $ \ \hat{\Omega}(y,\varphi) = \varphi + \Omega(y,\varphi) \ $, and assume that
	\begin{equation}
	\{\hat{N}^{(\ell)},\hat{\Omega}^{(m)}\} = \delta_{\ell m}, \qquad \{\hat{N}^{(\ell)},\hat{N}^{(m)}\} = 0, \qquad \{\hat{\Omega}^{(\ell)},\hat{\Omega}^{(m)}\} = 0 \ ,
	\end{equation}
	for all $  \ \ell,m \in \N \ $. Then, $ \ \Phi \ $ is a symplectomorphism.
	\end{prop}

\begin{proof}
Let $ \ m\in \N \ $. Since the power series defining the exponential function converges in any ball, and given that $ \ \Omega^{(m)} \ $, together with its derivatives, is bounded on $ \ \mathscr{U}_{r_0}(\ell^\infty_w)\times \ell^\infty \ $, the series
\begin{equation*}
\sum_{n=0}^\infty \frac{\im^n}{n!}(\hat{\Omega}^{(m)})^n \ , \qquad \sum_{n=0}^\infty \frac{\im^n}{n!}\frac{\partial(\hat{\Omega}^{(m)})^n}{\partial y_k} \qquad \text{and} \qquad \sum_{n=0}^\infty \frac{\im^n}{n!}\frac{\partial(\hat{\Omega}^{(m)})^n}{\partial \varphi_k} \ ,
\end{equation*}
converge uniformly on $ \ \mathscr{U}_{r_0}(\ell^\infty_w)\times \ell^\infty \ $ for all $ \ k \in \N \ $. Moreover, for all $ \ \ell,m \in \N \ $, both the series
\begin{equation*}
	\sum_{k \in \N}\sum_{n=0}^\infty \frac{\partial \hat{N}^{(\ell)}}{\partial y_k}\frac{1}{n!}(\hat{\Omega}^{(m)})^n\frac{\partial \hat{\Omega}^{(m)}}{\partial \varphi_k} \ \qquad \text{and} \qquad \sum_{k \in \N}\sum_{n=0}^\infty \frac{\partial \hat{N}^{(\ell)}}{\partial \varphi_k}\frac{1}{n!}(\hat{\Omega}^{(m)})^n\frac{\partial \hat{\Omega}^{(m)}}{\partial y_k}
\end{equation*}
converge absolutely for all $ \ (y,\varphi) \in \mathscr{U}_{r_0}(\ell^\infty_w)\times \ell^\infty \ $. Indeed, consider for e.g. the first one. We have that
\begin{gather*}
	\sum_{k \in \N}\sum_{n=0}^\infty \bigg|\frac{\partial \hat{N}^{(\ell)}}{\partial y_k}(y,\varphi)\bigg|\frac{1}{n!}|\hat{\Omega}^{(m)}(y,\varphi)|^n\bigg|\frac{\partial \hat{\Omega}^{(m)}}{\partial \varphi_k}(y,\varphi)\bigg| 
	\\
	=e^{|\hat{\Omega}^{(m)}(y,\varphi)|}\sum_{k \in \N} \bigg|\delta_{\ell k}+\frac{\partial N^{(\ell)}}{\partial y_k}(y,\varphi)\bigg|\bigg|\delta_{m k}+\frac{\partial \Omega^{(m)}}{\partial \varphi_k}(y,\varphi)\bigg|
	\\
	\leq e^{ |\hat{\Omega}^{(m)}(y,\varphi)|}\bigg( \delta_{\ell m}+\bigg|\frac{\partial N^{(\ell)}}{\partial y_m}(y,\varphi)\bigg|+\bigg|\frac{\partial \Omega^{(m)}}{\partial \varphi_\ell}(y,\varphi)\bigg|+\sum_{k \in \N}\bigg|\frac{\partial N^{(\ell)}}{\partial y_k}(y,\varphi)\bigg|\bigg|\frac{\partial \Omega^{(m)}}{\partial \varphi_k}(y,\varphi)\bigg|\bigg) \ .
	\end{gather*}
The latter quantity is finite because $ \ N^{(\ell)} \ $ and $ \ \Omega^{(m)}$ are regular for all $ \ \ell,m \in \N \ $. Similar estimates give that the other series is also absolutely convergent. Consequently, the order of the summations over $n$ and $k$ may be interchanged, yielding
\begin{equation*}
\{\hat{N}^{(\ell)}, e^{i \hat{\Omega}^{(m)}}\} = \sum_{n=0}^{\infty} \frac{i^n}{n!} \{\hat{N}^{(\ell)}, (\hat{\Omega}^{(m)})^n\} \ ,
\end{equation*}
for all $\ell, m \in \mathbb{N}$. Then we have
\begin{equation}\label{NOM}
\begin{aligned}
	\{\Phi^*y_\ell,\Phi^*e^{\im \varphi_m}\}&=\{\hat{N}^{(\ell)},e^{\im \hat{\Omega}^{(m)}}\}=\sum_{n=0}^{\infty}\frac{\im^n}{n!}\{\hat{N}^{(\ell)}, (\hat{\Omega}^{(m)})^n\} 
	\\
	&=\sum_{n=0}^\infty \frac{\im^n}{n!}n(\hat{\Omega}^{(m)})^{n-1}\{\hat{N}^{(\ell)}, \hat{\Omega}^{(m)}\}=\im\delta_{\ell m}\sum_{n=0}\frac{i^n}{n!}(\hat{\Omega}^{(m)})^{n}
	\\
	&=	\im \delta_{\ell m} e^{\im \hat{\Omega}^{(m)}}=\im\delta_{\ell m}\Phi^*e^{\im \varphi_m}=	\Phi^*\{y_\ell,e^{\im\varphi_m}\} \ .
	\end{aligned}
	\end{equation}
	Following an analogous argument, one can prove that
	\begin{equation}\label{OMOM}
		\{\Phi^*e^{\im \varphi_\ell},\Phi^*e^{\im \varphi_m}\}=\{e^{\im \hat{\Omega}^{(\ell)}},e^{\im \hat{\Omega}^{(m)}}\}=\Phi^*\{e^{\im \varphi_\ell},e^{\im \varphi_m}\}=0 \ .
	\end{equation}
	By a direct calculation, and by using the Leibniz rule, it is easy to prove that the results in \eqref{NOM} and \eqref{OMOM} imply that
	\begin{equation*}
		\{\Phi^* (y^\alpha e^{\im \nu \cdot \varphi}),\Phi^*(y^{\alpha'} e^{\im \nu' \cdot \varphi})\}=\Phi^*\{y^\alpha e^{\im \nu \cdot \varphi},y^{\alpha'} e^{\im \nu' \cdot \varphi}\}
	\end{equation*}
	holds for all $ \ (\nu,\alpha) \in \Z^{(\N)}\times\N_0^{(\N)} \ $. Then, item $ \ (i) \ $ of Definition \ref{sympl} follows by the total convergence of the involved series. Moreover,
	\begin{equation*}
		\begin{aligned}
		\{\Phi^*H,\Phi^*\varphi_m\}&=\{\Phi^*\sum_{\nu,\alpha}H_{\nu\alpha}y^\alpha e^{\im \nu \cdot \varphi}, \Phi^*\varphi_m\}=\sum_{\nu,\alpha}H_{\nu\alpha}\{\Phi^*y^\alpha e^{\im \nu \cdot \varphi}, \Phi^*\varphi_m\}
		\\
		&=\sum_{\nu,\alpha}H_{\nu\alpha}\{\hat{N}^\alpha e^{\im \nu \cdot \hat{\Omega}}, \hat{\Omega}^{(m)}\}=\sum_{\nu,\alpha}H_{\nu\alpha}\{\hat{N}^\alpha, \hat{\Omega}^{(m)}\}e^{\im \nu \cdot \hat{\Omega}} \ .
		\end{aligned}
	\end{equation*}
	By the Leibniz rule, and by the fact that $ \ \{\hat{N}^{(\ell)},\hat{\Omega}^{(m)}\} = \delta_{\ell m} \ $, it follow that $ \ \{\hat{N}^\alpha, \hat{\Omega}^{(m)}\}=\alpha_m \hat{N}^{\alpha-e_m} \ $. Therefore,
	\begin{equation*}
		\sum_{\nu,\alpha}H_{\nu\alpha}\{\hat{N}^\alpha, \hat{\Omega}^{(m)}\}e^{\im \nu \cdot \hat{\Omega}}=\sum_{\nu \alpha}H_{\nu \alpha}\alpha_m \hat{N}^{\alpha-e_m}e^{\im \nu \cdot \hat{\Omega}}=\Phi^*\frac{\partial H}{\partial y_k}=\Phi^*\{H,\varphi_m\} \ ,
	\end{equation*}
	which yields formula $ \ (ii) \ $ of Definition \ref{sympl}. This concludes the proof.
\end{proof}
\subsection{Hamilton's equations and Hamiltonian flow}\label{HAMEQSEC}
In this subsection, we extend the classical Hamiltonian formalism to our infinite-dimensional framework, defining the equations of motion and investigating their behavior under symplectic transformations.

 \smallskip

 Given $ \ r,s>0 \ $ and a function $ \ H \in \mathcal{A}_{r,s} \ $, we say that $ \ H \ $ induces a \emph{Hamiltonian flow} $ \ \phi_H \ $ on an open neighbourhood $ \ \mathcal{N}_0 \subseteq \mathscr{U}_{r}(\ell^\infty_w) \times \mathbb{T}^\infty \ $ of $ \ \{0\}\times \T^\infty \ $ if there exists a group action $ \ \phi_H \ $ of $ \ \mathbb{R} \ $ on $ \ \mathscr{N}_0 \ $
\begin{equation*}
\mathbb{R}\times\mathscr{N}_0 \ni (t,y,\varphi) \mapsto \phi^t_H(y,\varphi)\in \mathscr{N}_0 
\end{equation*}
such that, denoting $\phi^t_H(y,\varphi)=(y_1(t),\ldots,y_k(t),\ldots,\varphi_1(t),\ldots,\varphi_k(t),\ldots)$, the maps $t \mapsto y_k(t)$ and $t \mapsto \varphi_k(t)$ are differentiable functions for all $k \in \mathbb{N}$ and satisfy the \emph{Hamilton's equations}
\begin{equation*}
\begin{cases}
\dot{y}_k=\{H,y_k\}
\\
\dot{\varphi}_k=\{H,\varphi_k\}
\end{cases}
, \quad k \in \mathbb{N} \ .
\end{equation*}
In the following, we shall refer to $H$ as the \emph{Hamiltonian function}. 

Note that, in general, a Hamiltonian flow is not even continuous in time as a function of the time $t \in \R$ taking values in the Banach manifold. This is essentially due to the fact that, as already pointed out in Remark \ref{hamvfrmk}, the gradient of an analytic function defined on an infinite-dimensional Banach space $ \ X \ $ does not take values in $ \ X \ $, in general. Let us see this with a simple example. Given an unbounded sequence of real numbers $ \ \om=(\om_j)_{j \in \N} \ $ such that $ \ \sum_{j \in \N}\om_j w_j^{-1} \ $ is finite, consider the Hamiltonian function
\begin{equation}\label{hamex}
	H=\sum_{j \in \N}\om_j y_j \ , \qquad y \in \mathscr{U}_r(\ell^\infty_w) \ .
\end{equation}
Clearly, $ \ H \ $ belongs to $ \ \mathcal{A}_{r,s} \ $ for all $ \ r,s>0 \ $. The Hamilton's equations associated with the Hamiltonian $ \ H \ $ are given by
\begin{equation*}
\begin{cases}
	\dot{y}_k=0 
	\\
	\dot{\varphi}_k=\om_k
	\end{cases}
	, \quad k \in \mathbb{N} \ ,
\end{equation*}
which can be easily integrated yielding the solutions $ \ y_k(t)=\text{const} \ $, $ \ \varphi_k(t)=\om_k t+\text{const} \ $. Therefore, for all $ \ r>0 \ $, $ \ H \ $ induces the Hamiltonian flow on $ \ \mathscr{U}_{r}(\ell^\infty_w) \times \mathbb{T}^\infty \ $ defined by
\begin{equation*}
	\phi^t_H(y,\varphi)=(y,\varphi_1+\om_1t,\ldots,\varphi_k+\om_kt,\ldots) \ .
\end{equation*}
Since $ \ (\omega_j)_{j \in \mathbb{N}} \ $ is unbounded, this flow cannot be continuous in time. Indeed, for any $ \ (y,\varphi) \in \ell^\infty_w \times \ell^\infty \ $, any two distinct points on the orbit of $ \ (y,\varphi) \ $ under the action of $\phi_H$ are infinitely separated. Specifically, for any $ \ t_1 \neq t_2 \ $, one has
\begin{equation*}
\|\varphi+\om t_1-\varphi-\om t_2\|_{\ell^\infty} = \sup_{j \in \mathbb{N}} |\omega_j (t_1 - t_2)| = \infty.
\end{equation*}
For the sake of clarity, we emphasize that, in accordance with standard arguments from the theory of differentiable functions between Banach spaces, the gradient of the Hamiltonian defined in \eqref{hamex} belongs to the topological dual of $ \ \ell^\infty_w \times \ell^\infty \ $. This follows directly from the fact that the sum $ \ \sum_{j \in \mathbb{N}} \omega_j w_j^{-1} \ $ is finite.

\subsubsection{Invariance of the Hamilton's equations under symplectomorphisms} We now show that the Hamilton's equations are invariant under symplectic changes of coordinates. Namely, we prove that for any $ \ H \in \mathcal{A} \ $ inducing a Hamiltonian flow and any symplectomorphism $ \ \Phi \ $, the Hamilton's equations in the coordinate frame $ \ (x,\theta) := \Phi(y,\varphi) \ $ take the form
\begin{equation*}
\begin{cases}
\dot{x}_k = \{H \circ \Phi^{-1}, x_k\} \\
\dot{\theta}_k = \{H \circ \Phi^{-1}, \theta_k\}
\end{cases}
, \quad k \in \mathbb{N}.
\end{equation*}
To this end, we first establish the following lemma, which is central to our argument.

\begin{lemma}\label{lemmaflow}
	Let $ \ r,s>0 \ $ and let $ \ H \in \mathcal{A}_{r,s} \ $ be a Hamiltonian function that induces a Hamiltonian flow. Then, for all $ \ f \in \mathcal{R}_{r,s} \ $,
	\begin{equation*}
		\frac{d}{dt}(f\circ \phi^t_H)=\{H,f\}\circ \phi^t_H \ .
	\end{equation*}
	\begin{proof}
		The above equality follows easily for any power-trigonometric polynomial. The claim then follows for all $ \ f \in \mathcal{R}_{r,s} \ $ by an argument strictly analogous to the proof of Proposition \ref{proplie}.
	\end{proof}
\end{lemma}

With the help of the result established in the above lemma, let us prove the invariance of the Hamilton's equations under symplectic changes of coordinates. More precisely, we prove the following proposition, which is sufficient and crucial for our purposes.

\begin{prop}\label{propham}
Given $ \ r,s>0 \ $, let $ \ H \in \mathcal{A}_{r,s} \ $ and consider a symplectomorphism $ \ \Phi \in \mathcal{C}_{r,s}(\mathcal{R}) \ $. Assume that $ \ H\circ \Phi^{-1} \ $ induces a Hamiltonian flow $ \ \phi_{H \circ \Phi^{-1}} \ $ on $ \ \mathscr{U}_{r-\rho}(\ell^\infty_w)\times \T^\infty \ $, with $ \ \rho>0 \ $ as in Definition \ref{regchange}. Then, $ \ H \ $ induces a Hamiltonian flow on the neighbourhood $ \ \mathscr{N}_0:=\Phi^{-1}(\mathscr{U}_{r-\rho}(\ell^\infty_w)\times \T^\infty) \ $ of $ \ \{0\}\times \T^\infty \ $. In particular, we have that
\begin{equation*}
	\phi^t_H=\Phi^{-1}\circ \phi^t_{H \circ \Phi^{-1}}\circ \Phi 
\end{equation*}
holds for all $ \ t \in \R \ $.
\begin{proof}
First, note that $ \ H \circ \Phi^{-1} \ $ is a well defined analytic function belonging to $ \ \mathcal{A}_{r-\rho,s-\sigma} \ $, with $ \ \sigma>0 \ $ as in Definition \ref{regchange}. Moreover, since $ \ \Phi \ $ is in particular a change of coordinates, as already noted in the proof of Theorem \ref{FFoliation}, the open set $ \ \Phi^{-1}(\mathscr{U}_{r-\rho}(\ell^\infty_w)\times \T^\infty) \ $ is indeed a neighbourhood of $ \ \{0\}\times \T^\infty \ $. Now, by assumption, $ \ H\circ \Phi^{-1} \ $ induces a Hamiltonian flow. Hence, the above Lemma implies 
\begin{equation*}
	\frac{d}{dt}(y_k \circ \Phi^{-1} \circ \phi^t_{H\circ\Phi^{-1}})=\{H\circ \Phi^{-1},y_k \circ \Phi^{-1}\}\circ \phi^t_{H\circ\Phi^{-1}} \ .
\end{equation*}
Since $ \ \Phi \ $ is a symplectomorphism, by Remark \ref{symplrmk}, the right-hand side of the above equation is equal to
\begin{equation*}
\{H,y_k\}\circ \Phi^{-1}\circ \phi^t_{H\circ\Phi^{-1}}
\end{equation*}
By a similar, slightly refined argument one can prove that
\begin{equation*}
\frac{d}{dt}(\varphi_k \circ \Phi^{-1} \circ \phi^t_{H\circ\Phi^{-1}})=\{H,\varphi_k\}\circ \Phi^{-1}\circ \phi^t_{H\circ\Phi^{-1}} \ .
\end{equation*}
We remark that the argument requires refinement because, in general, $ \ \varphi_k \circ \Phi^{-1} \ $ does not belong to $ \ \mathcal{R} \ $. This precludes the direct application of Lemma \ref{lemmaflow}. Nevertheless, this is not a real obstacle: by the definition of change of coordinates and denoting the new angular variables by $ \ \theta_k \ $ (with $ \ k \in \mathbb{N} \ $), it follows that $ \ \varphi_k \circ \Phi^{-1} \ $ is simply the projection onto the $ \ \theta_k \ $ coordinate plus a function in $ \ \mathcal{R}_{r,s} \ $. Therefore, we conclude that $ \ \Phi^{-1} \circ \phi^t_{H \circ \Phi^{-1}}\circ \Phi \ $ is the Hamiltonian flow induced by $ \ H \ $.
\end{proof}
\end{prop}

Let us make some comments on the above proposition. The fact that $ \ H \circ \Phi^{-1} \ $ induces a Hamiltonian flow on $ \ \mathscr{U}_{r-\rho}\times \T^\infty \ $ means that there exists a group action $ \ \phi_{H\circ \Phi^{-1}} \ $
\begin{equation*}
\mathbb{R}\times \mathscr{U}_{r-\rho}\times \T^\infty \ni (t,x,\theta) \mapsto \phi^t_{H\circ\Phi^{-1}}(x,\theta)\in \mathscr{U}_{r-\rho}\times \T^\infty
\end{equation*}
such that, denoting $\phi^t_{H \circ \Phi^{-1}}(x,\theta)=(x_1(t),\ldots,x_k(t),\ldots,\theta_1(t),\ldots,\theta_k(t),\ldots)$, the maps $t \mapsto x_k(t)$ and $t \mapsto \theta_k(t)$ are differentiable functions for all $k \in \mathbb{N}$ and satisfy the Hamilton's equations
\begin{equation}\label{ODE}
\begin{cases}
\dot{x}_k=\{H\circ \Phi^{-1},x_k\}
\\
\dot{\theta}_k=\{H\circ\Phi^{-1},\theta_k\}
\end{cases}
, \quad k \in \mathbb{N} \ .
\end{equation}
Then, by Proposition \ref{propham}, the Hamilton's equations associated with the Hamiltonian $ \ H \ $, namely, 
\begin{equation*}
\begin{cases}
\dot{y}_k=\{H,y_k\}
\\
\dot{\varphi}_k=\{H,\varphi_k\}
\end{cases}
, \quad k \in \mathbb{N} \ ,
\end{equation*}
admit a solution $ \ (y_1(t),\ldots,y_k(t),\ldots,\varphi_1(t),\ldots,\varphi_k(t),\ldots) \ $ for any initial datum $ \ (y,\varphi) \ $ in the open neighbourhood $ \ \mathscr{N}_0:= \Phi^{-1}(\mathscr{U}_{r-\rho}(\ell^\infty_w)\times \T^\infty) \ $ of $ \ \{0\}\times \T^\infty \ $. Specifically, the solution is given by 
\begin{equation*}
\phi^t_H(y,\varphi)=\Phi^{-1}(\phi^t_{H \circ \Phi^{-1}}(\Phi(y,\varphi))) \ , \qquad \text{for all} \ \ (y,\varphi)\in \mathscr{N}_0 \ ,
\end{equation*}
and $ \ \phi_H \ $ is a group action of $ \ \R \ $ on $ \ \mathscr{N}_0 \ $.

\subsubsection{Action-angle coordinates}\label{AAsub} The above argument suggests that a good strategy to solve Hamilton's equations for a given Hamiltonian $ \ H \in \mathcal{A} \ $ consists in trying to find a symplectomorphism $ \ \Phi \ $ such that the Hamilton's equations for the new Hamiltonian $ \ H \circ \Phi^{-1} \ $ are easily solvable. For instance, if $ \ H \circ \Phi^{-1} \ $ does not depend on the variables $ \ \theta_k \ $, $ \ k \in \N \ $, the infinite system of ODE in \eqref{ODE} admit the unique solution
\begin{equation*}
\phi^t_{H \circ \Phi^{-1}}(x,\theta)=(x,\theta_1+\om_1(x)t,\ldots,\theta_k+\om_k(x)t,\ldots) \ , \quad \text{where} \ \ \om_k(x):=\frac{\partial H\circ \Phi^{-1}}{\partial x_k}(x) \ , \ k \in \N \ .
\end{equation*}
Namely, $ \ H\circ\Phi^{-1} \ $ induces a Hamiltonian flow $ \ \phi_{H \circ \Phi^{-1}} \ $ such that, for all $ \ x \ $, the leaf $ \ \{x\}\times \T^\infty \ $ is invariant and the restriction of the Hamiltonian flow on such invariant manifold, namely
\begin{equation*}
	\{x\}\times \T^\infty \in \theta \mapsto \phi^{t}_{H \circ \Phi^{-1}}(x,\theta) \in \{x\}\times \T^\infty \ ,
\end{equation*}
is simply a linear flow on the infinite torus. Therefore, $ \ \mathscr{N}_0 \ $ is foliated by infinite tori of the form $ \ \Phi^{-1}(\{x\}\times \T^\infty) \ $ and, by virtue of Proposition \ref{propham}, each leaf is invariant for the Hamiltonian flow $ \ \phi_H \ $, which is conjugated by $ \ \Phi \ $ to a linear flow on the infinite torus. The set of coordinates $ \ (x,\theta) \ $ defined by such symplectomorphism are called \emph{action-angle coordinates}. 

A final observation is in order: as previously observed, a Hamiltonian flow is not even continuous in time, in general. Thus, at this stage, we cannot talk about topological conjugacy, even if the conjugating map has very strong regularity. Such kind of dynamics is always continuous if one endows the infinite torus with the product topology.

In the next section, we provide conditions for a Hamiltonian $ \ H \in \mathcal{A} \ $ to admit a symplectomorphism $ \ \Phi \ $ such that the new Hamiltonian 
$ \ H \circ \Phi^{-1} \ $ is independent of the variables $ \ \theta_k \ $ for $ \ k \in \mathbb{N} \ $.

\section{Symplectic straightening and Liouville-Arnold integrability in infinite dimension}
\noindent We shall identify the elements of the torus $ \ \mathbb{T}^\infty$ with their representatives in the universal cover $\ell^\infty$. With a slight abuse of notation, we shall use the same symbol $\varphi$ for both, as well as for functions defined on the torus and their associated multi-periodic counterparts defined on $\ell^\infty$, provided that this does not lead to confusion.
 
 The primary goal of this section is to show that, given $P \in \mathcal{R}^{\infty,w}_{r,s}$ sufficiently small for some $r,s>0$, if the components of the map $\Gamma_P(y,\varphi)=y+P(y,\varphi)$ Poisson-commute, namely
\begin{equation}\label{commutano}
\{\Gamma_P^{(m)},\Gamma_P^{(n)}\}=\sum_{k \in\N}\bigg(\frac{\partial \Gamma^{(m)}_P}{\partial y_k}\frac{\partial \Gamma^{(n)}_P}{\partial \varphi_k}-\frac{\partial \Gamma^{(m)}_P}{\partial \varphi_k}\frac{\partial \Gamma^{(n)}_P}{\partial y_k} \bigg)=0, \qquad \text{for all} \ m,n \in \N,
\end{equation}
there exists a symplectic straightening of the foliation provided by Theorem \ref{FFoliation}. As a dynamic counterpart we shall prove in Theorem \ref{completethm} that any analytic Hamiltonian that Poisson commutes with every component of $\Gamma_P$ is completely integrable and admits action-angle variables. 

More precisely, we construct a coordinate transformation $\Delta$ such that the composition 
\begin{equation}\label{delta}
\Phi := \Delta \circ \Phi_P, \quad \text{where} \quad \Phi_P(y,\varphi)=(y+P(y,\varphi),\varphi),
\end{equation}
defines a symplectomorphism (cf. Definition \ref{sympl}). This map induces a new coordinate system $(x, \theta)$ where $x$ depends on the original variables $(y, \varphi)$ exclusively through $\Gamma_P$. Such a coordinate frame, which we refer to as \emph{action-angle} coordinates, is adapted to the foliation: the coordinates $x = (x_k)_{k \in \mathbb{N}}$ act as labels for the leaves, while $\theta = (\theta_k)_{k \in \mathbb{N}}$ serves as the coordinate system along each leaf. To this end, we lift the classical generating function framework to an infinite-dimensional setting by adapting the standard Hamiltonian machinery to our context. It is clear that in infinite dimension one has to carefully check that all the algebraic manipulations make sense at a quantitative level. The key technical result is the following.

\begin{thm}[Symplectic Straightening]\label{SYMPL}
Let $r_0, s_0 > 0$, $P \in \mathcal{R}^{\infty,w}_{r_0,s_0}$ and let
\begin{equation*}
\Gamma_P : \mathscr{U}_{r_0}(\ell^\infty_w) \times \mathbb{T}^\infty \to \ell^\infty_w \qquad (y,\varphi) \mapsto \Gamma_P(y, \varphi) = y + P(y, \varphi)
\end{equation*}
 be such that $ \ \{\Gamma_P^{(m)},\Gamma^{(n)}_P\}=0 \ $ for all $m,n \in \N$. Then, there is a constant $C > 1$ such that, the smallness condition
\begin{equation*}
\|P\|^{\infty,w}_{r_0,s_0} < \frac{1}{4C} \min\{r_0, s_0, 1\}
\end{equation*}
implies that, for all $\rho, \sigma$ such that
$$ \ C \|P\|^{\infty,w}_{r_0,s_0} \leq \rho < r_0/4 \qquad \text{and}  \qquad 2C \|P\|^{\infty,w}_{r_0,s_0} \leq \sigma < s_0/2,$$
there exists a symplectomorphism 
\begin{equation*}
\Phi : \mathscr{U}_{r_0-4\rho}(\ell^\infty_w) \times \mathbb{T}^\infty \to \ell^\infty_w \times \mathbb{T}^\infty, \quad (y,\varphi) \mapsto \Phi(y,\varphi)=(x,\theta),
\end{equation*}
that straightens the foliation induced by the map $ \ \Gamma_P \ $. Namely for all $\gamma \in \mathscr{U}_{r_0-\rho}(\ell^\infty_w)$ $$\Phi(\Gamma^{-1}_P(\gamma))=\{x\}\times \T^\infty, \quad x=\gamma-A_0(\gamma),$$ 
where $A_0$ is the average over the angles of the fixed point solution  $A$ defined in equation \eqref{A}.
\noindent Note that $\Phi$ is a regular change of coordinates of class $ \  \mathcal{C}^{\rho,\sigma/2}_{r_0-4\rho, s_0}(\mathcal{R}) \ $, with~$ \ \Phi^{-1}\in\mathcal{C}^{\rho,\sigma/2}_{r_0-3\rho, s_0-\sigma}(\mathcal{R})$.
\end{thm}

\subsection{An informal presentation of the proof}\label{heuristics}

To provide some intuition, let us sketch the underlying heuristic, which follows the ideas developed for the finite-dimensional setting. For the moment we disregard any issues concerning the regularity of the maps and the convergence of the series involved, as they will be addressed rigorously later in Subsection \ref{SubSs}.

\subsubsection*{Action coordinates}\label{CAA} Using the same notation as in Theorem \ref{FFoliation}, all $ \ \gamma \in \mathscr{U}_r(\ell^\infty_w) \ $, we endow the infinite torus $ \ \Gamma^{-1}_P(\gamma) \ $ with the pullback measure $ \ \lambda:=p_\varphi^*\mu \ $, where $ \ p_\varphi \ $ is the projection onto the second factor $ \ p_\varphi :\ell^\infty_w \times \T^\infty \to \T^\infty \ $ and $ \ \mu \ $ is the Lebesgue product measure of the infinite torus, namely
\begin{equation*}
\int_{\T^\infty}d\mu=\lim_{N\to\infty}\frac{1}{(2\pi)^N}\int_{\T^N}d\varphi_1 \ldots d\varphi_N \ .
\end{equation*}
In this way, for all $ \ \gamma \in \mathscr{U}_{r}(\ell^\infty_w) \ $ and for all $ \ f \in \mathcal{A}_{r+\rho,s_0} \ $, with $ \ \rho:=3|P|_{r_0,s_0}^{\infty,w} \ $, one has that
\begin{equation*}
\int_{\Gamma_P^{-1}(\gamma)}f \ d\lambda := \int_{\T^\infty}f\circ (\left. \Phi_{P} \right|_{\Gamma_P^{-1}(\gamma)})^{-1} \ d\mu = \lim_{N \to \infty}\frac{1}{(2\pi)^N}\int_{\T^N}f(\gamma-A(\gamma,\varphi),\varphi) \ d\varphi_1 \ldots d\varphi_N \ ,
\end{equation*}
where the map $ \ A\in \mathcal{A}_{r,s_0} \ $ is uniquely defined by (see formula \eqref{A} in the proof of Theorem \ref{FFoliation})
\begin{equation}\label{GA}
\Gamma_P(\gamma-A(\gamma,\varphi),\varphi)=\gamma \ \quad \text{for all} \ (\gamma,\varphi) \in \mathscr{U}_r(\ell^\infty_w)\times \ell^\infty \ .
\end{equation}
Drawing an analogy with statistical mechanics, we refer to the quantity above, which is well defined by Proposition \ref{manalytic}, as the average of $f$ over the infinite torus $\Gamma_P^{-1}(\gamma)$. 

Consider the map $\gamma \mapsto x(\gamma)=(x_k(\gamma))_{k \in \mathbb{N}}$ defined by
\begin{equation}\label{actionintegral}
\mathscr{U}_r(\ell^\infty_w) \ni \gamma \mapsto x(\gamma) := \bigg( \int_{\Gamma_P^{-1}(\gamma)} y_k \, d\lambda \bigg)_{k \in \mathbb{N}} \in \ell^\infty_w
\end{equation}
and note that, since for all $k \in \mathbb{N}$ we have
\begin{equation}\label{ipsilon}
y_k \circ (\left. \Phi_{P} \right|_{\Gamma_P^{-1}(\gamma)})^{-1}(\gamma,\varphi)=\gamma_k-A^{(k)}(\gamma,\varphi),
\end{equation}
the definition of the integral over $\Gamma^{-1}_P(\gamma)$ with respect to the pullback measure $\lambda$ yields the explicit formula
\begin{equation}\label{ics}
x_k(\gamma)=\gamma_k-A^{(k)}_0(\gamma), \quad \forall \gamma \in \mathscr{U}_r(\ell^\infty_w), \ k \in \mathbb{N}\,,
\end{equation}
where $A^{(k)}_0$ is the average of $A^{(k)}$ over the angles.
Notice that definition of the map $x(\gamma)$ in \eqref{actionintegral} provides the infinite-dimensional analogue of the integral of the Liouville one-form $\sum_n y_n d\varphi_n$ along the $k^\mathrm{th}$-cycle of the torus defined by the level set~$\Gamma_P^{-1}(\gamma)$. It is important to remark that the consistency of this identification relies on the commuting property \eqref{commutano}. Indeed, as a consequence of the relations \eqref{Acommutes}, the integral of the Liouville form along the full cycle parametrized by $\varphi_k$ results in a constant value. Since the torus is endowed with a probability measure, this constant must coincide with the average over the torus, justifying the expression for the candidate action coordinates $x_k$ as the mean values of the momenta $y_k$.

We shall prove that, given a suffciently small $P \in \mathcal{R}^{\infty,w}_{r,s}$, the map
\begin{equation*}
\mathscr{U}_{r}(\ell^\infty_w)\times \ell^\infty \ni (\gamma,\varphi) \mapsto (x(\gamma),\varphi) \in \ell^\infty_w \times \ell^\infty 
\end{equation*}
indeed defines a regular change of coordinates, in the sense of definition \ref{regchange}. We claim that the coordinates $ \ x=(x_k)_{k \in \N} \ $ defined above will play the role of \emph{action} coordinates for the given foliation. However, the definition of action variables is inherently incomplete without their canonical counterparts; indeed, their characterization as actions is only fully realized when they are paired with the corresponding angle variables $ \ \theta=(\theta_k)_{k \in \N} \ $ that satisfy the standard Poisson bracket relations. To complete this framework, we now briefly illustrate the construction of the associated angles. We achieve this by adapting the classical generating function method to the present infinite-dimensional setting.

\subsubsection*{Angle coordinates.}\label{Cac} We seek an analytic generating function $ \ S=S(x,\varphi) \ $ that satisfies
\begin{equation}\label{yS}
y_k(x,\varphi)=\dfrac{\partial S}{\partial \varphi_k}(x,\varphi) \quad \text{for all } k \in \mathbb{N} \ 
\end{equation}
and such that the relations 
\begin{equation}\label{thetaS}
\theta_k(x,\varphi)=\dfrac{\partial S}{\partial x_k}(x,\varphi) \quad \text{for all } k \in \mathbb{N}
\end{equation}
are invertible in $ \ \ell^\infty \ $ with respect to $ \ \varphi \ $. Assuming that, for all $ \ k \in \N \ $, the functions $ \ \theta_k(x,\varphi)-\varphi_k \ $ and $ \ y_k(x,\varphi) \ $ are sufficiently regular and $ \ 2\pi$-periodic in each $ \ \varphi_j$-variable, $ \ j \in \N \ $, equations \eqref{yS} and \eqref{thetaS} define a set of coordinates $(x,\theta)$ such that, for all $m,n\in\N$,
\begin{equation}
    \{x_m,x_n\}=0, \quad \{x_m,\theta_n\}=\delta_{n,m}, \quad \{\theta_m,\theta_n\}=0,
\end{equation}
and the list $ \ (\theta_k)_{k \in \mathbb{N}} \ $ constitutes a system of angle coordinates on each leaf of the foliation. 

We assume that the generating function $S$ is an analytic perturbation of the generator of the identity,
\begin{equation}\label{SW}
S(x, \varphi) = \sum_{k \in \mathbb{N}} x_k \varphi_k - W(x, \varphi),
\end{equation}
The perturbation $ \ W \ $ is assumed to be $ \ 2\pi$-periodic in each $ \ \varphi_k$-variable, $ \ k \in \N \ $. Then, Equation \eqref{yS} reads
\begin{equation}\label{W}
y_k(x,\varphi)=x_k-\frac{\partial W}{\partial \varphi_k}(x,\varphi) \quad \text{for all} \ k \in \N \ .
\end{equation}
Given $ \ m \in \N \ $, let $ \ \Lambda_m \ $ be the subset of $ \ \Z^{(\N)} \ $ defined by
\begin{equation*}
\Lambda_m:=\{ \nu \in \Z^{(\N)} \ : \ \nu_m\neq 0 \ , \nu_j=0 \ \ \text{for all} \ j<m \} \ .
\end{equation*}
We claim that the function 
\begin{equation}\label{WW}
W(x,\varphi)=\sum_{m \in \N}\sum_{\nu \in \Lambda_m}\dfrac{A^{(m)}_\nu(\gamma(x))}{i \nu_m}e^{i \nu \cdot \varphi}
\end{equation}
satisfies equation \eqref{W}. To justify this assertion at a formal level, deferring analytic technicalities to Subsection \ref{SubSs}, we first show that if the components of $ \ \Gamma_P \ $ Poisson commute, the function $ \ A \in \mathcal{A}_{r,s} \ $, whose definition is recalled in \eqref{GA}, satisfies the relations
\begin{equation}\label{Acommutes}
\nu_\ell A^{(m)}_\nu = \nu_m A^{(\ell)}_\nu \quad \text{for all} \ \ell,m \in \mathbb{N} \ .
\end{equation}
Note that, by the definition of the function $ \ A \ $, see \eqref{GA}, it follows that
\begin{equation*}
\left. \frac{\partial}{\partial \varphi_k} \right|_{\gamma=\Gamma_P(y,\varphi)} \Gamma_P(\gamma-A(\gamma,\varphi),\varphi)=0 \ ,
\end{equation*}
which implies
\begin{equation}\label{oo}
\sum_{\ell \in \mathbb{N}}\frac{\partial \Gamma_P^{(m)}}{\partial y_\ell}(y,\varphi)\frac{\partial A^{(\ell)}}{\partial \varphi_k}(\Gamma_P(y,\varphi),\varphi)=\frac{\partial \Gamma_P^{(m)}}{\partial \varphi_k}(y,\varphi) \ , \quad \text{for all} \ \ m \in \N \ .
\end{equation}
To ease the notation, in the following computations it is understood that the derivatives of $A$ are pre-composed with the map $(y,\varphi) \mapsto (\Gamma_P(y,\varphi),\varphi) \ $, even if not explicitly indicated. Using the assumption that the components of $ \ \Gamma_P \ $ Poisson commute, we have:
\begin{equation}\label{oa}
\begin{aligned}
0&=\{\Gamma_P^{(m)},\Gamma_P^{(n)}\}=\sum_{k \in \mathbb{N}} \left( \frac{\partial \Gamma_P^{(m)}}{\partial y_k}\frac{\partial \Gamma_P^{(n)}}{\partial \varphi_k}-\frac{\partial \Gamma_P^{(m)}}{\partial \varphi_k}\frac{\partial \Gamma_P^{(n)}}{\partial y_k} \right)
\\
&\stackrel{\eqref{oo}}{=}\sum_{k \in \mathbb{N}}\frac{\partial \Gamma_P^{(m)}}{\partial y_k}\sum_{\ell \in \mathbb{N}}\frac{\partial \Gamma_P^{(n)}}{\partial y_\ell}\frac{\partial A^{(\ell)}}{\partial \varphi_k}-	\sum_{k \in \mathbb{N}}\sum_{\ell \in \mathbb{N}}\frac{\partial \Gamma_P^{(m)}}{\partial y_\ell}\frac{\partial A^{(\ell)}}{\partial \varphi_k}\frac{\partial \Gamma_P^{(n)}}{\partial y_k} 
\\
&=\sum_{k \in \mathbb{N}}\sum_{\ell\in \mathbb{N}}\frac{\partial \Gamma_P^{(m)}}{\partial y_k}\frac{\partial \Gamma_P^{(n)}}{\partial y_\ell}\bigg(\frac{\partial A^{(\ell)}}{\partial \varphi_k}-\frac{\partial A^{(k)}}{\partial \varphi_\ell}\bigg) \ .
\end{aligned}
\end{equation}
\noindent Notice that the operator $ \ \partial_y \Gamma_P \ $ is invertible, being a small bounded perturbation of the identity. Then the above equalities yield
\begin{equation*}
\frac{\partial A^{(\ell)}}{\partial \varphi_k}=\frac{\partial A^{(k)}}{\partial \varphi_\ell} \ ,
\end{equation*}
which proves formula \eqref{Acommutes}. Consequently, the following equalities hold for all $ \ k \in \mathbb{N} \ $:
\begin{equation*}
\begin{aligned}
\frac{\partial W}{\partial \varphi_k}(x,\varphi)&=\sum_{m \in \N}\sum_{\nu \in \Lambda_m} i \nu_k \dfrac{A^{(m)}_\nu(\gamma(x))}{i \nu_m}e^{i \nu \cdot \varphi} \stackrel{\eqref{Acommutes}}{=} \sum_{m \in \N}\sum_{\nu \in \Lambda_m} A^{(k)}_\nu(\gamma(x))e^{i \nu \cdot \varphi}
\\
&=\sum_{\nu \neq 0}A_\nu^{(k)}(\gamma(x))e^{i \nu \cdot \varphi} = A^{(k)}(\gamma(x),\varphi)-A^{(k)}_0(\gamma(x))
\\
&=\gamma_k(x)-A^{(k)}_0(\gamma(x))-(\gamma_k(x)-A^{(k)}(\gamma(x),\varphi))
\\
&=x_k-y_k(x,\varphi) \ .
\end{aligned}
\end{equation*}
The last equality follows from Equations \eqref{ipsilon} and \eqref{ics}. Hence, the claim is formally proved. According to \eqref{thetaS}, and recalling \eqref{SW} and \eqref{WW}, we can assert that the candidate angle variables are formally defined by
\begin{equation}\label{angleS}
\theta_k(x,\varphi)=\varphi_k-\sum_{m\in \N}\sum_{\nu \in \Lambda_m}\frac{1}{\im \nu_m}\frac{\partial A^{(m)}_\nu (\gamma(x))}{\partial x_k}e^{\im \nu \cdot \varphi} \ , \quad k \in \N \ .
\end{equation}

\subsection{Proof of the Symplectic straightening theorem}\label{SubSs} 
We proceed with a series of Lemmas. \\ 
We recall that, given $ \ G \in \mathcal{A}^{\infty,w}_{r,s} \ $, we denote by $ \ \Phi_G \ $ the map
\begin{equation*}
\Phi_G(y,\varphi)=(y+G(y,\varphi),\varphi) \ , \qquad \forall \ y \in \mathscr{U}_r(\ell^\infty_w) \ , \ \varphi \in \ell^\infty \ .
\end{equation*}
The following result  ensures that the fixed-point argument developed in Section \ref{sectionstability} (see Theorem~\ref{fix}) is closed within the framework of regular Hamiltonians. 

\begin{lemma}\label{ao}
		Given $r, s > 0$, let $0 \le \rho < r$. Consider $F \in \mathscr{U}_{r/2}(\mathcal{R}^{\infty,w}_{r,s})$ and let $G \in \overline{\mathscr{U}}_{\rho}(\mathcal{A}^{\infty}_{r-\rho,s}) \ $. Assume that $G = \Phi^*_G F:=F\circ\Phi_G \ .$ Then, $G \in \mathcal{R}^{\infty,w}_{r-\rho,s} \ $ and satisfies the estimate
	\begin{equation*}
	\|G\|^{\infty,w}_{r-\rho,s} \leq \frac{2}{1-(2\|F\|_{r,s}^{\infty,w}/r)}\|F\|_{r,s}^{\infty,w} \ ,
	\end{equation*}
\end{lemma}
	\begin{proof}
		Let us begin by noting that the bound 
		\begin{equation*}
		 |G|^{\infty,w}_{r-\rho,s}= |\Phi^*_GF|^{\infty,w}_{r-\rho,s}\leq 2 |F|^{\infty,w}_{r,s} 
		\end{equation*}
		 follows by Proposition \ref{composition}. We now give suitable estimates for the derivatives. We adopt the notation $ \ (x,\theta)=\Phi_G(y,\varphi) \ $. Then, for all $ \ m,k \in \N \ $, by the chain rule one has that
		\begin{gather*}
		 \ \frac{\partial G^{(m)}}{\partial y_k}=\frac{\partial  \ \Phi^*_G F^{(m)}}{\partial y_k}=  \sum_{\ell \in \N}\Phi^*_G\frac{\partial F^{(m)}}{\partial x_\ell} \ \bigg(\delta_{\ell k}+\frac{\partial G^{(\ell)}}{\partial y_k}\bigg)
\\
		 \qquad \frac{\partial G^{(m)}}{\partial \varphi_k}=\frac{\partial  \ \Phi^*_G F^{(m)}}{\partial \varphi_k}=  \sum_{\ell \in \N}\Phi^*_G\frac{\partial F^{(m)}}{\partial x_\ell} \ \frac{\partial G^{(\ell)}}{\partial \varphi_k}+\Phi^*_G\frac{\partial F^{(m)}}{\partial \varphi_k}
		\end{gather*}
		The above equations can also be written as
		\begin{gather*}
			\sum_{\ell \in \N} \bigg( \delta_{m \ell}-\Phi^*_G\frac{\partial F^{(m)}}{\partial x_\ell}\bigg)\frac{\partial G^{(\ell)}}{\partial y_k}=\Phi^*_G\frac{\partial F^{(m)}}{\partial x_k}
			\\
			 \ \ \sum_{\ell \in \N} \bigg( \delta_{m \ell}-\Phi^*_G\frac{\partial F^{(m)}}{\partial x_\ell} \bigg)\frac{\partial G^{(\ell)}}{\partial \varphi_k}=\Phi^*_G\frac{\partial F^{(m)}}{\partial \varphi_k} \ ,
		\end{gather*}
		and they can be inverted by Neumann series, namely
		\begin{gather*}
		 \ \frac{\partial G}{\partial y}=\sum_{n=0}^\infty\bigg(\Phi^*_G\frac{\partial F}{\partial x}\bigg)^n\Phi^*_G\frac{\partial F}{\partial x}
		\\
		\quad \frac{\partial G}{\partial \varphi}=\sum_{n=0}^\infty\bigg(\Phi^*_G\frac{\partial F}{\partial x}\bigg)^n\Phi^*_G\frac{\partial F}{\partial \varphi} \ . \label{secondphi}
		\end{gather*}
		Indeed, by the triangular inequality and by the algebra property established in Proposition \ref{algebra}, one has
		\begin{equation*}
		\begin{aligned}
			\sup_{k \in \N}\sum_{m \in \N}\bigg|\frac{\partial G^{(m)}}{\partial y_k}\bigg|_{r-\rho,s} 
			&\leq
		    \sup_{k \in \N}\sum_{n=0}^\infty\sum_{\stackrel{m \in \N}{\ell_1,\ldots,\ell_n \in \N}}\bigg|\Phi^*_G\frac{\partial F^{(m)}}{\partial x_{\ell_1}}\ldots \Phi^*_G\frac{\partial F^{(\ell_{n-1})}}{\partial x_{\ell_n}} \  \Phi^*_G\frac{\partial F^{\ell_n}}{\partial x_k}\bigg|_{r-\rho,s}
		    \\
		    &\leq \sum_{n=0}^\infty \sup_{\ell_{n+1} \in \N} \sum_{\ell_0,\ell_1,\ldots,\ell_n \in \N}\prod_{j=0}^{n}\bigg|\Phi^*_G\frac{\partial F^{(\ell_j)}}{\partial x_{\ell_{j+1}}}\bigg|_{r-\rho,s} \ .
			\end{aligned}
		\end{equation*}
		By Proposition \ref{composition}, above quantity is bounded by
		\begin{equation*}
		\begin{split}
		&\sum_{n=0}^\infty 2^{n+1} \sup_{\ell_{n+1} \in \mathbb{N}} \left( \sum_{\ell_0,\dots,\ell_n \in \mathbb{N}} \prod_{j=0}^{n} \bigg|\frac{\partial F^{(\ell_j)}}{\partial x_{\ell_{j+1}}}\bigg|_{r,s} \right) \\
		&\qquad \leq \sum_{n=0}^\infty 2^{n+1}\prod_{j=0}^{n} \sup_{\ell_{j+1}\in \mathbb{N}} \sum_{\ell_j \in \mathbb{N}} \bigg|\frac{\partial F^{(\ell_j)}}{\partial x_{\ell_{j+1}}}\bigg|_{r,s} \\
		&\qquad \leq \sum_{n=0}^\infty \bigg(\frac{2\|F\|_{r,s}^{\infty,w}}{r}\bigg)^{n+1} 
		= \frac{2\|F\|_{r,s}^{\infty,w}}{r-2\|F\|_{r,s}^{\infty,w}} \ .
		\end{split}
		\end{equation*}
		Therefore,
		\begin{equation*}
        r\sup_{k \in \N}\sum_{m \in \N}\bigg|\frac{\partial G^{(m)}}{\partial y_k}\bigg|_{r-\rho,s}\leq \frac{2}{1-(2\|F\|_{r,s}^{\infty,w}/r)}\|F\|_{r,s}^{\infty,w} \ .
		\end{equation*}
		By analogous computations one can prove the same bound for the derivatives with respect to the periodic variables, yielding the claim.
	\end{proof}

The next step is to prove that the map 
\begin{equation*}
\mathscr{U}_r(\ell^\infty_w) \ni \gamma \mapsto x(\gamma) := \bigg( \int_{\Gamma_P^{-1}(\gamma)} y_k \, d\lambda \bigg)_{k \in \mathbb{N}} \in \ell^\infty_w
\end{equation*}
induces a regular change of coordinates. Recalling that, for all $ \ k \in \N \ $, we have the explicit formula
\begin{equation*}
	\int_{\Gamma_P^{-1}(\gamma)} y_k \, d\lambda = \gamma_k-A_0^{(k)}(\gamma) \ ,
\end{equation*}
where $ \ A_0 \ $ is the average of the function $ \ A \ $ defined by $ \ A= \Phi^*_{-A}P \ $ (see Theorem \ref{FFoliation}). The following lemma gives the desired claim.

\begin{lemma}\label{A_0,B}
	Given $r_0, s_0 > 0$, let $P \in \mathcal{R}^{\infty,w}_{r_0,s_0}$ satisfy $\|P\|^{\infty,w}_{r_0,s_0} < r_0/24$ and let $A_0$ be the average defined above. Then the map
	\begin{equation*}
	\Phi_{-A_0}:\mathscr{U}_{r_0-\rho}(\ell^\infty_w)\times \ell^\infty \to \ell^\infty_w\times \ell^\infty, \quad (\gamma,\varphi) \mapsto (\gamma-A_0(\gamma),\varphi),
	\end{equation*}
	is a regular change of coordinates for all $\rho$ such that $7\|P\|^{\infty,w}_{r_0,s_0}\leq\rho<r_0/2 $ and $0\leq \sigma <s_0/2$. Its inverse is given by
	\begin{equation*}
	\Phi_{B}:\mathscr{U}_{r_0-2\rho}(\ell^\infty_w)\times \ell^\infty \to \ell^\infty_w\times \ell^\infty, \quad (x,\varphi) \mapsto (x+B(x),\varphi),
	\end{equation*}
	where $B \in \mathcal{R}^{\infty,w}_{r_0-2\rho,s_0}$ is defined by $B = \Phi_{B}^* A_0$. Furthermore, the following bounds hold:
	\begin{equation}\label{bounAB}
	\|A_0\|^{\infty,w}_{r_0-\rho,s_0} \leq \|A\|^{\infty,w}_{r_0-\rho,s_0} \leq \frac{24}{11}\|P\|^{\infty,w}_{r_0,s_0}, \qquad \|B\|^{\infty,w}_{r_0-2\rho,s_0} \leq \frac{48}{7}\|P\|^{\infty,w}_{r_0,s_0}.
	\end{equation}
\begin{proof}
	We first verify that $\Phi_{-A_0}$ is a well-defined change of coordinates. Recalling that $A = \Phi^*_{-A} P$, Proposition \ref{composition} implies that $|A|^{\infty,w}_{r_0-\rho,s_0}$ is bounded by $2|P|^{\infty,w}_{r_0,s_0}$, which leads to the chain of inequalities
	\begin{equation*}
	|A_0|^{\infty,w}_{r_0-\rho,s_0} \leq |A|^{\infty,w}_{r_0-\rho,s_0} \leq 2|P|^{\infty,w}_{r_0,s_0}.
	\end{equation*}
	Under the smallness assumption $\|P\|^{\infty,w}_{r_0,s_0}<r_0/24$, it follows directly from the above inequalities that
	\begin{equation*}
	|A_0|^{\infty,w}_{r_0-\rho,s_0} < \frac{r_0-\rho}{6}.
	\end{equation*}
	Thus, the hypotheses of Theorem \ref{fix}, item \ref{fix-2}, are satisfied. Consequently, the map
	\begin{equation*}
	\Phi_{-A_0}:\mathscr{U}_{r_0-\rho}(\ell^\infty_w)\times \ell^\infty \to \ell^\infty_w\times \ell^\infty
	\end{equation*}
	is a change of coordinates with $\rho' = 3|A_0|^{\infty,w}_{r_0-\rho,s_0}$ and $\sigma = 0$\footnote{Note that there is no loss of regularity with respect to the variables $ \ \varphi_k \ $, $ \ k \in \N \ $, because the transformation does not change the angles.}. Its inverse is
	\begin{equation*}
	\Phi_{B}:\mathscr{U}_{r_0-\rho-\rho'}(\ell^\infty_w)\times \ell^\infty \to \ell^\infty_w\times \ell^\infty,
	\end{equation*}
	where $B \in \mathcal{A}^{\infty,w}_{r_0-\rho-\rho',s_0}$ is defined by $B = \Phi_{B}^* A_0$ and satisfies the bound $|B|_{r_0-\rho-\rho',s_0}^{\infty,w} \leq 2|A_0|^{\infty,w}_{r_0-\rho,s_0}$.
	
Finally, since $\rho' \leq \rho$, to complete the proof, it suffices to show that $\Phi_{-A_0}$ is also a regular change of coordinates with $  7\|P\|^{\infty,w}_{r_0,s_0} \leq \rho <r_0/2$. To this end, we first note that the smallness condition on $P$ and Lemma \ref{ao} directly imply the bounds in \eqref{bounAB}. In particular, the norms $ \ \|A_0\|^{\infty,w}_{r_0-\rho,s_0} \ $ and $ \ \|B\|^{\infty,w}_{r_0-2\rho,s_0} \ $ are smaller than $ \ \rho \ $. Then Remark \ref{regrmk} ensures that $\Phi_{-A_0}$ is indeed a regular change of coordinates with inverse $\Phi_B$. This completes the proof.
\end{proof}
\end{lemma}

It is worth highlighting that the estimates stated in above theorem imply that the proposed generating function $S$ defined by
\begin{equation}\label{SSS}
S(x,\varphi)=\sum_{k \in \N}x_k\varphi_k-W(x,\varphi), \qquad W(x,\varphi):=\sum_{m \in \N}\sum_{\nu \in \Lambda_m}\dfrac{A^{(m)}_\nu(\gamma(x))}{i \nu_m}e^{i \nu \cdot \varphi},
\end{equation}
is indeed a real majorant analytic function $S:\mathscr{U}_{r}(\ell^\infty_w)\times \ell^\infty \to \R$, $r\leq r_0-2\rho$: In fact, by the nature of the elements involved, the first sum is totally convergent, and $|W|_{r,s}$ can be bounded by $\|A\|^{\infty,w}_{r,s}$, up to a multiplicative constant, for any $r\leq r_0-2\rho$.

The next step concludes the construction of the action coordinates.

\begin{lemma}[Action coordinates]\label{Delta_1}
	Given $r_0, s_0 > 0$, let $P \in \mathcal{R}^{\infty,w}_{r_0,s_0}$ satisfy $\|P\|^{\infty,w}_{r_0,s_0} < r_0/27$. Then the map
	\begin{equation*}
	\Phi_{-A_0}\circ \Phi_P:\mathscr{U}_{r_0-2\rho}(\ell^\infty_w)\times \ell^\infty \to \ell^\infty_w\times \ell^\infty
	\end{equation*}
	is a regular change of coordinates for all $\rho$ such that $ 9\|P\|^{\infty,w}_{r_0,s_0}\leq\rho <r_0/3$ and $0 \leq \sigma <s_0/2$. Its inverse is given by
	\begin{equation*}
	\Phi_{-A}\circ\Phi_{B}:\mathscr{U}_{r_0-2\rho}(\ell^\infty_w)\times \ell^\infty \to \ell^\infty_w\times \ell^\infty.
	\end{equation*}
	Moreover, the set $(\Phi_{-A}\circ \Phi_B)(\mathscr{U}_{r_0-3\rho}(\ell^\infty_w) \times \ell^\infty)$ is an open neighbourhood of $\{0\}\times \ell^\infty$.
	
	\noindent Finally, assuming that 
	\begin{equation*}
			\{\Gamma^{(m)}_P,\Gamma^{(n)}_P\}=0 \qquad \text{for all} \ \ m,n \in \N \ ,
	\end{equation*}
	and denoting $(x,\varphi) = (\Phi_{-A_0}\circ \Phi_P) \ (y,\varphi)$, the following hold
	\begin{equation*}
	\{x_n,x_m\}=0 \qquad \text{for all } n,m \in \mathbb{N}.
	\end{equation*}
	\begin{proof}
First, note that the maps defined above are indeed the inverse of each other. In fact, by Lemma \ref{A_0,B}, we have:
\begin{equation*}
\begin{aligned}
\Phi_{-A_0}\circ\Phi_P \circ \Phi_{-A} \circ \Phi_B &\stackrel{\text{Thm}\ref{FFoliation}}{=} \Phi_{-A_0}\circ \id_{\mathscr{U}_{r_0-\rho}(\ell^{\infty}_w)\times \ell^\infty} \circ \Phi_{B} = \Phi_{-A_0} \circ \Phi_{B} 
\\
&\stackrel{\text{Lm}\ref{A_0,B}}{=} \id_{\mathscr{U}_{r_0-2\rho}(\ell^{\infty}_w)\times \ell^\infty}
\end{aligned}
\end{equation*}
\begin{equation*}
\begin{aligned}
\Phi_{-A}\circ\Phi_{B} \circ \Phi_{-A_0} \circ \Phi_P &\stackrel{\text{Lm}\ref{A_0,B}}{=} \Phi_{-A}\circ \id_{\mathscr{U}_{r_0-2\rho}(\ell^{\infty}_w)\times \ell^\infty} \circ \Phi_{P} 
= \Phi_{-A} \circ \Phi_{P} 
\\
&\stackrel{\text{Thm}\ref{FFoliation}}{=} \id_{\mathscr{U}_{r_0-3\rho}(\ell^{\infty}_w)\times \ell^\infty} \ .
\end{aligned}
\end{equation*}
We now verify that the composition
\begin{equation*}
\Phi_{-A_0}\circ \Phi_P:\mathscr{U}_{r_0-2\rho}(\ell^\infty_w)\times \ell^\infty \to \ell^\infty_w\times \ell^\infty
\end{equation*}
is indeed a regular change of coordinates. We write the explicit action of the maps involved:
\begin{equation*}
\begin{aligned}
&(\Phi_{-A_0} \circ \Phi_P)(y,\varphi)=(y+P(y,\varphi)-A_0(y+P(y,\varphi)),\varphi) \ , 
\\
&(\Phi_{-A} \circ \Phi_B)(x,\varphi)=(x+B(x)-A(x+B(x),\varphi),\varphi) \ ,
\end{aligned}
\end{equation*}
for all $(y,\varphi) \in \mathscr{U}_{r_0-2\rho}(\ell^\infty_w)\times \ell^\infty$. Therefore, by defining the maps $N$ and $M$ as
\begin{equation*}
N(y,\varphi)=P(y,\varphi)-A_0(y+P(y,\varphi)) \ , \qquad M(x,\varphi)=B(x)-A(x+B(x),\varphi) \ ,
\end{equation*}
we have that $\Phi_{-A_0} \circ \Phi_P = \id +(N,0)$ and $\Phi_{-A} \circ \Phi_B = \id +(M,0)$. By using the bounds in \eqref{bounAB} and estimating the composition of functions with respect to the $\|\cdot\|^{\infty,w}_{r,s}$-norm, it follows that both $\|N\|^{\infty,w}_{r_0-2\rho,s_0}$ and $\|M\|^{\infty,w}_{r_0-2\rho,s_0}$ are smaller than $\rho$. Consequently, the first claim follows from Remark~\ref{regrmk}. 

To prove the second claim, note that, given a change of coordinates $ \ \Phi \in \mathcal{C}^{\rho,0}_{r_0-2\rho,s_0}(\mathcal{A}) \ $ with inverse $ \ \Psi \in \mathcal{C}^{\rho,0}_{r_0-2\rho,s_0-\sigma}(\mathcal{A})$, one has that
\begin{equation}\label{nbdhof0T}
\{0\}\times \ell^\infty \subset \mathscr{U}_{r_0-3\rho}(\ell^\infty_w)\times \ell^\infty = \Psi\circ\Phi(\mathscr{U}_{r_0-3\rho}(\ell^\infty_w)\times \ell^\infty) \ \subset \Psi(\mathscr{U}_{r_0-2\rho}\times \ell^\infty) \ .
\end{equation}
The last inclusion follows because $ \ \Phi:\mathscr{U}_{r}(\ell^\infty_w)\times \ell^\infty \to \mathscr{U}_{r+\rho}(\ell^\infty_w)\times \ell^\infty \ $ for all $ \ r\leq r_0-2\rho \ $.

We conclude the proof by showing that the coordinate functions $(x_k)_{k \in \N}$ induced by the map $\Phi_{-A_0}\circ \Phi_P$ Poisson commute. To simplify the notation, in the following computations we assume that the derivatives of $x_k$ with respect to $\gamma_\ell$ (for $k, \ell \in \mathbb{N}$) are evaluated at $(\Gamma_P(y,\varphi),\varphi)$, even if not explicitly indicated. Then we have:
\begin{equation*}
\begin{aligned}
\{x_n,x_m\}&=\sum_{k \in \mathbb{N}}\frac{\partial x_n}{\partial y_k}\frac{\partial x_m}{\partial \varphi_k}-\frac{\partial x_n}{\partial \varphi_k}\frac{\partial x_m}{\partial y_k}
\\
&=\sum_{k \in \mathbb{N}}\sum_{\ell \in \mathbb{N}}\sum_{h \in \mathbb{N}}\frac{\partial x_n}{\partial \gamma_\ell}\frac{\partial x_m}{\partial \gamma_h}\bigg(\frac{\partial \Gamma^{(\ell)}_P}{\partial y_k}\frac{\partial \Gamma^{(h)}_P}{\partial \varphi_k}-\frac{\partial \Gamma^{(\ell)}_P}{\partial \varphi_k}\frac{\partial \Gamma^{(h)}_P}{\partial y_k}\bigg) \ .
\end{aligned}
\end{equation*}
Since, by assumption, the components of $\Gamma_P$ Poisson commute, the latter expression is formally equal to $0$. By performing Cauchy estimates and exploiting the strong regularity of the involved maps (expressed in terms of the $\|\cdot\|_{r,s}^{\infty,w}$-norms), one can prove that, for every $(y,\varphi)\in \mathscr{U}_{r_0-2\rho}(\ell^\infty_w)\times\ell^\infty$, there exists a neighbourhood in which the series is totally convergent. This concludes the proof.
\end{proof}
\end{lemma}

We now turn to the issue of proving that the angles defined formally by (see Equation \eqref{angleS} and the previous discussion)
\begin{equation}\label{th}
\theta_k(x,\varphi)=\varphi_k-\sum_{m\in \N}\sum_{\nu \in \Lambda_m}\frac{1}{\im \nu_m}\frac{\partial A^{(m)}_\nu (\gamma(x))}{\partial x_k}e^{\im \nu \cdot \varphi} \ , \quad k \in \N \ 
\end{equation}
define indeed a regular change of coordinates. Let us start by justifying the algebraic manipulations that we performed for proving the relations in \eqref{Acommutes}, which we recall below
\begin{equation}\label{hol}
\nu_\ell A^{(m)}_\nu = \nu_m A^{(\ell)}_\nu \quad \text{for all} \ \ell,m \in \mathbb{N} \ .
\end{equation}
Based on the construction of the candidate angle coordinates presented in Subsection \ref{heuristics}, to prove \eqref{Acommutes}, we must show that involved formal series are in fact totally convergent.

\begin{lemma}
	Given $ \ r,s >0 \ $, and $ \ 0 < \rho < r \ $, let $ \ P\in \mathcal{R}^{\infty,w}_{r,s} \ $ be such that $ \ |P|^{\infty,w}_{r,s}<\rho \ $ and consider $ \ A \in \mathcal{A}^{\infty,w}_{r-\rho,s} \ $. Then, for all $ \ (y,\varphi)\in \mathscr{U}_{r-2\rho}(\ell^\infty_w)\times \ell^\infty \ $, there exists a neighbourhood such that the formal series
	\begin{gather*}
		\sum_{\ell \in \mathbb{N}}\frac{\partial \Gamma_P^{(m)}}{\partial y_\ell}(y,\varphi)\frac{\partial A^{(\ell)}}{\partial \varphi_k}(\Gamma_P(y,\varphi),\varphi) \ , \qquad m \in \N \ ,
		\\
		\sum_{k \in \mathbb{N}}\sum_{\ell\in \mathbb{N}}\frac{\partial \Gamma_P^{(m)}}{\partial y_k}(y,\varphi)\frac{\partial \Gamma_P^{(n)}}{\partial y_\ell}(y,\varphi)\bigg(\frac{\partial A^{(\ell)}}{\partial \varphi_k}(\Gamma_P(y,\varphi),\varphi)-\frac{\partial A^{(k)}}{\partial \varphi_\ell}(\Gamma_P(y,\varphi),\varphi)\bigg) \ , \quad m,n \in \N \ ,
	\end{gather*}
	are totally convergent.
	\begin{proof}
Regarding the first series, the statement actually holds for all $P \in \mathscr{U}_{\rho}(\mathcal{A}_{r,s}^{\infty,w})$; that is, the regularity of $P$ is not required. In fact, total convergence in a neighbourhood of any point in the given domain can be shown as follows: by performing a Cauchy estimate on the derivative of $A$ and subsequently providing a bound for the composition (see Proposition \ref{composition}), the analyticity of $\Gamma_P$ allows one to complete the argument with a further Cauchy estimate. For the second series, the reasoning is similar, with the difference that one must exploit the regularity of $P$ to achieve the desired estimate.
	\end{proof}
\end{lemma}

Adopting the same notations as in Lemmas~\ref{A_0,B} and~\ref{Delta_1}, we now show that the function $Q: \mathscr{U}_{r_0-3\rho}(\ell^\infty_w) \times \ell^\infty \to \ell^\infty_w \times \ell^\infty$, defined for all $ \ (x,\varphi) \in \mathscr{U}_{r_0-2\rho}(\ell^\infty_w) \times \ell^\infty \ $ by
\begin{equation}\label{Omega}
Q^{(k)}(x,\varphi) = \sum_{m\in \mathbb{N}} \sum_{\nu \in \Lambda_m} \frac{1}{i \nu_m} \left. \frac{\partial A^{(m)}_\nu \circ \Phi_B}{\partial x_k} \right|_x e^{i \nu \cdot \varphi}, \quad k \in \mathbb{N},
\end{equation}
belongs to $\mathcal{R}^{\infty}_{r_0-3\rho,s_0}$. More precisely, we have the following result:

\begin{lemma}\label{QP}
	Given $ \ r_0, s_0 > 0 \ $, let $ \ P \in \mathcal{R}^{\infty,w}_{r_0,s_0} \ $ satisfy $ \ \|P\|^{\infty,w}_{r_0,s_0} < r_0/28 \ $. Assume that the components of $ \ \Gamma_P \ $ Poisson commute. Then the function $ \ Q \ $ defined in \eqref{Omega} belongs to $ \ \mathcal{R}^{\infty}_{r_0-3\rho,s_0} \ $, for all $ \ \rho \ $ such that $ \ 7\|P\|^{\infty,w}_{r_0,s_0} \leq \rho < r_0/4 \ $, and satisfies the bound
	\begin{equation}
		\|Q\|^{\infty}_{r_0-3\rho,s_0} \leq  c\|P\|^{\infty,w}_{r_0,s_0} \ ,
	\end{equation}
	for some constant $c>1$.
	\begin{proof}
First, note that all assumptions of Lemma \ref{A_0,B} are satisfied; thus, the function $B$ is well-defined and belongs to $\overline{\mathscr{U}}_{\rho}(\mathcal{R}^{\infty,w}_{r_0-2\rho,s_0})$. Consequently, the composition $A\circ \Phi_B \in \mathcal{R}^{\infty,w}_{r_0-2\rho,s_0}$ is also well-defined. The additional condition $\rho < r_0/4$ ensures sufficient analyticity width to perform Cauchy estimates that are uniform with respect to $\rho$. Without this condition, the constant $c$ would depend on $\rho$ such that $\lim_{\rho \to r_0/3} c(\rho) = \infty$. By the aforementioned lemma, the fact that the components of $\Gamma_P$ Poisson commute implies that the formulas in \eqref{hol} hold. Consequently, for any $k \in \mathbb{N}$, the derivative of $Q$ with respect to $\varphi_k$ can be expressed as
\begin{equation*}
	\sum_{m\in \mathbb{N}} \sum_{\nu \in \Lambda_m} \frac{\im \nu_k}{i \nu_m} \left. \frac{\partial A^{(m)}_\nu \circ \Phi_B}{\partial x_k} \right|_x e^{i \nu \cdot \varphi}= \sum_{\nu\neq 0} \left. \frac{\partial A^{(k)}_\nu \circ \Phi_B}{\partial x_k} \right|_x e^{i \nu \cdot \varphi} \ .
\end{equation*}
Therefore, the claim then follows by applying the properties of the $\|\cdot\|_{r,s}^{\infty,w}$-norm under composition, following the same idea as in the proof of Lemma~\ref{ao}.
	\end{proof}
\end{lemma}

We now discuss the invertibility of the map $ \ (x,\varphi) \mapsto (x,\theta) \ $, where the coordinates $ \ \theta=(\theta_k)_{k \in \N} \ $, according to \eqref{th} and using the notation \eqref{Omega}, are defined by
\begin{equation*}
\theta_k(x,\varphi)=\varphi_k-Q(x,\varphi) \ , \qquad \text{for all} \ \ k \in \N \ .
\end{equation*}
Then, parametrizing $ \ \varphi \ $ as $ \ \varphi=\theta+R(x,\theta) \ $, we want to solve
\begin{equation}\label{bi}
R(x,\theta)=Q(x,\theta+R(x,\theta)) \ .
\end{equation}
Let us introduce a notation:

\begin{notation}\label{phiggi2}
	Given $ \ G \in \mathcal{A}^{\infty}_{r,s} \ $, we denote by $ \ \Psi_G \ $ the map
	\begin{equation*}
	\Psi_G(x,\varphi)=(x,\varphi+G(x,\varphi)) \ , \qquad \forall \ x \in \mathscr{U}_r(\ell^\infty_w) \ , \ \varphi \in \ell^\infty \ .
	\end{equation*}
\end{notation}

\noindent Therefore, Equation \eqref{bi} can be written as $ \ R=\Psi_R^*Q \ $, where 
\begin{equation*}
\Psi^*_R Q (x,\theta) := \left( \Psi^*_R Q^{(m)} (x,\theta) \right)_{m \in \mathbb{N}}=( Q^{(m)}(x,\theta+R(x,\theta))_{m \in \N} \ .
\end{equation*}
As in the previous section, the strategy consists in formalizing our question as a fixed-point problem for the nonlinear operator $ \ G \mapsto \Psi^*_G Q \ $ acting on a closed ball of the Banach space $ \ \mathcal{A}^{\infty}_{r,s} \ $. In order to do this, we show that it is Fréchet differentiable and provide a suitable estimate for the differential. These technical results are contained in Proposition \ref{gat2} and Lemma \ref{contraction2}. Once the existence and uniqueness of the fixed point in the regularity class $\mathcal{A}^\infty$ have been established, we shall prove, by applying Lemma \ref{ao2} (which is analogous to Lemma \ref{ao}), that the fixed point inherits the stronger regularity properties of $Q$. Namely, we show that if $ \ Q \ $ is bounded with respect to the stronger norm $ \ \|\cdot\|^{\infty}_{r,s} \ $, the same property holds for the fixed point.

Our first step is to establish the differentiability of this operator. We shall omit the details of the proof because it follows by the same argument as in Proposition \ref{gat}.

\begin{prop}\label{gat2}
	Given $ \ r,s>0 \ $, let $ \ F \in \mathcal{A}^{\infty}_{r,s}$. Then, for any $\sigma \in (0,s)$, the map
	\begin{equation*}
	\mathscr{F}:\mathscr{U}_{\sigma/2}(\mathcal{A}^{\infty}_{r,s-\sigma}) \ni G \mapsto \Psi^*_GF \in \mathcal{A}^{\infty}_{r,s-\sigma}
	\end{equation*}
	is differentiable in the sense of Fréchet. Moreover, the following estimate holds
	\begin{equation}\label{boundgat2}
	|d\mathscr{F}(G)|_{\text{op}}\leq \frac{2}{e(\sigma-2|G|^{\infty}_{r,s-\sigma})}|F|^{\infty}_{r,s} \ .
	\end{equation}
\end{prop}

The estimate \eqref{boundgat2} suggests that the operator $ \ \mathscr{F} \ $ is a contraction provided that $ \ F \ $ is sufficiently small. This ensure the existence and uniqueness of the fixed point, as proved in the following lemma.

\begin{lemma}[Second Fixed Point Lemma]\label{contraction2}
	Given $r,s>0$, let $F \in \mathscr{U}_{s/8}(\mathcal{A}^{\infty}_{r,s})$ and define $\sigma:=8|F|^{\infty,w}_{r,s}$. Consider the map
	\begin{equation*}
	\mathscr{F}:\overline{\mathscr{U}}_{\sigma/4}(\mathcal{A}^{\infty}_{r,s-\sigma}) \ni G \mapsto \Psi^*_GF \in \mathcal{A}^{\infty}_{r,s-\sigma} 
	\end{equation*}
	with $\Psi_G$ as in \ref{phiggi2}. Then, there exists a unique $\hat{G} \in \overline{\mathscr{U}}_{\sigma/4}(\mathcal{A}^{\infty}_{r,s-\sigma})$ such that
	\begin{equation*}
	\hat{G}=\Psi^*_{\hat{G}}F.
	\end{equation*}
	
	\begin{proof}
		We follow the same strategy pursued for proving Theorem \ref{fix}. Namely, we first we show that
		\begin{equation*}
		\sup_{G \in  \overline{\mathscr{U}}_{\sigma/4}(\mathcal{A}^{\infty}_{r,s-\sigma})}|\Psi^*_GF|^{\infty}_{r,s-\sigma}\leq \sigma/4 \ ;
		\end{equation*}
		then, we prove that $\mathscr{F}$ acts as a contraction on $\overline{\mathscr{U}}_{\sigma/4}(\mathcal{A}^{\infty}_{r,s-\sigma})$ and the thesis will follow by the fact that a closed ball in a Banach space is indeed a complete metric space.
		
		The first point follows as an immediate consequence of Proposition \ref{composition} by setting, with the same notations of Proposition \ref{composition}, $ \ \Psi=\Psi_G \ $, i.e., $ \ N=0 \ $, $ \ \Omega\equiv G \ $.
		Indeed,
		\begin{equation*}
		|G|^{\infty}_{r,s-\sigma} \leq \sigma/4 \leq \sigma/2 \ .
		\end{equation*}
		Then, for all $G \in \overline{\mathscr{U}}_{\sigma/4}(\mathcal{A}^{\infty}_{r,s-\sigma}) $, the sequence $(G_n)_{n \in \N}$ defined by
		\begin{equation}\label{seq2}
		G_{n+1}:=\mathscr{F}(G_n), \qquad G_1:=G
		\end{equation}
		is contained in $\overline{\mathscr{U}}_{\sigma/4}(\mathcal{A}^{\infty}_{r,s-\sigma})$.
		
		Since $\mathscr{F}$ is Fréchet differentiable, if we show that its Fréchet differential $d\mathscr{F}$ satisfies
		\begin{equation*}
		\sup_{G \in \overline{\mathscr{U}}_{\sigma/4}(\mathcal{A}^{\infty}_{r,s-\sigma})}|d\mathscr{F}(G)|_{op}\leq K<1,
		\end{equation*}
		by Lagrange's inequality, one has
		\begin{equation*}
		|\mathscr{F}(G)-\mathscr{F}(G')|^{\infty}_{r,s-\sigma}\leq \sup_{G \in \overline{\mathscr{U}}_{\sigma/4}(\mathcal{A}^{\infty}_{r,s-\sigma})}|d\mathscr{F}(G)|_{op} |G-G'|^{\infty}_{r,s-\sigma}\leq K|G-G'|^{\infty}_{r,s-\sigma}.
		\end{equation*}
		Namely, $ \ \mathscr{F} \ $ would turn out to be a contraction with a Lipschitz constant equal to $ \ K<1 \ $.
		Note that, by \eqref{boundgat2} we have that
		\begin{equation*}
		\sup_{G \in \overline{\mathscr{U}}_{\sigma/4}(\mathcal{A}^{\infty}_{r,s-\sigma})}|d\mathscr{F}(G)|_{\text{op}} \leq \frac{4}{e\sigma}|F|^{\infty}_{r,s}=\frac{1}{2e}<1 \ .
		\end{equation*}
		
		Therefore, there exists $ \ \hat{G}\in \overline{\mathscr{U}}_{\sigma/4}(\mathcal{A}^{\infty}_{r,s-\sigma}) \ $ such that, for all $ \ G \in \overline{\mathscr{U}}_{\sigma/4}(\mathcal{A}^{\infty}_{r,s-\sigma}) \ $, the sequence $ \ (G_n)_{n \in \N} \ $ defined in \eqref{seq2} satisfies
		\begin{equation*}
		\lim_{n \to \infty}|G_n-\hat{G}|^{\infty}_{r,s-\sigma}=0.
		\end{equation*}
		Of course, $ \ \hat{G} \ $ satisfies $ \ \hat{G}=\Psi^*_{\hat{G}}F \ $.
	\end{proof}
\end{lemma}

We now focus on the regularity of the fixed point $ \ \hat{G} \ $. The next lemma shows that if $ \ F \ $ is bounded and sufficiently small with respect to the $\|\cdot\|^{\infty,w}_{r,s}$-norm, then $ \ \hat{G} \ $ is also bounded with respect to the same norm. Since this result is analogous to Lemma \ref{ao} and follows the same logic, we omit the proof.

\begin{lemma}\label{ao2}
	Given $r, s > 0$, let $0 \le \sigma < s$. Consider $F \in \mathscr{U}_{1/2}(\mathcal{R}^{\infty}_{r,s})$ and let $G \in \overline{\mathscr{U}}_{\sigma/2}(\mathcal{A}^{\infty}_{r,s-\sigma})$. Assume that $G = \Psi^*_G F$. Then, $G \in \mathcal{R}^{\infty}_{r,s-\sigma}$ and satisfies the estimate
	\begin{equation*}
	\|G\|^{\infty}_{r,s-\sigma} \leq \frac{2}{1-2\|F\|^{\infty}_{r,s}}\|F\|^{\infty}_{r,s} \ .
	\end{equation*}
\end{lemma}

With these technical results at hand, we now have all the necessary ingredients to prove the invertibility of the map $(x,\varphi) \mapsto (x,\theta)$.

\begin{lemma}\label{Q,R}
	Given $r_0, s_0 > 0$, let $P \in \mathcal{R}^{\infty,w}_{r_0,s_0}$ satisfy 
	\begin{equation*}
	\|P\|^{\infty,w}_{r_0,s_0} < \min\bigg\{\frac{r_0}{28} \ , \ \frac{s_0}{16c} \ , \ \frac{1}{4c}\bigg\}
	\end{equation*}
	where $c>1$ is the constant defined in Lemma \ref{QP}. Let us define $r_1:=r_0-27\|P\|^{\infty,w}_{r_0,s_0}$ and assume that the components of $\Gamma_P$ Poisson-commute. Then the map
	\begin{equation*}
	\Psi_{-Q}:\mathscr{U}_{r_1}(\ell^\infty_w)\times \ell^\infty \to \ell^\infty_w\times \ell^\infty, \quad (x,\varphi) \mapsto (x,\varphi-Q(x,\varphi)),
	\end{equation*}
	is a regular change of coordinates with $\rho = 0$ and $8c\|P\|^{\infty,w}_{r_0,s_0} \leq \sigma <s_0/2 $. Its inverse is given by
	\begin{equation*}
	\Psi_{R}:\mathscr{U}_{r_1}(\ell^\infty_w)\times \ell^\infty \to \ell^\infty_w\times \ell^\infty, \quad (x,\theta) \mapsto (x,\theta+R(x,\theta)),
	\end{equation*}
	where $R \in \mathcal{R}^{\infty}_{r_1,s_0-\sigma}$ is defined by $R = \Psi_{R}^* Q$. Furthermore, the following bounds hold:
	\begin{equation}\label{boundQR}
	\|Q\|^{\infty}_{r_1,s_0} \leq c\|P\|^{\infty,w}_{r_0,s_0}, \qquad \|R\|^{\infty}_{r_1,s_0-\sigma} \leq 4c\|P\|^{\infty,w}_{r_0,s_0}.
	\end{equation}
	\begin{proof}
		Note that all the assumptions of Lemma \ref{A_0,B} and Proposition \ref{QP} are satisfied. Thereferore, using the notation established previously, the maps $A$ and $B$ are well-defined together with $Q \in \mathcal{R}^\infty_{r_1,s_0}$, ensuring that $\Psi_{-Q}$ is also well-defined. 
		
		\noindent Moreover, the condition $\|P\|^{\infty,w}_{r_0,s_0} < s_0/(16c)$ together with the bound $\|Q\|_{r_1,s_0} \leq c\|P\|^{\infty,w}_{r_0,s_0}$ implies that $\|Q\|_{r_1,s_0} < s_0/16$. Therefore, we can apply the Second Fixed Point Lemma (see \ref{contraction2}, setting $F=Q$ and $\hat{G}=R$) to conclude that $\Psi_{-Q}$ is a change of coordinates with inverse $\Psi_{R}$, where $R \in \mathcal{A}^\infty_{r_1,s_0-\sigma}$ satisfies the bound $|R|^{\infty}_{r_1,s_0-\sigma}\leq 2|Q|^{\infty}_{r_1,s_0}$. Note that while $|Q|_{r_1,s_0} < s_0/8$ would suffice to apply Lemma \ref{contraction2}, the stronger requirement $\|P\|^{\infty,w}_{r_0,s_0} < s_0/(16c)$ ensures that $2\sigma < s_0$, so that the domain of the inverse is well-defined.
		
		The first bound in \eqref{boundQR} was established in Proposition \ref{QP}. The second bound follows from Lemma \ref{ao2}, which is applicable since
		\begin{equation*}
		|R|^{\infty}_{r_1,s_0-\sigma}\leq 2|Q|^{\infty}_{r_1,s_0} \leq 2\|Q\|^{\infty}_{r_1,s_0}\leq 2c\|P\|^{\infty,w}_{r_0,s_0}=\frac{2c\sigma}{8c}=\frac{\sigma}{4}\leq \frac{\sigma}{2} \ .
		\end{equation*}
		Indeed, the assumption $\|P\|^{\infty,w}_{r_0,s_0} \leq 1/(4c)$ combined with the first bound of \eqref{boundQR} implies $\|Q\|^\infty_{r_1,s_0} \leq 1/4$. Thus, Lemma \ref{ao2} yields
		\begin{equation*}
		\|R\|^\infty_{r_1,s_0-\sigma}\leq \frac{2}{1-2\|Q\|^{\infty}_{r_1,s_0}}\|Q\|^{\infty}_{r_1,s_0} \leq \frac{2}{1-1/2}\|Q\|^{\infty}_{r_1,s_0}=4\|Q\|^{\infty}_{r_1,s_0}\leq 4c\|P\|^{\infty,w}_{r_0,s_0} \ .
		\end{equation*}
		
		Since the right-hand side of these inequalities is equal to $\sigma/2$, both $\|R\|^\infty_{r_1,s_0-\sigma}$ and $\|Q\|^\infty_{r_1,s_0}$ are bounded by $\sigma/2$. This implies that $\Psi_{-Q} \in \mathcal{C}_{r_1,s_0}^{0,\sigma/2}(\mathcal{R})$ and $\Psi_{R} \in \mathcal{C}_{r_1,s_0-\sigma}^{0,\sigma/2}(\mathcal{R})$, and the claim follows from Remark \ref{regrmk}.
	\end{proof}
\end{lemma}

To make a connection with the discussion at the beginning of this section, we emphasize that the composition $\Psi_{-Q} \circ \Phi_{-A_0}$ provides the explicit construction of the map $\Delta$ (see \eqref{delta} and the related discussion). This transformation is precisely the regular change of coordinates required to straighten the foliation, ensuring that the resulting map $\Phi = \Delta \circ \Phi_P$ defines a symplectomorphism. By gathering the arguments developed so far, we are now in a position to prove the following theorem; from this result, the \emph{Symplectic Straightening Theorem}, which was informally stated in the introduction, will follow as a direct corollary.

\begin{prop}\label{tecmain}
Given $r_0, s_0 > 0$, let $P \in \mathcal{R}^{\infty,w}_{r_0,s_0}$ and assume that the components of $\Gamma_P$ pairwise Poisson-commute. Then, there exists a constant $C > 1$ such that, provided
\begin{equation*}
\|P\|^{\infty,w}_{r_0,s_0} < \frac{1}{4C} \min\{r_0, s_0, 1\},
\end{equation*}
the map
\begin{equation}\label{symplecto}
\Psi_{-Q} \circ \Phi_{-A_0} \circ \Phi_P : \mathscr{U}_{r_0-4\rho}(\ell^\infty_w) \times \ell^\infty \to \ell^\infty_w \times \ell^\infty
\end{equation}
is a symplectomorphism with inverse given by
\begin{equation*}
\Phi_{-A} \circ \Phi_{B} \circ \Psi_R : \mathscr{U}_{r_0-3\rho}(\ell^\infty_w) \times \ell^\infty \to \ell^\infty_w \times \ell^\infty.
\end{equation*}
for every $\rho$ and $\sigma$ satisfying $ \ C\|P\|^{\infty,w}_{r_0,s_0} \leq \rho < r_0/4 \ \ \text{and} \ \ 2C\|P\|^{\infty,w}_{r_0,s_0} \leq \sigma < s_0/2.$
\begin{proof}
Note that, by assuming $C \geq \max\{9, 4c\}$, all the hypotheses of Lemmas~\ref{A_0,B}, \ref{Delta_1}, and \ref{Q,R} are simultaneously satisfied. Consequently, both maps introduced in the statement are well-defined. To show they are the inverse of each other, we note that
	\begin{equation*}
	\begin{aligned}
	\Psi_{-Q} \circ \Phi_{-A_0} \circ \Phi_P \circ \Phi_{-A} \circ \Phi_{B} \circ \Psi_R &\stackrel{Thm \ref{FFoliation}}{=}\Psi_{-Q} \circ \Phi_{-A_0} \circ \id_{\mathscr{U}_{r_0-2\rho}(\ell^\infty_w)\times \ell^\infty}\circ \Phi_{B} \circ \Psi_R
	\\
	&=\Psi_{-Q} \circ \Phi_{-A_0} \circ \Phi_B\circ\Psi_R
	\\
	&\stackrel{Lm \ref{A_0,B}}{=}\Psi_{-Q} \circ \id_{\mathscr{U}_{r_0-3\rho}(\ell^\infty_w)\times \ell^\infty} \circ\Psi_R
	\\
	&=\Psi_{-Q} \circ \Psi_R\stackrel{Lm \ref{Q,R}}{=}\id_{\mathscr{U}_{r_0-3\rho}(\ell^\infty_w)\times \ell^\infty} \ ,
	\end{aligned}
	\end{equation*}
	and
	\begin{equation*}
	\begin{aligned}
	\Phi_{-A} \circ \Phi_{B} \circ \Psi_R \circ \Psi_{-Q} \circ \Phi_{-A_0} \circ \Phi_P &\stackrel{Lm \ref{Q,R}}{=}\Phi_{-A} \circ \Phi_{B} \circ \id_{\mathscr{U}_{r_0-3\rho}(\ell^\infty_w)\times \ell^\infty}\circ \Phi_{-A_0} \circ \Phi_P
	\\
	&=\Phi_{-A} \circ \Phi_{B} \circ \Phi_{-A_0}\circ\Phi_P
	\\
	&\stackrel{Lm \ref{A_0,B}}{=}\Phi_{-A} \circ \id_{\mathscr{U}_{r_0-3\rho}(\ell^\infty_w)\times \ell^\infty} \circ\Phi_P
	\\
	&=\Phi_{-A} \circ \Phi_P\stackrel{Thm \ref{FFoliation}}{=}\id_{\mathscr{U}_{r_0-4\rho}(\ell^\infty_w)\times \ell^\infty} \ .
	\end{aligned}
	\end{equation*}
	
	Next, recall the maps $N$ and $M$ introduced in the proof of Theorem~\ref{Delta_1}:
	\begin{equation*}
	N(y,\varphi)=P(y,\varphi)-A_0(y+P(y,\varphi)) \ , \qquad M(x,\varphi)=B(x)-A(x+B(x),\varphi) \ .
	\end{equation*}
	In terms of these maps, we can write
	\begin{equation*}
	\Phi_{-Q}\circ\Phi_{-A_0} \circ \Phi_P=\id+(N,Q\circ \Phi_{-A_0}\circ \Phi_P), \qquad \text{and} \qquad \Phi_{-A}\circ \Phi_B \circ \Psi_R=\id+(M\circ\Psi_R,R) \ .
	\end{equation*}
	To establish that $\Psi_{-Q} \circ \Phi_{-A_0} \circ \Phi_P$ is a regular change of coordinates, it remains to verify the following smallness conditions:
	\begin{equation*}
	\max\{\|N\|^{\infty,w}_{r_0-4\rho},\|M\circ\Psi_R\|^{\infty,w}_{r_0-4\rho}\} \leq \rho
	\end{equation*}
	\begin{equation*}
	\max\{\|Q\circ \Phi_{-A_0}\circ \Phi_P\|^{\infty}_{r_0-3\rho}, \|R\|^{\infty}_{r_0-3\rho}\}\leq \frac{\sigma}{2} \ .
	\end{equation*}
	If these hold, the thesis follows from Remark~\ref{regrmk}. Clearly, by using bounds established previously, and performing estimates for the norms of the compositions, one can prove that there exist constants $C_1, C_2 > 1$ such that
	\begin{equation*}
	\max\{\|N\|^{\infty,w}_{r_0-4\rho},\|M\circ\Psi_R\|^{\infty,w}_{r_0-4\rho}\} \leq C_1\|P\|^{\infty,w}_{r_0,s_0}
	\end{equation*}
	\begin{equation*}
	\max\{\|Q\circ \Phi_{-A_0}\circ \Phi_P\|^{\infty}_{r_0-3\rho}, \|R\|^{\infty}_{r_0-3\rho}\}\leq C_2\|P\|^{\infty,w}_{r_0,s_0} \ .
	\end{equation*}
	Thus, by defining
	\begin{equation}\label{Constant}
	C := \max\{9, 4c, C_1, C_2\},
	\end{equation}
	we deduce that $\Psi_{-Q} \circ \Phi_{-A_0} \circ \Phi_P$ is indeed a regular change of coordinates.
	
To conclude, we observe that the map preserves the Poisson brackets by construction. Indeed, this change of variables is induced by the generating function $S$ in \eqref{SSS}, which formally defines a symplectomorphism. The estimates provided above ensure that these algebraic manipulations are well-defined on the respective domains, thereby completing the proof.
\end{proof}
\end{prop}

	\begin{proof}[Proof of Theorem \ref{SYMPL}]
			The result follows directly from the construction developed in the previous paragraphs. By Proposition~\ref{tecmain}, the composition $\Delta:=\Psi_{-Q} \circ \Phi_{-A_0}$ provides the regular change of coordinates required to define the symplectomorphism $\Phi = \Delta \circ \Phi_P$. The straightening property is a consequence of the fact that $\Phi_P$ maps the level sets of $\Gamma_P$ into the standard horizontal leaves $\{\gamma = \text{const}\} \times \mathbb{T}^\infty$, while the transformation $\Delta$ preserves this foliation while restoring the symplectic structure.
	\end{proof}

 In the following subsections, we shall see how the existence of such a straightening transformation allows us to integrate the equations of motion for a class of infinite-dimensional systems, leading to a natural generalization of the Liouville-Arnold Theorem.

\subsection{Liouville-Arnold integrability} 
We now provide the dynamical counterpart to the geometric results established in Subsection \ref{SubSs}. In the classical finite-dimensional setting, the Liouville-Arnold Theorem ensures that a Hamiltonian system admitting a complete set of Poisson-commuting integrals is integrable by quadratures and can be described in terms of action-angle variables. Here, we prove that under suitable smallness conditions on the perturbation, an analogous conclusion holds for systems defined on $\ell^\infty_w \times \mathbb{T}^\infty$. Specifically, we show that if a Hamiltonian $H$ commutes with a set of integrals that define a Poisson-commuting foliation, the straightening symplectomorphism $\Phi$ constructed in Theorem \ref{SYMPL} acts as the linearizing transformation for the Hamilton's equations.

\begin{thm}[Complete Integrability]\label{completethm}
	Let $r,s>0$ and consider a majorant real analytic Hamiltonian 
	\begin{equation*}
		H:\mathscr{U}_r(\ell^\infty_w)\times\T^\infty \to \R
	\end{equation*}
	with $ H\in \mathcal{A}_{r,s}$. Assume that there exists $P \in \mathcal{R}^{\infty,w}_{r_0,s_0}$, with $r_0 \in (0,r]$ and $s_0\in(0,s]$, satisfying the following properties:
	\begin{itemize}
		\item[(a.1)] the following bound holds
		\begin{equation}\label{cazzo}
		\|P\|^{\infty,w}_{r_0,s_0} < \frac{1}{4C} \min\{r_0, s_0, 1\},
		\end{equation}
		where the constant $C$ is defined in \eqref{Constant};

		\item[(a.2)] the map $\Gamma_P$ defined by 
		\begin{equation*}
		\Gamma_P : \mathscr{U}_{r_0}(\ell^\infty_w) \times \mathbb{T}^\infty \to \ell^\infty_w \qquad \Gamma_P(y, \varphi) = y + P(y, \varphi)
		\end{equation*}
		is such that, for all $m,n \in \N$, $ \ \{\Gamma_P^{(m)},\Gamma^{(n)}_P\}=0 \ $ and $ \ \{\Gamma_P^{(m)},H\}=0 \ $ for all $m \in \N$.
		 
		\end{itemize}
Then,
	\begin{itemize}
	 
	\item[(i)] $H$ induces a Hamiltonian flow $\phi_H$ on a neighbourhood $ \ \mathscr{N}_0 \ $ of $\ \{0\}\times \T^\infty \ $;

	\item[(ii)] there is a neighbourhood $ \ \mathscr{U}_0 \ $ of $ \ \gamma=0 \ $ in  \ $\ell^\infty_w \ $ such that, for all $\gamma \in \mathscr{U}_0$, the level set $\Gamma_P^{-1}(\gamma)$ is an analytic Banach manifold, diffeomorphic to the infinite torus $\T^\infty$, which is invariant for the Hamiltonian flow $\phi_H$;

	\item[(iii)] there exists a set of \emph{action-angle} coordinates $(x,\theta)$ in which Hamilton's equations read
	
\begin{equation}\label{hameqnuove}
\begin{cases}
\dot{x}_k = 0 \\
\dot{\theta}_k = \dfrac{H\circ \Phi^{-1}}{\partial x_k}
\end{cases}
, \quad k \in \mathbb{N},
\end{equation}
where $\Phi$ is the symplectomorphism defined in \eqref{symplecto} and $(x,\theta)=\Phi(y,\varphi)$.
\end{itemize}

\noindent In particular the Hamilton's equations can be integrated by quadratures and the Hamiltonian flow reduces to a linear flow on the infinite torus $\T^\infty$.

	\begin{proof}
		To prove (i), we first address point (iii). In the new coordinates $(x, \theta) = \Phi(y, \varphi)$ provided by Theorem \ref{SYMPL}, the Hamiltonian takes the form $K := H \circ \Phi^{-1}$. The fact that $\{H, \Gamma_P^{(m)}\} = 0$ implies, via the symplecticity of $\Phi$, that $\{K, x_m\} = 0$ for all $m \in \mathbb{N}$. Consequently, $K$ is independent of the angular variables $\theta$. As discussed in the dedicated section, the Hamilton's equations for $K$ take the form \eqref{hameqnuove} and such a Hamiltonian $K$ induces a flow $\phi_K$ on $\mathscr{U}_{r_0-3\rho}(\ell^\infty_w) \times \mathbb{T}^\infty$, with $\rho=C\|P\|^{\infty,w}_{r_0,s_0}$, which is easily solvable and linear in the angles. Therefore, by applying Proposition \ref{propham}, the symplectomorphism $\Phi$ conjugates the flow of $K$ to the flow of $H$ on the neighbourhood of $\{0\}\times \T^\infty$ given by $\mathscr{N}_0:=\Phi^{-1}(\mathscr{U}_{r_0-3\rho}(\ell^\infty_w) \times \mathbb{T}^\infty)$.
		
		Point (ii) follows from the regularity of $\Phi$ and the straightening property: since $\Phi$ is a diffeomorphism and, for all $\gamma \in \mathscr{U}_0:= \Delta^{-1}(\mathscr{U}_{r_0-3\rho}(\ell^\infty_w) \times \mathbb{T}^\infty)$, maps $\Gamma_P^{-1}(\gamma)$ to the torus $\{x = \text{const}\} \times \mathbb{T}^\infty$, the level sets are indeed analytic Banach manifolds diffeomorphic to $\mathbb{T}^\infty$. The invariance of these manifolds under $\phi_H$ is a direct consequence of the fact that $\dot{x}_k = 0$ for all $k \in \N$.
	\end{proof}
\end{thm}

Despite the formal analogy with the classical Liouville-Arnold Theorem, some crucial distinctions must be emphasized. In the finite-dimensional case, integrability is often a global or semi-global geometric property derived from the algebraic structure of the symmetry groups involved. In this infinite-dimensional setting, however, our result relies heavily on the specific functional form of the momentum map $\Gamma_P$. The requirement that $P$ belongs to the class of majorant analytic functions $\mathcal{R}^{\infty,w}_{r_0,s_0}$ and satisfies the smallness condition (a.1) indicates that our approach is inherently perturbative, reflecting the technical challenges posed by the infinite-dimensional setting. Nonetheless, this version of the Liouville-Arnold Theorem provides a rigorous framework for studying several physically relevant models. As we shall see in the next sections, while the $\ell^2$ topology presents obstructions to the existence of foliations by infinite tori, our setting allows us to recover a form of complete integrability for systems defined on infinite lattices and for certain classes of Hamiltonian PDEs on the circle.

\section[Hamiltonian Systems on Infinite Lattices and PDEs on the Circle]{Hamiltonian Systems on Infinite Lattices and PDEs on the Circle}\label{latticiekdv}
A broad class of physical models characterized by an infinite number of degrees of freedom admits a Hamiltonian formulation. During the 1960s, it was discovered that certain nonlinear evolution equations possess an infinite number of conserved quantities (\cite{Gardner71}); consequently, their integrability has been rigorously investigated through diverse mathematical frameworks. Depending on the topological properties of the underlying manifolds, these approaches include algebraic geometry, inverse scattering transforms, and Riemann-Hilbert problems (\cite{Ablowitz81}). Such integrable structures arise both in continuous models—notably in the context of shallow water waves described by the \textit{Korteweg–de Vries} (KdV) equation—and in discrete lattices, such as the \textit{Toda lattice}, which models the dynamics of an infinite chain of interacting particles \cite{Toda89}. The case of Hamiltonian PDEs on the circle is particularly interesting because, in an appropriate set of coordinates, these infinite dimensional Hamiltonian systems look like infinitely many coupled harmonic oscillators, namely systems on an infinite lattice. In this section, we exploit this correspondence by embedding these systems into a weighted complex Banach space of sequences framework and provide conditions that allow to construct a set of action-angle coordinates.

\subsection{Systems on Infinite Lattices}
Building upon the results of the previous section, we provide here a class of Hamiltonian systems on infinite lattices that admit action-angle coordinates on an open subset of the phase space. Let us introduce the setting.

\subsubsection{The phase space}\label{Subps} Given a sequence of positive real numbers $(\eta_j)_{j \in \mathbb{N}}$, we consider as a phase space the complex weighted Banach space $\ell^\infty_\eta(\mathbb{C})$ defined by
\begin{equation*}
\ell^{\infty}_{\eta}(\mathbb{C}) := \left\{ u \in \mathbb{C}^\mathbb{N} : \|u\|_{\ell^{\infty}_\eta(\mathbb{C})} := \sup_{j \in \mathbb{N}} \eta_j |u_j| < \infty \right\}.
\end{equation*}
Note that $\ell^\infty_\eta(\mathbb{C})$ can be embedded into $\ell^2(\mathbb{C})$ if the weight $\eta$ satisfies the condition 
\begin{equation}\label{etaell2}
\tc_\eta:=\sum_{j \in \mathbb{N}} \eta_j^{-2} < \infty.
\end{equation}
From this embedding, one can inherit the standard symplectic form of $\ell^2(\mathbb{C})$, namely
\begin{equation*}
\omega := \im \sum_{j \in \mathbb{N}} du_j \wedge d\bar{u}_j,
\end{equation*}
by pulling it back onto $\ell^\infty_\eta(\mathbb{C})$.

\begin{rmk}\label{weak2}
	The embedding $\iota:\ell^\infty_\eta(\C)\hookrightarrow \ell^2(\C)$ is compact, and hence the pull back of the standard symplectic form of $\ell^2(\C)$ to $\ell^\infty_\eta(\C)$ is represented by a bounded operator $(\iota^*\om)^\#:\ell^\infty_\eta(\C)\to (\ell^\infty_\eta(\C))^*$, $\om^\#(V)=\om(\cdot,V)$, which is compact. Therefore $(\iota^*\om)^\#$ is injective but not surjective, making $\iota^*\om$ a weak symplectic form on $\ell^\infty_\eta(\C)$.
\end{rmk}

Let us begin to draw a connection with the phase space $\ell^\infty_w \times \mathbb{T}^\infty$ introduced in the previous sections. Let $\xi = (\xi_k)_{k \in \mathbb{N}}$ be a sequence of real numbers such that
\begin{equation}\label{xiinf}
0 < \inf_{k \in \mathbb{N}} \eta^2_k \xi_k \leq \sup_{k \in \mathbb{N}} \eta^2_k \xi_k < \infty.
\end{equation}
It then follows that the set
\begin{equation}\label{Txi}
\mathcal{T}_{\xi} := \left\{ u \in \ell^\infty_\eta(\mathbb{C}) : |u_k|^2 = \xi_k \text{ for all } k \in \mathbb{N} \right\}
\end{equation}
is analytically diffeomorphic to $\mathbb{T}^\infty$ via the map $\mathbb{T}^\infty \ni \varphi \mapsto (\sqrt{\xi_k}e^{\im \varphi_k})_{k \in \mathbb{N}}$. Indeed, by Theorem \ref{uniformanalytic}, this map is analytic as it is defined by a uniformly bounded family of analytic functions. The condition on the infimum in \eqref{xiinf} ensures both the invertibility of the map and the analyticity of its inverse, $(u_k \mapsto -\im \ln(u_k/\sqrt{\xi_k}))_{k \in \mathbb{N}}$, which again follows from Theorem \ref{uniformanalytic}. Consequently, $\mathcal{T}_\xi$ is an analytic Banach submanifold.\footnote{We refer the reader to \cite{PoTr} for a comprehensive treatment of analytic submanifolds in Banach spaces.}

Furthermore, the condition $0 < \inf_{k \in \mathbb{N}} \eta_k |u_k|$ defines an open subset of $\ell^\infty_\eta(\mathbb{C})$, which allows for the possibility of a foliation whose leaves are analytically diffeomorphic to $\mathbb{T}^\infty$. Specifically, any $\xi$ satisfying \eqref{xiinf} possesses an open neighborhood foliated by infinite tori. To see this, let $w$ be the weight defined by $w_k = \eta_k^2$ for all $k \in \mathbb{N}$ and note that, for any $r < \inf_{k \in \mathbb{N}} \eta^2_k \xi_k$, the map
\begin{equation}\label{phixi}
\Phi_\xi : \mathscr{U}_{r}(\ell^\infty_w) \times \mathbb{T}^\infty \ni (y, \varphi) \mapsto (\sqrt{\xi_k + y_k}e^{\im \varphi_k})_{k \in \mathbb{N}} \in \ell^\infty_\eta(\mathbb{C})
\end{equation}
is bi-analytic onto its image, providing the desired local coordinates for the foliation. This property is remarkable for the following reason: while $\ell^2(\mathbb{C})$ appears to be the most natural candidate for an infinite-dimensional Hamiltonian phase space, as it inherently provides a symplectic form via its Hermitian product, its topology does not allow for a foliation by infinite tori. Indeed, in $\ell^2(\mathbb{C})$, the set of finite-dimensional tori of the form \eqref{Txi} (with $\xi_k=0$ for all but finitely many $k$) is dense. Equivalently, as pointed out in \cite{Guan_2014}, the locus 
\begin{equation*}
\Game := \{u \in \ell^2(\mathbb{C}) : u_k = 0 \text{ for some } k \in \mathbb{N}\}
\end{equation*}
is dense in the strong topology. This density constitutes a topological obstruction that prevents Liouville-Arnold integrability, directing research in recent years towards the closely related concept of Birkhoff integrability for infinite-dimensional Hamiltonian systems.

Finally, provided that $\eta$ is such that $\tc_\eta$ is finite, the symplectic form $\sum_k dy_k \wedge d\varphi_k$ can be obtained by pulling back $\im\sum_{k}du_k\wedge d\overline{u}_k$ through $\Phi_\xi$, where $\xi$ satisfies \eqref{xiinf}, leading to the conclusion that $\Phi_\xi$ is a symplectomorphism.

\subsubsection{Space of Hamiltonians}\label{subSoH}
	For $ R>0$, let $\mathcal{H}_R$ be the space of Hamiltonians 
	$$
	H : \mathscr{U}_R(\ell^\infty_\eta(\C)) \to \R
	$$ 
	such that there exists a pointwise  absolutely convergent power series expansion\footnote{For the sake of brevity, the index set in summations will be omitted in the following.}
	\begin{equation*}
	H(u)=\sum_{a,b \in \N_0^{(\N)}}H_{ab}u^a\overline{u}^b
	\end{equation*}
	satisfying the following properties:

	\begin{itemize}
	\item[(i)] Reality condition:
		\begin{equation*}
		H_{a,b}= \overline{H_{b,a}}
		\end{equation*}
	\item[(ii)] Majorant analyticity:
	Given $H$ as above, we define its majorant as
	\begin{equation*}
		H_M(u)=\sum_{a,b}|H_{ab}|u^a\overline{u}^b \ .
	\end{equation*}
	Then we require that
	\begin{equation*}
	\lfloor H \rfloor_R:=	\sup_{u \in \mathscr{U}_R(\ell^\infty_\eta(\C))}|H_M(u)|<\infty \ .
	\end{equation*}
\end{itemize}
Note that, given $R > 0$, if we let $v_{\eta,R} \in \ell^\infty_\eta(\mathbb{C})$ be the vector defined by $(v_{\eta,R})_j = R \eta_j^{-1}$, then
\begin{equation*}
\floornorm{H}_R = \sum_{a,b} |H_{ab}| v_{\eta,R}^{a+b}.
\end{equation*}
It is worth stressing that the finiteness of the $\floornorm{\cdot}_R$-norm guarantees that the function is majorant-analytic in the ball $\mathscr{U}_R(\ell^\infty_\eta(\mathbb{C}))$, meaning that its majorant is itself an analytic function, as it is the uniform limit of polynomials (see, e.g., \cite{mujica}).

We now show that, given $\xi=(\xi_k)_{k \in \N}$, $\xi_k \in \R$ for all $k \in \N$, satisfying \eqref{xiinf}, the pullback of any majorant analytic Hamiltonian $H \in \mathcal{H}_R$ by $\Phi_\xi$ admits a Taylor-Fourier series which belongs to the class $\mathcal{A}$. More precisely, we fix a sequence of positive real numbers $\eta=(\eta_k)_{k \in \N}$ and we prove the following

\begin{prop}\label{Rrs}
Let $\xi_-,\xi_+ \in \R$ satisfy $0<\xi_-\leq \xi_+$. Then, for every sequence of real numbers $\xi=(\xi_k)_{k \in \N}$ such that
\begin{equation*}
	0<\xi_- \leq \inf_{k \in \N}\eta^2_k\xi_k \leq \sup_{k \in \N}\eta^2_k\xi_k \leq \xi_+ \ ,
\end{equation*}
and for all $R,r,s>0$ satisfying
\begin{equation}\label{ass}
\xi_+^{1/2}e^s \sum_{m=0}^{\infty}|a_m|\bigg(\frac{r}{\xi_-}\bigg)^m\leq R,
\end{equation}
where $(a_m)_{m=0}^\infty$ are the Taylor coefficients of the function $x\mapsto\sqrt{1+x}$, the pullback of any majorant analytic Hamiltonian $H\in\mathcal{H}_R$ by the map 
\begin{equation*}
	\Phi_\xi : \mathscr{U}_{r}(\ell^\infty_w) \times \mathbb{T}^\infty \ni (y, \varphi) \mapsto (\sqrt{\xi_k + y_k}e^{\im \varphi_k})_{k \in \mathbb{N}} \in \ell^\infty_\eta(\mathbb{C}), \qquad w_k=\eta_k^2 \ \ \text{for all} \ \ k \in \N, 
\end{equation*}
belongs to $\mathcal{A}_{r,s}$ and satisfies the bound
\begin{equation*}
|\Phi^*_{\xi}H|_{r,s}\leq	\floornorm{H}_R \ .
\end{equation*}
\begin{proof}
At a purely formal level, one has
\begin{equation*}
\begin{aligned}
	(H\circ\Phi_\xi)(y,\varphi)&=\sum_{a,b}H_{ab}\prod_{k \in \N}(\xi_k+y_k)^{(a_k+b_k)/2} \ e^{\im (a-b)\cdot\varphi}
	\\
	&=\sum_{\nu \in \Z}\sum_{a-b=\nu}H_{ab}\prod_{k \in \N}(\xi_k+y_k)^{(a_k+b_k)/2} \ e^{\im \nu\cdot\varphi}
	\\
    &=\sum_{\nu \in \Z}\sum_{a-b=\nu}H_{ab}\prod_{k \in \N}\bigg(\xi_k^{1/2} \ \sqrt{1+\frac{y_k}{\xi_k}} \ \bigg)^{(a_k+b_k)}e^{\im \nu\cdot\varphi}
    	\\
    &=\sum_{\nu \in \Z}\sum_{a-b=\nu}H_{ab}\prod_{k \in \N}\bigg(\xi_k^{1/2} \ \sum_{m=0}^{\infty}a_m\bigg(\frac{y_k}{\xi_k}\bigg)^m \bigg)^{(a_k+b_k)}e^{\im \nu\cdot\varphi}.
	\end{aligned}
\end{equation*}
The latter series is formally equal to
\begin{equation*}
	\sum_{\nu,\alpha}K_{\nu\alpha}y^\alpha e^{\im \nu \cdot \varphi}
\end{equation*}
for suitable coefficients $K_{\nu\alpha}$, with $\alpha \in \N_0^{(\N)}$ and $\nu \in \Z^{(\N)}$. Clearly, by the triangular inequality, and by using the assumptions on $R,r$ and $s$, it follows that
\begin{equation*}
\begin{aligned}
	\sum_{\nu,\alpha}|K_{\nu\alpha}|v_{w,r}^\alpha e^{s|\nu|}&\leq \sum_{\nu \in \Z}\sum_{a-b=\nu}|H_{ab}|\prod_{k \in \N}\bigg(\xi_k^{1/2} \ \sum_{m=0}^{\infty}|a_m|\bigg(\frac{rw_k^{-1}}{\xi_k}\bigg)^m \bigg)^{(a_k+b_k)}e^{s|\nu|}
	\\
	&\leq \sum_{a,b}|H_{ab}|\prod_{k \in \N}\bigg(\xi_+^{1/2}\eta_k^{-1}e^s \ \sum_{m=0}^{\infty}|a_m|\bigg(\frac{r}{\xi_-}\bigg)^m \bigg)^{(a_k+b_k)}
	\\
	&\leq \sum_{a,b}|H_{ab}|\prod_{k \in \N}v_{\eta,R}^{a+b}=\floornorm{H}_R.
	\end{aligned}
\end{equation*}
Therefore, $H\circ\Phi_\xi$ admits a power-trigonometric series representation, namely $(H\circ\Phi_\xi)_{\nu\alpha}=K_{\nu\alpha}$, and belongs to $\mathcal{A}_{r,s}$ with $|\Phi^*_\xi H|_{r,s}\leq\floornorm{H}_R$.
\end{proof}
\end{prop}

\subsubsection{Hamiltonian vector fields}\label{SubHvf} In this paragraph we assume that the constant $\tc_\eta$ is finite. A vector field $V: \mathscr{U}_R(\ell^\infty_\eta(\mathbb{C})) \to \ell^\infty_\eta(\mathbb{C})$ is Hamiltonian if there exists $H \in \mathcal{H}_R$ such that
\begin{equation}\label{HAMVF}
\omega(\cdot, V) = dH[\cdot].
\end{equation}
In this case, we say that the vector field is associated with the Hamiltonian $H$ and denote it by $V_H$. It should be noted that, in general, not all Hamiltonians in $\mathcal{H}_R$ induce a Hamiltonian vector field. Indeed, the majorant-analyticity property only ensures that the differential belongs to the topological dual of the phase space. In the following, in line with the , we introduce a class of analytic Hamiltonians that do induce Hamiltonian vector fields.

\subsubsection{Space of regular Hamiltonians}\label{SorH} We say that a Hamiltonian $f \in \mathcal{H}_R$ is regular if the norm
\begin{equation*}
\floorNorm{f}_R := \max \left\{ \lfloor f \rfloor_R, \, R \sup_{k \in \mathbb{N}} \eta_k \Floornorm{\frac{\partial f}{\partial \overline{u}_k}}_R \right\}
\end{equation*}
is finite. Note that, if $\tc_\eta < \infty$, the finiteness of the $\floorNorm{\cdot}_R$-norm ensures that the function defines a Hamiltonian vector field in the sense of \eqref{HAMVF}, as the gradient of any such regular Hamiltonian is necessarily an element of $\ell^2(\mathbb{C})$. Explicitly, the Hamiltonian vector field associated with a regular Hamiltonian $f$ is given by
\begin{equation}\label{hamvf2}
V^{(k)}_f = \im \frac{\partial f}{\partial \overline{u}_k}, \quad k \in \mathbb{N}.
\end{equation}
It is worth noting that, even when $\tc_\eta = \infty$, the regularity of $f$ ensures that the associated Hamiltonian vector field is analytic as a map $V_f:\mathscr{U}_R(\ell^\infty_\eta(\mathbb{C})) \to \ell^\infty_\eta(\mathbb{C})$. This is in contrast to the case where $f$ is merely analytic; in the latter situation, the Hamiltonian vector field, defined as the sequence of derivatives as in \eqref{hamvf2}, is only guaranteed to be bounded when viewed as a map taking values into the dual space $(\ell^\infty_\eta(\C))^*$. 

With the aim of exploiting the analytical background developed in the previous section, we now prove that given any regular Hamiltonian $f$, the composition $f \circ \Phi_\xi$ defines indeed a regular function in the sense of the previous section.  Let us first recall that given $r,s>0$ the $|\cdot|_{r,s}$-norm of a regular function $f\in\mathcal{R}_{r,s}$ is defined as
\begin{equation*}
\|f\|_{r,s}:=\max\bigg\{|f|_{r,s}, \ \sup_{k\in \N}w_k\bigg|\frac{\partial f}{\partial \varphi_k}\bigg|_{r,s}, \ r\sup_{k\in \N}\bigg|\frac{\partial f}{\partial y_k}
\bigg|_{r,s}\bigg\}<\infty.
\end{equation*}

\begin{prop}\label{Rrs2}
	With the same notations of Proposition \ref{Rrs}, and under the same assumptions, we have that the pullback of any regular Hamiltonian $ \ f:\mathscr{U}_R(\ell^\infty_\eta)\to\R \ $ by the map $\Phi_\xi$ belongs to $\mathcal{R}_{r,s}$ and satisfies the bound
	\begin{equation*}
		\|\Phi^*_{\xi}f\|_{r,s}\leq 2\floorNorm{f}_R.
	\end{equation*}
	\begin{proof}
		It follows by Proposition \ref{Rrs} that $|\Phi^*_{\xi}f|_{r,s}\leq	\floornorm{f}_R\leq \floorNorm{f}_R$. Furthermore, by using the relations
		\begin{equation*}
			\frac{\partial}{\partial y_k}=\frac{1}{2}\bigg(\frac{1}{\overline{u}_k}\frac{\partial}{\partial u_k}+\frac{1}{u_k}\frac{\partial}{\partial \overline{u}_k}\bigg), \qquad 
			\frac{\partial}{\partial \varphi_k}=\im\bigg(u_k\frac{\partial}{\partial u_k}+\overline{u}_k\frac{\partial}{\partial \overline{u}_k}\bigg), \qquad \text{for all} \ \ k \in \N,
		\end{equation*}
		and the reality condition $f_{ab}=\overline{f_{ab}}$, one can follows the same argument as in the proof of Proposition \ref{Rrs} to prove the following bounds
		\begin{gather*}
			r\sup_{k\in \N}\bigg|\frac{\partial \Phi^*_\xi f}{\partial y_k}
			\bigg|_{r,s} \leq \frac{r}{R^2} R \sup_{k \in \mathbb{N}} \eta_k \Floornorm{\frac{\partial f}{\partial \overline{u}_k}}_R \leq \frac{r}{R^2}\floorNorm{f}_R
			\\
			\sup_{k\in \N}w_k\bigg|\frac{\partial \Phi^*_\xi f}{\partial \varphi_k}
			\bigg|_{r,s} \leq 2R\sup_{k \in \mathbb{N}} \eta_k \Floornorm{\frac{\partial f}{\partial \overline{u}_k}}_R \leq 2\floorNorm{f}_R.
		\end{gather*}
		Since the assumption \eqref{ass} implies that $r/R^2<1$, the claim follows.
	\end{proof}
\end{prop}

In particular, the Hamiltonian vector field associated with a regular Hamiltonian $H\in\mathcal{H}_R$ is pushed forward to $\ell^\infty_w\times\T^\infty$ by $\Phi^{-1}_\xi$, and the resulting vector field is analytic and bounded. More precisely, it belongs to the space of vector field $\mathfrak{X}_{r,s}$ (see Definition \ref{analvect}).

\subsubsection{Poisson bracket}\label{subpb2} As usual, the symplectic form induces a Poisson bracket on the space of \emph{regular} functions by setting
\begin{equation*}
\{f,g\}:=\om(V_f,V_g) \ .
\end{equation*}
By a direct calculation, one discovers that
\begin{equation}\label{pbuubar}
	\{f,g\}=\im\sum_{k \in \N}\frac{\partial f}{\partial \overline{u}_k}\frac{\partial g}{\partial u_k}-\frac{\partial f}{\partial u_k}\frac{\partial g}{\partial \overline{u}_k} \ .
\end{equation}
According to the identification between the symplectic structures $\im \sum_{k}du_k\wedge d\overline{u}_k$ and $\sum_k dy_k \wedge d\varphi_k$ just established, we have that
\begin{equation}\label{pbuubaryphi}
\Phi^*_\xi\bigg(\im\sum_{k \in \N}\frac{\partial f}{\partial \overline{u}_k}\frac{\partial g}{\partial u_k}-\frac{\partial f}{\partial u_k}\frac{\partial g}{\partial \overline{u}_k}\bigg)=\sum_{k \in \N}\frac{\partial\Phi^*_\xi f}{\partial y_k}\frac{\partial \Phi^*_\xi g}{\partial \varphi_k}-\frac{\partial \Phi^*_\xi f}{\partial \varphi_k}\frac{\partial\Phi^*_\xi g}{\partial y_k} \ .
\end{equation}
Therefore, in line with the observations in Section \ref{SecRHam}, one can take \eqref{pbuubar} as a definition on the Poisson bracket, overcoming then the necessity of the requirement $c_\eta<\infty$. From now on, we shall dispense with the embedding into $\ell^2(\mathbb{C})$, as our results do not depend on such an assumption. Another important property that we inherit from the analysis developed in Section \ref{SecRHam} is the possibility of extending the operation induced by the Poisson bracket between pairs of functions in which one of the two is regular and the other one belongs to $\mathcal{H}_R$ for some $R>0$. This can be easily seen by interpreting this operation as the Lie derivative of an analytic function in the direction of the Hamiltonian vector field associated with a regular function.

\subsubsection{Hamilton's equations}\label{SubHameq2} Given $ \ R>0 \ $ and a function $ \ H \in \mathcal{H}_{R} \ $, we say that $ \ H \ $ induces a \emph{Hamiltonian flow} $ \ \phi_H \ $ on an open set $\mathscr{U}\subseteq \ell^\infty_\eta(\C)$ if there exists a group action $ \ \phi_H \ $ of $ \ \mathbb{R} \ $ on $ \ \mathscr{U} \ $
\begin{equation*}
\mathbb{R}\times\mathscr{U}\ni (t,u) \mapsto \phi^t_H(u)\in \mathscr{U}
\end{equation*}
such that, denoting $\phi^t_H(u)=(u_1(t),\ldots,u_k(t),\ldots)$, the maps $t \mapsto u_k(t)$ are differentiable functions for all $k \in \mathbb{N}$ and satisfy the \emph{Hamilton's equations}
\begin{equation*}
\dot{u}_k=\{H,u_k\}, \qquad k \in \N.
\end{equation*}
According to the notation previously introduced, we denote by $V_H:=(\{H,u_k\})_{k \in\N}$ the (typically unbounded) Hamiltonian vector field associated with $H$. Clearly,
\begin{equation*}
V_H^{(k)}=\im \frac{\partial H}{\partial \overline{u}_k}.
\end{equation*}
The same discussion and remarks of Section \ref{HAMEQSEC} on the regularity and on the existence of the Hamiltonian flow hold literally in this setting as well.

\subsubsection{The main theorem: Liouville-Arnold integrability}\label{subAL} 
Here we prove Theorem \ref{COMPL} on complete integrability (à la Liouville-Arnold) of Hamiltonian systems on infinite lattices. The idea is the following. In the phase space $\ell^\infty_\eta(\C)$, any torus of the form $\mathcal{T}_\xi$ (see \eqref{Txi}) has a tubular neighbourhood that can be analytically identified with the product space $\mathscr{U}_r(\ell^\infty_w)\times \T^\infty$ via the symplectomorphism $\Phi_\xi$ defined in \eqref{phixi}. Hence, by Propositions \ref{Rrs} and \ref{Rrs2}, $\Phi^{-1}_\xi$ pushes forward the Hamiltonian vector field $\im\partial_{\overline{u}}H$, restricted to this neighbourhood, to a vector field on $\mathscr{U}_r(\ell^\infty_w)\times\T^\infty$. Hamilton's equations retain their standard form in these coordinates.

	\begin{thm}[Complete Integrability for Systems on Infinite Lattices]\label{COMPL}
		Let $R>0$ and consider a majorant real analytic Hamiltonian 
		\begin{equation*}
		H:\mathscr{U}_R(\ell^\infty_\eta(\C))\to \R
		\end{equation*}
		with $ H\in \mathcal{H}_{R}$. Suppose that there are functions $(f^{(m)})_{m \in \N}$, where $f^{(m)}\in\mathcal{H}_R \ $ for all $m \in \N$, satisfying the following properties:
		\begin{itemize}
			 
			\item[(a.1)] the Taylor series of every $f^{(m)}$ starts at order larger or equal than three.

			\item[(a.2)] denoting $f=(f^{(m)})_{m \in \N}$, the following quantity is bounded
			\begin{equation*}
				\norm f \norm_R:=\max\bigg\{\sup_{m \in \N}\eta^2_m\floornorm{f^{(m)}}_R, \ R\sup_{k \in \N}\eta_k\sum_{m \in \N} \Floornorm{\frac{\partial f^{(m)}}{\partial \overline{u}_k}}_R\bigg\};
			\end{equation*}

			\item[(a.3)] let $w_k=\eta^2_k \ $ for all $ k \in \N$, then the map $\Gamma : \mathscr{U}_{R}(\ell^\infty_\eta(\C)) \to \ell^\infty_w$  defined by 
			\begin{equation*}
			\Gamma^{(m)}(u) = |u_m|^2 + f^{(m)}(u)
			\end{equation*}
			is such that, for all $m,n \in \N$, $ \ \{\Gamma^{(m)},\Gamma^{(n)}\}=0 \ $ and $ \ \{\Gamma^{(m)},H\}=0 \ $ for all $m \in \N$.
			 
		\end{itemize}
		Then, there exists $\xi_-,\xi_+\in\R$, with $0<\xi_-\leq\xi_+<R^2$, such that for every sequence of real numbers $\xi = (\xi_k)_{k \in \mathbb{N}}$ satisfying
		\begin{equation}\label{xiinftheo}
		\xi_-\leq \inf_{k \in \mathbb{N}} \eta^2_k \xi_k \leq \sup_{k \in \mathbb{N}} \eta^2_k \xi_k \leq \xi_+
		\end{equation}
		\begin{itemize}
			 
			\item[(i)] $H$ induces a Hamiltonian flow $\phi_H$ on a neighbourhood $ \ \mathscr{M}_\xi \ $ of $\ \mathcal{T}_\xi \ $;

			\item[(ii)] there is a neighbourhood $ \ \mathscr{U}_\xi \ $ of $ \ \gamma=\xi \ $ in  \ $\ell^\infty_w \ $ such that, for all $\gamma \in \mathscr{U}_\xi$, the level sets $\Gamma^{-1}(\gamma)$ are invariant analytic Banach submanifold diffeomorphic to the infinite torus $\T^\infty$, and form an analytic foliation of $\mathscr M_\xi$.

			\item[(iii)] there exists a set of \emph{action-angle} coordinates $(x,\theta)$ in which Hamilton's equations read
			
			\begin{equation}\label{hameq3}
			\begin{cases}
			\dot{x}_k = 0 \\
			\dot{\theta}_k = \dfrac{H \circ \Phi_\xi \circ \Phi^{-1}}{\partial x_k}
			\end{cases}
			, \quad k \in \mathbb{N},
			\end{equation}
			where $\Phi$ and $\Phi_\xi$ are the symplectomorphisms defined in \eqref{symplecto} and \ref{phixi}, respectively, and $(x,\theta)=\Phi \circ \Phi^{-1}_\xi(u)$.
		\end{itemize}

\noindent In particular the Hamilton's equations can be integrated by quadratures and the Hamiltonian flow reduces to a linear flow on the infinite torus $\T^\infty$.
\begin{proof}
	 Note that the requirement that the perturbations $f^{(m)}$, for all $m \in \N$, are of order three or higher implies that the quantity $R^{-2}\norm f\norm_R$ can be made arbitrarily small by simply reducing the radius $R$. Then, there exists $R_0>0$ such that
	 \begin{gather}
	 R_0\leq \min\{R,2^{5/2}\sqrt{\ln(2)}\} \label{laprima}
	 \\
	 R_0^{-2}\norm f \norm_{R_0}\leq \frac{1}{2^9C} \ , \label{laseconda}
	 \end{gather}
	 where the constant $C$ is defined in \eqref{Constant}. Let $\xi_-:=32C\norm f \norm_{R_0}$ and $\xi_+:=R^2_0/16$. Clearly, \eqref{laseconda} implies that $\xi_-\leq\xi_+$. Finally, define $r_0=\xi_-/2$ and let $s_0:=\ln(2)$. It is immediate to verify that $R_0,r_0,s_0$ satisfies \eqref{ass}. Therefore, by setting $P=\Phi^*_\xi f$ and $\Gamma_P=\Phi^*_\xi \Gamma$, it follows directly from Proposition \ref{Rrs2} that
	 \begin{equation*}
	 \|P\|^{\infty,w}_{r_0,s_0} \leq 2\norm f\norm_{R_0}.
	 \end{equation*}
	 This bound allows the smallness assumption of Theorem \ref{SYMPL} to be met. Indeed, one has the following chains on inequalities
	 \begin{equation*}
	 \begin{aligned}
	 	&\|P\|^{\infty,w}_{r_0,s_0} \leq 2\norm f \norm_{R_0}=\frac{\xi_-}{16C}=\frac{r_0}{8C}<\frac{r_0}{4C}
	 	\\
	 	&\|P\|^{\infty,w}_{r_0,s_0} \leq 2\norm f \norm_{R_0} \stackrel{\eqref{laseconda}}{\leq} \frac{R_0^2}{2^8C}\stackrel{\eqref{laprima}}{\leq}\frac{\ln(2)}{8C}<\frac{\ln(2)}{4C}<\frac{1}{4C},
	 	\end{aligned}
	 \end{equation*}
	 which imply
\begin{equation*}
\|P\|^{\infty,w}_{r_0,s_0} < \frac{1}{4C} \min\{r_0, s_0, 1\}.
\end{equation*}
Thus, since $\Phi_\xi$ is a symplectomorphism, the claim follows by Theorem \ref{completethm} and Proposition \ref{propham}  by setting, with the same notations of Theorem \ref{completethm}, $\mathscr{M}_\xi:=\Phi_\xi(\mathscr{N}_0)$ and $\mathscr{U}_\xi=\xi+\mathscr{U}_0$.
\end{proof}
\end{thm}

It is worth stressing that, besides providing a connection (established in Proposition \ref{Rrs2}) with the theory developed in the previous section, the boundedness of the $\norm \cdot\norm_R$ norm ensures that the momentum map is a majorant analytic map taking values in the space of actions, and that the associated pullback map preseves the Poisson algebra of regular Hamiltonians. This property, together with result established in the above theorem, should be compared with \cite{BaSt}, where he authors show how a family of commuting analytic vector fields can be reduced to a normal form near the origin via an analytic coordinate transformation. Their approach provides an abstract framework to construct Birkhoff coordinates for infinite-dimensional systems in a small neighbourhood around the origin, with direct applications to the KdV, the NLS equation, and the Toda lattice. Similarly to their approach, our result is local in nature; however, it applies to an open region of the phase space that excludes the origin, whose existence is fundamentally linked to the $\ell^\infty$-based topology. On the other hand, our construction provides a more detailed geometric description in terms of foliations into infinite-dimensional invariant submanifolds.

Another remarkable point is that, by interpreting the functions $u \mapsto |u_m|^2$ for all $m \in \mathbb{N}$ as Hamiltonians of harmonic oscillators, the first integrals $f^{(m)}$ can be seen as small (with respect to the $\norm \cdot\norm_R$ norm) perturbations of these systems. Finally, we point out that one has the explicit expression for the actions
\begin{equation}
x_m=\int_{\Gamma^{-1}(\gamma)}(|u_m|^2-\xi_m) d\lambda', \qquad m \in \N,
\end{equation}
where $\lambda'$ is the pushforward by $\Phi_\xi$ of the measure $\lambda$ defined in Subsection \ref{heuristics}. This directly follows by the definition of the actions (see \eqref{actionintegral}) and by the fact that $|u_m|^2=\xi_m+y_m$ for all $m \in \N$. Since, by \eqref{hameq3}, these actions are preserved by the Hamiltonian flow, we find that, although the energies of the harmonic oscillators $u \mapsto |u_m|^2$ are not preserved, the averages of such energies on the invariant tori are constant in time.

\section{The KdV Equation}\label{kdvsec}
In this section we shall apply the results established thus far for showing that the KdV equation on the circle is a completely integrable infinite-dimensional Hamiltonian system in the sense of Liouville-Arnold (see Theorem \ref{LA_KdV}). By the nature of our results, complete integrability is to be intended in a region of small data bounded away from the origin, structured as a product of thin complex annuli of the form $0 < r < |u_m| < R$ for each Fourier component. 
First we place the KdV equation within our $\ell^\infty$-weighted space setting. We show that, when viewed as a system on an infinite lattice via the Fourier transform, it defines a completely integrable infinite-dimensional Hamiltonian system in the sense of Theorem \ref{completethm}. To this end, it is necessary to investigate the structure and regularity of the spectral gaps of the associated periodic Schr\"odinger operator. We achieve this by interpreting the spectral problem as a first-order dynamical system, describing a perturbed harmonic oscillator where the potential $u$ acts as an external periodic force. We then apply a reduction scheme to transform this system into one with constant coefficients. Such a reduction significantly simplifies the analysis of the spectral gaps, which, in this setting, correspond precisely to the instability intervals of the associated dynamical system.

For the sake of notational simplicity and to avoid the proliferation of factors of $\pi$, we prefer to work with $2\pi$-periodic functions rather than $1$-periodic ones. This choice allows us to identify the Fourier lattice directly with $\mathbb{Z}$ instead of $2\pi\mathbb{Z}$. To this end, we perform a rescaling of the parameters by setting $t \to (2\pi)^3 t$, $u \to (2\pi)^{-2} u$, and $x \to 2\pi x$. One can easily verify that this transformation leaves the form of the KdV equation, including all its coefficients, invariant. Therefore, in what follows, we shall adopt this scaling without changing the notation for the variables.

Our goal is to prove that the square of the spectral gaps satisfy the hypotheses of Theorem~\ref{completethm}. We deduce the structure of the gaps via an implicit function theorem argument.

\subsection{Spectral problem for the periodic Schr\"odinger operator}\label{subfloquet} Consider the Schr\"odinger operator $L_u = -\partial_{xx} + u$, where $u$ is $2\pi$-periodic, with zero mean and analytic. The eigenvalue equation $$-\psi_{xx} + u\psi = E\psi$$ can be reformulated as a first-order dynamical system. Setting $t=x$, $q = \psi$ and $p = \psi_t$, we obtain:
	\begin{equation}\label{dyna}
	\frac{d}{dt} \begin{pmatrix} q \\ p \end{pmatrix} = \left[ \begin{pmatrix} 0 & 1 \\ -E & 0 \end{pmatrix} + \begin{pmatrix} 0 & 0 \\ u & 0 \end{pmatrix} \right] \begin{pmatrix} q \\ p \end{pmatrix},
	\end{equation}
	namely the Hamilton's equation for a harmonic oscillator forced by a periodic perturbation. From Floquet theory, it is well known that the spectrum is composed of bands and gaps, that correspond to bounded and unbounded solution for the associated dynamical system \eqref{dyna}. For the unperturbed case $u=0$, the periodic and anti-periodic eigenvalues are doubly degenerate and equal to $E_n = n^2/4$, with associated eigenfunctions $e^{\pm \im nt/2}$, for $n \in \mathbb{N}$. Clearly, for even $n \in \N$, one has $2\pi$-periodic eigenfunctions, while for odd $n$, the eigenfunctions are $2\pi$-anti-periodic. There is no spectrum for $E<0$, while $n=1$ yields the first anti-periodic eigenvalue $E_1 = 1/4$. For small $u$, the effect of the perturbation is that of open gaps between the unperturbed doubly degenerate eigenvalues. Moreover, $E_n$ remains close to its unperturbed value. Then the study of the spectrum amounts to identifying the regions of $\lambda>0$, where $E = \lambda^2$, that yield bounded eigenfunctions. Clearly, the study of the spectral gaps coincides with the identification of the complementary region.

Let us introduce the complex coordinates
\begin{equation*}
z=\frac{1}{\sqrt{2}}\bigg(\frac{p}{\sqrt{\lambda}}+\im\sqrt{\lambda}q\bigg) \qquad \text{and} \qquad \bar{z}=\frac{1}{\sqrt{2}}\bigg(\frac{p}{\sqrt{\lambda}}-\im\sqrt{\lambda}q\bigg).
\end{equation*}
Under this change of variables, the eigenvalue equation $L_u\psi = \lambda^2 \psi$ recasts as:
\begin{equation}\label{hamzzbar}
\frac{d}{dt} \begin{pmatrix} z \\ \bar{z} \end{pmatrix} = \bigg[ \lambda \begin{pmatrix} \im & 0 \\ 0 & -\im \end{pmatrix} + \frac{u}{2\lambda} \begin{pmatrix} -\im & \im \\ -\im & \im \end{pmatrix} \bigg] \begin{pmatrix} z \\ \bar{z} \end{pmatrix}.
\end{equation}
This is a parameter family of Hamiltonian dynamical systems w.r. to the standard  symplectic form $\omega = \im dz \wedge d\bar{z}$. 
Accordingly the flow is represented by an element of the group $SU(1,1)$\footnote{The notation is motivated by the fact that, from a complex geometry perspective, $\om$ represents a $(1,1)$-differential form.} 
\begin{equation*}
SU(1,1):=\bigg\{ \begin{pmatrix} 
\alpha & \beta \\
\bar{\beta} & \bar{\alpha}
\end{pmatrix} \ : \ \alpha,\beta \in \C, \ |\alpha|^2-|\beta|^2=1 \bigg\}.
\end{equation*} 
The infinitesimal generators of such a flow belong to the corresponding Lie algebra $ \mathfrak{su}(1,1)$:
\begin{equation*}
\mathfrak{su}(1,1):=\bigg\{ \begin{pmatrix} 
\im a & b \\
\bar{b} & -\im a
\end{pmatrix} \ : \ a \in \R, \ b\in\C \bigg\}\,,
\end{equation*}
with (real) basis  $\{J_1, J_2, J_3\}$
\begin{equation*}
J_1 = \begin{pmatrix} \im & 0 \\ 0 & -\im \end{pmatrix}, \quad J_2 = \begin{pmatrix} 0 & 1 \\ 1 & 0 \end{pmatrix}, \quad J_3 = \begin{pmatrix} 0 & \im \\ -\im & 0 \end{pmatrix}.
\end{equation*}
In this basis, the Hamilton's equations \eqref{hamzzbar} take the form
\begin{equation*}
\frac{d}{dt} \begin{pmatrix} z \\ \bar{z} \end{pmatrix} = D(u,\lambda) \begin{pmatrix} z \\ \bar{z} \end{pmatrix},\quad D(u,\lambda) =  \left( \lambda - \frac{u}{2\lambda} \right) J_1 + \frac{u}{2\lambda} J_3=: \lambda J_1 + P(u,\lambda).
\end{equation*}
We seek a time-dependent change of coordinates in $SU(1,1)$ that reduces the Hamilton's equations \eqref{hamzzbar} to a system with constant coefficients. Assuming such a transformation exists, the new equations of motion will take the form:
\begin{equation*}
\frac{d}{dt} \begin{pmatrix} w \\ \bar{w} \end{pmatrix} = \mathcal{D}(u,\lambda) \begin{pmatrix} w \\ \bar{w} \end{pmatrix},
\end{equation*}
where $\mathcal{D}(u,\lambda)\in \mathfrak{su}(1,1)$ has constant coeffiecients. In particular, $\mathcal{D}$ is traceless, and its eigenvalues are given by:
\begin{equation*}
\mu_{\pm}(u,\lambda)=\pm \sqrt{-\det(\mathcal{D}(u,\lambda))}.
\end{equation*}
Therefore, if $\det(\mathcal{D}(u,\lambda)) > 0$, the dynamics is equivalent to that of a harmonic oscillator, with the frequency determined by the eigenvalues. Conversely, if $\det(\mathcal{D}(u,\lambda)) < 0$, the eigenvalues are real, leading to unbounded dynamics. The intervals of $E=\lambda^2$ for which this occurs correspond to the spectral gaps.

A differentiable, $2\pi$-periodic map $\mathbb{R} \ni t \mapsto B(t) \in SU(1,1)$, conjugates  the equations of motion to
\begin{equation}\label{D_1}
\frac{d}{dt}
\begin{pmatrix}
\xi \\ \bar{\xi}
\end{pmatrix}
= D_1(\lambda)\begin{pmatrix}
\xi \\ \bar{\xi}
\end{pmatrix}, \qquad D_1(\lambda)=B \ (\lambda J_1+P(\lambda))B^{-1}-B\frac{d}{dt}B^{-1}.
\end{equation}

Assuming that $B=\exp(-A)$, for some $A \in \mathfrak{su}(1,1)$ and setting  $\mathrm{ad}_{A}M=[A,M]$ 
equation \eqref{D_1} reads
\begin{equation}\label{xixibar}
\frac{d}{dt}
\begin{pmatrix}
\xi \\ \bar{\xi}
\end{pmatrix}
= 
\bigg(\sum_{k=0}^\infty \frac{(\mathrm{ad}_{-A})^k}{k!}(\lambda J_1 +P) - \sum_{k=1}^\infty \frac{(\mathrm{ad}_{-A})^{k-1}}{k!}\dot{A}\bigg)
\begin{pmatrix}
\xi \\ \bar{\xi}
\end{pmatrix}.
\end{equation}

Assuming that $A = O(|P|)$, the transformation intended to eliminate the time dependence yields the following homological equation at order $O(|P|)$:
\begin{equation}\label{homo}
\dot{A} + \lambda [A, J_1] = P - \mathcal{Z}_0,
\end{equation}
where the matrix $\mathcal{Z}_0$ is introduced to account for the (possibly) non-trivial kernel of the operator
\begin{equation}\label{Llamb}
\mathcal{L}_\lambda := \frac{d}{dt} + \lambda [ \ \cdot \ , J_1].
\end{equation}
{It is well known that, if 
\begin{equation}\label{cond0}
|2\lambda - n| \geq  \gamma_0, \quad  \gamma_0\sim O(|P|)
\end{equation}
then by standard perturbation theory $\lambda$ belongs to the continuous spectrum. Thus for a fixed $n \in \mathbb{N}$, we restrict our study to the region:
\begin{equation}\label{cond}
|2\lambda - n| < \gamma, \quad \frac{1}{4} \leq \gamma \leq \frac{1}{2}.
\end{equation}
This condition ensures that the operator $\mathcal{L}_\lambda$ encounters no singularities other than the resonance at $\lambda = n/2$, providing the necessary information to characterize the $n$-th spectral gap. }

Our strategy it to find, via a contraction argument, a change of  coordinates which eliminates the time-dependent terms lying in the range of the operator $\mathcal{L}_{n/2}$. 
Recall that
\begin{equation*}
\mathrm{Ker}(\mathcal{L}_{n/2}) = 
\bigg\{ \begin{pmatrix} 
\im a & be^{\im nt} \\
\bar{b}e^{-\im nt} & -\im a
\end{pmatrix} \ : \ a \in \mathbb{R}, \ b\in\mathbb{C} \bigg\}.
\end{equation*}

The matrix $A\in \mathfrak{su}(1,1)$ must solve the equation
\begin{equation}
\Pi_{rg} \bigg( -\sum_{k=1}^\infty\frac{\mathrm{ad}_{-A}^{k-1}}{k!}\mathcal{L}_\lambda A+ e^{\mathrm{ad}_{-A}}P\bigg)=0,
\end{equation}
where $\Pi_{rg}$ denotes the projector onto the range of $\mathcal{L}_{n/2}$. Denoting by $\mathcal{B}(A)$ the inverse of the operator $X \mapsto \sum_{k=1}^\infty(\mathrm{ad}_{-A}^{k-1})/k!X$ (see Lemma \ref{invB} for details), and noting that the operator $\mathcal{L}_\lambda$ is invertible on the range of $\mathcal L_{n/2}$ we reduce to proving that the map
\begin{equation} \label{eq:phi_P_final}
\phi_P(A) := \mathcal{L}^{-1}_\lambda \Pi_{rg} \left( \mathcal{B}(A) \Pi_{rg} \left( e^{\mathrm{ad}_{-A}} P \right) \right)\,,
\end{equation}
admits a fixed point in an appropriate space of majorant analytic functions from some ball $B_r(..)$ into the space of analytic time-periodic  matrices in $\mathfrak{su}(1,1) $.

This process leads to a vector field  represented by a matrix of the form:
\begin{equation}
\lambda J_1 + \tilde{\mathcal{Z}}(u,\lambda,t),
\end{equation}
where the residual perturbation $\tilde{\mathcal{Z}}$ belongs to $\mathrm{Ker}(\mathcal{L}_{n/2})$. Applying the standard 
\begin{equation}\label{S_n}
S_n(t) = \begin{pmatrix} e^{-\im nt/2} & 0 \\ 0 & e^{\im nt/2} \end{pmatrix}\in SU(1,1).
\end{equation}
we complete the reduction to constant coefficients as desired. 

\begin{thm}[Reducibility to constant coefficients]\label{final_reducibility_thm}
	Given $p>1$ and $s\geq 0$, let $|P|_{s,p,R} \leq (48 K)^{-1}$, where the constant $K$ is defined in Lemma \ref{banachliealgA}, and let $A_P$ be the fixed point provided by Theorem \ref{fpt}. Then, the transformation defined by the composition $\eta = S_ne^{-A_P} z$ reduces the Hamiltonian system \eqref{questo} to a constant coefficient system of the form:
	\begin{equation}\label{constant_final}
	\frac{d}{dt} \begin{pmatrix} \eta \\ \bar{\eta} \end{pmatrix} = \left[ (\lambda - n/2) J_1 + \mathcal{Z}_0(u, \lambda) + \mathcal{P}_\infty(u, \lambda) \right] \begin{pmatrix} \eta \\ \bar{\eta} \end{pmatrix},
	\end{equation}
	where,
	\begin{equation}\label{gapmatrix}
	\mathcal{Z}_0 = \begin{pmatrix} 0 & P_n^{+-} \\ P_n^{+-} & 0  \end{pmatrix}=\frac{1}{2\lambda}\begin{pmatrix} 0 & \im u_n \\ -\im \bar{u}_n & 0  \end{pmatrix} \quad \text{and} \quad	\mathcal{P}_\infty = \begin{pmatrix} \im a & b \\ \bar{b} & -\im a \end{pmatrix}.
	\end{equation}
	The coefficients $a(u,\lambda) \in \mathbb{R}$ and $b(u,\lambda) \in \mathbb{C}$ are majorant analytic in $u$ and $\lambda$, more precisely they belong to $\mathcal{H}_{s,p,R}$, and satisfy the following estimates:
	\begin{equation*}
	\|a\|_{s,p,R}, \|b\|_{s,p,R} , \|\partial_\lambda a\|_{s,p,R}, \|\partial_\lambda b\|_{s,p,R} \leq C |P|_{s,p,R}^2,
	\end{equation*}
	for a suitable constant $C>0$, uniformly in the region $|2\lambda - n| \leq 1/2$.
\end{thm}
See Appendix \ref{AAE} for the proof.

\subsubsection{Homological equation and spectral gaps at the leading order}\label{leao}

While the full reduction is established through the convergence of the contractive scheme, it is highly instructive to examine the result at the first order in $|P|$. 
For any matrix $W\in\text{Mat}(2,\mathbb{C})$, we write
\begin{equation*}
W=
\begin{pmatrix}
W^{++} & W^{+-} \\ W^{-+} & W^{--}
\end{pmatrix},
\end{equation*}
or $W=(W^{\sigma\sigma'})_{\sigma\sigma'=\pm1}$. 
 Since $u$ is analytic, the matrix-valued function $t \mapsto P(t)$ can be expanded in a Fourier series:
\begin{equation*}
P(t)=\sum_{\ell \in \mathbb{Z}}e^{\im \ell t}
\begin{pmatrix}
P^{++}_\ell & P^{+-}_\ell \\ P^{-+}_\ell & P^{--}_\ell
\end{pmatrix},\quad P^{++}_\ell=P^{-+}_\ell=\frac{-\im u_\ell}{2\lambda} \qquad \text{and} \qquad P^{+-}_\ell=P^{--}_\ell=\frac{\im u_\ell}{2\lambda}.
\end{equation*}
By substituting this expansion into the homological equation
$\mathcal{L}_\lambda A = P - \mathcal{Z}_0,$
we obtain a system of algebraic equations for the Fourier coefficients of the generator $A$:
\begin{equation*}
\im( \ell-(\sigma-\sigma')\lambda)A^{\sigma\sigma'}_\ell= P^{\sigma\sigma'}_\ell - \tilde{\mathcal{Z}}_{0,\ell}^{\sigma\sigma'}.
\end{equation*}
Consequently, the matrix $\tilde{\mathcal{Z}}_0(t)$ representing the projection of $P(t)$ onto the kernel of $\mathcal{L}_\lambda$ takes the form:
\begin{equation*}
\tilde{\mathcal{Z}}_0(t)=
\begin{pmatrix}
0 & P^{+-}_n e^{\im n t} \\ P^{-+}_{-n} e^{-\im n t} & 0
\end{pmatrix}.
\end{equation*}
Under this transformation, the vector field is mapped to $\lambda J_1 + \tilde{\mathcal{Z}}_0 + \mathcal{R}$, where the remainder $\mathcal{R}$ collects all higher-order terms.
Finally, by moving to the rotating frame as defined in \eqref{S_n}, we obtain a matrix of the form
\begin{equation}\label{leadi}
    D_1(\lambda)+ O(|P|^2)\,,\quad D_1(\lambda):=(\lambda - n/2) J_1+
\frac{1}{2\lambda}
\begin{pmatrix}
0 & \im u_n \\ -\im \bar{u}_n & 0
\end{pmatrix}.
\end{equation}
By imposing the instability condition $\det(\mathcal{D}(u, \lambda)) \leq 0$, we can determine the endpoints $\lambda_{\pm}$ of the $n$-th spectral gap, which at first order are
\begin{equation*}
\lambda_+ = \frac{n + \sqrt{n^2 + 8|u_n|}}{4}, \quad \lambda_- = \frac{n + \sqrt{n^2 - 8|u_n|}}{4}.
\end{equation*}
Defining $E_+ = \lambda_+^2$ and $E_- = \lambda_-^2$, the squared gap at first order is 
\begin{equation}\label{gn2}
g^2_n= 4|u_n|^2+\sum_{k=2}^{\infty} \frac{c_{2k}}{n^{4(k-1)}} |u_n|^{2k}\,,\quad |c_{h}| \leq 8^h\,,
\end{equation}
and hence it is a majorant analytic function in the whole ball of radius $R = n^2/8$ centered at the origin.

Since at higher orders the $n$-th gap may depend on all the Fourier coefficients a careful analysis of the corrections coming form $\mathcal P_\infty$ is needed.

\subsubsection{Higher order corrections and proof of Theorem \ref{LA_KdV}}\label{SG} Given the structure of the reduced matrix in Theorem \ref{final_reducibility_thm}, the determinant is given by
\begin{equation}
\det(\mathcal{D}(u, \lambda)) = \left( \lambda - \frac{n}{2} + a(u, \lambda) \right)^2 - \left| \frac{\im u_n}{2\lambda} + b(u, \lambda) \right|^2\,,
\end{equation}
so that
\begin{equation}\label{gap_roots}
\lambda_{\pm} - \frac{n}{2} + a(u, \lambda_{\pm}) = \pm \left| \frac{\im u_n}{2\lambda_{\pm}} + b(u, \lambda_{\pm}) \right|.
\end{equation}
The presence of the term $a(u, \lambda)$ induces a shift in the center of the gap, while the modulus on the right-hand side, which includes the higher-order corrections $b(u, \lambda)$, determines the width of the instability interval. 

Therefore, by looking at the structure of the equation \eqref{gap_roots} defining the end points of the spectral gaps, it is clear that the first and most fundamental step is to investigate the regularity of the coefficients $\alpha$ and $\beta$ defined by:
\begin{equation}
\begin{pmatrix} \im \alpha & \beta \\ \bar{\beta} & -\im \alpha \end{pmatrix} := 
\mathcal{P}_\infty \circ |\partial|^{1/2}.
\end{equation}

From now on, the proof is articulated in the following main steps:
\begin{enumerate}
	\item Analysis of the regularity of the coefficients $\alpha$ e $\beta$.
	\item Majorant analyticity of the spectral gaps.
	\item Derivatives of the spectral gaps and proof of the main theorem.
\end{enumerate}
\bigskip 

\noindent \textbf{Step 1.} Let us start by proving the following proposition, whose technical ingredients are provided in Appendix \ref{AAE}.

\begin{prop}[Regularity of the coefficients $\alpha$ and $\beta$]\label{regcoeff}
Given $p>3/2$ and $s\geq0$, for sufficiently small $R>0$, the maps $\alpha,\beta:\mathscr{U}_R(\ell^\infty_{s,p})\to \mathbb{C}$ (where $\alpha$ is real-valued) are majorant analytic and depends analytically on $\lambda$. Specifically, $\alpha,\beta \in \mathcal{H}_{s,p,R}$ and satisfy the bounds:
	\begin{equation}\label{albet}
	\|\alpha\|_{s,p,R}\leq \frac{n}{\lambda^2}R^2, \quad \|\beta\|_{s,p,R}\leq \frac{n^{3/2}}{\lambda^2}R^2\frac{e^{-s|n|}}{\jap{n}^p}
	\end{equation}
	and
	\begin{equation}\label{dalbet}
	\|\partial_{z_j}\alpha\|_{s,p,R}\leq \frac{n}{\lambda^2}R\frac{e^{-s|j|}}{\jap{j}^p}, \quad \|\partial_{z_j}\beta\|_{s,p,R}\leq \frac{n^{3/2}}{\lambda^2}R\frac{e^{-s|j-n|}}{\jap{j-n}^p}.
	\end{equation}
\end{prop}

\begin{proof}
The majorant analyticity of $\alpha$ and $\beta$ follows from Theorem \ref{final_reducibility_thm} and from the fact that if $f \in \mathcal{H}_{s,p-1/2,R}$, then $f \circ |\partial|^{1/2} \in \mathcal{H}_{s,p}$. Analyticity in $\lambda$ is also inherited from Theorem \ref{final_reducibility_thm}.
	
	Using the structure of $\mathrm{ker}(\mathcal{L}_{n/2})$ given in \eqref{kernel}, and the fact that $\mathcal{P}_\infty=S_n \tilde{\mathcal{P}}_\infty S_n^{-1}$ (see \eqref{PinftySn}), one immediately deduces that 
	\begin{equation}\label{alfabeta}
	\im \alpha = \tilde{\mathcal{P}}^{++}_{\infty,0}\circ |\partial|^{1/2} \quad \text{and} \quad \beta = \tilde{\mathcal{P}}^{+-}_{\infty,n}\circ |\partial|^{1/2}.
	\end{equation}
	Moreover, by Lemma \ref{vittoria}, one can deduce the following bound:
	\begin{equation}\label{vitt1}
		\sup_{\ell \in \Z}\jap{\ell}^{p-1/2}e^{s|\ell|}\max_{\sigma\sigma'=\pm 1}\|\tilde{\mathcal{P}}_{\infty,\ell}^{\sigma\sigma'}\circ|\partial|^{1/2}\|_{s,p,R}\leq \frac{n}{\lambda^2}R^2,
	\end{equation}
	for $R>0$ sufficiently small. Therefore, combining \eqref{vitt1} and \eqref{alfabeta}, we get the bounds \eqref{albet}.
	
	 Lemma \ref{vittoria}, states also that 
	 \begin{equation}\label{vitt2}
	 \sup_{\ell \in \Z}\jap{\ell}^{p-1/2}e^{s|\ell|}\max_{\sigma\sigma'=\pm 1}\|\partial_{z_j}(\tilde{\mathcal{P}}_{\infty,\ell}^{\sigma\sigma'}\circ|\partial|^{1/2})\|_{s,p,R}\leq \frac{n}{\lambda^2}R,
	 \end{equation}
	 from which one can easily derive the bounds in \eqref{dalbet}.
\end{proof}
\bigskip

\noindent \textbf{Step 2.} We now have all the necessary ingredients to study the regularity of the spectral gaps. 

\begin{rmk}\label{constantC}
	Throughout this subsection, $K$ denotes a generic positive constant, independent of $n, z$, and $R$, whose value may change at each occurrence. Furthermore, since we need to recall the majorant norm $|\cdot|_{r,s}$ for the power-trigonometric series (introduced in the second section), we denote the majorant analytic norm in the $z$ variables by $\| \cdot \|_{a,p,R}$ and the corresponding weights by $\langle j \rangle^p e^{a|j|}$, where $a \geq 0$. This choice of notation is intended to avoid any confusion between the decay of the $z$ variables and the regularity parameter $s$.
\end{rmk}

Let us rewrite formula \eqref{gap_roots} as
\begin{equation}\label{fez0}
\phi(z, \lambda_{\pm}) = \pm \left| \frac{i \sqrt{n} z_n}{2} + \lambda_{\pm} \beta(z, \lambda_{\pm}) \right|,
\end{equation}
where $\phi$ is the function defined by
\begin{equation}\label{fez}
\phi(z, \lambda) := \lambda (\lambda - n/2 + \alpha(z, \lambda)).
\end{equation}
Since, by Proposition \ref{regcoeff}, $\alpha$ is analytic as a function of $\lambda$ within the region $|2\lambda - n| \leq 1/2$, it follows that $\phi$ inherits the same regularity properties. In particular, $\phi$ is differentiable with respect to $\lambda$ and its derivative is
\begin{equation*}
\partial_\lambda \phi(z, \lambda) = 2\lambda - n/2 + \alpha(z, \lambda) + \lambda \partial_\lambda \alpha(z, \lambda).
\end{equation*}
This derivative can be bounded from below by a constant proportional to $n$, provided the domain is suitably restricted. For instance, assuming $|2\lambda - n| \leq 1/8$, one ensures that $2\lambda - n/2 \geq n/4$ for all $n \geq 1$. Furthermore, from Proposition \ref{regcoeff}, we have $\|\alpha\|_{a,p,R} \leq n \lambda^{-2} R^2$. Since in the considered region $\lambda \sim n/2$, this term is $O(R^2 n^{-1})$. The derivative $\partial_\lambda \alpha$ can be estimated via a standard Cauchy estimate, yielding a bound proportional to $R^2$. Consequently, for $R$ sufficiently small, the perturbation terms $\alpha + \lambda \partial_\lambda \alpha$ are small compared to $n/4$. Thus, we conclude that
\begin{equation*}
\inf_{\substack{\|z\|_{a,p} \leq R \\ |2\lambda - n| \leq 1/8}} |\partial_\lambda \phi(z, \lambda)| \geq K n,
\end{equation*}
ensuring the invertibility of the map. 

This lower bound on the derivative allows us to control the separation between the roots $\lambda_+(z)$ and $\lambda_-(z)$. By the Mean Value Theorem applied to $\phi$ with respect to $\lambda$, there exists a value $\bar{\lambda}$ in the segment $[\lambda_-(z), \lambda_+(z)]$ such that:
\begin{equation*}
\phi(z,\lambda_+(z)) - \phi(z,\lambda_-(z)) = \partial_\lambda \phi(z,\bar{\lambda}) (\lambda_+(z) - \lambda_-(z)).
\end{equation*}

Recalling the definition of $\phi(z, \lambda_{\pm})$ and the lower bound $|\partial_\lambda \phi(z, \bar{\lambda})| \geq K n$, we have:
\begin{align*}
|\lambda_+(z) - \lambda_-(z)| &= \frac{|\phi(z,\lambda_+(z)) - \phi(z,\lambda_-(z))|}{|\partial_\lambda \phi(z,\bar{\lambda})|} 
\\
&\leq \frac{1}{K n} \left( \left| \frac{i \sqrt{n} z_n}{2} + \lambda_+ \beta(z, \lambda_+) \right| + \left| \frac{i \sqrt{n} z_n}{2} + \lambda_- \beta(z, \lambda_-) \right| \right)
\\
&\leq \frac{1}{K n} \left(\sqrt{n} |z_n| + |\lambda_+ \beta(z, \lambda_+)|+|\lambda_- \beta(z, \lambda_-)|  \right).
\end{align*}
From Proposition \ref{regcoeff}, we have $\|\beta\|_{a,p,R} \leq K n^{3/2} \lambda^{-2} R^2 e^{-a|n|} \jap{n}^{-p}$. Since $\lambda \sim n/2$, the terms $\lambda_{\pm} \beta(z, \lambda_{\pm})$ are of order $\sqrt{n} R^2 e^{-a|n|} \jap{n}^{-p}$. Similarly, the contribution from the potential is given by $\sqrt{n} |z_n| \leq \sqrt{n} R e^{-a|n|} \jap{n}^{-p}$.

By combining these contributions, the numerator is bounded by $K \sqrt{n} (R + R^2) e^{-a|n|} \jap{n}^{-p}$. We then conclude that the following estimate holds uniformly for all $z \in \mathscr{U}_R(\ell^\infty_{a,p})$:
\begin{equation*}
|\lambda_+(z) - \lambda_-(z)| \leq K R\frac{e^{-a|n|}}{\sqrt{n} \jap{n}^p}.
\end{equation*}

Finally, the estimate for the sum $\lambda_+(z) + \lambda_-(z)$ follows directly from the fact that both roots belong to the region $|2\lambda - n| \leq 1/8$. This condition implies $|\lambda_{\pm}| \leq n/2 + 1/16$, and summing these contributions yields:
\begin{equation*}
|\lambda_+(z) + \lambda_-(z)| \leq |\lambda_+| + |\lambda_-| \leq n + 1/8 \leq K n,
\end{equation*}
which again holds uniformly for all $z \in \mathscr{U}_R(\ell^\infty_{a,p})$ and $n \geq 1$.

By combining the previous results, we can now estimate the term $(\lambda_+(z)^2 - \lambda_-(z)^2)$. Using the identity $a^2 - b^2 = (a-b)(a+b)$, we can write:
\begin{equation*}
|\lambda_+(z)^2 - \lambda_-(z)^2| = |\lambda_+(z) - \lambda_-(z)| \cdot |\lambda_+(z) + \lambda_-(z)|.
\end{equation*}
By substituting the uniform estimates obtained above, we get:
\begin{equation*}
|\lambda_+(z)^2 - \lambda_-(z)^2| \leq K R \sqrt{n} \frac{e^{-a|n|}}{\jap{n}^p}.
\end{equation*}
These quantities will be used to construct the action variables for the Korteweg-de Vries equation. In order for the actions to have the correct physical dimension (consistent with the order of $|z|^2$) and to map properly into the considered sequence spaces, it is necessary to introduce a rescaling factor $n^{-1}$. Indeed, by squaring the previous estimate and dividing by $n$, we obtain the bound for rescaled squared gap lengths:
\begin{equation}\label{gammino}
\gamma_n(z):=\frac{(\lambda_+(z)^2 - \lambda_-(z)^2)^2}{n} \leq K R^2 \frac{e^{-2a|n|}}{\jap{n}^{2p}},
\end{equation}
which holds uniformly for all $z \in \mathscr{U}_R(\ell^\infty_{a,p})$.

However, a significant challenge arises from the implicit definition of the gaps via the modulus $|\cdot|$ in \eqref{gap_roots}, which introduces a singularity at the origin ($z=0$). To establish the regularity of the gaps without explicitly addressing the singularity, we exploit the possibility of constructing local symplectic polar coordinates $(y, \varphi)$ around any infinite torus which is sufficiently far form the origin, as explained in Subsection \ref{Subps}. Namely, let $w = (w_k)_{k \in \mathbb{N}}$ be the weight defined by $w_k = e^{2a|k|} \jap{k}^{2p}$. For a fixed sequence $\xi = (\xi_k)_{k \in \mathbb{N}}$ of strictly positive real numbers, we assume that there exist constants $0 < \xi_- \leq \xi_+$ such that:
\begin{equation}\label{xi_decay}
0 < \xi_- \leq \inf_{k \in \mathbb{N}} w_k \xi_k \leq \sup_{k \in \mathbb{N}} w_k \xi_k \leq \xi_+ < \frac{R^2}{4}.
\end{equation}
Under this assumption, for small $r>0$, we define the local chart $\Phi_\xi$ as in \eqref{phixi}
\begin{equation*}
\Phi_\xi : \mathscr{U}_{r}(\ell^\infty_w) \times \mathbb{T}^\infty \ni (y, \varphi) \mapsto z = (\sqrt{\xi_k + y_k}e^{\im \varphi_k})_{k \in \mathbb{N}} \in \ell^\infty_{a,p}
\end{equation*}
in a neighbourhood of the torus
\begin{equation*}
\mathcal{T}_{\xi} := \left\{ z \in \ell^\infty_{a,p} : |z_k|^2 = \xi_k \text{ for all } k \in \mathbb{N} \right\}.
\end{equation*}

The regularity of the spectral gaps in these coordinates is then established through the norm defining the Banach space $\mathcal{A}_{r,s}$ of power-trigonometric series introduced in Definition \ref{matcalA}. Recall that for $f \in \mathcal{A}_{r,s}$, the majorant norm is given by:
\begin{equation*}
|f|_{r,s} := \sum_{\nu \in \mathbb{Z}^{(\mathbb{N})}} \sum_{\alpha \in \mathbb{N}_0^{(\mathbb{N})}} |f_{\nu \alpha}| v_{w,r}^\alpha e^{s|\nu|}, \quad (v_{w,r})_j = r w_j^{-1}.
\end{equation*}
Crucially, the right-hand sides of our estimates for $\gamma_n$ are composed of functions (such as $z_n$ and $\beta(z, \lambda)$) that are majorant analytic in the $z$ variables even at the origin. This allows us to bound their $|\cdot|_{r,s}$ norm by further estimating them with the $\| \cdot \|_{a,p,R}$ norm. 

The transition between these two norms is governed by Proposition \ref{Rrs}. Specifically, for any Hamiltonian $H(z)$ that is majorant real analytic on $\mathscr{U}_R(\ell^\infty_{a,p})$, its pullback $\Phi_\xi^* H$ belongs to $\mathcal{A}_{r,s}$ and satisfies:
\begin{equation*}
|\Phi^*_{\xi}H|_{r,s} \leq \|H\|_{a,p,R},
\end{equation*}
provided the parameters $r, s, \xi, R$ satisfy the consistency condition \eqref{ass}. This framework leads to the following regularity result:

\begin{thm}\label{gap1}
	Given $p>3/2$ and $a\geq 0$, let $\xi\in\ell^{\infty}_{2a,2p}$ be a sequence satisfying the condition \eqref{xi_decay}, and let $R>0$ be the analyticity radius of the coefficients $\alpha$ and $\beta$ (see Proposition \ref{regcoeff}). Then, for $r>0$ and $s\geq 0$ sufficiently small (depending on $R$, as explained in Proposition \ref{Rrs}), the rescaled squared gaps in the coordinate system $(y,\varphi)$, namely
	\begin{equation*}
	\Phi^*_\xi\gamma_n := \frac{((\Phi^*_\xi\lambda_+)^2 - (\Phi^*_\xi\lambda_-)^2)^2}{n},
	\end{equation*}
	belong to the space $\mathcal{A}_{r,s}$ (see Definition \ref{matcalA}) and satisfy the following bound:
	\begin{equation*}
	\sup_{n \in \mathbb{N}} e^{2a|n|} \jap{n}^{2p} |\Phi^*_\xi\gamma_n|_{r,s} \leq K R^2,
	\end{equation*}
	for a suitable positive constant $K$.
\end{thm}

The Theorem above ensures that the rescaled squared gaps are well-behaved as elements of $\mathcal{A}_{r,s}$ and the map $\Gamma := (\Phi^*_\xi\gamma_n)_{n \in \mathbb{N}}$ take values in $\ell^\infty_{2a,2p}$. 
\bigskip

\noindent \textbf{Step 3.} In order to satisfy the smallness and regularity requirements of Theorem \ref{completethm}, we must also ensure a suitable bound for the derivatives of $\Gamma$. 

Specifically, the applicability of the construction of action-angle coordinates requires the finiteness of the following quantity:
\begin{equation}\label{eqB}
\mathcal{K}(\gamma) := \max\bigg\{ r\sup_{j \in \mathbb{N}} \sum_{n=1}^\infty \bigg|\frac{\partial \Phi_\xi^*\gamma_n}{\partial y_j} \bigg|_{r,s}, \ \sup_{j \in \mathbb{N}} e^{2a|j|} \jap{j}^{2p} \sum_{n=1}^\infty \bigg|\frac{\partial \Phi_\xi^*\gamma_n}{\partial \varphi_j} \bigg|_{r,s} \bigg\}.
\end{equation}

To estimate \eqref{eqB}, we relate the derivatives in the local chart $(y,\varphi)$ to those in the canonical variables $(z_j, \bar{z}_j)$ via the relations:\footnote{We recall that, being the potential $u$ a real-valued function, we have that $z_{-\ell}=\bar{z}_\ell$.}
\begin{equation}\label{eqC}
\frac{\partial}{\partial y_j}=\frac{1}{2}\bigg(\frac{1}{\bar{z}_j}\frac{\partial}{\partial z_j}+\frac{1}{z_j}\frac{\partial}{\partial \bar{z}_j}\bigg) \quad \text{and} \quad   
\frac{\partial}{\partial \varphi_j}=\im\bigg(z_j\frac{\partial}{\partial z_j} - \bar{z}_j\frac{\partial}{\partial \bar{z}_j}\bigg),
\end{equation}
and we perform our calculations in the variables $z$. As already observed, for the rescaled squared gaps $\gamma_n$ the majorant analytic norm $\|\cdot\|_{a,p,R}$ in the $z$ variables is not well defined due to the implicit presence of the modulus in \eqref{gap_roots}, which prevents the Taylor expansion around the origin. Nevertheless, as we shall see from the explicit calculations, the derivatives of $\gamma_n$ can be bounded by functions involving the coefficients $\alpha, \beta$ and their derivatives $\partial_{z_j}\alpha, \partial_{z_j}\beta$, for which the majorant norm $\| \cdot \|_{a,p,R}$ is well defined. 

Specifically, we shall prove that the terms in \eqref{eqB} are dominated by the majorant norms of these bounding functions. Namely, by virtue of Proposition \ref{Rrs2}, if we denote such functions by $f^{(n)}$, the quantity $\mathcal{K}(\gamma)$ satisfies:
\begin{equation}\label{eqD}
\mathcal{K}(\gamma) \leq 2 R \sup_{j \in \mathbb{N}} \jap{j}^{p}e^{a|j|} \sum_{n=1}^\infty \left\| \frac{\partial f^{(n)}}{\partial z_j} \right\|_{a,p,R}.
\end{equation}
By substituting the bounds from Proposition \ref{regcoeff} into the expressions for $f^{(n)}$, we shall see that this sum is finite, thus ensuring the applicability of Theorem \ref{completethm}.

Let us proceed with the explicit calculation of the gradient of the rescaled squared gaps with respect to the canonical variables $z_j$. From the definition of $\gamma_n(z)$ (see \eqref{gammino}), we have:
\begin{equation}\label{eqF}
\partial_{z_j}\gamma_n = \frac{(\lambda_+^2-\lambda_-^2)}{n} \left( 2\lambda_+ \partial_{z_j}\lambda_+ - 2\lambda_- \partial_{z_j}\lambda_- \right).
\end{equation}
The regularity of $\partial_{z_j}\gamma_n$ depends on the behavior of the gradients $\partial_{z_j}\lambda_{\pm}$. To compute these, we differentiate the relation \eqref{gap_roots} squared:
\begin{equation}\label{eqE}
\lambda_{\pm}^2\left( \lambda_{\pm} - \frac{n}{2} + \alpha(z, \lambda_{\pm}) \right)^2 = \left| \frac{i \sqrt{n} z_n}{2} + \lambda_{\pm} \beta(z, \lambda_{\pm}) \right|^2.
\end{equation}
By differentiating the left-hand side of \eqref{eqE} with respect to $z_j$, and recalling the definition of $\phi(z, \lambda)$ from \eqref{fez}, we obtain:
\begin{equation}\label{lhsderE}
2\phi(z,\lambda_{\pm}) \left[ \left( 2\lambda_{\pm} - \frac{n}{2} + \alpha(z,\lambda_{\pm}) + \lambda_{\pm}\partial_\lambda \alpha(z,\lambda_{\pm}) \right) \partial_{z_j}\lambda_{\pm} + \lambda_{\pm}\partial_{z_j}\alpha(z,\lambda_{\pm}) \right].
\end{equation}
Similarly, differentiating the right-hand side of \eqref{eqE} yields:
\begin{align}\label{rhsderE}
&\bigg(-\im\sqrt{n}\bar{z}_n+2\lambda_{\pm}\bar{\beta}(z,\lambda_{\pm})\bigg)
\bigg[ \frac{\im\sqrt{n}}{2}\delta_{jn}+\bigg(\beta(z,\lambda_{\pm})+\lambda_{\pm}\partial_{\lambda}\beta(z,\lambda_{\pm})\bigg)\partial_{z_j}\lambda_{\pm}+\lambda_{\pm}\partial_{z_j}\beta(z,\lambda_{\pm})
\bigg] \nonumber \\
&+\bigg(\im\sqrt{n}z_n+2\lambda_{\pm}\beta(z,\lambda_{\pm})\bigg)
\bigg[
\bigg(\bar{\beta}(z,\lambda_{\pm})+\lambda_{\pm}\partial_{\lambda}	\bar{\beta}(z,\lambda_{\pm})\bigg)\partial_{z_j}\lambda_{\pm}+\lambda_{\pm}\partial_{z_j}\bar{\beta}(z,\lambda_{\pm})
\bigg].
\end{align}
Let us introduce the phase factor
\begin{equation}\label{eqH}
e^{\im \theta_{\pm}} := \frac{\im\sqrt{n}z_n + 2\lambda_{\pm}\beta(z,\lambda_{\pm})}{|\im\sqrt{n}z_n + 2\lambda_{\pm}\beta(z,\lambda_{\pm})|},
\end{equation}
Note that $e^{\im \theta_{\pm}}$ is well-defined as long as the denominator is non-zero; this condition is guaranteed within the neighborhood of the torus $\mathcal{T}_\xi$ defined by the chart $\Phi_\xi$. Dividing both sides of \eqref{rhsderE} by the factor $\phi(z, \lambda_{\pm}) = \pm |i \sqrt{n} z_n/2 + \lambda_{\pm} \beta(z, \lambda_{\pm})|$, we can isolate the derivative $\partial_{z_j}\lambda_{\pm}$. The resulting equation for the gradient is:
\begin{align}\label{eqI}
&\left[ \left( 2\lambda_{\pm} - \frac{n}{2} + \alpha + \lambda_{\pm}\partial_\lambda\alpha \right) - \mathrm{Re}\left( e^{-\im \theta_{\pm}} (\beta + \lambda_{\pm}\partial_{\lambda}\beta) \right) \right] \partial_{z_j}\lambda_{\pm} \nonumber \\
&= -\lambda_{\pm}\partial_{z_j}\alpha + \frac{e^{-\im \theta_{\pm}}}{2}\sqrt{n}\delta_{jn} + \mathrm{Re}\left( e^{\im \theta_{\pm}}\lambda_{\pm}\partial_{z_j}\beta \right).
\end{align}
Following the same argument used for the invertibility of $\phi$, the coefficient in the square brackets on the left-hand side of \eqref{eqI} is bounded from below by $Kn$. Consequently, the derivative $\partial_{z_j}\lambda_{\pm}$ can be expressed as:
\begin{equation}\label{eqL_expl}
\partial_{z_j}\lambda_{\pm} = \frac{-\lambda_{\pm}\partial_{z_j}\alpha + \frac{e^{-\im \theta_{\pm}}}{2}\sqrt{n}\delta_{jn} + \mathrm{Re}\left( e^{\im \theta_{\pm}}\lambda_{\pm}\partial_{z_j}\beta \right)}{\left( 2\lambda_{\pm} - \frac{n}{2} + \alpha + \lambda_{\pm}\partial_\lambda\alpha \right) - \mathrm{Re}\left( e^{-\im \theta_{\pm}} (\beta + \lambda_{\pm}\partial_{\lambda}\beta) \right)}.
\end{equation}
By applying the estimates for $\alpha$ and $\beta$ from Proposition \ref{regcoeff} to \eqref{eqL_expl}, we obtain:
\begin{equation}\label{eqL}
|\partial_{z_j}\lambda_{\pm}| \leq \frac{K}{n} \left[ \sqrt{n}\delta_{jn} + R \left( \frac{e^{-a|j|}}{\jap{j}^p} + \sqrt{n}\frac{e^{-a|j-n|}}{\jap{j-n}^p} \right) \right].
\end{equation}
Substituting this bound and the estimate for the gap length $|\lambda_+^2 - \lambda_-^2| \leq K R \sqrt{n} e^{-a|n|} \jap{n}^{-p} $ into \eqref{eqF}, we get:
\begin{equation}\label{eqM}
|\partial_{z_j}\gamma_n(z)| \leq K R \frac{e^{-a|n|}}{\jap{n}^p}\delta_{jn}
+ K R^2 \left( \frac{1}{\sqrt{n}}\frac{e^{-a(|j|+|n|)}}{\jap{n}^p\jap{j}^p} + \frac{e^{-a(|j-n|+|n|)}}{\jap{n}^p\jap{j-n}^p} \right).
\end{equation}
As previously discussed, the right-hand side of \eqref{eqL_expl} consists of functions that are majorant analytic in the $z$ variables, thus providing the required bounding functions $f^{(n)}(z)$ for the estimate \eqref{eqD}. To verify the convergence of the quantity that gives the bound for $\mathcal{K}(\gamma)$, we observe that:
\begin{equation}\label{eqN}
\sup_{j \in \mathbb{N}}\jap{j}^pe^{a|j|} \sum_{n \in \mathbb{N}} \frac{e^{-a(|n|+|j|)}}{\jap{n}^p} \delta_{jn} \leq 1,
\end{equation}
and the other terms satisfy:
\begin{equation}\label{eqZ}
\sup_{j \in \mathbb{N}}\jap{j}^pe^{a|j|}\sum_{n \in \mathbb{N}} \left( \frac{1}{\sqrt{n}}\frac{e^{-a(|j|+|n|)}}{\jap{n}^p\jap{j}^p} + \frac{e^{-a(|j-n|+|n|)}}{\jap{n}^p\jap{j-n}^p} \right) \leq \sum_{n \in \mathbb{N}} \frac{1}{\jap{n}^{p+1/2}}+\sum_{n \in \N}\dfrac{\jap{j}^p}{\jap{j-n}^p\jap{n}^p} < \infty.
\end{equation}

These bounds ensure that $\mathcal{K}(\gamma) < \infty$, proving the desired regularity requirement for the map $\Gamma=(\Phi^*_\xi \gamma_n)_{n \in \N}$ in the local coordinates $(y, \varphi)$. 

From the estimate \eqref{eqM}, we observe that the leading contribution (of order $R$) is purely diagonal and is given by the term $K R e^{-a|n|} \langle n \rangle^{-p} \delta_{jn}$. This term corresponds to the derivative of $4|z_n|^2$, which is consistent with the linear approximation derived in Subsection \ref{leao}. Indeed, as shown in \eqref{gn2} when studying the spectral gaps $g_n$ of the approximated system, we found
\begin{equation}\label{gn2_bis}
g^2_n = 4|u_n|^2 + O(|u|^4),
\end{equation}
where $u_k = \langle k \rangle^{1/2} z_k$ and, in our rescaled setting, $\gamma_n \sim g_n^2 / n$. This identification allows us to decompose the map $\Gamma$ as $\Gamma(y, \varphi) = 4(\xi + y) + P(y, \varphi)$, where $4(\xi + y)$ is the pullback of the leading term $4|z|^2$. The perturbation $P(y, \varphi)$ collects all higher-order terms and, by \eqref{eqD}, \eqref{eqM}, \eqref{eqN} and \eqref{eqZ}, satisfies:
\begin{equation}\label{eqP_small}
\max\bigg\{ r\sup_{j \in \mathbb{N}} \sum_{n=1}^\infty \bigg|\frac{\partial P^{(n)}}{\partial y_j} \bigg|_{r,s}, \ \sup_{j \in \mathbb{N}} e^{2a|j|} \jap{j}^{2p} \sum_{n=1}^\infty \bigg|\frac{\partial P^{(n)}}{\partial \varphi_j} \bigg|_{r,s} \bigg\} \leq K R^3.
\end{equation}

This bound with $R^3$ is crucial for ensuring the smallness condition needed for the applicability of Theorem \ref{completethm}.

We can now summarize these results in the following Theorem.

\begin{thm}\label{gap2}
	Let $\xi$ be a sequence satisfying \eqref{xi_decay}. For $r, s > 0$ sufficiently small, the rescaled squared gaps $(\gamma_n)_{n \in \N}$ satisfy
	\begin{equation}
	\sup_{j \in \mathbb{N}} \sum_{n=1}^\infty \bigg|\frac{\partial \Phi_\xi^*\gamma_n}{\partial y_j} \bigg|_{r,s} < \infty 
	\quad 
	\text{and} 
	\quad 
	\sup_{j \in \mathbb{N}} e^{2a|j|} \jap{j}^{2p} \sum_{n=1}^\infty \bigg|\frac{\partial \Phi_\xi^*\gamma_n}{\partial \varphi_j} \bigg|_{r,s}<\infty.
	\end{equation}
	Moreover, the map $\Gamma = (\Phi^*_\xi\gamma_n)_{n \in \mathbb{N}}$ can be written as
	\begin{equation*}
	\Gamma(y, \varphi) = 4(\xi + y) + P(y, \varphi), 
	\end{equation*}
	where the perturbation $P$ is majorant analytic and satisfies the bound
	\begin{equation}\label{figa}
		\|P\|^{\infty,w}_{r,s} \leq K R^3
	\end{equation}
	for a suitable positive constant $K$, where $w_j=\jap{j}^{2p}e^{2a|j|}$, for all $j \in \N$, and the norm $\| \cdot \|^{\infty,w}_{r,s}$ is defined in Definition~\ref{definitionnorm}.
\end{thm}

Theorems \ref{gap1} and \ref{gap2} provide the final verification of the assumptions of Theorem \ref{completethm}. Specifically, the decomposition of the map $\Gamma$ into a leading diagonal part plus a perturbation $P$ of order $R^3$ ensures that the smallness condition \textnormal{(a.1)} is satisfied for a sufficiently small but finite radius $R$. Furthermore, it is well-known that the squared gap lengths are integrals of the motion in involution, i.e., $\{\gamma_m, \gamma_n\} = 0$ and $\{\gamma_n, H_{KdV}\} = 0$ for all $m, n \in \mathbb{N}$. Since the pullback via $\Phi_\xi$ preserves the Poisson brackets (see the discussion in \ref{subpb2}), the map $\Gamma$ provides a complete set of integrals in involution for the Hamiltonian $H_{KdV}\circ \Phi_\xi$, fulfilling condition \textnormal{(a.2)} of Theorem \ref{completethm}.

\begin{proof}[Proof of Theorem \ref{LA_KdV}]
The proof consists in verifying that the Hamiltonian $H_{KdV}$ and the map $\Gamma = (\Phi^*_\xi \gamma_n)_{n \in \mathbb{N}}$ satisfy the hypotheses of Theorem \ref{completethm}. 
	
First, we observe that the Hamiltonian $H_{KdV}$ is majorant analytic on $\ell^\infty_{a,p}$ if $p>2$. Its pullback via the local chart, $\Phi_\xi^* H_{KdV}$, thus belongs to the class $\mathcal{A}$.
	
By Theorems \ref{gap1} and \ref{gap2} the map can be decomposed as $\Gamma(y, \varphi) = 4(\xi + y) + P(y, \varphi)$, where $P$ satisfies the cubic bound $\|P\|^{\infty,w}_{r,s} \leq K R^3$ (see \eqref{eqP_small}), with $R>0$ being the analyticity radius of the coefficients $\alpha$ and $\beta$ determined in Proposition \ref{regcoeff}. 
	
Recalling the conditions \eqref{cazzo}, together with \eqref{figa}, and \eqref{ass}, one can find $0<\xi_- \leq\xi_+$ such that, there exist $R'\leq R$, $r'\leq r$ and $s' \leq s$ satisfying
	\begin{equation*}
		\xi_+^{1/2}e^{s'} \sum_{m=0}^{\infty}|a_m|\bigg(\frac{r'}{\xi_-}\bigg)^m\leq R' 
		\quad \text{and} \quad
		K(R')^3\leq \frac{1}{4C}\min\{r',s',1\},
	\end{equation*}
	ensuring the condition $(a.1)$ of Theorem \ref{completethm}.
	Finally, the commutativity condition \textnormal{(a.2)} follows from the standard spectral theory: the squared gap lengths $\gamma_n$ are integrals of motion in involution. Since the map $\Phi_\xi$ is a symplectomorphism, it preserves the Poisson brackets, ensuring that the components of $\Gamma$ commute among themselves and with the Hamiltonian.
\end{proof}

\appendix

\section{Properties of majorant analytic functions}\label{appanalytic}

\begin{prop}\label{uniform}	
	If $f \in \mathcal{A}_{r,s}$ for some $r,s>0$, then it is absolutely convergent and defines 
	a function
	\begin{equation*}
	\mathscr{U}_r(\ell^\infty_w)\times \ell^\infty \to \mathbb{R} \ , \qquad (y,\varphi) \mapsto f(y,\varphi) \ ,
	\end{equation*}
	whose complex extension $f_\mathbb{C}$, obtained by considering $(y,\varphi)\in D_{r,s}$ in formula \eqref{fourier0}, is totally convergent. 
	In particular, for all exhaustions $(S_n)_{n \in \mathbb{N}}$ of $\, \mathbb{N}_0^{(\mathbb{N})} \times \mathbb{Z}^{(\mathbb{N})}$ by finite subsets, one has
	\begin{equation}\label{lim}
	\lim_{n \to \infty} \sup_{(y,\varphi) \in D_{r,s}}\bigg|f_\mathbb{C}(y,\varphi)-\sum_{(\alpha,\nu) \in S_n}f_{\alpha,\nu}y^\alpha \ e^{\im \nu \cdot \varphi} \bigg|=0.
	\end{equation}
	\begin{proof}
		The complex extension of $f$ is a totally convergent series because
		\begin{equation*}
		\sum_{\alpha \in \mathbb{N}_0^{(\mathbb{N})}}\sum_{\nu \in \mathbb{Z}^{(\mathbb{N})}}\sup_{(y,\varphi)\in D_{r,s}} |f_{\alpha,\nu}y^\alpha e^{\im \nu \cdot \varphi}| = |f|_{r,s}<\infty.
		\end{equation*}
		In particular, $f$ is absolutely convergent and hence the attained values do not depend on the order of summation. Namely, $f$ defines a function
		\begin{equation*}
		\mathscr{U}_r(\ell^\infty_w)\times \ell^\infty \to \mathbb{R} \ , \qquad (y,\varphi) \mapsto f(y,\varphi),
		\end{equation*} 
		whose complex extension $f_\mathbb{C}$ is represented by a totally convergent series. Formula \eqref{lim} follows by the total convergence.
	\end{proof}
\end{prop}

\bigskip
\begin{proof}[Proof of Proposition \ref{manalytic}]
Let us start by proving that any $f\in\cA$ is majorant real analytic. To this purpose, let $(S_n)_{n \in \mathbb{N}}$ be an exhaustion of $\mathbb{Z}^{(\mathbb{N})}\times \mathbb{N}_0^{(\mathbb{N})}$ by finite subsets.
		Define the sequence of functions $(f_{M,n})_{n \in \mathbb{N}}\subset \mathcal{A}_{r,s}$ of finitely many variables as
		\begin{equation*}
		f_{M,n}(y,\varphi):=\sum_{(\alpha,\nu)\in S_n}|f_{\alpha,\nu}|y^\alpha e^{\im \nu \cdot \varphi}.
		\end{equation*}
		Now, being each $f_{M,n}$ a trigonometric polynomial with polynomial functions as coefficients, it extends to a complex valued holomorphic function $f_{M,n,\C}$ by considering $(y,\varphi)\in D_{r,s}$ in formula \eqref{fourier0}. Moreover, by Proposition \ref{uniform}, the sequence $(f_{M,n,\C})_{n \in \mathbb{N}}$ converges uniformly in $D_{r,s}$ to $f_{M,\C}$. The conclusion follows by the fact that the uniform limit of a sequence of holomorphic functions on complex Banach spaces is a holomorphic function. Of course, the limit is real on reals.

        Let us prove equality \eqref{ansia1}. First note that, since for all $y \in \mathscr{U}_r(\ell^\infty_w)$ the series $\sum_{\nu}f_\nu(y) e^{\im \nu \cdot \varphi}$ converges uniformly for $\varphi \in \ell^\infty$, one has 
		\begin{equation*}
		\int_{\mathbb{T}^N}f(y,\varphi)d\varphi_1 \ldots d\varphi_N=\sum_{\nu}f_\nu(y)\int_{\mathbb{T}^N}e^{\im \nu \cdot \varphi}d\varphi_1 \ldots d\varphi_N 
		\end{equation*}
		for all $N \in \mathbb{N}$.
		Moreover, for any $N \in \mathbb{N}$, let $S_N$ denote the subset of $\mathbb{Z}^{(\mathbb{N})}$ given by
		\begin{equation*}
		S_N:=\{ \nu \in \mathbb{Z}^{(\mathbb{N})} \ : \ \nu \neq 0 \ , \ \nu_j=0 \quad \forall j>N \} \ .
		\end{equation*}
		Thus, for all $N\in \mathbb{N}$ and any $\nu \in S_N$, one has
		\begin{equation*}
		\frac{1}{(2 \pi)^{N}} \int_{\mathbb{T}^{N}} e^{\im \nu \cdot \varphi}\, d\varphi_1 \ldots d\varphi_N = 0 \ .
		\end{equation*}
		Hence,
	\[
	\begin{aligned}
	& \frac{1}{(2 \pi)^{N}}\int_{\mathbb{T}^{N}}f(y,\varphi) \, d\varphi_1 \ldots d\varphi_N \\
	& \qquad \qquad = \frac{1}{(2 \pi)^{N}}\int_{\mathbb{T}^{N}} \bigg(f_0(y) + \sum_{\nu \in S_N} f_\nu(y) e^{\im \nu \cdot \varphi}   
	+ \sum_{0\neq\nu \notin S_N} f_{\nu}(y) e^{\im \nu \cdot \varphi} \bigg) d\varphi_1 \ldots d\varphi_N  
	\end{aligned}
	\]
	\[
	\begin{aligned}
	& \qquad \qquad = f_0(y) +  \sum_{0\neq\nu \notin S_N} f_\nu(y) \frac{1}{(2 \pi)^{N}}\int_{\mathbb{T}^{N}} e^{\im \nu \cdot \varphi} \, d\varphi_1 \ldots d\varphi_N \ .
	\end{aligned}
	\]
		The latter integral can only take the values \(0\) or \(\pm1\). Moreover, the tail of the series with $0\neq\nu \notin S_N$ goes to zero as $N \to \infty$. This proves \eqref{ansia1} for $\nu=0$.
		
		Now, take $\nu' \in \mathbb{Z}^{(\mathbb{N})} \setminus \{ 0 \}$ and set
		\[
		f(y,\varphi,\nu') := f(y,\varphi) e^{- \im \nu' \cdot \varphi} = \sum_{\nu} f_{\nu}(y) e^{\im (\nu - \nu') \cdot \varphi} = \sum_{\nu} f_{\nu+\nu'}(y) e^{\im \nu \cdot \varphi}\,.
		\]
		By applying \eqref{ansia1} with $\nu=0$ to $f(y,\varphi,\nu')$ and using the fact that $\nu'$ is arbitrary, the equality \eqref{ansia1} follows for $\nu\neq0$ as well.
        \end{proof}

\bigskip
        \begin{defn}
	Let $X,Y$ be Banach spaces and $\mathscr{U}\subseteq X$ an open subset. A countable family of functions $(f^{(m)})_{m \in \N}$, $ \ f^{(m)}:\mathscr{U}\to Y$, is uniformly bounded if there exists $C<\infty$ such that
	\begin{equation*}
	\sup_{m \in \N}\sup_{x \in \mathscr{U}}|f^{(m)}(x)|_Y \leq C.
	\end{equation*}
\end{defn}

\begin{lemma}\label{uniformanalytic}
	Let $X$ be a complex Banach space and $\mathscr{U}\subseteq X$ an open subset. Let $(f^{(m)})_{m \in \N}$ be a uniformly bounded countable family of analytic functions $f^{(m)}:\mathscr{U}\to\mathbb{C}$. Then, the function $f$ defined by $f(x):=(f^{(m)}(x))_{m \in \N}$ defines an analytic function 
	\[ f:\mathscr{U}\to \ell^\infty_\C \ . \]
    \end{lemma}
	\begin{proof}
		For any $n \in \mathbb{N}$, let $L_S^n(X,\ell^\infty_\C)$ denote the space of all bounded $n$-linear symmetric maps from $X^n$, i.e., the direct product of $n$ copies of $X$, taking values in $\ell^\infty_\C$, and let us define the map
		\begin{equation*}
		d^nf:\mathscr{U} \to L_S^n(X,\ell^\infty_\C), \qquad x \mapsto d^nf(x):=(d^nf^{(m)}(x))_{m \in \N},
		\end{equation*}
		where, for any $m \in \N$, $d^nf^{(m)}$ denotes the $n^{\text{th}}$-differential of $f^{(m)}$. Namely, for any $n\in \mathbb{N}$ and for all $(h^{(j)})_{j=1}^n \in X^n$ and all $x \in \mathscr{U}$,
		\begin{equation*}
		d^nf(x)[h^{(1)},\dots,h^{(n)}]=(d^nf^{(m)}(x)[h^{(1)},\dots,h^{(n)}])_{m \in \N}.
		\end{equation*}
		Note that, for all $x \in \mathscr{U}$, by the $n$-linearity and the symmetry property of each $d^nf^{(m)}(x)$, $d^nf(x)$ is indeed an $n$-linear and symmetric map. Moreover, given $x \in \mathscr{U}$ and $h^{(1)},\dots,h^{(n)} \in X$, let $r>0$ be such that $\overline{\mathscr{U}}_{r,x}(X)\subset \mathscr{U}$ and let $(\rho_j)_{j=1}^n$ be a list of real numbers such that, for all $j\in \{1,\dots,n\}$, $0<\rho_j<r/(n|h^{(j)}|_X)$. Then, for all $m \in \N$, since each $f^{(m)}$ is analytic (and hence weak-analytic), by the Cauchy's integral formula applied to the complex valued functions of $n$ complex variables
		\begin{equation*}
		\prod_{j=1}^n\mathscr{U}_{r/(n|h^{(j)}|_X)} \ni (\xi_1,\dots,\xi_n)\mapsto f^{(m)}\left(x+\sum_{j=1}^n\xi_jh^{(j)}\right)\in \mathbb{C},
		\end{equation*}
		we have
		\begin{equation*}\label{cauchyintegraln}
		d^nf^{(m)}(x)[h^{(1)},\dots,h^{(n)}]=\frac{1}{(2\pi \im)^n}\int_{|\xi_1|=\rho_1} \dots \int_{|\xi_n|=\rho_n}\frac{f^{(m)}(x+\sum_{j=1}^n \xi_jh^{(j)})}{\prod_{j=1}^n\xi^2_j}d\xi_1\dots d\xi_n,
		\end{equation*}
		and the integrals do not depend on $(\rho_j)_{j=1}^n$ as long as, $\forall j=1,\dots,n$, $0<\rho_j< r/(n |h^{(j)}|_X)$. Therefore, one immediately has the estimate 
		\begin{equation*}
		|d^nf^{(m)}(x)[h^{(1)},\dots,h^{(n)}]| \leq C\prod_{j=1}^n\frac{1}{2\pi}\int_{|\xi_j|=\rho_j}\frac{|d\xi_j|}{|\xi_j|^2}=C\prod_{j=1}^n\frac{1}{\rho_j},
		\end{equation*}
		where $C:=\sup_{m \in \N}\sup_{x \in \mathscr{U}}|f^{(m)}(x)|$, from which, by taking the limits $\rho_j \to r/(n|h^{(j)}|_X)$, one derives the boundedness of $d^nf(x)$, namely,
		\begin{equation*}
		|d^nf(x)[h^{(1)},\dots,h^{(n)}]| _{\ell^\infty_\C}\leq C\frac{n^n}{r^n}\prod_{j=1}^n|h^{(j)}|_X.
		\end{equation*}
		
		Now we show that, for all $n \in \mathbb{N}$, $d^n f$ is indeed the $n^{\text{th}}$-differential of $f$. Of course, we set $d^0f:=f$. Now, given $x \in \mathscr{U}$ and $(h^{(j)})_{j=1}^n\in X^n$, let $r$ and $(\rho_j)_{j=1}^n$ be as before. Then, one has
		\begin{equation}\label{n-1}
		\begin{aligned}
		&d^{n-1}f^{(m)}(x+h^{(n)})[h^{(1)},\dots,h^{(n-1)}] \\
		&\quad =\frac{1}{(2\pi \im)^{n-1}}\int_{|\xi_1|=\rho_1} \dots \int_{|\xi_{n-1}|=\rho_{n-1}}\frac{f^{(m)}(x+\sum_{j=1}^{n-1} \xi_jh^{(j)}+h^{(n)})}{\prod_{j=1}^{n-1}\xi^2_j}d\xi_1\dots d\xi_{n-1}.
		\end{aligned}
		\end{equation}
		By implementing again Cauchy's integral formula applied to the analytic function
		\begin{equation*}
		\mathscr{U}_{r/(n|h^{(n)}|_X)}\ni z \mapsto f^{(m)}\left(x+\sum_{j=1}^{n-1}\xi_jh^{(j)}+zh^{(n)}\right)\in \mathbb{C},
		\end{equation*}
		one gets
		\begin{equation*}
		\eqref{n-1}=\frac{1}{(2\pi \im)^{n}}\int_{|\xi_1|=\rho_1} \dots \int_{|\xi_{n}|=\rho_{n}}\frac{f^{(m)}(x+\sum_{j=1}^{n} \xi_jh^{(j)})}{\prod_{j=1}^{n-1}\xi^2_j(\xi_n-1)}d\xi_1\dots d\xi_{n}.
		\end{equation*}
		Furthermore, the expansion
		\begin{equation*}
		\frac{1}{\xi_n-1}=\frac{1}{\xi_n}+\frac{1}{\xi_n^2}+\frac{1}{\xi_n^2(\xi_n-1)} \ ,
		\end{equation*}
		and further applications of Cauchy's integral formula, yield
		\begin{equation*}
		\begin{aligned}
		\eqref{n-1} &= d^{n-1}f^{(m)}(x)[h^{(1)},\dots,h^{(n-1)}]+d^nf^{(m)}(x)[h^{(1)},\dots,h^{(n)}] \\
		&\quad +\frac{1}{(2\pi \im)^{n}}\int_{|\xi_1|=\rho_1} \dots \int_{|\xi_{n}|=\rho_{n}}\frac{f^{(m)}(x+\sum_{j=1}^{n} \xi_jh^{(j)})}{\prod_{j=1}^{n}\xi^2_j(\xi_n-1)}d\xi_1\dots d\xi_{n}.
		\end{aligned}
		\end{equation*}
		Therefore,
		\begin{equation*}
		\begin{aligned}
		&\left|(d_{x+h^{(n)}}^{n-1}f^{(m)}-d^{n-1}f^{(m)}(x))[h^{(1)},\dots,h^{(n-1)}]-d^nf^{(m)}(x)[h^{(1)},\dots,h^{(n)}]\right| \\
		&\quad \leq C\frac{1}{\rho_n(\rho_n-1)} \prod_{j=1}^{n-1}\frac{1}{\rho_j},
		\end{aligned}
		\end{equation*}
		from which, by taking the limits $\rho_j \to r/(n|h^{(j)}|_X)$ one gets
		\begin{equation*}
		\begin{aligned}
		&\sup_{h^{(1)},\dots,h^{(n-1)}\in\mathscr{U}_1(X)}|(d^{n-1}f(x+h^{(n)})-d^{n-1}f(x))[h^{(1)},\dots,h^{(n-1)}]-d^nf(x)[h^{(1)},\dots,h^{(n)}]| _{\ell^\infty_\C} \\
		&\quad \leq C \frac{n^n}{r^n}\frac{|h^{(n)}|^2_X}{(r/n-1)}=o(|h^{(n)}|_X).
		\end{aligned}
		\end{equation*}
		
		We conclude the proof by showing that, for any $x \in \mathscr{U}$, there exists $r_0>0$ such that the Taylor series
		\begin{equation*}
		\sum_{n=0}^\infty\frac{1}{n!}d^nf(x)[h,\dots,h]
		\end{equation*}
		converges to $f(x+h)$ uniformly for $h \in \mathscr{U}_{r_0}(X)$. The thesis follows by P\"oschel-Trubowitz, Appendix A, Theorem $1$. First of all, let $r>0$ be such that $\overline{\mathscr{U}}_{r,x}(X)\subset \mathscr{U}$, $h \in X$, and take any $\rho \in (1,r/|h|_{X})$. Then, Cauchy's integral formula gives
		\begin{equation*}
		d^nf^{(m)}(x)[h,\dots,h]=\frac{n!}{2\pi\im}\int_{|\xi|=\rho}\frac{f^{(m)}(x+\xi h)}{\xi^{n+1}}d\xi.
		\end{equation*}
		This formula, together with the expansion
		\begin{equation*}
		\frac{1}{(\xi-1)}=\sum_{n=0}^N\frac{1}{\xi^{n+1}}+\frac{1}{\xi^{N+1}(\xi-1)} \ ,
		\end{equation*}
		yields
		\begin{equation*}
		f^{(m)}(x+h)-\sum_{n=0}^N\frac{1}{n!}d^nf^{(m)}(x)[h,\dots,h]=\frac{1}{2\pi\im}\int_{|\xi|=\rho}\frac{f^{(m)}(x+\xi h)}{\xi^{N+1}(\xi-1)}d\xi.
		\end{equation*}
		Therefore,
		\begin{equation}\label{tail}
		\left|f(x+h)-\sum_{n=0}^N\frac{1}{n!}d^nf(x)[h,\dots,h]\right| _{\ell^\infty_\C}
		\end{equation}
		is bounded by
		\begin{equation*}
		C\frac{|h|^{N+1}_X}{r^N(r-|h|_X)}.
		\end{equation*}
		From this estimate, we conclude that
		\begin{equation*}
		\lim_{N\to\infty} \eqref{tail}=0
		\end{equation*}
		uniformly for $h \in \mathscr{U}_{r_0}(X)$, for any $r_0<r$.
	\end{proof}

\section[A Fixed Point Theorem in Banach Lie Algebras]{A Fixed Point Theorem in Banach-Lie Algebras}\label{AAD}

In this appendix, we investigate the existence and uniqueness of a solution for a nonlinear equation of the form $X = \phi(X)$. We operate within a Banach Lie algebra $(\mathcal{V}, [\cdot,\cdot] |\cdot|)$ characterized by a sub-multiplicative norm such that $|[X,Y]| \leq K |X| |Y|$ for all $X, Y \in \mathcal{V}$. Within this framework, we consider two analytic operator-valued maps $\mathcal{E}$ and $\mathcal{B}$, and we construct formally the map $\phi_P: \mathscr{U}\subset \mathcal{V} \to \mathcal{V}$ in a suitable neighbourhood $\mathscr{U}$ of the origin as:
\begin{equation*}
\phi_P(X) = T \Pi \left( \mathcal{B}(X) \Pi (\mathcal{E}(X) P) \right)
\end{equation*}
where $P \in \mathcal{V}$ is a fixed element, and $\Pi \in \mathcal{L}(\mathcal{V})$ is a bounded linear operator satisfying $\Pi^2 = \Pi$, while $T: \mathrm{range}(\Pi) \to \mathcal{V}$ is a bounded linear operator defined on the range of $\Pi$. Our objective is to prove that $\phi_P$ admits a unique fixed point $X^* \in \mathcal{V}$ by showing that it is a contraction. To achieve this, we first establish the analyticity of the constituent operators; second, we derive their Fréchet derivatives; and third, we provide quantitative norm estimates to satisfy the hypotheses of the Banach Fixed Point Theorem.

\subsubsection{Functional framework}

Let $(\mathcal{V}, [\cdot,\cdot], |\cdot|)$ be a Banach Lie algebra with sub-multiplicative norm: $|[X,Y]| \leq K|X| |Y|$ for all $X, Y \in \mathcal{V}$. We denote by $\mathcal{L}(\mathcal{V})$ the space of bounded linear operators on $\mathcal{V}$ with the operator norm $\|\cdot\|$. For $X \in \mathcal{V}$, the adjoint operator $\mathrm{ad}_{X} \in \mathcal{L}(\mathcal{V})$ is defined by $\mathrm{ad}_{X}(Y) = [X, Y]$. Its norm satisfies:
\begin{equation} \label{eq:ad_norm}
\|\mathrm{ad}_{X}\| \leq K|X|.
\end{equation}

\subsubsection{Definition and analyticity of map $\mathcal{E}$ and $\mathcal{B}$}\label{EandB}

For $X \in \mathcal{V}$, we formally define two operators $\mathcal{E}(X)$ and $\mathcal{B}(X)$ via the Taylor series of two complex functions $f(z) = e^z$ and $g(z) = \frac{z}{e^z - 1}$. 
\bigskip

It is well known that:
\begin{enumerate}
	\item $f(z) = \sum_{k=0}^{\infty} \frac{1}{k!} z^k$ has radius of convergence $R_f = \infty$;
	\item let $g(z) = \sum_{k=0}^{\infty} \beta_k z^k$ be the Taylor expansion of the function $g$ around $z=0$. Then the radius of convergence is $R_g = 2\pi$.
\end{enumerate}
The operators are then defined as $\mathcal{E}(X) = f(\mathrm{ad}_{-X})$ and $\mathcal{B}(X) = g(\mathrm{ad}_{-X})$.

\begin{prop}[Analyticity of $\mathcal{E}$ and $\mathcal{B}$]\label{propcattiva}
	The maps $\mathcal{E}: \mathcal{V} \to \mathcal{L}(\mathcal{V})$ and $\mathcal{B}: \mathscr{U}_{2\pi/K}(\mathcal{V}) \to \mathcal{L}(\mathcal{V})$ are analytic
	\end{prop}

\begin{proof}
	Let $h(z) = \sum_{k=0}^{\infty} c_k z^k$ be a power series with radius of convergence $R$. Define the sequence of partial sums $S_n(X) = \sum_{k=0}^{n} c_k (\mathrm{ad}_{-X})^k$. 
	\begin{enumerate}
		\item \textbf{Analyticity of polynomials:} Since $X \mapsto \mathrm{ad}_{-X}$ is a bounded linear map, each term $c_k (\mathrm{ad}_{-X})^k$ is a continuous polynomial in $X$. In any Banach space, continuous polynomials are analytic functions. Thus, $S_n(X)$ is analytic for every $n \in \mathbb{N}$.
		\item \textbf{Uniform convergence:} Let $r < R$ and consider $X \in \mathscr{U}_{r/K}(\mathcal{V})$. From the bound $\|\mathrm{ad}_{-X}\| \leq K|X|$ established in \eqref{eq:ad_norm}, we have $K|X| \leq r$, which implies:
		\[ \sum_{k=0}^{\infty} \|c_k (\mathrm{ad}_{-X})^k\| \leq \sum_{k=0}^{\infty} |c_k| \|\mathrm{ad}_{-X}\|^k \leq \sum_{k=0}^{\infty} |c_k| r^k. \]
		Since $r < R$, the numerical series $\sum |c_k| r^k$ converges. By the Weierstrass M-test, the series of operators converges absolutely and uniformly in the operator norm $\|\cdot\|$ on $\mathscr{U}_{r/K}(\mathcal{V})$. This proves that the sequence of partial sums $S_n(X)$ converges to a limit operator $S(X)$ for every $X$ in the open ball $\mathscr{U}_{R/K}(\mathcal{V})$.
		\item \textbf{Limit property:} According to a fundamental result in complex analysis (cf. Mujica, \textit{Complex Analysis in Banach Spaces}), if a sequence of analytic functions $\{S_n\}$ on an open set $\Omega$ converges uniformly to a function $S$, then the limit $S$ is analytic on $\Omega$. 
	\end{enumerate}
	Since the analyticity radii of $f(z)$ and $g(z)$ are $R_f = \infty$ and $R_g = 2\pi$ respectively, the maps $\mathcal{E}$ (defined on the whole $\mathcal{V}$) and $\mathcal{B}$ (defined on $\mathscr{U}_{2\pi/K}(\mathcal{V})$) are analytic.
\end{proof}

\subsubsection{Definition and analyticity of the map $\phi_P$}

Let $P \in \mathcal{V}$. Consider a bounded linear operator $\Pi \in \mathcal{L}(\mathcal{V})$ such that $\Pi^2 = \Pi$ (a projection), and a linear operator $T: \mathrm{range}(\Pi) \to \mathcal{V}$ which satisfies the norm bound $\|T\|_{\mathcal{L}(\mathrm{range}(\Pi), \mathcal{V})} \leq 2$. We also assume $\|\Pi\| \leq 1$. 
\bigskip  

\begin{defn}
	We define the map $\phi_P: \mathscr{U}_{2\pi/K}(\mathcal{V}) \to \mathcal{V}$ as follows:
	\begin{equation*}
	\phi_P(X) = T \Pi \left( \mathcal{B}(X) \Pi (\mathcal{E}(X) P) \right),
	\end{equation*}
	where the non-linear analytic maps $\mathcal{E}$ and $\mathcal{B}$ are defined in Subsection \ref{EandB}.
\end{defn}
The analyticity of $\phi_P$ follows directly from its construction as a composition of analytic operator-valued maps and constant, bounded linear operators applied to an element $P\in\mathcal{V}$.

\subsubsection{Differentiation of $\mathcal{E}$ and $\mathcal{B}$}

To differentiate the map $\phi_P, $ we must first establish the precise form of the Fr\'echet derivative of the operator-valued power series.
\bigskip 

\noindent \textbf{Derivative of the adjoint operator.}
The map $\alpha: \mathcal{V} \to \mathcal{L}(\mathcal{V})$ defined by $\alpha(X) = \mathrm{ad}_{-X}$ is a bounded linear operator. By the definition of the Fr\'echet derivative for linear maps, we have:
\begin{equation*}
\mathrm{d} \alpha(X)[H] = \mathrm{ad}_{-H}, \quad \forall H \in \mathcal{V}.
\end{equation*}

\noindent \textbf{The Leibniz rule for operator powers.}
Consider the $k$-th power of the adjoint operator as a product of $k$ linear maps. Applying the Leibniz rule for linear maps in a Banach algebra of operators, the derivative in the direction $H$ is the sum of $k$ terms, where the derivative $\mathrm{ad}_{-H}$ occupies each possible position in the product:
\begin{equation*}
\mathrm{d} \alpha^k (X) [H] = \sum_{j=0}^{k-1} (\mathrm{ad}_{-X})^j \circ (\mathrm{ad}_{-H}) \circ (\mathrm{ad}_{-X})^{k-1-j}.
\end{equation*}
This generalized Leibniz formula accounts for the non-commutativity of the operators within $\mathcal{L}(\mathcal{V})$.
\bigskip

\noindent \textbf{Term-wise differentiation of the series.}
Assuming that we can differentiate term-by-term, the derivative of $\mathcal{E}$ and $\mathcal{B}$ are
\begin{equation} \label{eq:deriv_E}
\mathrm{d} \mathcal{E}(X)[H] = \sum_{k=1}^{\infty} \frac{1}{k!} \sum_{j=0}^{k-1} (\mathrm{ad}_{-X})^j \circ (\mathrm{ad}_{-H}) \circ (\mathrm{ad}_{-X})^{k-1-j},
\end{equation}
\begin{equation} \label{eq:deriv_B}
\mathrm{d} \mathcal{B}(X)[H] = \sum_{k=1}^{\infty} \beta_k \sum_{j=0}^{k-1} (\mathrm{ad}_{-X})^j \circ (\mathrm{ad}_{-H}) \circ (\mathrm{ad}_{-X})^{k-1-j},
\end{equation}
where $\beta_k$ are the Taylor coefficients of $g(z) = z/(e^z - 1)$. In the next subsection, we provide quantitative bounds to establish that these operators are indeed the Fréchet derivatives of $\mathcal{E}$ and $\mathcal{B}$.

\subsubsection{Estimates}

In this subsection, we derive explicit bounds for the maps $\mathcal{E}$ and $\mathcal{B}$ and their differentials. These estimates are fundamental to establishing the contraction property of the map $\phi_P$.
\bigskip

\noindent \textbf{Estimates for $\mathcal{E}$ and $\mathcal{B}$.}
Let $\rho = \|\mathrm{ad}_{-X}\| \leq K|X|$. For $\rho < 1/2$, the series defining $\mathcal{E}(X)$ and $\mathcal{B}(X)$ are totally convergent in $\mathcal{L}(\mathcal{V})$. This implies that the following norm estimates hold uniformly for all $X$ in the ball $\mathscr{U}_{1/(2K)}(\mathcal{V})$:
\begin{itemize}
	\item  By the sub-multiplicative property and the triangle inequality, the norm of the power series is bounded by the power series of the norms:
	\[ \|\mathcal{E}(X)\| = \left\| \sum_{k=0}^{\infty} \frac{1}{k!} (\mathrm{ad}_{-X})^k \right\| \leq \sum_{k=0}^{\infty} \frac{1}{k!} \|\mathrm{ad}_{-X}\|^k = \sum_{k=0}^{\infty} \frac{\rho^k}{k!} = e^\rho. \]
	For $\rho < 1/2$, we have $e^\rho < e^{1/2}$, thus $\|\mathcal{E}(X)\| < 2$.
	\item  The Taylor coefficients $\beta_k$ satisfy $|\beta_k| \leq 1/(2\pi)^k \leq 1$. To obtain a robust upper bound, we compare the series with the dominant geometric series:
	\[ \|\mathcal{B}(X)\| = \left\| \sum_{k=0}^{\infty} \beta_k (\mathrm{ad}_{-X})^k \right\| \leq \sum_{k=0}^{\infty} |\beta_k| \rho^k \leq \sum_{k=0}^{\infty} \rho^k = \frac{1}{1-\rho}. \]
	For $\rho < 1/2$, the term $1/(1-\rho) < 2$.
\end{itemize}
\bigskip

\noindent \textbf{Estimates for the differential of the $k$-th power.}
To estimate the norm of the operator $\mathrm{d}\alpha^k(X)$, we rely on the Banach algebra structure of the space of bounded operators $\mathcal{L}(\mathcal{V})$. A key property of such an algebra is that the norm is sub-multiplicative: $\|S \circ T\| \leq \|S\| \|T\|$. 

Applying this property to each individual term in the Leibniz expansion of the $k$-th power, and using the bound $\|\mathrm{ad}_{-H}\| \leq K|H|$ from \eqref{eq:ad_norm}, we have:
\[ \| (\mathrm{ad}_{-X})^j \circ (\mathrm{ad}_{-H}) \circ (\mathrm{ad}_{-X})^{k-1-j} \| \leq \|\mathrm{ad}_{-X}\|^j \|\mathrm{ad}_{-H}\| \|\mathrm{ad}_{-X}\|^{k-1-j} \leq K|H| \|\mathrm{ad}_{-X}\|^{k-1}. \]
Then, by applying the triangle inequality, we obtain the fundamental bound:
\begin{equation} \label{eq:bound_k_power}
\| \mathrm{d}\alpha^k (X) [H] \| \leq \sum_{j=0}^{k-1} \| (\mathrm{ad}_{-X})^j \circ (\mathrm{ad}_{-H}) \circ (\mathrm{ad}_{-X})^{k-1-j} \| \leq k K |H| \|\mathrm{ad}_{-X}\|^{k-1}.
\end{equation}
\bigskip

\noindent \textbf{Bounds for the differentials $\mathrm{d}\mathcal{E}$ and $\mathrm{d}\mathcal{B}$.}
Recalling that the derivative of $\mathrm{ad}_{-X}$ is $\mathrm{ad}_{-H}$, and substituting \eqref{eq:bound_k_power} into \eqref{eq:deriv_E} and \eqref{eq:deriv_B}, we obtain:
\begin{itemize}
	\item Substituting the bound $k K |H| \rho^{k-1}$ for each term $k \geq 1$:
	\[ \| \mathrm{d}\mathcal{E}(X)[H]\| \leq \sum_{k=1}^{\infty} \frac{1}{k!} (k K |H| \rho^{k-1}) = K|H| \sum_{k=1}^{\infty} \frac{k}{k!} \rho^{k-1}. \]
	Since $k/k! = 1/(k-1)!$, we shift the index $m = k-1$ to recover the exponential series:
	\[ \| \mathrm{d}\mathcal{E}(X)[H]\| \leq K|H| \sum_{m=0}^{\infty} \frac{\rho^m}{m!} = K|H| e^\rho. \]
	For $\rho < 1/2$, we have $K|H|e^\rho < K|H|e^{1/2} < 2K|H|$.
	
	\item Using the same term-wise estimation and the property $|\beta_k| \leq 1$:
	\[ \| \mathrm{d} \mathcal{B}(X)[H]\| \leq \sum_{k=1}^{\infty} |\beta_k| (k K |H| \rho^{k-1}) \leq K|H| \sum_{k=1}^{\infty} k \rho^{k-1}. \]
	The sum $\sum_{k=1}^{\infty} k \rho^{k-1}$ is the derivative of the geometric series $\frac{d}{d\rho} (1-\rho)^{-1} = (1-\rho)^{-2}$. Therefore:
	\[ \| \mathrm{d} \mathcal{B}(X)[H]\| \leq K|H| \frac{1}{(1-\rho)^2}. \]
	For $\rho < 1/2$, the factor $1/(1-\rho)^2 < 4$, yielding the bound $4K|H|$.
\end{itemize}
\medskip 

\noindent\textbf{Invertibility and Estimates for $\mathcal{B}(X)^{-1}$.} Let us show that $\mathcal{B}(X)$ is invertible and provide quantitative bounds for its norm and its Fréchet differential.

\begin{lemma}\label{invB}
	For any fixed $X \in \mathscr{U}_{2\pi/K}(\mathcal{V})$, the operator $\mathcal{B}(X) \in \mathcal{L}(\mathcal{V})$ is invertible. Its inverse is given by the power series:
	\begin{equation} \label{eq:inverse_B_integrated}
	\mathcal{B}(X)^{-1} = \sum_{k=0}^{\infty} \frac{1}{(k+1)!} (\mathrm{ad}_{-X})^k,
	\end{equation}
	which defines an analytic map $X \mapsto \mathcal{B}(X)^{-1}$ on the whole space $\mathcal{V}$. Moreover, let $\rho = \|\mathrm{ad}_{-X}\| \leq K|X|$; for $\rho < 1/2$, the following estimates hold:
	\begin{enumerate}
		\item $\|\mathcal{B}(X)^{-1}\| \leq \frac{1}{1-\rho} < 2$;
		\item $\|\mathrm{d}(\mathcal{B}(X)^{-1})[H]\| \leq \frac{K|H|}{(1-\rho)^2} < 4K|H|$.
	\end{enumerate}
\end{lemma}

\begin{proof}
	The existence of the inverse follows from the functional identity $g(z)q(z) = 1$, where $g(z) = z/(e^z - 1)$ and $q(z) = (e^z - 1)/z$. The Taylor expansion of $q(z)$ is $q(z) = \sum_{k=0}^{\infty} \frac{1}{(k+1)!} z^k$, which represents an entire function with infinite radius of convergence.
	
	Let us derive the quantitative estimates:
	\[ \|\mathcal{B}(X)^{-1}\| \leq \sum_{k=0}^{\infty} \frac{1}{(k+1)!} \rho^k \leq \sum_{k=0}^{\infty} \rho^k = (1-\rho)^{-1}. \]
	The estimate for the Fréchet derivative is obtained by term-wise differentiation of the series \eqref{eq:inverse_B_integrated}. Using the bound $\|\mathrm{d}\alpha^k(X)[H]\| \leq k K |H| \rho^{k-1}$ previously established, we obtain:
	\[ \|\mathrm{d}(\mathcal{B}(X)^{-1})[H]\| \leq K|H| \sum_{k=1}^{\infty} \frac{k}{(k+1)!} \rho^{k-1} \leq K|H| \sum_{k=1}^{\infty} k \rho^{k-1} = K|H| (1-\rho)^{-2}. \]
	Under the assumption $\rho < 1/2$, we immediately recover the bounds $\|\mathcal{B}(X)^{-1}\| < 2$ and $\|\mathrm{d}(\mathcal{B}(X)^{-1})[H]\| < 4K|H|$.
\end{proof}

\subsubsection{The Fixed Point Theorem}

\noindent \textbf{Differentiation of the map $\phi_P$.}
The Fréchet derivative $\mathrm{d} \phi_P$ is computed by applying the Leibniz rule to the expression of $\phi_P$. Since $T$, $\Pi$, and $P$ are constant, we only differentiate $\mathcal{B}(X)$ and $\mathcal{E}(X)$ according to their position in the product:
\begin{equation} \label{eq:deriv_phi_final}
\mathrm{d} \phi_P(X)[H] = T \Pi \left[ \Big( \mathrm{d}\mathcal{B}(X)[H] \Big) \Pi (\mathcal{E}(X) P) + \mathcal{B}(X) \Pi \Big( \mathrm{d} \mathcal{E}(X)[H] P \Big) \right].
\end{equation}

\noindent \textbf{Convergence via the Banach Fixed Point Theorem.}
To prove that $\phi_P$ is a contraction, we employ the Mean Value Inequality for Banach spaces: given any two points $X_1, X_2$ in the closed ball $\overline{\mathscr{U}}_{8|P|}(\mathcal{V})$, the segment connecting them is entirely contained within the ball due to convexity. The inequality states:
\begin{equation} \label{eq:MVI}
|\phi_P(X_1) - \phi_P(X_2)| \leq \left( \sup_{X \in \overline{\mathscr{U}}_{8|P|}(\mathcal{V})} \|\mathrm{d}\phi_P(X)\| \right) |X_1 - X_2|.
\end{equation}
By finding a suitable uniform bound for $\|\mathrm{d}\phi_P(X)\|$ that holds for every $X$ in the ball, we establish that $\phi_P$ is a contraction.

\begin{thm}[Fixed Point Theorem]\label{FPT}
	Let $|P| < (24 K)^{-1}$. Then the map $\phi_P$ has a unique fixed point $X(P)$ in the ball $\overline{\mathscr{U}}_{8|P|}(\mathcal{V})$. Moreover, the fixed point is such that $X(0) = 0$ and satisfies the following asymptotic expansion for small $P$:
	\begin{equation*} 
	X(P) = T \Pi P + O(|P|^2).
	\end{equation*}
\end{thm}

\begin{proof}
	For any $X \in \overline{\mathscr{U}}_{8|P|}(\mathcal{V})$, we have $\|X\| \leq 8|P|$, which implies $\rho = \|\mathrm{ad}_{-X}\| \leq 8K|P|$. Under the smallness condition on $P$, we have $\rho < 1/3 < 1/2$. Therefore, we can apply the uniform norm bounds for $\mathcal{B}$ and $\mathcal{E}$ derived in the previous subsection:
	\[ |\phi_P(X)| \leq \|T\|_{\mathcal{L}(\mathrm{range}(\Pi), \mathcal{V})} \|\Pi\| \|\mathcal{B}(X)\| \|\Pi\| \|\mathcal{E}(X)\| |P| \leq 2 \cdot 1 \cdot 2 \cdot 1 \cdot 2 |P| = 8|P|. \]
	This confirms that $\phi_P$ maps the ball $\overline{\mathscr{U}}_{8|P|}(\mathcal{V})$ into itself.
	
	We now derive a uniform bound for the norm of the operator $\mathrm{d}\phi_P(X)$ for an arbitrary $X \in \overline{\mathscr{U}}_{8|P|}(\mathcal{V})$. Applying the triangle inequality and the sub-multiplicative property to the Leibniz expansion \eqref{eq:deriv_phi_final}, and recalling that $\|T\|_{\mathcal{L}(\mathrm{range}(\Pi), \mathcal{V})} \leq 2$ and $\|\Pi\| \leq 1$:
	\[ |\mathrm{d} \phi_P(X)[H]| \leq \|T\|_{\mathcal{L}(\mathrm{range}(\Pi), \mathcal{V})} \|\Pi\| \left( \|\mathrm{d} \mathcal{B}(X)[H]\| \cdot \|\Pi\| \cdot |\mathcal{E}(X) P| + \|\mathcal{B}(X)\| \cdot \|\Pi\| \cdot |\mathrm{d} \mathcal{E}(X)[H] P| \right). \]
	Substituting the quantitative bounds for the differentials and operators found in the previous subsection, and noting that these bounds are valid throughout the ball since $\rho < 1/3 < 1/2$:
	\begin{equation}\label{eq: bounddifferential}
	|\mathrm{d} \phi_P(X)[H]| \leq 2 \cdot 1 \cdot \left( [4K|H|] \cdot 1 \cdot (2 |P|) + 2 \cdot 1 \cdot ([2K|H|] \cdot |P|) \right) = 24K|P||H|.
	\end{equation}
	Since this bound holds for all $X \in \overline{\mathscr{U}}_{8|P|}(\mathcal{V})$, the Lipschitz constant in \eqref{eq:MVI} is bounded by $\gamma = 24K|P|$. The condition $\gamma < 1$ is satisfied by the hypothesis on $P$. Therefore, by the Banach Fixed Point Theorem, $\phi_P$ is a contraction on a complete metric space, ensuring the existence of a unique fixed point $X(P)$.
	
	To establish the asymptotic expansion, we first observe that for $P=0$, the map $\phi_0(X)$ vanishes identically for any $X$, since $P$ enters the definition as a linear factor. This immediately implies that $X(0) = \phi_0(X(0)) = 0$. 
	
	For $P \neq 0$, we expand $\phi_P$ around the origin $X=0$ using Taylor's theorem:
	\[ X(P) = \phi_P(0) + \mathcal{R}_1(X(P)), \quad \text{where} \quad \mathcal{R}_1(X(P)) = \int_0^1 \mathrm{d}\phi_P(t X(P))[X(P)] \mathrm{d}t. \]
	Evaluating $\phi_P$ at $X=0$, and noting that $\mathrm{ad}_0 = 0$ implies $\mathcal{B}(0) = \mathcal{E}(0) = \id$, we obtain the first-order term:
	\[ \phi_P(0) = T \Pi (\id \ \Pi (\id P)) = T \Pi^2 P = T \Pi P. \]
	The remainder $\mathcal{R}_1$ accounts for higher-order terms. Using the previously derived bound $\|\mathrm{d}\phi_P\| \leq 24K|P|$ and the fact that $|X(P)| \leq 8|P|$, we estimate the norm of the integral:
	\[ |\mathcal{R}_1(X(P))| \leq \sup_{t \in [0,1]} \|\mathrm{d}\phi_P(t X(P))\| |X(P)| \leq (24K|P|)(8|P|) = 192K|P|^2. \]
	Thus, $\mathcal{R}_1 = O(|P|^2)$, confirming that $X(P) = T \Pi P + O(|P|^2)$.
\end{proof}

\begin{cor}\label{corapp}
	Suppose that $T = T(\lambda)$ and $P = P(\lambda)$ depend analytically on a real parameter $\lambda$ for $\lambda$ in a closed interval $\Lambda \subset \mathbb{R}$ and satisfy the bounds $\|T(\lambda)\|_{\mathcal{L}(\mathrm{range}(\Pi), \mathcal{V})} \leq 2$, $|P(\lambda)| < (24 K)^{-1}$ for all $\lambda \in \Lambda$. Then the unique fixed point $X(P, \lambda)$ is analytic in $\lambda$. Furthermore, its derivative satisfies the bound:
	\begin{equation*}
	\left| \frac{dX}{d\lambda} \right| \leq \frac{4 \|T'\|_{\mathcal{L}(\mathrm{range}(\Pi), \mathcal{V})} |P| + 8 |P'|}{1 - 24K|P|}
	\end{equation*}
	where $T' = \frac{d}{d\lambda} T$ and $P' = \frac{d}{d\lambda} P$.
\end{cor}

\begin{proof}
	The proof starts by ensuring that the fixed point is analytic in $\lambda$. Since $T(\lambda)$ and $P(\lambda)$ are analytic, the operator $\phi_{P,\lambda}$ depends analytically on $\lambda$. Moreover, the condition $|P(\lambda)| < (24K)^{-1}$ on the closed interval $\Lambda$ ensures that $\phi_{P,\lambda}$ is always a contraction with a constant $\gamma = 24K|P| < 1$, uniformly in $\lambda$. Therefore, the Implicit Function Theorem guarantees that the fixed point $X(P,\lambda)$ is analytic in $\lambda$.
	
	Since $T$ and $P$ are analytic on the closed interval $\Lambda$, their derivatives $T'$ and $P'$ are continuous and their norms remain finite for all $\lambda \in \Lambda$. To find the bound for the derivative of the fixed point, we differentiate both sides of the identity $X(\lambda) = \phi_{P,\lambda}(X(\lambda))$ with respect to $\lambda$ using the chain rule:
	\[ \frac{dX}{d\lambda} = \frac{\partial \phi_{P, \lambda}}{\partial \lambda} + \mathrm{d}\phi_{P, \lambda}(X) \left[ \frac{dX}{d\lambda} \right]. \]
	Here, $\partial_\lambda \phi_{P, \lambda}$ represents the derivative with respect to the external parameter $\lambda$, while $\mathrm{d}\phi_{P, \lambda}$ is the Fréchet differential with respect to $X$. Moving the terms involving $dX/d\lambda$ to the left, we get:
	\[ \left( \id - \mathrm{d}\phi_{P, \lambda}(X) \right) \frac{dX}{d\lambda} = \frac{\partial \phi_{P, \lambda}}{\partial \lambda}. \]
	Since the norm of the operator $\mathrm{d}\phi_{P, \lambda}$ is bounded by $24K|P|<1$, the linear operator $\id - \mathrm{d}\phi_{P, \lambda}(X)$ is invertible. Using the usual geometric series estimate for the Neumann series that defines the inverse, we obtain the bound
	\[ \left| \frac{dX}{d\lambda} \right| \leq \frac{1}{1 - 24K|P|} \left| \frac{\partial \phi_{P, \lambda}}{\partial \lambda} \right|. \]
	The partial derivative $\partial_\lambda \phi_{P, \lambda}$ is calculated by differentiating $\phi_P(X) = T \Pi (\mathcal{B}(X) \Pi (\mathcal{E}(X) P))$ while keeping $X$ fixed:
	\[ \frac{\partial \phi_{P, \lambda}}{\partial \lambda} = T' \Pi \left( \mathcal{B}(X) \Pi (\mathcal{E}(X) P) \right) + T \Pi \left( \mathcal{B}(X) \Pi (\mathcal{E}(X) P') \right). \]
	Using the bounds $\|\Pi\| \leq 1$, $\|T\|_{\mathcal{L}(\mathrm{range}(\Pi), \mathcal{V})} \leq 2$, $\|\mathcal{B}\| < 2$, and $\|\mathcal{E}\| < 2$, we find:
	\[ \left| \frac{\partial \phi_{P, \lambda}}{\partial \lambda} \right| \leq \|T'\|_{\mathcal{L}(\mathrm{range}(\Pi), \mathcal{V})} \cdot (1 \cdot 2 \cdot 1 \cdot 2 |P|) + 2 \cdot (1 \cdot 2 \cdot 1 \cdot |P'|) = 4 \|T'\|_{\mathcal{L}(\mathrm{range}(\Pi), \mathcal{V})} |P| + 8 |P'|. \]
	Plugging this bound into the inequality for $|dX/d\lambda|$ gives the final estimate.
\end{proof}

\subsubsection{Additional Estimates for Applications}

Throughout this subsection, we assume that $T$ maps the range of $\Pi$ into itself and that there exists a left inverse $U$ such that $U \circ T = \id$ on $\mathrm{range}(\Pi)$. We do not require $U$ to be a bounded operator on $\mathcal{V}$. 

If $X$ is the unique fixed point established in Theorem \ref{FPT}, then $X \in \mathrm{range}(T)$, which ensures that the expression $UX$ is always well-defined. We can therefore characterize the fixed point through the identity:
\begin{equation*}
UX = \Pi \left( \mathcal{B}(X) \Pi \mathcal{E}(X)P \right).
\end{equation*}

\begin{lemma}[Estimates for $UX$]\label{UX}
	Let $|P| < (24K)^{-1}$ and let $X$ be the unique fixed point of $\phi_{P}$ in $\overline{\mathscr{U}}_{8|P|}$. Then, $|UX| \leq 4|P|$ and satisfies the following asymptotic expansion for small $P$:
	\begin{equation}
	UX = \Pi P + O(|P|^2).
	\end{equation}
\end{lemma}

\begin{proof}
	To prove the bound, we use the property $\| \Pi \| \leq 1$ and the bounds $\|\mathcal{B}(X)\| \leq 2$, $\|\mathcal{E}(X)\| \leq 2$, which are valid since $K|X| < 1/2$ for any $X$ in the ball $\overline{\mathscr{U}}_{8|P|}$. We have:
	\begin{equation*}
	|UX| \leq \|\Pi\|^2 \|\mathcal{B}(X)\| \|\mathcal{E}(X)\| |P| \leq 1 \cdot 2 \cdot 2 |P| = 4|P|.
	\end{equation*}
	
	For the asymptotic expansion, we observe that if $P=0$, then $X(0)=0$ and $UX$ vanishes. For $P \neq 0$, we expand $UX$ around $X=0$ using the fundamental theorem of calculus:
	\begin{equation*}
	UX = (UX)\big|_{X=0} + \mathcal{R}_1(X, P), \quad \text{where} \quad \mathcal{R}_1(X, P) = \int_0^1 \frac{d}{dt} \left[ \Pi \left( \mathcal{B}(tX) \Pi \mathcal{E}(tX)P \right) \right] \mathrm{d}t.
	\end{equation*}
	Since $\mathcal{B}(0) = \mathcal{E}(0) = \id$, the first-order term is $(UX)\big|_{X=0} = \Pi^2 P = \Pi P$. 
	By applying the Leibniz rule, the integrand becomes:
	\begin{equation*}
	\frac{d}{dt} \left[ \Pi \mathcal{B}(tX) \Pi \mathcal{E}(tX) P \right] = \Pi \left( \mathrm{d}\mathcal{B}(tX)[X] \right) \Pi \mathcal{E}(tX) P + \Pi \mathcal{B}(tX) \Pi \left( \mathrm{d}\mathcal{E}(tX)[X] \right) P.
	\end{equation*}
	Recalling the estimates from the previous section, for $X \in \overline{\mathscr{U}}_{8|P|}$ and $t \in [0,1]$, we have the bounds $\|\mathrm{d}\mathcal{B}(tX)[X]\| \leq 4K|X|$ and $\|\mathrm{d}\mathcal{E}(tX)[X]\| \leq 2K|X|$, as well as $\|\mathcal{B}(tX)\| \leq 2$ and $\|\mathcal{E}(tX)\| \leq 2$. It follows that:
	\begin{equation*}
	\begin{aligned}
	\left| \frac{d}{dt} \left[ \Pi \mathcal{B}(tX) \Pi \mathcal{E}(tX) P \right] \right| &\leq \|\Pi\|^2 \|\mathrm{d}\mathcal{B}(tX)[X]\| \cdot \|\mathcal{E}(tX)\| |P| + \|\Pi\|^2 \|\mathcal{B}(tX)\| \cdot \|\mathrm{d}\mathcal{E}(tX)[X]\| |P| \\
	&\leq 1 \cdot (4K|X|) \cdot 2 |P| + 1 \cdot 2 \cdot (2K|X| |P|) \\
	&= 8K|X||P| + 4K|X||P| = 12K|X||P|.
	\end{aligned}
	\end{equation*}
	Integrating over $t \in [0,1]$ and using the bound $|X| \leq 8|P|$ from the fixed point theorem, we obtain:
	\begin{equation*}
	|\mathcal{R}_1(X, P)| \leq \int_0^1 12K|X||P| \, \mathrm{d}t \leq 12K(8|P|)|P| = 96K|P|^2.
	\end{equation*}
	This confirms the expansion $UX = \Pi P + O(|P|^2)$.
\end{proof}

We conclude with the following technical lemma:

\begin{lemma}[Estimates for $\Delta$]\label{Deltaestimates}
	Let $|P| < (24K)^{-1}$ and $X$ be the unique fixed point. Under the assumption $\|\id - \Pi\| \leq 1$, the quantity 
	\begin{equation*}
	\Delta := (\id - \Pi) \left( (\mathcal{E}(X) - \id)P - (\mathcal{B}(X)^{-1} - \id)UX \right) 
	\end{equation*}
	satisfies $|\Delta| \leq 64 K |P|^2$. Furthermore, $\Delta$ satisfies the following asymptotic expansion for small $P$:
	\begin{equation*}
	\Delta = (\id - \Pi) \left( [(\id - \Pi)P, T \Pi P] + \frac{1}{2} [\Pi P, T \Pi P] \right) + O(|P|^3).
	\end{equation*}
\end{lemma}

\begin{proof}
	First, we substitute the fixed point identity $UX = \Pi \left( \mathcal{B}(X) \Pi \mathcal{E}(X)P \right)$ into the definition of $\Delta$:
	\begin{equation*}
	\Delta = (\id - \Pi) (\mathcal{E}(X) - \id)P - (\id - \Pi) (\mathcal{B}(X)^{-1} - \id) \Pi \mathcal{B}(X) \Pi \mathcal{E}(X)P.
	\end{equation*}
	We now use the triangle inequality to split the norm into two pieces:
	\begin{equation*}
	|\Delta| \leq |\underbrace{(\id - \Pi)(\mathcal{E}(X) - \id)P}_{\Delta_1}| + |\underbrace{(\id - \Pi) (\mathcal{B}(X)^{-1} - \id) \Pi \mathcal{B}(X) \Pi \mathcal{E}(X)P}_{\Delta_2}|.
	\end{equation*}
	Let $\rho = \|\mathrm{ad}_{-X}\| \leq 8K|P|$. Since $|P| < (24K)^{-1}$, we have $\rho < 1/3$ and, as established in Lemma \ref{invB}, the bound $\|(\mathcal{B}(X))^{-1}\| \leq \frac{1}{1-\rho} < 2$ holds.
	
	\medskip
	\noindent \textbf{Term $\Delta_1$.} The operator $\mathcal{E}(X) - \id$ is the series $\sum_{k=1}^\infty \frac{1}{k!}(\mathrm{ad}_{-X})^k$. Taking the norm and using the bound $\|\id - \Pi\| \leq 1$:
	\begin{equation*}
	|\Delta_1| \leq 1 \cdot \left( \sum_{k=1}^\infty \rho^k \right) |P| = \frac{\rho}{1-\rho} |P|.
	\end{equation*}
	Substituting $\rho \leq 8K|P| $ and $\frac{1}{1-\rho} < 2$, we get:
	\begin{equation*}
	|\Delta_1| \leq 2 \cdot (8K|P|) \cdot |P| = 16 K |P|^2.
	\end{equation*}
	
	\medskip
	\noindent \textbf{Term $\Delta_2$.} To estimate this piece, we use the property from Lemma \ref{invB} where $\|(\mathcal{B}(X))^{-1} - \id\| \leq \frac{\rho}{1-\rho}$. Combining this with $\|\mathcal{B}(X)\| \leq 2$, and $|\mathcal{E}(X)P| \leq 2|P|$:
	\begin{equation*}
	|\Delta_2| \leq \frac{\rho}{1-\rho} \cdot 1 \cdot 2 \cdot 1 \cdot 2|P| = 4 \frac{\rho}{1-\rho} |P| \leq 48 K |P|^2.
	\end{equation*}
	Summing the estimates yields $|\Delta| \leq 64 K |P|^2$.
	
	\medskip
	\noindent \textbf{Asymptotic Expansion.} To derive the leading order for small $P$, we perform a Taylor expansion of $\Delta$ around $X=0$. Since $\Delta(0)=0$, we have $\Delta(X) = \mathrm{d}\Delta(0)[X] + \mathcal{R}_2(X)$. We differentiate the two components $\Delta_1$ and $\Delta_2$ separately:
	\begin{align*}
	\Delta_1(X) &= (\id - \Pi) (\mathcal{E}(X) - \id) P, \\
	\Delta_2(X) &= (\id - \Pi) (\mathcal{B}(X)^{-1} - \id) \Pi \mathcal{B}(X) \Pi \mathcal{E}(X) P.
	\end{align*}
	From Section \ref{EandB}, we have $\mathrm{d}\mathcal{E}(0)[H] = \mathrm{ad}_{-H}$. Then,
	
	\medskip
	\noindent \textit{Differential of $\Delta_1$:}
	Differentiating $\Delta_1$ at $X=0$ in the direction $H$:
	\begin{equation*}
	\mathrm{d}\Delta_1(0)[H] = (\id - \Pi) \left( \mathrm{d}\mathcal{E}(0)[H] \right) P = (\id - \Pi) \mathrm{ad}_{-H} P = (\id - \Pi) [P, H].
	\end{equation*}
	
	\medskip
	\noindent \textit{Differential of $\Delta_2$:}
	We apply the Leibniz rule to the product of operators in $\Delta_2$ and evaluate at $X=0$. Recalling that $\mathcal{B}(0) = \mathcal{E}(0) = \id$, we compute the differential. Since $(\mathcal{B}(0)^{-1} - \id) = 0$, only the term where the derivative acts on the first factor remains non-zero:
	\begin{align*}
	\mathrm{d}\Delta_2(0)[H] &= (\id - \Pi) \left[ \left(\mathrm{d}(\mathcal{B}(0)^{-1})[H]\right) \Pi \mathcal{B}(0) \Pi \mathcal{E}(0) P \right] \\
	&= (\id - \Pi) \left(\mathrm{d}(\mathcal{B}(0)^{-1})[H]\right) \Pi P.
	\end{align*}
	Now, from the analytic expansion in Lemma \ref{invB}, we substitute $\mathrm{d}(\mathcal{B}(0)^{-1})[H] = \tfrac{1}{2} \mathrm{ad}_{-H}$ to obtain
	\begin{equation*}
	\mathrm{d}\Delta_2(0)[H] = (\id - \Pi) \left( \tfrac{1}{2} \mathrm{ad}_{-H} \Pi P \right)=\tfrac{1}{2} (\id - \Pi) [\Pi P, H].
	\end{equation*}
	
	\medskip
	\noindent \textit{Differential of $\Delta$:}
	Combining the differentials of $\Delta_1$ and $\Delta_2$, we have $\mathrm{d}\Delta(0)[H] = \mathrm{d}\Delta_1(0)[H] - \mathrm{d}\Delta_2(0)[H]$. Substituting $X = T \Pi P + O(|P|^2)$ into the linear expressions:
	\begin{align*}
	\mathrm{d}\Delta(0)[X] = (\id - \Pi) \left( \frac{1}{2} [\Pi P, T \Pi P] + [(\id - \Pi)P, T \Pi P] \right) + O(|P|^3).
	\end{align*}
	
	\medskip
	\noindent \textit{Estimate for the reminder:}
	The remainder $\mathcal{R}_2(X)$ is defined by the integral formula:
	\begin{equation*}
	\mathcal{R}_2(X) = \int_0^1 (1-t) \mathrm{d}^2 \Delta(tX)[X,X] \mathrm{d}t.
	\end{equation*}
	The second Fréchet derivative $\mathrm{d}^2 \Delta(X)[X,X]$ is a sum of terms involving double commutators with $P$, of the form $[X, [X, P]]$, multiplied by operators $\mathcal{E}, \mathcal{B}, \mathcal{B}(X)^{-1}$ and their derivatives. 
	Given the algebra property and the fact that these operators are analytic and bounded for $\rho < 1/3$, the norm of the second derivative satisfies:
	\begin{equation*}
	\| \mathrm{d}^2 \Delta(X)[X,X] \| \leq C K^2 |X|^2 |P|,
	\end{equation*}
	where $C$ is a constant. Substituting this bound into the integral:
	\begin{align*}
	|\mathcal{R}_2(X)| &\leq \int_0^1 (1-t) \| \mathrm{d}^2 \Delta(tX)[X,X] \| \mathrm{d}t \\
	&\leq \left( \int_0^1 (1-t) \mathrm{d}t \right) \cdot \sup_{t \in [0,1]} \left( C K^2 |X|^2 |P| \right) \\
	&= \frac{1}{2} \cdot C K^2 |X|^2 |P|.
	\end{align*}
	Using the bound for the fixed point $|X| \leq 8|P|$ we get
	\begin{equation*}
	|\mathcal{R}_2(X)| \leq \frac{1}{2} \cdot C K^2 (8 |P|)^2 |P| = 32 C K^2 |P|^3.
	\end{equation*}
	This shows that $|\mathcal{R}_2(X)| \leq M |P|^3$ for a constant $M$, which confirms $\mathcal{R}_2 = O(|P|^3)$.
\end{proof}

\begin{rmk}
	If the operators $T$ and $P$ depend analytically on a real parameter $\lambda \in \Lambda$ (where $\Lambda$ is a closed interval of the real line), the fixed point $X(\lambda)$ inherits the analyticity, as proved in Corollary \ref{corapp}. Consequently, $UX$ and $\Delta$ are analytic in $\lambda$ since they are obtained by composing analytic maps (the operators $\mathcal{B}$ and $\mathcal{E}$) with the analytic functions $X(\lambda)$ and $P(\lambda)$. Moreover, the constants appearing in the estimates depend only on the uniform bounds of the operators for $\lambda \in \Lambda$, ensuring that all the asymptotic expansions hold uniformly in the parameter $\lambda$.
\end{rmk}

\section{Technical Lemmas for the Analysis of the Spectral Gaps}\label{AAE}
To ensure a self-contained treatment, we detail here the functional framework and provide the proofs for several auxiliary results. These technical lemmas are fundamental to the study of the spectral gaps developed in Subsection \ref{SG}.
\subsubsection{Functional setting}\label{FS}

Let $p > 1$ and $s \ge 0$. For any $\ell \in \mathbb{Z}$, we denote the standard weight as $\langle \ell \rangle := \max\{1, |\ell|\}$.

\paraga \textbf{The spaces $\ell^\infty_{s,p}(\mathscr{B})$.} Let $(\mathscr{B},| \cdot |_{\mathscr{B}}, \cdot)$ be a Banach algebra with constant $C_{\mathrm{alg}}(\mathscr{B})\in\R^+$, that is, for all $x,y \in \mathscr{B}$, one has
\begin{equation*}
	|x\cdot y|_{\mathscr{B}}\leq C_{\mathrm{alg}}(\mathscr{B})|x|_\mathscr{B}|y|_\mathscr{B}.
\end{equation*}
We define $\ell^\infty_{s,p}(\mathscr{B})$ as the Banach space of sequences $X = (X_\ell)_{\ell \in \mathbb{Z}}$, $X_\ell \in \mathscr{B}$ for all $\ell \in \Z$, equipped with the norm:
\begin{equation*}
\|X\|^{\mathscr{B}}_{s,p} := \sup_{\ell \in \mathbb{Z}} \langle \ell \rangle^p e^{s|\ell|} |X_\ell|_\mathscr{B} < \infty.
\end{equation*}
This space serves as the space of Fourier coefficients for the classes of periodic functions considered herein.

The convolution of two sequences $X,Y \in \ell^\infty_{s,p}(\mathscr{B})$ is the sequence $Z=X\star Y$ defined componentwise as
\begin{equation}\label{convprod}
	Z_\ell=\sum_{k \in \Z}X_{\ell-k}Y_k.
\end{equation}
The definition is well posed because
\begin{align}
	\sup_{\ell \in \mathbb{Z}} \langle \ell \rangle^p e^{s|\ell|} |Z_\ell|&\leq C_{\mathrm{alg}}(\mathscr{B}) \sup_{\ell \in \mathbb{Z}} \langle \ell \rangle^p e^{s|\ell|}\sum_{k \in \Z}|X_{\ell-k}|_{\mathscr{B}}|Y_k|_\mathscr{B} \notag
	\\
	&\leq C_{\mathrm{alg}}(\mathscr{B}) \|X\|^{\mathscr{B}}_{s,p}\|Y\|^{\mathscr{B}}_{s,p}\sup_{\ell \in \Z}\sum_{k \in \Z}\frac{\jap{\ell}^p}{\jap{\ell-k}^p\jap{k}^p}e^{s(|\ell|-|\ell-k|-|k|)} \notag
	\\
	&\leq C_{\mathrm{alg}}(\mathscr{B}) \|X\|^{\mathscr{B}}_{s,p}\|Y\|^{\mathscr{B}}_{s,p}\sup_{\ell \in \Z}\sum_{k \in \Z}\frac{\jap{\ell}^p}{\jap{\ell-k}^p\jap{k}^p} \notag
	\\
	&\leq C_\star(p)C_{\mathrm{alg}}(\mathscr{B}) \|X\|^{\mathscr{B}}_{s,p}\|Y\|^{\mathscr{B}}_{s,p}, \label{cpcalg}
\end{align}
where
\begin{equation}\label{C(p)}
	C_\star(p):=2^p(1+2\zeta(p))
\end{equation}
denotes the convolution algebra constant of $\ell^\infty_{s,p}$ and $\zeta(p) = \sum_{n \geq 1} n^{-p}$ is the \emph{Riemann function}.

Concretely, in the following discussion we shall consider the cases $\mathscr{B}=\C$, where $C_{\mathrm{alg}}(\C)=1$, and $\mathscr{B}=\text{Mat}(2,\C)$, endowed with the norm
\begin{equation*}
		\normS{W}:=\max_{\sigma,\sigma'=\pm1}|W^{\sigma\sigma'}|,
\end{equation*}
for which $C_{\mathrm{alg}}(\text{Mat}(2,\C))=2$.

For the case $\mathscr{B}=\C$, we adopt the shortened notation $\ell^\infty_{s,p}(\C)=\ell^{\infty}_{s,p}$.

\paraga \textbf{The algebra $\mathscr{S}_{s,p} \ $.} We define $\mathscr{S}_{s,p}$ as the Banach space of $2\pi$-periodic functions $S: \R \to \mathfrak{su}(1,1)$, $S(t) = \sum_{\ell \in \mathbb{Z}} S_\ell e^{\im \ell t}$, with $S_\ell=(S_\ell^{\sigma \sigma'})_{\sigma,\sigma'=\pm 1} \in \mathfrak{sl}(2,\C)$, satisfying
\begin{equation*}
|S|_{s,p} := \sup_{\ell \in \mathbb{Z}} \langle \ell \rangle^p e^{s|\ell|}\max_{\sigma,\sigma'=\pm 1}|S^{\sigma\sigma'}_\ell|.
\end{equation*}
It is clear that the list of matrices $(S_\ell)_{\ell\in\Z}$ must necessarily satisfy the constrains
\begin{equation}
	S^{--}_\ell=-S_{\ell}^{++}, \quad \text{and} \quad \overline{S^{-+}_{-\ell}}=S^{+-}_\ell, \qquad \text{for all} \ \ell \in \Z,
\end{equation}
in order to define a $\mathfrak{su}(1,1)$ valued function via the Fourier series.

\paraga \textbf{The algebra $\hat{\mathscr{S}}_{s,p} \ $.}
We denote by $\hat{\mathscr{S}}_{s,p}$ the subspace of $\ell^\infty_{s,p}(\text{Mat}(2,\mathbb{C}))$ consisting of matrix-valued sequences $X = (X_\ell)_{\ell \in \mathbb{Z}}$ such that, for every $\ell \in \mathbb{Z}$, the matrix $X_\ell$ satisfies:
\begin{enumerate}
	\item $X^{--}_\ell = -X^{++}_\ell$ (Trace-free condition);
	\item $X^{++}_\ell = -\overline{X^{++}_{-\ell}}$ (Diagonal antisymmetry);
	\item $X^{-+}_\ell = \ \ \ \overline{X^{+-}_{-\ell}}$ (Off-diagonal conjugation).
\end{enumerate}
These relations ensure that $\hat{\mathscr{S}}_{s,p}$ is the space of Fourier coefficients for maps taking values in the Lie algebra $\mathfrak{su}(1,1)$. The following lemma implies that the Fourier transform defines an isomorphism of Banach-Lie algebras between $\hat{\mathscr{S}}_{s,p}$ and $\mathscr{S}_{s,p}$.

\begin{lemma}\label{FSsp}
	The space $\hat{\mathscr{S}}_{s,p}$, equipped with the norm 
	\begin{equation*}
	\normS{X}_{s,p} :=  \sup_{\ell \in \mathbb{Z}} \langle \ell \rangle^p e^{s|\ell|} \normS{X_\ell}= \|X\|^{\text{Mat}(2,\mathbb{C})}_{s,p},
	\end{equation*}
	is a Banach-Lie algebra with respect to the Lie bracket $[X, Y]_{\star} = X \star Y - Y \star X$. In particular, the bracket is continuous and satisfies the estimate:
	\begin{equation}\label{liebound}
	\normS{[X, Y]_\star}_{s,p} \leq 4 C_\star(p) \normS{X}_{s,p} \normS{Y}_{s,p},
	\end{equation}
	where $C_\star(p)$ the constant defined in \eqref{C(p)}.
\end{lemma}

\begin{proof}
	The space $\hat{\mathscr{S}}_{s,p}$ is a closed subspace of the Banach space $\ell^\infty_{s,p}(\text{Mat}(2,\mathbb{C}))$. To show it is a Banach-Lie algebra, we verify the algebraic closure and the continuity of the bracket. 
	\smallskip
	
	\noindent \textit{Algebraic closure:} For any $X, Y \in \hat{\mathscr{S}}_{s,p}$, let $Z = [X, Y]_{\star}$. The condition $\mathrm{Tr}(Z_\ell) = 0$ follows from the linearity of the trace and the identity $\mathrm{Tr}([A,B])=0$. Regarding the diagonal antisymmetry, a direct computation on the $(++)$ component of the convolution yields:
	\begin{align*}
	-\overline{(X \star Y)^{++}_{-\ell}} &= -\sum_{k \in \mathbb{Z}} \left( \overline{X^{++}_{-\ell-k}} \overline{Y^{++}_k} + \overline{X^{+-}_{-\ell-k}} \overline{Y^{-+}_k} \right) \\
	&= -\sum_{k \in \mathbb{Z}} \left( (-X^{++}_{\ell+k})(-Y^{++}_{-k}) + X^{-+}_{\ell+k} Y^{+-}_{-k} \right) = -(Y \star X)^{++}_\ell.
	\end{align*}
	It follows that $-\overline{Z^{++}_{-\ell}} = -(Y \star X)^{++}_\ell + (X \star Y)^{++}_\ell = Z^{++}_\ell$. Similarly, evaluating the conjugation of the $(+-)$ component:
	\begin{align*}
	\overline{(X \star Y)^{+-}_{-\ell}} &= \sum_{k \in \mathbb{Z}} \left( \overline{X^{++}_{-\ell-k}} \overline{Y^{+-}_k} + \overline{X^{+-}_{-\ell-k}} \overline{Y^{--}_k} \right) \\
	&= \sum_{j \in \mathbb{Z}} \left( X^{++}_{\ell-j} Y^{-+}_j + X^{-+}_{\ell-j} Y^{--}_j \right) = (Y \star X)^{-+}_\ell.
	\end{align*}
	Thus $\overline{Z^{+-}_{-\ell}} = (Y \star X)^{-+}_\ell - (X \star Y)^{-+}_\ell = -Z^{-+}_\ell$, confirming that $Z^{-+}_\ell = \overline{Z^{+-}_{-\ell}}$.
	
	\noindent \textit{Continuity of the bracket:} The sub-multiplicativity of the convolution product established in \eqref{cpcalg} implies:
	\begin{align*}
	\normS{[X, Y]_\star}_{s,p} &\leq \normS{X \star Y}_{s,p} + \normS{Y \star X}_{s,p} \\
	&\leq 2 C_\star(p) C_{\mathrm{alg}}(\text{Mat}(2,\mathbb{C})) \normS{X}_{s,p} \normS{Y}_{s,p}.
	\end{align*}
	Substituting $C_{\mathrm{alg}}(\text{Mat}(2,\mathbb{C})) = 2$, we obtain the bound \eqref{liebound}. The boundedness of the bracket ensures that $\hat{\mathscr{S}}_{s,p}$ is a Banach-Lie algebra.
\end{proof}

\paraga \textbf{The functional space $\mathcal{H}_{s,p,R}$.}
Given $R>0$, we denote by $\mathcal{H}_{s,p,R}$ the Banach space of functions
\begin{equation*}
f:\mathscr{U}_R(\ell^\infty_{s,p})\to \mathbb{C}
\end{equation*}
of the form $f(u)=\sum_{\alpha \in \mathbb{N}_0^\mathbb{Z}}f_\alpha u^\alpha$, with $f_\alpha \in \mathbb{C}$, satisfying
\begin{equation}\label{normpsR}
\|f\|_{s,p,R}=\sum_{\alpha \in \mathbb{N}^\mathbb{Z}}|f_\alpha|v^\alpha_{s,p,R} < \infty,
\end{equation}
where $v_{s,p,R}:=(\langle \ell \rangle^{-p}e^{-s|\ell|}R)_{\ell \in \mathbb{Z}}$. Note that, if one introduces the \emph{majorant} $\underline{f}$ of $f$ as the function defined by
\begin{equation}
\underline{f}(u)=\sum_{\alpha \in \mathbb{N}_0^\mathbb{Z}}|f_\alpha|u^\alpha,
\end{equation}
the condition \eqref{normpsR} can be rewritten as
\begin{equation*}
\|f\|_{s,p,R}=\sup_{\|u\|_{s,p}\leq R}|\underline{f}(u)|<\infty.
\end{equation*}
Crucially, the majorant of a product of functions precedes the product of the respective majorants., i.e., $\underline{fg} \preceq \underline{f} \underline{g}$. This notation means that the inequality holds coefficient-wise for the corresponding formal power series. Consequently, for any $u$ with non-negative components (such as the weight vector $v_{s,p,R}$), one has $\underline{fg}(u) \leq \underline{f}(u)\underline{g}(u)$. This property ensures that $\mathcal{H}_{s,p,R}$ is a Banach algebra, namely 
\begin{equation*}
\|fg\|_{s,p,R} \leq \|f\|_{s,p,R}\|g\|_{s,p,R}.
\end{equation*}
It is worth stressing that any $f \in \mathcal{H}_{s,p,R}$ is majorant analytic, as its majorant $\underline{f}$ is the uniform limit of continuous polynomials.

\paraga \textbf{The Banach-Lie Algebra $\underline{\mathcal{A}}_{R}^{\pi}(\ell^\infty_{s,p},\hat{\mathscr{S}}_{s,p})$.} Given $R>0$, we define the space $\underline{\mathcal{A}}_{R}^{\pi}(\ell^\infty_{s,p},\hat{\mathscr{S}}_{s,p})$ as the space of analytic maps $X : \mathscr{U}_R(\ell^\infty_{s,p}) \to \hat{\mathscr{S}}_{s,p}$. An element $X$ in this space is identified with a sequence of matrix-valued power series $X(u) = (X_\ell(u))_{\ell \in \mathbb{Z}}$, where each component $X_\ell(u)$ satisfies a \textit{momentum conservation condition}. Specifically, the Taylor expansion of the entries $X_{\ell}^{\sigma\sigma'} \in \mathcal{H}_{s,p,R}$ involves only multi-indices $\alpha \in \mathbb{N}^\mathbb{Z}$ whose \emph{momentum} $\pi(\alpha)$ matches the Fourier index $\ell$:
\begin{equation*}
X_{\ell}^{\sigma\sigma'}(u) = \sum_{\alpha \in \mathscr{M}_\ell} X_{\ell,\alpha}^{\sigma\sigma'} u^\alpha, \quad \text{with} \quad \mathscr{M}_\ell := \bigg\{ \alpha \in \mathbb{N}_0^\mathbb{Z} : \pi(\alpha) := \sum_{m \in \mathbb{Z}} m \alpha_m = \ell \bigg\}.
\end{equation*}
The momentum conservation condition ensures that each Fourier component $X_\ell$ of the map transforms, with respect to a time shift $t$, with the same frequency as the phase factor $e^{\im\ell t}$. This implies that the map depends on $t$ only through $u$, as shown by the following chain of equalities:
\begin{equation*}
e^{\im \ell t} X_\ell(u) = \sum_{\alpha \in \mathscr{M}_\ell} X_{\ell,\alpha} u^\alpha e^{\im \pi(\alpha) t} = \sum_{\alpha \in \mathscr{M}_\ell} X_{\ell,\alpha} \prod_m (u_m e^{\im mt})^{\alpha_m} = X_\ell(u(t)).
\end{equation*}

To equip this space with a Banach structure, we introduce the \textit{majorant} $\underline{X}$ as the sequence of maps $\underline{X} = (\underline{X}_\ell)_{\ell \in \mathbb{Z}}$ obtained by replacing each entry with its majorant in $\mathcal{H}_{s,p,R}$, namely $\underline{X}_{\ell}^{\sigma\sigma'}(u) := \sum_{\alpha \in \mathscr{M}_\ell} |X_{\ell,\alpha}^{\sigma\sigma'}| u^\alpha$. The majorant $\underline{X}$ allows us to define the norm $|X|_{s,p,R}$ by considering the supremum of the $\hat{\mathscr{S}}_{s,p}$-norm over the ball of radius $R$:
\begin{equation}\label{majnorm_def}
|X|_{s,p,R} := \sup_{\|u\|_{s,p} \leq R} \normS{\underline{X}(u)}_{s,p} = \sup_{\ell \in \mathbb{Z}} \langle \ell \rangle^p e^{s|\ell|} \max_{\sigma,\sigma'=\pm1} \| X_{\ell}^{\sigma\sigma'} \|_{s,p,R},
\end{equation}
which, due to the non-negativity of the coefficients, coincides with evaluating the majorant on the weight vector $v_{s,p,R}$.

It should be observed that any $X \in \underline{\mathcal{A}}_{R}^{\pi}$ is majorant analytic. In fact, Theorem \ref{uniformanalytic} ensures that $\underline{X}$ is analytic as it is defined by the uniformly bounded family of (majorant) analytic functions $ X_\ell^{\sigma\sigma'} \in \mathcal{H}_{s,p,R}$.

\begin{lemma}\label{banachliealgA}
	The space $(\underline{\mathcal{A}}_{R}^{\pi}, |\cdot|_{s,p,R})$ is a Banach-Lie algebra with respect to the pointwise Lie bracket $[X, Y]_\star(u) = X(u) \star Y(u) - Y(u) \star X(u)$, satisfying:
	\begin{equation}
	|[X, Y]_\star|_{s,p,R} \leq K |X|_{s,p,R} |Y|_{s,p,R},
	\end{equation}
	with $K = 4 C_\star(p)$ (see \eqref{C(p)}).
\end{lemma}

\begin{proof}
The momentum condition is preserved by the convolution product: if $\pi(\alpha) = \ell-k$ and $\pi(\beta) = k$, then $\pi(\alpha+\beta) = \ell$, and hence $[X, Y]_\star(u) \in \hat{\mathscr{S}}_{s,p}$ by Lemma \ref{FSsp}.
To prove the continuity, we recall that for any $f, g \in \mathcal{H}_{s,p,R}$ one has $\underline{fg} \preceq \underline{f} \underline{g}$, which implies that $\underline{X \star Y} \preceq \underline{X} \star \underline{Y}$ in the sense of power series. Consequently, for any $u$ with non-negative components (and thus for the supremum over the ball), it follows that:
\begin{align*}
|[X, Y]_\star|_{s,p,R} &= \sup_{\|u\|_{s,p} \leq R} \normS{\underline{[X, Y]_\star}(u)}_{s,p} \\
&\leq \sup_{\|u\|_{s,p} \leq R} \left( \normS{(\underline{X} \star \underline{Y})(u)}_{s,p} + \normS{(\underline{Y} \star \underline{X})(u)}_{s,p} \right) \\
&\leq 4 C_\star(p) \sup_{\|u\|_{s,p} \leq R} \normS{\underline{X}(u)}_{s,p} \sup_{\|u\|_{s,p} \leq R} \normS{\underline{Y}(u)}_{s,p},
\end{align*}
where we used the algebra property of the $\hat{\mathscr{S}}_{s,p}$ norm established in Lemma \ref{FSsp} and the fact that the majorant of a product precedes the product of the majorants.
\end{proof}

\subsubsection{Reduction to constant coefficients}\label{sec:reduction}
In this subsection, we establish the reducibility of the linear Hamiltonian system
\begin{equation}\label{questo}
\frac{d}{dt} \begin{pmatrix} z \\ \bar{z} \end{pmatrix} = \left[ \lambda \begin{pmatrix} i & 0 \\ 0 & -i \end{pmatrix} + \frac{u}{2\lambda} \begin{pmatrix} -i & i \\ -i & i \end{pmatrix} \right] \begin{pmatrix} z \\ \bar{z} \end{pmatrix},
\end{equation}
where $u=u(t)\in\R$ is a zero-mean, analytic, $2\pi$-periodic function whose sequence of Fourier coefficients is small in the $\ell^\infty_{s,p}$ norm. We provide quantitative estimates for the reducing transformation; these results will be crucial in the next subsection to investigate the regularity of the spectral gaps with respect to the potential $u$. It is worth stressing that, in this dynamical setting, such gaps correspond precisely to the instability intervals of the periodic forcing. For any fixed $n \in \mathbb{N}$, the results hold for all spectral parameters $\lambda \in \mathbb{R}^+$ satisfying $|2\lambda - n| < \gamma \leq 1/2$. Accordingly, throughout the following discussion, we fix a natural number $n$ and we assume that $\lambda$ remains restricted to this region. 

Following the strategy outlined in Subsection \ref{subfloquet}, we aim to find a $2\pi$-periodic change of coordinates that reduces the Hamiltonian system \eqref{questo} to constant coefficients. To preserve the algebraic structure, this transformation must take values in the Lie group $SU(1,1)$. 

Given the identification between the space of $2\pi$-periodic functions $\mathscr{S}_{s,p}$ and its Fourier coefficients space $\hat{\mathscr{S}}_{s,p}$ established in Subsection \ref{FS} (see Lemma \ref{FSsp}), we seek the generator of this transformation within the functional space $\underline{\mathcal{A}}_{R}^{\pi}(\ell^\infty_{s,p},\hat{\mathscr{S}}_{s,p})$. This choice allows us to track the regularity of the transformation with respect to the potential $u$ and eventually derive estimates for the spectral gaps. In this regard, it is worth noting that the perturbation
\begin{equation}
P(u,\lambda)
=
\frac{u}{2\lambda}
\begin{pmatrix}
-\im & \im \\ -\im & \im
\end{pmatrix}
=
\frac{u}{2\lambda}(J_3-J_1).
\end{equation}
indeed belongs to $\underline{\mathcal{A}}_{R}^{\pi}(\ell^\infty_{s,p},\hat{\mathscr{S}}_{s,p})$, with norm $|P|_{s,p,R} = \frac{R}{2\lambda}$.

Given a map $A \in \underline{\mathcal{A}}_{R}^{\pi}$, we consider a change of variables $(\xi, \bar{\xi})^T = e^{-A(u,\lambda)} (z, \bar{z})^T$. Here, the exponential $e^{-A}$ and the transformation law are to be understood by identifying the sequence $A(u,\lambda) \in \hat{\mathscr{S}}_{s,p}$ with its corresponding matrix-valued periodic function in $\mathscr{S}_{s,p}$ via the Fourier isomorphism. With this identification, the time derivative $\dot{A}$ acts as a diagonal operator on the Fourier coefficients, namely $(\dot{A})_\ell = i\ell A_\ell$. Under this transformation, the matrix $D(u, \lambda)$ associated with the linear vector field inducing the Hamiltonian dynamics defined in \eqref{questo} is transformed into
\begin{equation}\label{D_1_transf}
\mathcal{D}(u, \lambda) = e^{-A(u,\lambda)} D(u, \lambda) e^{A(u,\lambda)} - e^{-A(u,\lambda)} \frac{d}{dt} e^{A(u,\lambda)}.
\end{equation}
By defining the adjoint action as $\mathrm{ad}_{A}(M) := [A, M]_\star$, where the Lie bracket $[\cdot, \cdot]_\star$ is defined in Lemma \ref{banachliealgA}, and recalling the standard expansion for the adjoint action and the derivative of the exponential map:
\begin{equation*}
e^{-A}Me^{A}=\sum_{k=0}^\infty \frac{(\mathrm{ad}_{-A})^k}{k!}M, \qquad e^{-A}\frac{d}{dt}e^{A}=\sum_{k=1}^\infty\frac{(\mathrm{ad}_{-A})^{k-1}}{k!}\dot{A},
\end{equation*}
the equation \eqref{D_1_transf} can be expressed as:
\begin{equation}\label{D_1_series}
\mathcal{D}(\lambda) = \sum_{k=0}^\infty \frac{(\mathrm{ad}_{-A})^k}{k!}(\lambda J_1 +P(\lambda)) - \sum_{k=1}^\infty \frac{(\mathrm{ad}_{-A})^{k-1}}{k!} \dot{A}.
\end{equation}
Thanks to the Banach-Lie algebra structure of $\underline{\mathcal{A}}_{R}^{\pi}$, the series in \eqref{D_1_series} is well-defined and convergent, provided that the norm $|A|_{s,p,R}$ is sufficiently small (see Appendix \ref{AAD} for details).

We remark that we restrict our search to generators $A$ belonging to the space $\underline{\mathcal{A}}_{R}^{\pi}$, thereby enforcing the momentum conservation condition. This assumption is not strictly necessary for the reduction to constant coefficients. However, since the initial perturbation trivially satisfies this property and the Lie bracket preserves it, the property remains valid at each step of the contractive iteration we shall implement to construct the solution. While not required for the reducibility itself, this selection rule ensures that the transformed system inherits the same symmetry, which will be crucial in the subsequent sections when studying the derivatives of the spectral gaps.

As already noted in Subsection \ref{subfloquet}, the reduction to constant coefficients for the component lying in the range of $\mathcal{L}_{n/2}$ can be recast as a fixed-point problem for the operator $\phi_P$, defined as
\begin{equation}
\phi_P(A) = \mathcal{L}_\lambda^{-1} \Pi_{rg} \left( \mathcal{B}(A) \Pi_{rg} \left( e^{\mathrm{ad}_{-A}} P \right) \right),
\end{equation}
where $\Pi_{rg}$ denotes the projector onto the range of the operator $\mathcal{L}_{n/2}$. Here, $\mathcal{B}(A)$ represents the inverse of the map $X \mapsto \sum_{k=1}^\infty \frac{1}{k!}(\mathrm{ad}_{-A})^{k-1} X$ (the details of which are provided in Lemma \ref{invB}), while $\mathcal{L}^{-1}_\lambda$ is the inverse of $\mathcal{L}_\lambda$, which is well defined on $\mathrm{range}(\mathcal{L}_{n/2})$.

Before stating and proving the fixed-point theorem, it is of crucial importance to make some comments on the operator $\mathcal{L}^{-1}_\lambda$: since $\mathcal{L}_\lambda$ acts diagonally on the components of the generator $A$, namely
\begin{equation}
(\mathcal{L}_{\lambda}A)^{\sigma\sigma'}_\ell=\im(\ell-\lambda(\sigma-\sigma'))A^{\sigma\sigma'}, \quad \text{for all} \ \sigma,\sigma'=\pm1, \ \ell \in \Z,
\end{equation}
$\mathcal{L}^{-1}_\lambda$ naturally preserves both the momentum conservation condition and the range of the operator $\mathcal{L}_{n/2}$. Furthermore, we have the following lemma:

\begin{lemma}\label{boundLinv}
	The inverse $\mathcal{L}^{-1}_\lambda$, viewed as an operator from $\Pi_{rg}(\underline{\mathcal{A}}_{R}^{\pi})$ into itself, is bounded. Specifically, for any $X \in \Pi_{rg}(\underline{\mathcal{A}}_{R}^{\pi})$, we have the estimate:
	\begin{equation}
	|\mathcal{L}^{-1}_\lambda X|_{s,p,R} \leq 2 |X|_{s,p,R}.
	\end{equation}
\end{lemma}

\begin{proof}
	Recall that for $X \in \Pi_{rg}(\underline{\mathcal{A}}_{R}^{\pi})$, the action of $\mathcal{L}^{-1}_\lambda$ on the Fourier components $X_\ell^{\sigma\sigma'}$ is given by the division by the factor $\im(\ell - (\sigma-\sigma')\lambda)$. Consequently, the majorant of the transformed component satisfies:
	\begin{equation*}
	\underline{\left( \frac{1}{\im(\ell - (\sigma-\sigma')\lambda)} X_{\ell}^{\sigma\sigma'} \right)} = \frac{1}{|\ell - (\sigma-\sigma')\lambda|} \underline{X}_{\ell}^{\sigma\sigma'}.
	\end{equation*}
	
	Noting that $|\ell - (\sigma-\sigma')\lambda| \geq 1/2$ for all terms in the range of $\mathcal{L}_{n/2}$ within the region $|2\lambda - n| < \gamma \leq 1/2$, we obtain:
	\begin{align*}
	|\mathcal{L}^{-1}_\lambda X|_{s,p,R} &= \sup_{\|u\|_{s,p} \leq R} \normS{\underline{\mathcal{L}^{-1}_\lambda X}(u)}_{s,p} \\
	&= \sup_{\ell \in \mathbb{Z}} \langle \ell \rangle^p e^{s|\ell|} \max_{\sigma,\sigma'=\pm1} \sup_{\|u\|_{s,p} \leq R} \bigg|\underline{\left( \frac{1}{\im(\ell - (\sigma-\sigma')\lambda)} X_{\ell}^{\sigma\sigma'} \right)}(u)\bigg| \\
	&= \sup_{\ell \in \mathbb{Z}} \langle \ell \rangle^p e^{s|\ell|} \max_{\sigma,\sigma'=\pm1} \frac{1}{|\ell - (\sigma-\sigma')\lambda|} \sup_{\|u\|_{s,p} \leq R} |\underline{X}_{\ell}^{\sigma\sigma'}(u)| \\
	&= \sup_{\ell \in \mathbb{Z}} \langle \ell \rangle^p e^{s|\ell|} \max_{\sigma,\sigma'=\pm1} \frac{1}{|\ell - (\sigma-\sigma')\lambda|} \| X_{\ell}^{\sigma\sigma'} \|_{s,p,R} \\
	&\leq 2 \sup_{\ell \in \mathbb{Z}} \langle \ell \rangle^p e^{s|\ell|} \max_{\sigma,\sigma'=\pm1} \| X_{\ell}^{\sigma\sigma'} \|_{s,p,R} = 2 |X|_{s,p,R},
	\end{align*}
	yielding the claim.
\end{proof}

Given $n \geq 1$, in the following theorem, we require the smallness condition on $P$ to be satisfied for all $\lambda \in \mathbb{R}^+$ such that $|2\lambda - n| \leq 1/2$. This ensures that both the existence of the fixed point and the derived estimates are valid uniformly throughout this spectral region.

\begin{thm}[Reduction to constant coefficients in the range of $\mathcal{L}_{n/2}$]\label{fpt}
Let $|P|_{s,p,R} \leq (48 K)^{-1}$, where the constant $K$ is defined in Lemma \ref{banachliealgA}. Then the map $\phi_P$ has a unique fixed point $A_P$ in the ball $\overline{\mathscr{U}}_{8|P|_{s,p,R}}(\Pi_{rg}(\underline{\mathcal{A}}_{R}^{\pi}))$. The fixed point is such that $A_{P=0} = 0$ and satisfies the following asymptotic expansion for small $P$:
\begin{equation*} 
A_P = \mathcal{L}^{-1}_\lambda \Pi_{rg} P + O(|P|_{s,p,R}^2).
\end{equation*}
Furthermore, $A_P$ is analytic in $\lambda$ and satisfies the estimate
\begin{equation*}
	|\partial_\lambda A_P|_{s,p,R}\leq 96 |P|_{s,p,R}.
\end{equation*}
Finally, $\mathcal{L}_\lambda A_P$, which is also analytic in $\lambda$, is estimated as $|\mathcal{L}_\lambda A_P|_{s,p,R} \leq 4|P|_{s,p,R}$ and satisfies the following asymptotic expansion for small $P$:
\begin{equation*}
\mathcal{L}_\lambda A_P = \Pi_{rg} P + O(|P|_{s,p,R}^2).
\end{equation*}
\end{thm}

\begin{proof}
	The proof follows from the application of the general fixed point framework developed in Appendix \ref{AAD} to the Banach-Lie algebra $(\underline{\mathcal{A}}_{R}^{\pi}, |\cdot|_{s,p,R})$. According to Lemma \ref{banachliealgA}, this space satisfies the sub-multiplicative property with constant $K$. We identify the operators in the appendix with our specific case by setting $\Pi = \Pi_{rg}$, $T = \mathcal{L}^{-1}_\lambda$, and $X = A$.
	\medskip
	
	\noindent \textit{Existence and uniqueness.} 
	By hypothesis, $|P|_{s,p,R} \leq (48K)^{-1} < (24K)^{-1}$. This satisfies the smallness condition required by the Fixed Point Theorem \ref{FPT}. Consequently, $\phi_P$ is a contraction on the ball $\overline{\mathscr{U}}_{8|P|_{s,p,R}}$, ensuring the existence of a unique fixed point $A_P$ such that $A_{P=0} = 0$. The first-order expansion follows directly from the same theorem:
	\begin{equation*}
	A_P = \mathcal{L}^{-1}_\lambda \Pi_{rg} P + O(|P|_{s,p,R}^2).
	\end{equation*}
	
	\noindent \textit{Analyticity in $\lambda$ and estimate for the derivative.}
	Since $\mathcal{L}^{-1}_\lambda$ and $P$ are analytic in $\lambda$, Corollary \ref{corapp} implies that $A_P$ is also analytic in $\lambda$. To derive the bound on $\partial_\lambda A_P$, we use the estimate provided in the same corollary:
	\begin{equation*}
	|\partial_\lambda A_P|_{s,p,R} \leq \frac{4 \|\partial_\lambda \mathcal{L}^{-1}_\lambda\| |P|_{s,p,R} + 8 |\partial_\lambda P|_{s,p,R}}{1 - 24K|P|_{s,p,R}}.
	\end{equation*}
	Given the condition $|P|_{s,p,R} \leq (48K)^{-1}$, the denominator is bounded as $(1 - 24K|P|)^{-1} \leq (1 - 1/2)^{-1} = 2$. Using the specific bounds for the derivatives of the operators in this functional setting, and the fact that $\partial_\lambda P=-\lambda^{-1}P$, where $\lambda \geq 1/4$, we obtain:
	\begin{equation*}
	|\partial_\lambda A_P|_{s,p,R} \leq 96 |P|_{s,p,R}.
	\end{equation*}
	\smallskip
	
	\noindent \textit{Estimates for $\mathcal{L}_\lambda A_P$.}
	Setting $U = \mathcal{L}_\lambda$ as the inverse of $T$ on $\text{range}(\mathcal{L}_{n/2})$, we apply Lemma \ref{UX}. Since $A_P$ is in the range of $\mathcal{L}^{-1}_\lambda$, the expression $\mathcal{L}_\lambda A_P$ is well-defined and satisfies:
	\begin{equation*}
	|\mathcal{L}_\lambda A_P|_{s,p,R} \leq 4 |P|_{s,p,R}.
	\end{equation*}
	Finally, the asymptotic expansion for $\mathcal{L}_\lambda A_P$ is given by:
	\begin{equation*}
	\mathcal{L}_\lambda A_P = \Pi_{rg} P + O(|P|_{s,p,R}^2),
	\end{equation*}
	which concludes the proof.
\end{proof}

As a consequence of Theorem \ref{fpt}, the matrix associated with the linear vector field \eqref{questo} is now reduced to the form
\begin{equation}\label{p1}
	\lambda J_1 +\tilde{\mathcal{Z}}(u,\lambda,t), \qquad \text{with} \quad \tilde{\mathcal{Z}} \in \mathrm{Ker}(\mathcal{L}_{n/2}).
\end{equation}
By construction, $\tilde{\mathcal{Z}}$ is defined by
\begin{equation}\label{p2}
	\tilde{\mathcal{Z}}= \Pi_{ker} \bigg( -\sum_{k=1}^\infty\frac{\mathrm{ad}_{-A}^{k-1}}{k!}\mathcal{L}_\lambda A+ e^{\mathrm{ad}_{-A}}P\bigg),
\end{equation}
where $\Pi_{ker}$ denotes the projector onto the kernel of $\mathcal{L}_{n/2}$, and $A$ denotes the fixed point provided by Theorem \ref{fpt}. Using the fact that $\mathcal{L}_\lambda A\in \mathrm{range}(\mathcal{L}_{n/2})$, and isolating the term of order $O(|P|_{s,p,R})$, we get
\begin{equation}\label{p3}
	\tilde{\mathcal{Z}}=\Pi_{ker}P+\tilde{\mathcal{P}}_\infty:= \Pi_{ker}P+\Pi_{ker} \bigg(  (e^{\mathrm{ad}_{-A}}-\id)P-\sum_{k=2}^\infty\frac{\mathrm{ad}_{-A}^{k-1}}{k!}\mathcal{L}_\lambda A\bigg).
\end{equation}

The following corollary, proved in Appendix \ref{AAD} (see Lemma \ref{Deltaestimates}), shows that $\tilde{\mathcal{P}}_\infty$ is indeed of order $O(|P|_{s,p,R}^2)$ and provides an asymptotic expansion near $P=0$.

\begin{cor}\label{questocor}
		The following bound holds: 
		\begin{equation*}
		|\tilde{\mathcal{P}}_\infty|_{s,p,R} \leq 64 K |P|_{s,p,R}^2.
		\end{equation*}
		Furthermore, $\tilde{\mathcal{P}}_\infty$ satisfies the following asymptotic expansion for small $P$:
		\begin{equation*}
		\tilde{\mathcal{P}}_\infty = \Pi_{ker}\left( [\Pi_{ker}P, \mathcal{L}^{-1}_\lambda \Pi_{rg} P]_{\star} + \frac{1}{2} [\Pi_{rg} P, \mathcal{L}^{-1}_\lambda \Pi_{rg} P]_{\star} \right) + O(|P|_{s,p,R}^3).
		\end{equation*}
\end{cor}

To understand the final step of the reduction, it is useful to recall the explicit structure of the kernel of $\mathcal{L}_{n/2}$, namely
\begin{equation}\label{kernel}
\mathrm{Ker}(\mathcal{L}_{n/2}) = 
\left\{ \begin{pmatrix} 
X^{++}_0 &  X^{+-}_ne^{\im nt} \\
X^{-+}_{-n} e^{-\im nt} & X^{--}_{0}
\end{pmatrix} \ : \ X^{--}_0 = -X^{++}_0 \in \im \R, \  X^{-+}_{-n} = \overline{X^{+-}_{n}} \in \C  \right\}.
\end{equation}
The elements in this kernel are still time-dependent, but we can then reach a frame where the system is stationary by applying an additional $SU(1,1)$ transformation, specifically a rotation $S_n(t)$ that counter-acts these oscillations:
\begin{equation*}
S_n(t) = \begin{pmatrix} e^{-\im nt/2} & 0 \\ 0 & e^{\im nt/2} \end{pmatrix}.
\end{equation*}
Under this final change of variables, the resonant terms are mapped to constant matrices, effectively completing the reduction to constant coefficients. As before, given $n \geq 1$, we assume the smallness condition on $P$ to hold for all $\lambda \in \mathbb{R}^+$ within the closed interval $|2\lambda - n| \leq 1/2$. 

\begin{proof}[Proof of Theorem \ref{final_reducibility_thm}]
	Let $\tilde{\mathcal{D}} = \lambda J_1 + \tilde{\mathcal{Z}}$ be the resonant form provided by Theorem \ref{fpt}, with $\tilde{\mathcal{Z}} = \Pi_{ker} P + \tilde{\mathcal{P}}_\infty$ defined in \eqref{p2} and \eqref{p3}. Following the change of variables $\eta = S_n e^{-A} z$, the reduced matrix is obtained as $$\mathcal{D} = S_n \tilde{\mathcal{D}} S_n^{-1} + \dot{S}_n S_n^{-1}.$$
	The diagonal part $\lambda J_1$ commutes with $S_n$, while the derivative of the rotation contributes the shift:
	\begin{equation*}
	\dot{S}_n S_n^{-1} = \begin{pmatrix} -\im n/2 & 0 \\ 0 & \im n/2 \end{pmatrix} = - \frac{n}{2} J_1,
	\end{equation*}
	yielding the term $(\lambda - n/2) J_1$. For the resonant part, the adjoint action of $S_n$ on any $M \in \mathrm{Ker}(\mathcal{L}_{n/2})$ cancels the time-dependence $e^{\pm \im nt}$. Specifically, 
	\begin{equation}\label{PinftySn}
		\mathcal{Z}_0 := S_n (\Pi_{ker} P) S_n^{-1} \quad \text{and} \quad \mathcal{P}_\infty := S_n \tilde{\mathcal{P}}_\infty S_n^{-1} \quad 
	\end{equation}
	are constant matrices. The explicit form of $\mathcal{Z}_0$ follows from the Fourier expansion of the perturbation $P$ in the kernel of $\mathcal{L}_{n/2}$. The majorant analyticity and the bounds for the coefficients of $\mathcal{P}_\infty$ are a direct consequence of the estimates for $\tilde{\mathcal{P}}_\infty$ established in Corollary \ref{questocor} and the analyticity of the fixed point $A_P$ with respect to $\lambda$.
\end{proof}

\subsubsection{Estimates for $P \circ |\partial|^{1/2}$ and $\partial_{z_j}(P \circ |\partial|^{1/2})$}
To ease the notation, we suppress the explicit dependence on the parameter $s$ in the Banach spaces we consider and their respective norms, since it remains constant and plays no role throughout this appendix; nevertheless, we emphasize that all subsequent results hold for any $s \geq 0$.

\begin{defn}
	Given $R>0$, $p>1$, and $s \geq 0$, for any $j \in \mathbb{Z}$ and $\tau = \pm 1/2$, we denote by $\mathcal{M}^{\tau,j}_{p,R}$ the space of maps $X: \mathscr{U}_R(\ell^\infty_{p}) \to \prod_{\ell \in \mathbb{Z}} \mathrm{Mat}(2,\mathbb{C})$ such that each entry $X_\ell^{\sigma\sigma'} \in \mathcal{H}_{p,R}$ satisfies the momentum conservation condition $\pi(\alpha) = \ell-j$. We equip this space with a Banach structure by introducing the majorant $\underline{X} = (\underline{X}_\ell)_{\ell \in \mathbb{Z}}$, where $\underline{X}_\ell^{\sigma\sigma'}(z) := \sum_{\alpha \in \mathscr{M}_{\ell-j}} |X_{\ell \alpha}^{\sigma\sigma'}| z^\alpha$. 
	
	To define the norm, we first generalize the norm on the space of sequences of $2 \times 2$ complex matrices $W=(W_\ell)_{\ell \in \Z}$ to the $(\tau, j)$-weighted case:
	\begin{equation*}
		\|W\|^{\tau,j}_{p} := \sup_{\ell \in \mathbb{Z}} \langle \ell-j \rangle^p e^{s|\ell-j|}\langle \ell \rangle^\tau \max_{\sigma,\sigma'=\pm1} |W_{\ell}^{\sigma\sigma'}|.
	\end{equation*}
	The norm on $\mathcal{M}^{\tau,j}_{p,R}$ is then given by:
	\begin{equation}\label{norm_shifted_weighted}
		|X|^{\tau,j}_{p,R} := \sup_{\|z\|_{p} \leq R} \|\underline{X}(z)\|^{\tau,j}_{p} = \sup_{\ell \in \mathbb{Z}} \langle \ell-j \rangle^p e^{s|\ell-j|} \langle \ell \rangle^\tau \max_{\sigma,\sigma'=\pm1} \| X_{\ell}^{\sigma\sigma'} \|_{p,R},
	\end{equation}
	which, due to the non-negativity of the coefficients, coincides with evaluating the majorant on the weight vector $v_{p,R}$. The space $\mathcal{M}^{\tau,j}_{p,R}$ consists of all such maps for which this norm is finite.
\end{defn}

\begin{rmk}\label{majz}
	Any $X \in \mathcal{M}^{\tau,0}_{p,R}$ is a majorant analytic function $X:\ell^\infty_{p} \to \ell^{\infty}_{p+\tau,R}(\mathrm{Mat}(2,\C))$. In fact, Theorem \ref{uniformanalytic} ensures that $\underline{X}$ is analytic as it is defined by the uniformly bounded family of analytic functions $\underline{X}_\ell^{\sigma\sigma'} \in \mathcal{H}_{p,R}$.
\end{rmk}

Given a positive real number $\lambda$, let $P$ be the map from $\ell^\infty_p$ to the space of sequences of $2 \times 2$ complex matrices, defined for all $u \in \ell^\infty_p$ by its components:
\begin{equation*}
	P_\ell(u)= \frac{u_\ell}{2\lambda}
	\begin{pmatrix}
		-\im & \im \\ -\im & \im
	\end{pmatrix}, \quad \text{for all } \ell \in \Z.
\end{equation*}

In the following, we denote by $|\partial|^{1/2} : \ell^\infty_{p} \to \ell^\infty_{p-1/2}$ the operator acting as a Fourier multiplier:
\begin{equation}
	u := |\partial|^{1/2}z, \quad \text{with} \quad u_\ell = \langle \ell \rangle^{1/2} z_\ell.
\end{equation}

\begin{lemma}\label{PeDP}
	For any $R>0$, $P \circ |\partial|^{1/2} \in \mathcal{M}^{-1/2,0}_{p,R}$ and $\partial_{z_j}(P \circ |\partial|^{1/2}) \in \mathcal{M}^{-1/2,j}_{p,R}$. Specifically, the following estimates hold:
	\begin{equation}
		|P \circ |\partial|^{1/2}|^{-1/2,0}_{p,R} = \frac{R}{2\lambda}, \quad |\partial_{z_j}(P \circ |\partial|^{1/2})|^{-1/2,j}_{p,R} \leq \frac{1}{2\lambda}.
	\end{equation}
\end{lemma}

\begin{proof}
	The first equality follows directly from the definition of $P$, noting that the weight of the norm compensates for the action of $|\partial|^{1/2}$. By the same reasoning, the second bound is immediate, as the derivative produces a $\delta_{j\ell}$ which is trivially bounded by $\jap{\ell-j}^{-p}e^{-s|\ell-j|}$.
\end{proof}

\subsubsection{Estimates for $A_P\circ|\partial|^{1/2}$ and $\mathcal{L}_\lambda A_P \circ |\partial|^{1/2}$}
Fix $n \in \mathbb{N}$ and let $\lambda \in \mathbb{R}^+$ be such that $|2\lambda-n| \leq 1/2$. Let $\mathcal{L}_\lambda$ be the operator defined in \eqref{Llamb}. By identifying periodic functions taking values in the Lie algebra $\mathrm{Mat}(2,\mathbb{C})$ with the sequences of their Fourier coefficients, the operator $\mathcal{L}_{\lambda}$ induces an action on the space of matrix sequences $W=(W_\ell)$ defined by
\begin{equation}
	(\mathcal{L}_\lambda W)^{\sigma\sigma'}_\ell = i(\ell-\lambda(\sigma-\sigma'))W^{\sigma\sigma'}_\ell.
\end{equation}
We have already observed that, if $\lambda=n/2$, this operator has a non-trivial kernel. Namely, a sequence $W$ belongs to $\mathrm{range}(\mathcal{L}_{n/2})$ if and only if its components $W_\ell^{\sigma\sigma'}$ are zero whenever $\ell - \frac{n}{2}(\sigma-\sigma') = 0$.

From the definition of the action of $\mathcal{L}_\lambda$ on the Fourier coefficients, it follows that the formal inverse $\mathcal{L}_\lambda^{-1}$ acts diagonally as:
\begin{equation}
	(\mathcal{L}_\lambda^{-1} W)^{\sigma\sigma'}_\ell = \frac{1}{i(\ell-\lambda(\sigma-\sigma'))} W^{\sigma\sigma'}_\ell,
\end{equation}
provided that the denominator does not vanish. For any $\lambda$ such that $|2\lambda - n| \leq 1/2$, we define $\mathcal{L}_\lambda^{-1}$ on $\mathrm{range}(\mathcal{L}_{n/2})$. Clearly, both $\mathcal{L}_\lambda$ and its inverse preserve $\mathrm{range}(\mathcal{L}_{n/2})$. 

For $\lambda$ in the specified region, the divisors for indices in the range of the operator satisfy the following estimate:
\begin{equation}\label{divisor_estimate}
	\frac{\langle \ell \rangle}{|\ell - \lambda(\sigma - \sigma')|} \leq 4 n.
\end{equation}
This bound shows that the inverse operator compensates for the loss of regularity, effectively shifting the weight from $\langle \ell \rangle^{-1/2}$ to $\langle \ell \rangle^{+1/2}$ at the cost of a factor $n$, as established in the following lemma.

\begin{lemma}\label{Lamenouno}
	Let $n \in \mathbb{N}$ and $\lambda \in \mathbb{R}^+$ satisfy $|2\lambda - n| \leq 1/2$. Then, $\mathcal{L}_\lambda^{-1}$ is a bounded linear operator
	\begin{equation*}
		\mathcal{L}_\lambda:\mathcal{M}^{-1/2,j}_{p,R} \cap \mathrm{range}(\mathcal{L}_{n/2}) \to \mathcal{M}^{+1/2,j}_{p,R} \cap \mathrm{range}(\mathcal{L}_{n/2})
	\end{equation*}
	In particular, there exists a constant $C > 0$ such that
	\begin{equation}
		|\mathcal{L}_\lambda^{-1} X|^{+1/2,j}_{p,R} \leq 4n |X|^{-1/2,j}_{p,R}
	\end{equation}
	for every $X \in \mathcal{M}^{-1/2,j}_{p,R} \cap \mathrm{range}(\mathcal{L}_{n/2})$.
\end{lemma}

\begin{proof}
	By the definition of the norm on $\mathcal{M}^{\tau,j}_{p,R}$ given in \eqref{norm_shifted_weighted}, we have
	\begin{equation*}
		|\mathcal{L}_\lambda^{-1} X|^{+1/2,j}_{p,R} = \sup_{\ell \in \mathbb{Z}} \left( \langle \ell-j \rangle^p e^{s|\ell-j|} \langle \ell \rangle^{1/2} \max_{\sigma,\sigma'} \| (\mathcal{L}_\lambda^{-1} X)_{\ell}^{\sigma\sigma'} \|_{p,R} \right).
	\end{equation*}
	Recall that $(\mathcal{L}_\lambda^{-1} X)_{\ell}^{\sigma\sigma'} = \frac{1}{i(\ell-\lambda(\sigma-\sigma'))} X_{\ell}^{\sigma\sigma'}$. Substituting this into the expression, we obtain
	\begin{equation}\label{sopra}
		|\mathcal{L}_\lambda^{-1} X|^{+1/2,j}_{p,R} = \sup_{\ell \in \mathbb{Z}} \left( \langle \ell-j \rangle^p e^{s|\ell-j|} \langle \ell \rangle^{1/2} \max_{\sigma,\sigma'} \frac{1}{|\ell-\lambda(\sigma-\sigma')|} \| X_{\ell}^{\sigma\sigma'} \|_{p,R} \right).
	\end{equation}
	We can rewrite the terms inside the supremum by multiplying and dividing by $\langle \ell \rangle^{-1/2}$ to recover the norm of $X$:
	\begin{equation*}
		\eqref{sopra} = \sup_{\ell \in \mathbb{Z}} \left( \frac{\langle \ell \rangle}{|\ell-\lambda(\sigma-\sigma')|} \cdot \langle \ell-j \rangle^p e^{s|\ell-j|} \langle \ell \rangle^{-1/2} \max_{\sigma,\sigma'} \| X_{\ell}^{\sigma\sigma'} \|_{p,R} \right).
	\end{equation*}
	Using the estimate \eqref{divisor_estimate}, we get
	\begin{equation*}
		|\mathcal{L}_\lambda^{-1} X|^{+1/2,j}_{p,R} \leq 4 n \sup_{\ell \in \mathbb{Z}} \left( \langle \ell-j \rangle^p e^{s|\ell-j|} \langle \ell \rangle^{-1/2} \max_{\sigma,\sigma'} \| X_{\ell}^{\sigma\sigma'} \|_{p,R} \right).
	\end{equation*}
	The supremum on the right-hand side is precisely the norm $|X|^{-1/2,j}_{p,R}$. Therefore, we get
	\begin{equation*}
		|\mathcal{L}_\lambda^{-1} X|^{+1/2,j}_{p,R} \leq 4n |X|^{-1/2,j}_{p,R},
	\end{equation*}
	which ends the proof.
\end{proof}

Let $A_P$ be the fixed point provided by Theorem \ref{fpt}, which we recall to be defined by
\begin{equation}\label{phiB}
	A_P=\phi_P(A_P) = \mathcal{L}_\lambda^{-1} \Pi_{rg} \left( \mathcal{B}(A_P) \Pi_{rg} \left( e^{\mathrm{ad}_{-A_P}} P \right) \right),
\end{equation}
where the operators $\mathcal{E}(X)$ and $\mathcal{B}(X)$ are formally defined via the Taylor series of the complex functions $f(z) = e^z = \sum_{k=0}^{\infty} \frac{1}{k!} z^k$ and $g(z) = \frac{z}{e^z - 1} = \sum_{k=0}^{\infty} \beta_k z^k$ (see Section \ref{EandB}). Specifically, $\mathcal{E}(X) = f(\mathrm{ad}_{-X})$ and $\mathcal{B}(X) = g(\mathrm{ad}_{-X})$.

In the following lemma, we first derive a bound for $\mathcal{L}_\lambda A_P \circ |\partial|^{1/2}$. We then apply the smoothing property of $\mathcal{L}_\lambda^{-1}$ established in Lemma \ref{Lamenouno} to prove that $A_P \circ |\partial|^{1/2}$ maps into a sequence space with stronger decay, with a resulting bound constant that grows linearly with $n$.

\begin{lemma}\label{AeLA}
	Let $n \in \mathbb{N}$ and $\lambda > 0$ such that $|2\lambda - n| \leq 1/2$. Then, for all $p>3/2$, we have that $\mathcal{L}_\lambda A_P \circ |\partial|^{1/2} \in \mathcal{M}^{-1/2,0}_{p,R}$ and $A_P \circ |\partial|^{1/2} \in \mathcal{M}^{+1/2,0}_{p,R}$, with the following estimates:
	\begin{equation}
		|\mathcal{L}_\lambda A_P \circ |\partial|^{1/2}|_{p,R}^{-1/2,0} \leq \frac{2R}{\lambda}, \quad |A_P \circ |\partial|^{1/2}|_{p,R}^{+1/2,0} \leq 8R\frac{n}{\lambda}.
	\end{equation}
\end{lemma}

\begin{proof}
	By Theorem \ref{fpt}, $|\mathcal{L}_\lambda A_P|_{p',R}\leq 4|P|_{p',R}$ for all $p'>1$. Then, the first bound follows by the fact that, if $X \in \underline{\mathcal{A}}_R^\pi$ satisfies $|X|_{p-1/2,R}<\infty$, then $X \circ |\partial|^{1/2} \in \mathcal{M}^{-1/2,0}_{p,R}$, with the identity $$|X \circ |\partial|^{1/2}|_{p,R}^{-1/2,0} = |X|_{p-1/2,R}$$ holding by definition, combined with Lemma \ref{PeDP}. The second estimate is obtained by writing $A_P= (\mathcal{L}_{\lambda}^{-1} \circ \mathcal{L}_{\lambda}) A_P$ and applying the first bound together with the one provided by Lemma  \ref{Lamenouno}.
\end{proof}

\subsubsection{Estimates for $\partial_{z_j}(A_P \circ |\partial|^{1/2})$ and $\partial_{z_j}(\mathcal{L}_{\lambda}A_P\circ|\partial|^{1/2})$}
In this section, we provide estimates for the maps $\partial_{z_j}(A_P \circ |\partial|^{1/2})$ and $\partial_{z_j}(\mathcal{L}_{\lambda}A_P\circ|\partial|^{1/2})$, where $A_P$ is the fixed point of the operator from Theorem \ref{fpt}, whose definition is recalled in \eqref{phiB}. Since it is  implicitly defined through commutators, the main task is to provide estimates for the latter by understanding how the convolution product works in our weighted spaces. The following lemmas are the key tool for all the calculations in this section.

\begin{lemma}[Convolution and Algebra Property]\label{MN}
	Let $X \in \mathcal{M}^{-1/2, j}_{p,R}$ and $Y \in \mathcal{M}^{+1/2, j'}_{p,R}$. The convolution product $Z = X \star Y$, defined by the sequence of matrices 
	$$Z_\ell(z) = \sum_{k \in \mathbb{Z}} X_{\ell-k}(z) Y_k(z),$$ 
	belongs to $\mathcal{M}^{-1/2,j+j'}_{p,R}$. More precisely, it satisfies the following properties:
	\begin{itemize}
		\item \textbf{Momentum Conservation:} For each $\ell \in \mathbb{Z}$, the entries of $Z_\ell(z)$ belong to $\mathcal{H}_{p,R}$ and their Taylor expansion involves only multi-indices $\alpha$ such that $\pi(\alpha) = \ell-(j+j')$.
		\item \textbf{Estimate:} The following estimates hold:
		\begin{equation}\label{conv_estimate}
			|X \star Y|^{-1/2, j+j'}_{p,R}, |Y \star X|^{-1/2, j'+j}_{p,R} \leq 2C_\star(p) |X|^{-1/2, j}_{p,R} \, |Y|^{+1/2, j'}_{p,R},
		\end{equation}
		where $C_\star(p) = 2^p(1+2\zeta(p))$ denotes the convolution algebra constant of $\ell^\infty_{s,p}$ and $\zeta(p) = \sum_{n \geq 1} n^{-p}$ is the Riemann zeta function.
	\end{itemize}
\end{lemma}

\begin{proof}
	\noindent \textit{Momentum Conservation.} Let $X \in \mathcal{M}^{-1/2, j}_{p,R}$ and $Y \in \mathcal{M}^{+1/2, j'}_{p,R}$. By definition, the entries of the convolution $Z = X \star Y$ are given by
	\begin{equation*}
		Z_{\ell}^{\sigma\sigma'}(z) = \sum_{k \in \mathbb{Z}} \sum_{\rho \in \{\pm 1\}} X_{\ell-k}^{\sigma\rho}(z) Y_{k}^{\rho\sigma'}(z).
	\end{equation*}
	Substituting the Taylor expansions of the entries, we obtain
	\begin{equation*}
		Z_{\ell}^{\sigma\sigma'}(z) = \sum_{k \in \mathbb{Z}} \sum_{\rho \in \{\pm 1\}} \left( \sum_{\alpha \in \mathscr{M}_{\ell-k-j}} X_{\ell-k,\alpha}^{\sigma\rho} z^\alpha \right) \left( \sum_{\beta \in \mathscr{M}_{k-j'}} Y_{k,\beta}^{\rho\sigma'} z^\beta \right).
	\end{equation*}
	The product of two monomials $z^\alpha z^\beta = z^{\alpha+\beta}$ satisfies the condition $\pi(\alpha+\beta) = \pi(\alpha) + \pi(\beta) = (\ell-k-j) + (k-j') = \ell - (j+j')$. Thus, each term in the expansion of $Z_\ell$ belongs to $\mathscr{M}_{\ell-(j+j')}$. Since $\mathcal{H}_{p,R}$ is an algebra, the series converges in $\mathcal{H}_{p,R}$, and the momentum condition is satisfied.
	\smallskip
	
	\noindent \textit{Estimate.} To estimate the norm $|X \star Y|^{-1/2, j+j'}_{p,R}$, we consider the majorant $\underline{Z}_{\ell}$. By the triangle inequality and the definition of the majorant norm, we have
	\begin{align*}
		\max_{\sigma,\sigma'} |\underline{Z}_{\ell}^{\sigma\sigma'}(z)| &\leq 2 \sum_{k \in \mathbb{Z}} \max_{\sigma,\rho} |\underline{X}_{\ell-k}^{\sigma\rho}(z)| \max_{\rho,\sigma'} |\underline{Y}_{k}^{\rho\sigma'}(z)| \\
		&\leq 2 \sum_{k \in \mathbb{Z}} \frac{\|\underline{X}(z)\|^{-1/2,j}_{p}}{\langle \ell-k-j \rangle^p e^{s|\ell-k-j|} \langle \ell-k \rangle^{-1/2}} \cdot \frac{\|\underline{Y}(z)\|^{+1/2,j'}_{p}}{\langle k-j' \rangle^p e^{s|k-j'|} \langle k \rangle^{1/2}}.
	\end{align*}
	Multiplying both sides by the weight $\langle \ell - (j+j') \rangle^p e^{s|\ell-(j+j')|} \langle \ell \rangle^{-1/2}$ and taking the supremum over $\ell$, we obtain:
	\begin{equation*}
		\|\underline{Z}(z)\|^{-1/2,j+j'}_{p} \leq 2\cdot2^{p} \left( \sup_{\ell \in \mathbb{Z}} \sum_{k \in \mathbb{Z}} \frac{1}{\langle \ell-k-j \rangle^p} \right) \|\underline{X}(z)\|^{-1/2,j}_{p} \|\underline{Y}(z)\|^{+1/2,j'}_{p}.
	\end{equation*}
	The sum in the parentheses is independent of $\ell$ and converges for any $p > 1$. In particular, we get:
	\begin{equation*}
		\|\underline{Z}(z)\|^{-1/2,j+j'}_{p} \leq 2C_\star(p) \|\underline{X}(z)\|^{-1/2,j}_{p} \|\underline{Y}(z)\|^{+1/2,j'}_{p},
	\end{equation*}
	Finally, taking the supremum over $\|z\|_{p} \leq R$ yields the desired estimate.
\end{proof}

\begin{lemma}\label{comm_estimate}
	Let $R>0$ and $X \in \mathcal{M}^{+1/2,0}_{p,R}$. The adjoint operator $\mathrm{ad}_X: Y \mapsto X \star Y - Y \star X$ is a bounded linear operator from $\mathcal{M}^{-1/2,j}_{p,R}$ to itself for every $j \in \mathbb{Z}$. Its operator norm, defined as
	\begin{equation}
		\|\mathrm{ad}_{X}\|_R := \sup_{j \in \mathbb{Z}} \sup_{|Y|^{-1/2,j}_{p,R} \leq 1} |X \star Y - Y \star X|_{p,R}^{-1/2,j},
	\end{equation}
	satisfies the estimate
	\begin{equation}
		\|\mathrm{ad}_{X}\|_R \leq 4C_\star(p) |X|^{+1/2,0}_{p,R}.
	\end{equation}
\end{lemma}

\begin{proof}
	Let $Y \in \mathcal{M}^{-1/2,j}_{p,R}$. By Lemma \ref{MN}, the star product satisfies $|X \star Y|^{-1/2,j}_{p,R} \leq C |X|^{+1/2,0}_{p,R} |Y|^{-1/2,j}_{p,R}$. Consequently, the commutator can be estimated as:
	\begin{equation*}
		|X \star Y - Y \star X|^{-1/2,j}_{p,R} \leq |X \star Y|^{-1/2,j}_{p,R} + |Y \star X|^{-1/2,j}_{p,R} \leq 4C_\star(p) |X|^{+1/2,0}_{p,R} |Y|^{-1/2,j}_{p,R}.
	\end{equation*}
	By dividing by the norm of $Y$ and taking the supremum over $Y \in \mathcal{M}^{-1/2,j}_{p,R} \setminus \{0\}$, we obtain:
	\begin{equation*}
		\sup_{|Y|^{-1/2,j}_{p,R} \leq 1} |X \star Y - Y \star X|^{-1/2,j}_{p,R} \leq 4C_\star(p) |X|^{+1/2,0}_{p,R}.
	\end{equation*}
	Since the right-hand side does not depend on the momentum index $j$, the inequality remains valid when taking the supremum over all $j \in \mathbb{Z}$. This concludes the proof that $\|\mathrm{ad}_X\|_R \leq 4C_\star(p) |X|^{+1/2,0}_{p,R}$.
\end{proof}

\begin{lemma}[Estimates for $\phi_Y$ and its derivative]\label{phi_estimates}
	Let $n \in \mathbb{N}$ and $\lambda > 0$ such that $|2\lambda - n| \leq 1/2$. Given $Y \in \bigcup_{j \in \Z}\mathcal{M}^{-1/2,j}_{p,R}$ and $X \in \mathcal{M}^{+1/2,0}_{p,R}$, consider the map:
	\begin{equation}
		\phi_Y(X) = \mathcal{L}_\lambda^{-1} \Pi_{rg} \left( \mathcal{B}(X) \Pi_{rg} \left( e^{{\mathrm{ad}}_{-X}} Y \right) \right).
	\end{equation}
	Then, the following properties hold:
	\begin{enumerate}
		\item For any fixed $X \in \mathcal{M}^{+1/2,0}_{p,R}$ satisfying $4C_\star(p)|X|^{+1/2,0}_{p,R} \leq 1/2$, the map $Y \mapsto \phi_Y(X)$ defines a bounded linear operator from $\mathcal{M}^{-1/2,j}_{p,R}$ to $\mathcal{M}^{+1/2,j}_{p,R}$ for all $j \in \Z$. In particular:
		\begin{equation}
			|\phi_Y(X)|^{+1/2,j}_{p,R} \leq 16n |Y|^{-1/2,j}_{p,R}.
		\end{equation}
		\item For a fixed $Y \in \mathcal{M}^{-1/2,0}_{p,R}$, the Fréchet derivative of $\phi_Y$ at $X \in \mathcal{M}^{+1/2,0}_{p,R}$, with $4C_\star(p)|X|^{+1/2,0}_{p,R} \leq 1/2$, denoted by $\mathrm{d}\phi_Y(X)$, satisfies:
		\begin{equation}
			\|\mathrm{d} \phi_Y(X)\|_R := \sup_{j \in \mathbb{Z}} \sup_{|H|^{+1/2,j}_{p,R} \leq 1} |\mathrm{d}\phi_Y(X)[H]|^{+1/2,j}_{p,R} \leq 256nC_\star(p) |Y|^{-1/2,0}_{p,R}.
		\end{equation}
	\end{enumerate}
\end{lemma}

\begin{proof}
	\textit{1. Estimate for $\phi_Y(X)$.} Let $X \in \mathcal{M}^{+1/2,0}_{p,R}$ be fixed and consider $Y \in \mathcal{M}^{-1/2,j}_{p,R}$. From Lemma \ref{comm_estimate}, the adjoint operator $\mathrm{ad}_{-X}$ is a bounded linear operator with norm:
	\begin{equation*}
		\|\mathrm{ad}_{-X}\|_R \leq 4C_\star(p) |X|^{+1/2,0}_{p,R} \leq \frac{1}{2}.
	\end{equation*}
	Consequently, the operators $\mathcal{E}(X) = \sum_{k=0}^\infty \frac{1}{k!} (\mathrm{ad}_{-X})^k$ and $\mathcal{B}(X) = \sum_{k=0}^\infty \beta_k (\mathrm{ad}_{-X})^k$ are well-defined and bounded. In fact, under the assumption $4C_\star(p) |X|^{+1/2,0}_{p,R} \leq 1/2$, we provide a coarse estimate for their norms by dominating the Taylor coefficients with $1$:
	\begin{equation}\label{geom_est}
		\|\mathcal{E}(X)\|_R, \ \|\mathcal{B}(X)\|_R \leq \sum_{k=0}^\infty \|\mathrm{ad}_{-X}\|^k_R \leq 2.
	\end{equation}
	Let $Z = \Pi_{rg}\mathcal{B}(X) \Pi_{rg} (\mathcal{E}(X) Y)$. Since both the adjoint action and the projectors preserve the momentum conservation condition, $Z \in \mathcal{M}^{-1/2,j}_{p,R} \cap \mathrm{range}(\mathcal{L}_{n/2})$. Using $\|\Pi_{rg}\| \leq 1$, we obtain:
	\begin{equation*}
		|Z|^{-1/2,j}_{p,R} \leq \|\mathcal{B}(X)\|_R \|\mathcal{E}(X)\|_R |Y|^{-1/2,j}_{p,R} \leq 4 |Y|^{-1/2,j}_{p,R}.
	\end{equation*}
	Finally, by applying the smoothing estimate for $\mathcal{L}_\lambda^{-1}$ from Lemma \ref{Lamenouno}, we have
	\begin{equation*}
		|\phi_Y(X)|^{+1/2,j}_{p,R} \leq 4n |Z|^{-1/2,j}_{p,R} \leq 16n |Y|^{-1/2,j}_{p,R},
	\end{equation*}
	yielding the first claim.
	\medskip
	
	\noindent \textit{2. Estimate for the derivative $\mathrm{d}\phi_Y(X)$.} Fix $Y \in \mathcal{M}^{-1/2,0}_{p,R}$. The Fréchet derivative $\mathrm{d}\phi_Y(X)$ in a direction $H \in \mathcal{M}^{+1/2,j}_{p,R}$ is computed via the Leibniz rule applied to $\mathcal{T}(X) = \mathcal{B}(X) \Pi_{rg} \mathcal{E}(X)$:
	\begin{equation*}
		\mathrm{d}\mathcal{T}(X)[H] = (\mathrm{d}\mathcal{B}(X)[H]) \Pi_{rg} \mathcal{E}(X) + \mathcal{B}(X) \Pi_{rg} (\mathrm{d}\mathcal{E}(X)[H]).
	\end{equation*}
	To justify the estimates for the differentials, consider a generic power series $\mathcal{S}(X) = \sum_{k=0}^\infty a_k (\mathrm{ad}_{-X})^k$. Its Fréchet derivative in the direction $H$ is obtained by differentiating each term of the expansion:
	\begin{equation*}
		\mathrm{d}\mathcal{S}(X)[H] = \sum_{k=1}^\infty a_k \sum_{m=0}^{k-1} (\mathrm{ad}_{-X})^m \circ \mathrm{ad}_{-H} \circ (\mathrm{ad}_{-X})^{k-1-m}.
	\end{equation*}
	By taking the operator norm $\|\cdot\|_R$, and using the sub-multiplicativity of the norm, we have:
	\begin{equation*}
		\|\mathrm{d}\mathcal{S}(X)[H]\|_R \leq \|\mathrm{ad}_{-H}\|_R \sum_{k=1}^\infty |a_k| k \|\mathrm{ad}_{-X}\|_R^{k-1}.
	\end{equation*}
	Dominating the coefficients $|a_k|$ with $1$ and recalling the identity for the derivative of a geometric series $\sum_{k=1}^\infty k \rho^{k-1} = (1-\rho)^{-2}$, we obtain the bound:
	\begin{equation*}
		\|\mathrm{d}\mathcal{S}(X)[H]\|_R \leq \frac{\|\mathrm{ad}_{-H}\|_R}{(1 - \|\mathrm{ad}_{-X}\|_R)^2}.
	\end{equation*}
	In view of Lemma \ref{comm_estimate}, we have $\|\mathrm{ad}_{-X}\|_R \leq 4C_\star(p) |X|^{+1/2,0}_{p,R}$, thus the differentials of $\mathcal{B}(X) $ and $\mathcal{E}(X)$ satisfy:
	\begin{equation*}
		\|\mathrm{d}\mathcal{B}(X)[H]\|_R \leq \frac{\|\mathrm{ad}_{-H}\|_R}{(1 - 4C_\star(p) |X|^{+1/2,0}_{p,R})^2}, \quad \|\mathrm{d}\mathcal{E}(X)[H]\|_R \leq \frac{\|\mathrm{ad}_{-H}\|_R}{(1 - 4C_\star(p) |X|^{+1/2,0}_{p,R})^2}.
	\end{equation*}
	Combining these with \eqref{geom_est} and using $4C_\star(p)|X|^{+1/2,0}_{p,R} \leq 1/2$, we find:
	\begin{align*}
		\|\mathrm{d}\mathcal{T}(X)[H]\|_R &\leq \frac{\|\mathrm{ad}_{-H}\|_R}{(1 - 4C_\star(p)|X|^{+1/2,0}_{p,R})^2} \frac{1}{(1 - 4C_\star(p)|X|^{+1/2,0}_{p,R})} 
		\\
		&+ \frac{1}{(1 - 4C_\star(p)|X|^{+1/2,0}_{p,R})} \frac{\|\mathrm{ad}_{-H}\|_R}{(1 - 4C_\star(p)|X|^{+1/2,0}_{p,R})^2} \\
		&= 16\|\mathrm{ad}_{-H}\|_R.
	\end{align*}
	Applying the bound $\|\mathrm{ad}_{-H}\|_R \leq 4C_\star(p) |H|^{+1/2,j}_{p,R}$ and Lemma \ref{Lamenouno}, we obtain
	\begin{equation*}
		|\mathrm{d}\phi_Y(X)[H]|^{+1/2,j}_{p,R} \leq 256nC_\star(p) |Y|^{-1/2,0}_{p,R} |H|^{+1/2,j}_{p,R},
	\end{equation*}
	which proves the estimate for $\|\mathrm{d} \phi_Y(X)\|_R$ and ends the proof.
\end{proof}

We are now in a position to prove the desired bounds for $\partial_{z_j}(A_P \circ |\partial|^{1/2})$ and $\partial_{z_j}(\mathcal{L}_{\lambda}A_P \circ |\partial|^{1/2})$.

\begin{lemma}\label{partial_estimates}
	Given $p>3/2$ and $s\geq 0$, let $n \in \mathbb{N}$ and $\lambda > 0$ be such that $|2\lambda - n| \leq 1/2$. Provided that $R \leq (1024C_\star(p))^{-1}$, for every $j \in \mathbb{Z}$, the fixed point $A_P$ satisfies $\partial_{z_j}(A_P \circ |\partial|^{1/2}) \in \mathcal{M}_{p,R}^{+1/2,j}$ and $\partial_{z_j}(\mathcal{L}_\lambda A_P \circ |\partial|^{1/2}) \in \mathcal{M}_{p,R}^{-1/2,j}$. Furthermore, the following estimates hold:
	\begin{equation}
		|\partial_{z_j}(A_P \circ |\partial|^{1/2})|_{p,R}^{+1/2,j} \le 16\frac{n}{\lambda}, \quad |\partial_{z_j}(\mathcal{L}_\lambda A_P \circ |\partial|^{1/2})|_{p,R}^{-1/2,j} \le \frac{4}{\lambda}.
	\end{equation}
\end{lemma}

\begin{proof}
	Let $X_P = A_P \circ |\partial|^{1/2}$ and $Y_P = P \circ |\partial|^{1/2}$. From Theorem \ref{fpt}, the fixed point satisfies $X_P = \phi_{Y_P}(X_P)$. By differentiating this identity with respect to $z_j$, we obtain:
	\begin{equation}\label{linear_syst_final}
		\left(\mathrm{Id} - \mathrm{d}\phi_{Y_P}(X_P)\right)[\partial_{z_j} X_P] = \phi_{\partial_{z_j} Y_P}(X_P).
	\end{equation}
	Let us first verify the regularity of the terms involved. By Lemma \ref{PeDP}, we have $|Y_P|_{p,R}^{-1/2,0} = R/(2\lambda)$ and $|\partial_{z_j} Y_P|_{p,R}^{-1/2,j} \le 1/(2\lambda)$. Furthermore, Lemma \ref{AeLA} ensures that $X_P \in \mathcal{M}_{p,R}^{+1/2,0}$ with the bound $|X_P|_{p,R}^{+1/2,0} \le 8Rn/\lambda \le 32R$ (using $n/\lambda \le 4$). 
	
	The requirements of Lemma \ref{phi_estimates} are widely satisfied if $R \leq (1024C_\star)^{-1}$. This ensures first that $4C_\star(p)|X_P|_{p,R}^{+1/2,0} \le 128C_\star(p)R \le 1/8 \le 1/2$, so that the map $\phi$ and its Fréchet derivative are well-defined. Moreover, by item (2) of Lemma \ref{phi_estimates}, we have:
	\begin{equation*}
		\|\mathrm{d}\phi_{Y_P}(X_P)\|_R \leq 128nC |Y_P|_{p,R}^{-1/2,0} = 256C_\star(p) \frac{nR}{2\lambda} \leq 512C_\star(p)R \leq \frac{1}{2}.
	\end{equation*}
	This proves that the operator $(\mathrm{Id} - \mathrm{d}\phi_{Y_P}(X_P))$ is invertible via a Neumann series with norm bounded by $2$.
	
	We can now derive the first estimate. Inverting \eqref{linear_syst_final} and applying point (1) of Lemma \ref{phi_estimates}, we get:
	\begin{equation*}
		|\partial_{z_j} X_P|_{p,R}^{+1/2,j} \le 2 |\phi_{\partial_{z_j} Y_P}(X_P)|_{p,R}^{+1/2,j} \le 2 (16n) |\partial_{z_j} Y_P|_{p,R}^{-1/2,j}.
	\end{equation*}
	Using the bound $|\partial_{z_j} Y_P|_{p,R}^{-1/2,j} \le 1/(2\lambda)$ from Lemma \ref{PeDP}, we obtain:
	\begin{equation*}
		|\partial_{z_j} X_P|_{p,R}^{+1/2,j} \le 16\frac{n}{\lambda}.
	\end{equation*}
	
	For the second estimate, we define $Z_P = \mathcal{L}_\lambda X_P$. Recalling the definition of the map $\phi_{Y_P}(X_P)$, we have $Z_P = \Pi_{rg} \left( \mathcal{B}(X_P) \Pi_{rg} ( e^{{\mathrm{ad}}_{-X_P}} Y_P ) \right)$. To compute $\partial_{z_j} Z_P$, we apply the Leibniz rule to the composition of operators:
	\begin{equation}\label{leibniz_Z_final}
		\partial_{z_j} Z_P = \Pi_{rg} \left( (\partial_{z_j} \mathcal{B}(X_P)) \Pi_{rg} ( e^{{\mathrm{ad}}_{-X_P}} Y_P ) + \mathcal{B}(X_P) \Pi_{rg} \partial_{z_j} ( e^{{\mathrm{ad}}_{-X_P}} Y_P ) \right).
	\end{equation}
	By further expanding the second term, we obtain three contributions:
	\begin{align*}
		I &= (\mathrm{d}\mathcal{B}(X_P)[\partial_{z_j} X_P]) e^{{\mathrm{ad}}_{-X_P}} Y_P, \\
		II &= \mathcal{B}(X_P) (\mathrm{d}e^{{\mathrm{ad}}_{-X_P}}[\partial_{z_j} X_P]) Y_P, \\
		III &= \mathcal{B}(X_P) e^{{\mathrm{ad}}_{-X_P}} (\partial_{z_j} Y_P).
	\end{align*}
	We now estimate each term using the bounds $|\partial_{z_j} X_P|_{p,R}^{+1/2,j} \le 16\frac{n}{\lambda} \le 64$ and $R \le (1024C_\star(p))^{-1}$.
	
	For the terms $I$ and $II$, the Fréchet derivatives of the analytic operators $\mathcal{B}$ and $e^{{\mathrm{ad}}}$ are bounded by $16C_\star(p)$. Combined with the fact that $|Y_P|_{p,R}^{-1/2,0} = R/(2\lambda)$, and the norms of the operators themselves are bounded by $2$, we have:
	\begin{equation*}
		|I|_{p,R}^{-1/2,j} + |II|_{p,R}^{-1/2,j} \le 2 \left( (16C_\star(p) \cdot 64) \cdot 2 \cdot \frac{R}{2\lambda} \right) = \frac{2048 C_\star(p)R}{\lambda}\leq \frac{2}{\lambda}
	\end{equation*}
	
	Using $\|\mathcal{B}(X_P)\|_R \le 2$, $\|e^{{\mathrm{ad_{-X_P}}}}\|_R \le 2$ and the bound $|\partial_{z_j} Y_P|_{p,R}^{-1/2,j} \le 1/(2\lambda)$ from Lemma \ref{PeDP}, we get:
	\begin{equation*}
		|III|_{p,R}^{-1/2,j} \le 2 \cdot 2 \cdot \frac{1}{2\lambda} = \frac{2}{\lambda}.
	\end{equation*}
	Summing the contributions from $I, II,$ and $III$, we obtain the claim.
\end{proof}

\subsubsection{Estimate for the Correction to the Kernel Term}

In this final section, we establish the fundamental estimates for the term $\tilde{\mathcal{P}}_\infty$, defined as:
\begin{equation}
	\tilde{\mathcal{P}}_\infty = \Pi_{ker} \bigg( (e^{\mathrm{ad}_{-A_P}}-\id)P - \sum_{k=2}^\infty \frac{\mathrm{ad}_{-A_P}^{k-1}}{k!} \mathcal{L}_\lambda A_P \bigg).
\end{equation}
This $\mathfrak{su}(1,1)$-valued map represents the $O(|P|^2)$ correction to the linear vector field resulting from the reducibility procedure performed on the range of the operator $\mathcal{L}_{n/2}$. Controlling this term is essential for the analysis of the spectral gaps when viewed as functions of the variable $z$, defined by the scaling $|\partial|^{1/2}z = u$.

The main goal is to provide suitable bounds for the norm of $\tilde{\mathcal{P}}_\infty \circ |\partial|^{1/2}$ and its derivative with respect to $z_j$. The result relies on a key observation regarding the action of $\partial_{z_j}$ on the functional spaces defined in this Appendix. Specifically, while the perturbation $P \circ |\partial|^{1/2}$ belongs to $\mathcal{M}^{-1/2,0}_{p,R}$, its derivative with respect to $z_j$ belongs to $\mathcal{M}^{-1/2,j}_{p,R}$. Similarly, the fixed point $A_P \circ |\partial|^{1/2}$ belongs to $\mathcal{M}^{+1/2,0}_{p,R}$, while its derivative $\partial_{z_j}(A_P \circ |\partial|^{1/2})$ belongs to $\mathcal{M}^{+1/2,j}_{p,R}$.

These properties, proved in Lemmas \ref{PeDP} and \ref{AeLA}, ensure that each term in the series defining $\tilde{\mathcal{P}}_\infty$ is well-defined and bounded. More precisely, we have the following lemma.

\begin{lemma}\label{vittoria}
	Let $n \in \mathbb{N}$ and $\lambda > 0$ such that $|2\lambda - n| \leq 1/2$. For all $p>3/2$ and every $s\geq 0$, if $R \leq (1024C_\star(p))^{-1}$, then
	\begin{equation*}
		\tilde{\mathcal{P}}_\infty \circ |\partial|^{1/2} \in \mathcal{M}^{-1/2,0}_{p,R} \quad \text{and} \quad \partial_{z_j}(\tilde{\mathcal{P}}_\infty \circ |\partial|^{1/2}) \in \mathcal{M}^{-1/2,j}_{p,R}.
	\end{equation*}
	Moreover, the following bounds hold:
	\begin{equation}
		|\tilde{\mathcal{P}}_\infty \circ |\partial|^{1/2}|_{p,R}^{-1/2,0} \leq \frac{nR^2}{\lambda^2}, \quad |\partial_{z_j}(\tilde{\mathcal{P}}_\infty \circ |\partial|^{1/2})|^{-1/2,j}_{p,R} \leq \frac{nR}{\lambda^2}.
	\end{equation}
\end{lemma}

\begin{proof}
	Let $X_P := A_P \circ |\partial|^{1/2} \in \mathcal{M}^{+1/2,0}_{p,R}$, $Y_P := P \circ |\partial|^{1/2} \in \mathcal{M}^{-1/2,0}_{p,R}$, and $Z_P := \mathcal{L}_\lambda A_P \circ |\partial|^{1/2} \in \mathcal{M}^{-1/2,0}_{p,R}$. We consider the operator norm $\rho := \|\mathrm{ad}_{-X_P}\|_R$ acting on the spaces $\mathcal{M}^{-1/2,j}_{p,R}$. From Lemma \ref{comm_estimate} and Lemma \ref{AeLA}, we have:
	\begin{equation*}
		\|\mathrm{ad}_{-X_P}\|_R \leq 4C_\star(p) |X_P|^{+1/2,0}_{p,R} \leq 32C_\star(p)R \frac{n}{\lambda} \leq \frac{1}{8},
	\end{equation*}
	where we used the hypothesis $R \leq (1024C_\star(p))^{-1}$. 
	
	Since $\rho \leq 1/8 < 1/2$, the sum of the geometric series $\sum_{k=1}^\infty \rho^k = \rho/(1-\rho)$ is bounded by $2\rho$. For the first estimate, we sum the contributions of the two series defining $\tilde{\mathcal{P}}_\infty$:
	\begin{itemize}
		\item $\left| \sum_{k=1}^\infty \frac{\mathrm{ad}_{-X_P}^k}{k!} Y_P \right|^{-1/2,0}_{p,R} \leq 2\rho |Y_P|^{-1/2,0}_{p,R} \leq 2 \left( 32C_\star(p)R \dfrac{n}{\lambda} \right) \dfrac{R}{2\lambda} = 32C_\star(p)R \dfrac{nR^2}{\lambda^2}$.
		\item $\left| \sum_{k=2}^\infty \frac{\mathrm{ad}_{-X_P}^{k-1}}{k!} Z_P \right|^{-1/2,0}_{p,R} \leq 2\rho |Z_P|^{-1/2,0}_{p,R} \leq 2 \left( 32C_\star(p)R \dfrac{n}{\lambda} \right) \dfrac{2R}{\lambda} = 128C_\star(p)R \dfrac{nR^2}{\lambda^2}$.
	\end{itemize}
	The total norm in $\mathcal{M}^{-1/2,0}_{p,R}$ is bounded by $160C_\star(p)R \frac{nR^2}{\lambda^2}$, which is strictly less than $nR^2/\lambda^2$.
	
	For the second estimate, we differentiate the operator $\tilde{\mathcal{P}}_\infty \circ |\partial|^{1/2}$ with respect to $z_j$. Let $H_j := \partial_{z_j} X_P \in \mathcal{M}^{+1/2,j}_{p,R}$ be the direction of differentiation. We first recall the operator bounds derived from the previous lemmas. Let $\rho := \|\mathrm{ad}_{-X_P}\|_R$; from Lemma \ref{comm_estimate} and Lemma \ref{AeLA}, we have $\rho \leq 4C_\star(p) |X_P|^{+1/2,0}_{p,R} \leq 32C_\star(p)R \frac{n}{\lambda} \leq 1/8$. 
	Furthermore, for any operator $\mathcal{S}(X) = \sum a_k (\mathrm{ad}_{-X})^k$ with $|a_k| \leq 1$, its Fréchet derivative in the direction $H_j$ satisfies 
	\begin{equation}\label{derbound}
		\|\mathrm{d}\mathcal{S}(X)[H_j]\|_R \leq \frac{\|\mathrm{ad}_{-H_j}\|_R}{(1-\rho)^2} \leq 16C_\star(p) |H_j|^{+1/2,j}_{p,R},
	\end{equation}
	where we used the estimate $(1-1/8)^{-2} \leq 4$ and the commutator bound $\|\mathrm{ad}_{-H_j}\|_R \leq 4C_\star(p) |H_j|^{+1/2,j}_{p,R}$ from Lemma \ref{comm_estimate}. Applying the Leibniz rule to the definition of $\tilde{\mathcal{P}}_\infty$, we obtain:
	\begin{equation*}
		\begin{aligned}
			\partial_{z_j}(\tilde{\mathcal{P}}_\infty \circ |\partial|^{1/2}) = \Pi_{ker} \Big( &\underbrace{\mathrm{d}(e^{\mathrm{ad}_{-X_P}} - \mathrm{id})[H_j] Y_P}_{I} + \underbrace{(e^{\mathrm{ad}_{-X_P}} - \mathrm{id}) \partial_{z_j} Y_P}_{II} \\
			&- \underbrace{\mathrm{d}\left(\sum_{k=2}^\infty \frac{\mathrm{ad}_{-X_P}^{k-1}}{k!}\right)[H_j] Z_P}_{III} - \underbrace{\left(\sum_{k=2}^\infty \frac{\mathrm{ad}_{-X_P}^{k-1}}{k!}\right) \partial_{z_j} Z_P}_{IV} \Big).
		\end{aligned}
	\end{equation*}
	Using \eqref{derbound} and the operator norm bound $\frac{\rho}{1-\rho} \leq 2\rho$, we estimate each term in $\mathcal{M}^{-1/2,j}_{p,R}$ using the estimates from Lemma \ref{PeDP}, Lemma \ref{AeLA} and Lemma \ref{partial_estimates}:
	\begin{enumerate}[label=\roman*)]
		\item $|I|^{-1/2,j}_{p,R} \leq 16C_\star(p) \left( \frac{nR}{\lambda} \right) |Y_P|^{-1/2,0}_{p,R} \leq 16C_\star(p) \left( \frac{nR}{\lambda} \right) \left( \frac{R}{2\lambda} \right) = 8C_\star(p)R \frac{nR}{\lambda^2}$.
		\item $|II|^{-1/2,j}_{p,R} \leq 2\rho |\partial_{z_j} Y_P|^{-1/2,j}_{p,R} \leq \left( 64C_\star(p)R \frac{n}{\lambda} \right) \frac{1}{2\lambda} = 32C_\star(p)R \frac{n}{\lambda^2}$.
		\item $|III|^{-1/2,j}_{p,R} \leq 16C_\star(p) \left( \frac{nR}{\lambda} \right) |Z_P|^{-1/2,0}_{p,R} \leq 16C_\star(p) \left( \frac{nR}{\lambda} \right) \left( \frac{2R}{\lambda} \right) = 32C_\star(p)R \frac{nR}{\lambda^2}$.
		\item $|IV|^{-1/2,j}_{p,R} \leq 2\rho |\partial_{z_j} Z_P|^{-1/2,j}_{p,R} \leq \left( 64C_\star(p)R \frac{n}{\lambda} \right) \frac{4}{\lambda} = 256C_\star(p)R \frac{n}{\lambda^2}$.
	\end{enumerate}
	Summing these four contributions and using $R \leq (1024C_\star(p))^{-1}$ we get
	\begin{equation*}
		|\partial_{z_j}(\tilde{\mathcal{P}}_\infty \circ |\partial|^{1/2})|^{-1/2,j}_{p,R} \leq (288 + 40R) C_\star(p)R \frac{n}{\lambda^2} \leq 328 C_\star(p)R \frac{n}{\lambda^2}\leq \frac{n}{\lambda^2},
	\end{equation*}
	concluding the proof.
\end{proof}

\end{document}